\documentclass[12pt,oneside]{amsbook}
\usepackage{LRTthesis}
\begin{document}

\thispagestyle{empty}
\centerline{Surgery on Paracompact Manifolds}
\centerline{by}
\centerline{Laurence R.~Taylor}
\bigskip

This is a \TeX'ed version of the author's 1971 thesis
from the University of California at Berkeley under the supervision 
of J.~B.~Wagoner.
It has never been published, but over the years I have received a
small but steady request for copies.
As I am running out of copies and the quality of the remaining ones is
not conducive to recopying or scanning I have decided to \TeX\ the manuscript
while I still have some copies in my possession. 
The result is a hyperlinked and searchable manuscript. 

I have resisted the temptation to rewrite the manuscript.
It was written in some haste in the summer of 1971 and the haste shows.
Forty five years of experience would certainly enable me to produce
a better manuscript, but might delay appearance for another forty five years.
No effort has been made to reproduce line breaks, kerning, fonts, etc. so
\TeX ing does produce a different looking manuscript from the original. 
The original has 231 pages (6 pages of front matter plus 225 pages of text ), this version only has 109. (This page and the last two bibliography pages are new.)
But, modulo new typing errors and minor corrections, this is a faithful copy.
Footnotes are used for added material unless indicated otherwise. 
The one original footnote is identified as such. 
In a few places outright errors have been corrected, both mine and others. 
For historical accuracy, the symbol \iff\ is retained for ``if and only if''.

The bibliography has been updated with some references to 
relevant works appearing after 1971. 
The last bibliography page lists some works citing this thesis. 

The references to theorems, lemmas, etc. have been updated. 
The old scheme was truly terrible (it is explained in the introduction). 
It was inspired by the scheme in Spanier's book \cite{bthirtyfive}. 
This sort of scheme made more sense in the pre{--}\TeX\ days when 
adding a lemma, theorem, etc. to an earlier section of the manuscript 
would involve locating and changing all references to subsequent 
lemmas, theorems, etc.
With this sort of scheme there was a good chance you could slip in an
additional item without having to change any current references. 
For example, at the end of a section you could record a result for later use.
The current scheme just numbers everything consecutively, but the old number
is listed in parentheses after the new one and citations include both the number and
a page number.

\bigskip
\hbox{
\vtop{\hsize=2in
\noindent
Department of Mathematics\hfill\\
University of Notre Dame\hfill\\
Notre Dame IN 46556, USA\vfill}
\hskip1in
\vtop{\hsize=2.5in
\noindent
taylor.2@nd.edu\hfill
\vskip14pt\noindent
https://www.nd.edu/\lower 4pt\hbox{\char`\~ taylor}\vfill}\hfill}

\vskip20pt
\noindent Partially supported by the N.S.F.
\newpage

\vbox{\title{Surgery on Paracompact Manifolds}
\author{Laurence R.~Taylor}
\begin{abstract}
In this thesis we solve the problem of surgery for an arbitrary, finite 
dimensional, paracompact manifold$^1$.
The problem of surgery is to decide whether, given a proper map
$f\colon M\to X$, a bundle (vector, PL{--}micro, or TOP{--}micro), $\nu$,
over $X$, and a stable bundle map $F\colon \nu_M\to\nu$ over $f$ 
($\nu_M$ is the normal bundle of $M$, so we must assume $M$ is respectively 
a differentiable, a PL, or a topological, finite dimensional, paracompact manifold),
we can find a cobordism $W$ with $\partial W= M\disjointunion N$, a proper map
$g\colon W\to X$ with $g\vert_M=f$, a stable bundle map 
$G\colon \nu_W\to\nu$ with $G\vert_{\nu_M}=F$, such that $g\vert_N$ is a
proper homotopy equivalence.

If this problem can be solved, we show this forces conditions on $X$, $\nu$ and $f$.
In particular, $X$ must be a Poincar\'e duality space (Chapter 2), $\nu$ must lift the
Spivak normal fibration of $X$, and $f$ must be degree $1$.

If $X$, $\nu$ and $f$ satisfy these conditions, there is a well{--}defined obstruction
to solving this problem if $m$, 
the dimension of $M$, is at least five (\fullRef{T.3.2.1}).
This obstruction lies in a naturally defined group, $L_m(X,w)$, and every element
of this group can be realized, in a specific fashion, as the obstruction to a surgery
problem, provided $m\geq6$ ( \fullRef{T.3.2.4}).
$L_m(X,w)$ depends only on the system of fundamental groups of $X$ 
(\fullRef{T.3.2.3}).

Finally, we have applications for paracompact manifolds along the same lines as
the compact case.
Perhaps the most interesting of these is the theoretical solution of the related questions
of when a Poincar\'e duality space has the proper homotopy type of a paracompact
manifold, and if a proper homotopy equivalence between paracompact manifolds can
be properly deformed to a homeomorphism, diffeomorphism, or PL{--}equivalence
(\fullRef{T.3.2.4}).

\noindent\vbox{\noindent\hrule width 1in depth 0pt height .4pt}

\noindent
$^1$ Of dimension at least $5$.
\end{abstract}

\maketitle
}
\newpage
\thispagestyle{empty}
\vbox{\tableofcontents\vfill}
\thispagestyle{empty}
\newpage
\thispagestyle{empty}
\setcounter{section}{0}
\section*{INTRODUCTION}
\newHead{INTRO}
\vskip20pt

The object of this work is to give an adequate theory
of surgery for paracompact manifolds and proper maps.
By adequate we mean first that it should contain the theory of surgery
for compact manifolds. Secondly, the theory should be general enough to
permit extensions of the theoretical results of compact surgery.

These objectives are largely realized.
We obtain surgery groups which characterize the problem in dimensions
greater than or equal to five.
These groups depend only on the proper $2${--}type of the problem.
Using these groups one can classify all paracompact manifolds of a 
given proper, simple homotopy type (see \cite{bthirtythree} or 
\cite{bten} for a definition of simple homotopy type).

The first chapter constitutes the chief technical results of this work. 
In \cite{bthirtythree}, Siebenmann gives a ``geometric'' characterization 
of proper homotopy equivalence (Proposition IV).
This characterization was also discovered by Farrell{--}Wagoner 
\cite{bnine} from whom I learned it.

In section 2 we develop an algebraic process to handle this 
characterization. 
In section 3 we apply this process to construct groups which are the
analogue of the homotopy and homology groups.
Thus we get actual groups measuring by how much a map fails to
be a proper homotopy equivalence.
These groups also satisfy a version of the Hurewicz and Namioka theorems,
so one can often use these homology groups, which satisfy a version of
excision, Mayer{--}Vietoris, etc.

In section 4 we construct a cohomology theory for our theory.
We get various products for this theory.
Section 5 is devoted to an analysis of simple homotopy type along the
lines set out by Milnor in \cite{btwentythree}.
Section 6 is devoted to constructing locally compact CW complexes 
with a given chain complex (see Wall \cite{bthirtyeight} for a
treatment of the compact case of this problem).

Chapter 2 is devoted to an analysis of Poincar\'e duality for paracompact
manifolds and its generalization to arbitrary locally compact, finite
dimensional CW complexes.

In Chapter 3 the actual surgeries are performed. 
It has been observed by several people (especially Quinn 
\cite{btwentynine} and \cite{bthirty}) that all the surgery one needs to
be able to do is the surgery for a pair $(X,\partial X)$ for which 
$\partial X\subseteq X$ is a proper $1${--}equivalence (in the compact
case this means the inclusion induces isomorphisms on components and
on $\pi_1$).
We do this in the first section.
In the second section, we sketch the general set up and applications 
of the theory of paracompact surgery.

A word or two is in order here about internal referencing. 
A reference reads from right to left, so that Corollary 3.4.1.5 is the fifth
corollary to the first theorem of section 4 in chapter 3.
If the reference is made from chapter 3 it would be Corollary 4.1.5,
and if from section 4, Corollary 1.5.
Theorem (Proposition, Lemma) 3.4.6 is the sixth theorem of section 4 of
chapter 3.

Perhaps we should also remark that our use of the term 
$n${--}ad agrees with the use of the term in 
Wall \cite{bfortyone} (see especially Chapter 0).
For an $n${--}ad, $K$, $\partial_iK$ denotes the $(n-1)${--}ad whose
total space is the $i^{\text{ th}}$ face of $K$ and with the $(n-1)${--}ad
structure induced by intersecting the other faces of $K$.
$\delta_iK$ is the $(n-1)${--}ad  obtained by deleting the $i^{\text{ th}}$
face.
$s_nK$ denotes the $(n+1)${--}ad obtained by making $K$ 
the $(n+1)^{\text{ st}}$ face (it can also be regarded as the $(n+1)${--}ad
$K\times I$, where $I$ has the usual pair structure).

Lastly, several acknowledgements are in order.
This thesis was written under the direction of J.~Wagoner,
to whom I am indebted for many suggestions during the preparation of
this work.
I am greatly indebted to him and to T.~Farrell for sharing their results and
intuition on proper homotopy with me at the very beginning.
Thanks are also due to G.~Cooke for many helpful discussions.
Many other friends likewise deserve thanks for their help.
The National Science Foundation should also be thanked for its support
during my graduate career.

\newpage
\chapter{The Proper Homotopy Category and Its Functors}
\section{Introduction, elementary results, and homogamous spaces}
\newHead{I.1}

The purpose of this chapter is to recall for the reader some of the basic
results we will need and to describe a ``good'' category in which to do
proper homotopy theory.

The notion of a proper map is clearly essential.
We define a map to be \emph{proper} \iff\ the inverse image of every closed
compact set is contained in a closed compact set.
We note that this definition is also found in Bredon \cite{btwo}, page 56.

With this definition of a proper map we immediately have the notions of
proper homotopy, proper homotopy equivalence, etc., 
and we can define the category of all topological spaces and proper maps.
Classically there are several functors which apply to this situation.
As examples we have sheaf cohomology with compact supports 
and Borel{--}Moore homology with closed supports 
(see Bredon \cite{btwo}).

We prefer to use singular theory whenever possible.
Here too we have cohomology with compact supports and homology with
locally finite chains.
Most of the results concerning such groups are scattered 
(or non{--}existent) in the literature.
As a partial remedy for this situation we will write out the definitions of 
these groups and at least indicate the results we need.

\begin{xDefinition}
A collection of subsets of $X$ is said to be \emph{locally finite} if every 
closed, compact subset of $X$ intersects only finitely many 
elements of this collection.
\end{xDefinition}

\begin{xDefinition}
$S_q^\locf(X;\Gamma)$, where $\Gamma$ is a local system of 
$R${--}modules on $X$ 
(see Spanier \cite{bthirtyfive} pages 58; 281{--}283),
is defined to be the $R${--}module which is the set of all formal sums
$\sum \alpha_\sigma \sigma$, where $\sigma$ is a singular $q${--}simplex
of $X$, and $\alpha_\sigma\in\Gamma\bigl(\sigma(V_0)\bigr)$ is
zero except for a set of $\sigma$ whose images in $X$ are locally finite.
\end{xDefinition}

$S^q(X;\Gamma)$ is the module of functions $\varphi$ assigning to every
singular $q${--}simplex $\sigma$ of $X$ an element 
$\varphi(\sigma)\in \Gamma\bigl(\sigma(V_0)\bigr)$.

For a family of supports $\psi$ on $X$ 
(see Bredon \cite{btwo} page 15 for a definition) let
$S^\psi_q(X;\Gamma)$ denote the submodule of $S_q^\locf(X;\Gamma)$
such that the union of all the images of the $\sigma$ 
occurring with non{--}zero coefficient in a chain lies in some element of
$\psi$.
$S^q_\psi(X;\Gamma)$ consists of the submodule of all functions 
$\varphi$ for which there exists an element $c\in\phi$ such that if
$\Image\, \sigma \cap c=\emptyset$, $\varphi(\sigma)=0$.

These modules become chain complexes in the usual fashion.
Note that for the family of compact supports, $\cmpsup$,
$S^\cmpsup_q(X;\Gamma)$ is just the ordinary singular chains with
local coefficients.

For a proper subspace $A\subseteq X$ (inclusion is a proper map)
we get relative chain groups $S^\psi_q(X, A;\Gamma)$ and
$S^q_\psi(X, A;\Gamma)$.
Actually proper subspace is sometimes stronger than we need;
i.e. $S^\cmpsup_q(X, A;\Gamma)$ and
$S^q(X, A;\Gamma)$ are defined for any $A\subseteq X$.
There is a similar definition for the chain groups of a (proper) $n${--}ad.

The homology of $S^\psi_\ast(X, A;\Gamma)$ will be denoted
$H^\psi_\ast(X, A;\Gamma)$ except when $\psi=\cmpsup$ when
we just write $H_\ast(X, A;\Gamma)$.
The homology of $S_\psi^\ast(X, A;\Gamma)$ will be written 
$H_\psi^\ast(X, A;\Gamma)$.

Now $S_\cmpsup^q(X, A;\Gamma)\subseteq S^q(X, A;\Gamma)$.
The quotient complex will be denoted 
$S_{\text{end}}^q(X, A;\Gamma)$ and its homology
$H_{\text{end}}^q(X, A;\Gamma)$.
We have similar definitions for proper $n${--}ads and also for homology.

We will next set out the properties of these groups we will use.
Some of the obvious properties such as naturality and long exact sequences
will be omitted.

\textbf{Cup products}: There is a natural cup product\hfill
\[H^q_\psi(X\Colon A_1,\cdots, A_n;\Gamma_1)\otimes
H^k(X\Colon A_{n+1},\cdots, A_m;\Gamma_2)\RA{\cup}
H^{q+k}(X\Colon A_{1}, \cdots, A_m;\Gamma_1\otimes\Gamma_2)
\]
for a proper $(m+1)${--}ad $(X\Colon A_1,\cdots A_m)$.
It is associative and commutative in the graded sense
(i.e. $a\cdot b=(-1)^{{\text{deg}}\>a\cdot{\text{deg}}\> b}b\cdot a$).

Since $S^q_\psi\subseteq S^q$, all this follows easily from the properties of the
ordinary cup product with local coefficients once one checks that
if a cochain was supported in $c\in\psi$, then its product with any other
cochain is supported in $c$ if one uses the Alexander{--}Whitney diagonal
approximation (Spanier \cite{bthirtyfive} page 250).
\footnote{I should have remarked here that the chain homotopies giving 
the associativity and the graded commutativity are correctly supported.}

{\bf Cross products}:
There are natural products
\[\begin{aligned}%
H^q_\psi(A\Colon A_1&,\cdots,A_n;\Gamma_1)\otimes
H^k(Y\Colon B_1,\cdots, B_m;\Gamma_2)\quad\RA{\qquad\times\qquad}\\
&H^{q+k}_{\pi_1^{-1}\psi}(X\times Y\Colon X\times B_1,\cdots,X\times B_m,
A_1\times Y,\cdots, A_m\times Y;\Gamma_1\otimes\Gamma_2)\\\end{aligned}\]
and
\[\begin{aligned}%
H_q^\psi(A\Colon A_1&,\cdots,A_n;\Gamma_1)\otimes
H_k^\locf(Y\Colon B_1,\cdots, B_m;\Gamma_2)\quad\RA{\qquad\times\qquad}\\
&H_{q+k}^{\psi\times Y}(X\times Y\Colon X\times B_1,\cdots,X\times B_m,
A_1\times Y,\cdots, A_m\times Y;\Gamma_1\otimes\Gamma_2)\\\end{aligned}\]
where 
$\pi_1^{-1}(\psi)=\{ K\subseteq X\times Y\ \vert\ \pi_1(K)\in\psi\}$
and $\psi\times Y=\{ K\times Y\subseteq X\times Y\ \vert\ K\in\psi\}$.
These satisfy the usual properties of the cross product.

We discuss this case in some detail.
Let us first define
\[\tau\colon S_n^\locf(X\times Y)\to \sum_{i+j=n}
S^{\cmpsup}_i(X)\ \widehat\otimes\ S^{\cmpsup}_j(Y)\]
where $\widehat\otimes$ is the completed tensor product,
i.e. infinite sums are allowed\footnote{As long as the resulting sum is locally finite.}.
If $\sigma\colon \Delta^n\to X\times Y$, and if $\pi_1$ and $\pi_2$
are the projections, 
$\displaystyle\tau(\sigma)=
\sum_{i+j=n}\ _i(\pi_1\sigma)\ \otimes\ (\pi_2\sigma)_j$
where $\ _i(\ )$ is the front $i${--}face and $(\ )_j$ is the back $j${--}face
(see Spanier \cite{bthirtyfive} page 250).
This extends over all of $S^\locf_n$ and is a natural chain map.

The cohomology cross product is then defined on the chain level by
$(c\times d)(\sigma)=c\bigl(\>_i(\pi_1\sigma)\bigr)\ \otimes
d\bigl((\pi_2\sigma)_j\bigr)$,
where $c$ is an $i${--}cochain, $d$ a $j${--}cochain, and $\sigma$ an
$(i+j)${--}chain.
One checks it has the usual properties.

We next define 
$\lambda\colon S^\locf_i(X)\ \otimes\ S^\locf_j(Y)\ \RA{\ }\ 
S^\locf_{i+j}(X\times Y)$ as follows.
Let $h_{i,j}\colon \Delta^{i+j}\ \to\ \Delta^i\times\Delta^j$
be a homeomorphism such that 
$\ _i(h_{i,j}) \colon \Delta^i\to\Delta^i\times\Delta^j$ is given by
$x\mapsto(x,0)$ and such that
$(h_{i,j})_j\colon \Delta^j\to\Delta^i\times\Delta^j$ is given by
$y\mapsto(0,y)$.
Define $\lambda(\sigma_X\otimes\sigma_Y)=
h_{i,j}\circ(\sigma_X\times\sigma_Y)$ and extend ``linearly'';
i.e. $\displaystyle\lambda\left(\sum \alpha\>\sigma_\alpha\otimes
\sum\beta\>\sigma_\beta\right)=
\sum_{\alpha,\beta}\ \alpha\otimes\beta\cdot
\lambda(\sigma_\alpha\otimes\sigma_\beta)$.
$\lambda$ then becomes a chain map, and the homology cross product
is then defined on the chain level as above.
It has the usual properties.

{\bf Slant product}:
There are natural products
\[\begin{aligned}%
H^q_\cmpsup(Y\Colon B_1,\cdots, B_m;&\Gamma_1)\otimes
H^\locf_{q+k}(X\times Y\Colon A_1\times Y,\cdots, A_n\times Y,
X\times B_1,\cdots, X\times B_m;\Gamma_2)
\\&\RA{\ \vert\ }\
H_k^\locf(X\Colon A_1,\cdots, A_n;\Gamma_1\otimes\Gamma_2)\\
\end{aligned}\]
and
\[\begin{aligned}%
H^q(Y\Colon B_1,\cdots, B_m;&\Gamma_1)\otimes
H_{q+k}(X\times Y\Colon A_1\times Y,\cdots, A_n\times Y,
X\times B_1,\cdots, X\times B_m;\Gamma_2)
\\&\RA{\ \vert\ }\
H_k(X\Colon A_1,\cdots, A_n;\Gamma_1\otimes\Gamma_2)\quad.\\\end{aligned}\]

The product is defined on the chain level by
\[c\vert \sigma = c\Big\vert {\sum \alpha\>\sigma_\alpha}=
\sum_\alpha\biggl(\sum_{i+j=q+k}
\ _i(\pi_1\sigma_\alpha)\ \otimes\ \Bigl(
c\bigl((\pi_2\sigma)_j\bigr)\otimes\alpha\Bigr)\biggr)\]
where $c$ applied to a chain is zero if the dimensions do not agree.
The slant product is natural on the chain level and has all the usual 
properties.

{\bf Cap product}:
There is a natural product
\[\begin{aligned}%
H^q_\psi(X\Colon A_1,\cdots,A_n;&\Gamma_1)\otimes 
H^\varphi_{q+k}(X\Colon A_1,\cdots,A_n,B_1,\cdots, B_m;\Gamma_2)\\
&\RA{\ \cap\ }\ 
H^{\varphi\cap\psi}_{k}(X\Colon B_1,\cdots,B_m;
\Gamma_1\otimes\Gamma_2)\quad .\\\end{aligned}\]
It is given by $u\cap z =u\vert d_\ast v$, where
$d\colon X\to X\times X$ is the diagonal map.
The cap product has all the usual properties.
We get better support conditions for our cap product than we did
for an arbitrary slant product because $d_\ast$ of a chain in $X\times X$
is ``locally finite'' with respect to sets of the form $c\times X$ and
$X\times c$ for any closed, compact $c\subseteq X$.

One of the most useful of the usual properties of the cap product is the

\setcounter{footnote}{0}
{\bf Browder Lemma}: 
(\cite{bthree}, \cite{bfour}).\footnote{Best reference is I.1.5 Theorem of \cite{Browder}.}
Let $(X,A)$ be a proper pair ($A$ is  a proper subspace), and let
$Z\in H^\psi_n(X,A;\Gamma_2)$.
Then $\partial Z\in H^\psi_{n-1}(A;\Gamma_2\vert_A)$ is defined.
The following diagram commutes, 
where $\Gamma_3=\Gamma_1\otimes\Gamma_2$
\[\begin{matrix}%
H^{\ast-1}_\varphi(A;\Gamma_1\vert_A)&\rsa&
H^\ast_\varphi(X,A;\Gamma_1)&\rsa&H^\ast_\varphi(X;\Gamma_1)
&\rsa&H^{\ast}_\varphi(A;\Gamma_1\vert_A)
\\
\downlabeledarrow[\bigg]{\cap(-1)^{n-1}\partial Z}{}&&\downlabeledarrow[\bigg]{\cap Z}{}
&&\downlabeledarrow[\bigg]{\cap Z}{}&&\downlabeledarrow[\bigg]{\cap \partial Z}{}\\
H^{\varphi\cap\psi}_{n-\ast}(A;\Gthree\vert_A)&\rsa&
H^{\varphi\cap\psi}_{n-\ast}(X;\Gthree)&\rsa&
H^{\varphi\cap\psi}_{n-\ast}(X,A;\Gthree)&\rsa&
H^{\varphi\cap\psi}_{n-\ast-1}(A;\Gthree\vert_A)\\
\end{matrix}\]

\LRTpageLabel{Universalcoefficientx}
In two cases, we also have a universal coefficient formula relating 
cohomology and homology.
We first have the ordinary universal coefficient formula; namely
\[0\to{\text{Ext}}\bigl(H_{\ast-1}(\_\ ;\Gamma),\Z\bigr)\to
H^\ast\bigl(\_\ ;{\Homx}(\Gamma,\Z)\bigr)\to
{\Homx}\bigl(H_\ast(\_\ ;\Gamma),\Z\bigr)\to 0\]
is split exact (see Spanier \cite{bthirtyfive}, page 283).

We have a natural chain map
\[\alpha\colon S^\locf_\ast\bigl(\_\ ,{\Homx}(\Gamma,\Z)\bigr)\ \to\ 
{\Homx}\bigl(S^\ast_\cmpsup(\_\ ;\Gamma),\Z\bigr)\]
given by $\alpha(c)(\varphi)=\varphi(c)$.
If the space $X$ is \HCLx\ Bredon \cite{btwo}*{II.9.23}, shows that $\alpha$ induces
a homology isomorphism, so we get
\[0\to{\text{Ext}}(H^{\ast+1}_\cmpsup(X ;\Gamma),\Z)\to
H_\ast^\locf(X ;{\Homx}(\Gamma,\Z))\to
{\Homx}(H^\ast_\cmpsup(X ;\Gamma),\Z)\to0\]
is split exact.

Write $\Bar{\Gamma}$ for ${\Homx}(\Gamma,\Z)$.
Then if $c\in H^k(\_\ ;\Gamma)$, the following diagram commutes
\[\begin{matrix}%
0\to&{\text{Ext}}(H_{\ast-1}(\_\ ;\Z),\Z)&\to&
H^\ast(\_\ ;\Z)&\to&{\Homx}(H_\ast(\_\ ;\Z),\Z)&\to0\\
&\downlabeledarrow[\bigg]{{\text{Ext}}(\cap c)}{}&&\downlabeledarrow[\bigg]{c\cup}{}
&&\downlabeledarrow[\bigg]{{\Homx}(\cap c)}{}\\
0\to&{\text{Ext}}(H_{\ast+k-1}(\_\ ;\Gamma),\Z)&\to&
H^{\ast+k}(\_\ ;\Bar\Gamma)&\to&
{\Homx}(H_{\ast+k}(\_\ ;\Gamma),\Z)&\to0\\
\end{matrix}\]

If $c\in H^\locf_k(\_\ ;\Gamma)$, and if the spaces in question are \HCLx,
the following diagram commutes
\[\begin{matrix}%
0\to&\text{Ext}(H_{k-\ast-1}(\_\ ;\Z),\Z)&\to&
H^{k-\ast}(\_\ ;\Z)&\to&{\Homx}(H_{k-\ast}(\_\ ;\Z),\Z)&\to0\\
&\downlabeledarrow[\bigg]{{\text{Ext}}(\cap c)}{}&&\downlabeledarrow[\bigg]{\cap c}{}
&&\downlabeledarrow[\bigg]{{\Homx}(\cap c)}{}\\
0\to&{\text{Ext}}(H^{\ast+1}_\cmpsup(\_\ ;\Gamma),\Z)&\to&
H^\locf_\ast(\_\ ;\Bar\Gamma)&\to&
{\Homx}(H^{\ast}_\cmpsup(\_\ ;\Gamma),\Z)&\to0\\
\end{matrix}\]
These formulas can actually be seen on the chain level by picking
representatives and using the Alexander{--}Whitney 
diagonal approximation.

These homology and cohomology groups enjoy other pleasant properties.
One which we shall exploit heavily throughout the remainder of this work 
is the existence of a transfer map for any arbitrary cover.
For particulars, let $\pi\colon \widetilde X\to X$ be a covering map.
Then we have homomorphisms
\[{\text{tr}}\colon H^\locf_\ast(X;\Gamma)\ \to\ 
H^\locf_\ast(\widetilde X;\pi^\ast\Gamma)\]
and 
\[{\text{tr}}\colon H_\cmpsup^\ast(\widetilde X;\pi^\ast\Gamma)
\ \to\  H_\cmpsup^\ast(X;\Gamma)
\ .\]
The first of these is given by defining ${\text{tr}}(\sigma)$ for a simplex
$\sigma$ and extending ``linearly.''
$\displaystyle{\text{tr}}(\sigma)=\sum_{p\in\pi^{-1}(v_0)}\sigma_p$,
where $p$ runs over all the points in $\pi^{-1}(v_0)$, where $v_0$
is a vertex of $\sigma$, and $\sigma_p$ is $\sigma$ 
lifted so that $v_0$ goes to $p$.
It is not hard to check ${\text{tr}}$ is a chain map.
For the cohomology trace define ${\text{tr}}(c)$ as the cochain whose value
on the simplex $\sigma$ in $X$ is $c\bigl({\text{tr}}(\sigma)\bigr)$; i.e.
$\bigl({\text{tr}}(c)\bigr)(\sigma)=c\bigl({\text{tr}}(\sigma)\bigr)$.
If $f\colon X\to Y$ is a proper map, and if $\pi\colon \widetilde Y\to Y$
is a cover, then, for the cover $\widetilde X\to X$ which is induced from
$\pi$ by $f$, $\tilde f_\ast({\text{tr}}Z)={\text{tr}}\tilde f_\ast Z$ and
${\text{tr}}(\tilde f^\ast c)=f^\ast({\text{tr}}\>c)$.

\insetitem{Warning}
The trace tends to be highly unnatural except in this one situation.

As an easy exercise, one may check that if 
$c\in H^k_\cmpsup(\widetilde X; \pi^\ast\Gamma_1)$
and if $Z\in H^\locf_{q+k}(X;\Gamma_2)$, then, in
$H_q(X;\Gamma_1\otimes\Gamma_2)$,
$\pi_\ast(c\cap {\text{tr}} Z)= {\text{tr}}\>c\ \cap Z$.

In the coming pages, we will want to study spherical fibrations and
paracompact manifolds.
For the former objects we have

{\bf Thom Isomorphism Theorem}:
Let $\xi$ be a spherical fibration of dimension $(q-1)$ over $B$.
Let $S(\xi)$ be its total space, and let $D(\xi)$ be the total space of the 
associated disc fibration.
Then there is a class 
$U_\xi\in H^q\bigl(D(\xi),S(\xi); p^\ast_\xi(\Gamma_\xi)\bigr)$
[where $p\colon D(\xi)\to B$ is the projection, and $\Gamma_\xi$
is the local system on $B$ given at $b\in B$ by
$H^q\bigl(p^{-1}(b),p^{-1}(b)\cap S(\xi);\Z) \bigr)$] such that
\[\cup U_\xi \ \colon H_\varphi^\ast(B;\Gamma)\ \to\
H^{\ast+q}_{p^{-1}(\varphi)}\bigl(
D(\xi), S(\xi); p^\ast(\Gamma\otimes\Gamma_\xi)\bigr)\]
is an isomorphism.
One also has
\[U_\xi\cap\ \colon H_\ast\bigl(D(\xi), S(\xi);p^\ast(\Gamma)\bigr)\ \to\ 
H_{\ast-q}(B;\Gamma_\xi\otimes\Gamma)\]
is an isomorphism.

Note that we have been (and will continue to be) a little sloppy.
If $c\in H^\ast_\varphi(B;\Gamma)$, $c\cup U_\xi$ should actually be
$p^\ast(c)\ \cup\ U_\xi$.
A similar notational amalgamation has occurred when we write $U_\xi\cap$.

This theorem is proved by a spectral sequence argument 
(see \cite{btwentysix}), so one need only check that we still have a Serre
spectral sequence with the appropriate supports.

For a paracompact manifold (i.e. a locally Euclidean, paracompact, Hausdorff
space), possibly with boundary, we have

{\bf Lefschetz Duality}: (\cite{btwenty}, \cite{bfortyfour}). 
If $(M,\partial M)$ is a paracompact manifold pair of dimension $n$, there is
a class $[M]\in H^\locf_n(M,\partial M;\Gamma_M)$ (where $\Gamma_M$
is the local system for the bundle $\nu$, the normal bundle of $M$)
such that the maps
\[\cap[M]\colon H^\ast_\psi(M,\partial M;\Gamma)\ \to\
H^\psi_{n-\ast}(M;\Gamma\otimes\Gamma_M)\]
and
\[\cap[M]\colon H^\ast_\psi(M;\Gamma)\ \to\
H^\psi_{n-\ast}(M,\partial M;\Gamma\otimes\Gamma_M)\]
are isomorphisms.

This completes the first objective of this section, so we turn to the second.
The functors above already give us much non{--}trivial information on
the category of all spaces and proper maps, but they are insufficient 
even to determine if a map is a proper homotopy equivalence on the
subcategory of locally compact, finite dimensional CW complexes, a
category in which we are surely going to be interested.
In fact, the next two sections will be concerned precisely with the
problem of constructing functors which will determine whether a map
is or is not a proper homotopy equivalence in this category.

If we restrict ourselves to finite complexes, the Whitehead Theorem
(\cite{bfortythree}) already provides the answer.
Notice that to solve the problem, even for finite complexes, we are
forced to consider homotopy, which means base points.
In order to solve the problem for locally finite complexes, we are going to
have to consider lots of base points simultaneously.
The category of spaces we are about to define is about the largest in
which we can place our points nicely.
It is also closed under proper homotopy equivalence.

\begin{xDefinition}
A set $B$ of points of $X$ is said to be a \emph{set of base points} for $X$
provided
\begin{enumerate}
\item[a)] every path component of $X$ contains a point of $B$
\item[b)] given any closed, compact set $c\subseteq X$, 
there is a closed compact set $D$ such that there is a point of $B$
in every path component of $X-c$ which is not contained in $D$.
\end{enumerate}
\end{xDefinition}

\begin{xDefinition}
A set of base points, $B$, for a path connected space $X$
is said to be \emph{irreducible} if,
for any set of base points $C$ for $X$ with $C\subseteq B$, the cardinality
of $C$ is equal to the cardinality of $B$.

A set of base points for any space $X$ is said to be \emph{irreducible} 
provided it is an irreducible set of base points for 
each path component of $X$.
\end{xDefinition}

\begin{xDefinition}
Two locally finite sets of points are said to be equivalent $(\sim)$
provided there is a 1{--}1 correspondence between the two sets which
is given by a locally finite set of paths.
\end{xDefinition}

\begin{xDefinition}
Consider the following two properties of a space $X$:
\begin{enumerate}
\item[]
\begin{enumerate}
\item[1)] Every set of base points for $X$ has an irreducible,
locally finite subset.
\item[2)] Any two irreducible, locally finite sets of base points for $X$
are equivalent.
\end{enumerate}
\item[] A space $X$ is said to be \emph{\prex{--}homogamous}
\footnote[1]{\text{In the original manuscript we started by calling this property 
homogamous and then}\\\text{redefined the term at the end of this section.}}
provided $X\times I$ satisfies 1) and 2).
\end{enumerate}
\end{xDefinition}

\bigskip
\BEGIN{P.1.1.1}
If $X$ has the proper homotopy type of a \hfill\\\prex{--}homogamous space,
then $X$ has properties 1) and 2).
\end{Proposition}
\medskip

\begin{proof}
We first prove two lemmas.

\BEGIN{L.1.1.1}
Let $f\colon X\to Y$ be a proper map which induces injections of
$H^0(Y)$ into $H^0(X)$ and 
of $H^0_{{\text{end}}}(Y)$ into $H^0_{{\text{end}}}(X)$.
Then if $\{p\}$ is a set of base points for $X$, then $\{f(p)\}$ is a set
of base points for $Y$.
\end{Lemma}

\begin{proof}
Since $f$ induces an injection on $H^0$, there is an $f(p)$ in
every path component of $Y$.

Now look at the path components of $Y-c$, where $c$ is some closed,
compact subset of $Y$.
Let $\{W_\alpha\}$ be the set of path components of $Y-c$ such that
$f^{-1}(W_\alpha)$ contains no point of $\{ p\}$.
Since $\{p\}$ is a set of base points for $X$, 
$\displaystyle\mathop{\cup}_\alpha f^{-1}(W_\alpha)\subseteq D$, 
where $D$ is some closed, compact subset of $X$. 
Then $f(X-D)\cap W_\alpha=\emptyset$ for all $\alpha$.

Define a cochain $\beta$ by
\[\beta(q)=\begin{cases}1& q\in W_\alpha\\ 0& q\notin W_\alpha\\\end{cases}\ .\]
Then $\delta\beta(\lambda)=
\beta\bigl(\lambda(1)\bigr)-\beta\bigl(\lambda(0)\bigr)=0$ 
if $\lambda\cap c = \emptyset$.
Hence $\delta\beta=0$ in $S^1_{{\text{end}}}(Y;\Z)$.
But since $f(X-D)\cap W_\alpha=\emptyset$, $f^\ast\beta=0$ in 
$S^0_{{\text{end}}}(X;\Z)$. 
Since $f^\ast$ is an injection on $H^0_{{\text{end}}}$, $\beta=0$ in
$H^0_{{\text{end}}}(Y;\Z)$.
But this implies $\displaystyle\mathop{\cup}_\alpha W_\alpha$
is contained in some compact set. 
\end{proof}

\medskip
\BEGIN{L.1.1.2}
Let $f$ be a map properly homotopic to the identity. 
Let $\{p\}$ be a locally finite set of base points.
Then $\{f(p)\}$ is equivalent to a subset of $\{p\}$.
\end{Lemma}
\begin{proof}
We have $F\colon X\times I\to X$ a proper map.
The set $\{p\times I\}$ is clearly locally finite.
Since $F$ is proper, $\{F(p\times I)\}$ is easily seen to be locally finite.
But $\{F(p\times I)\}$ provides an equivalence between $\{f(p)\}$
and some subset of $\{p\}$ [more than one $p$ may go to a given
$f(p)$]. 
\end{proof}

Now let $X$ have the proper homotopy type of $Y$, 
a \prex{--}homogamous space.
Hence we have proper maps $f\colon X\to Y$ and $g\colon Y\to X$
with the usual properties.

Let $\{p\}$ be a set of base points for $X$.
Then by \fullRef{L.1.1.1},
$\{f(p)\}$ is a set for $Y$ and 
$\{f(p)\times 0\}$ is a set for $Y\times I$.
Since $Y$ is \prex{--}homogamous, there is an irreducible, locally finite
subset $\{ f(p^\prime)\times 0\}$.
By 
\fullRef{L.1.1.1},
$\{g\circ f(p^\prime)\}$ is a locally finite set 
of base points for $X$.
But by 
\fullRef{L.1.1.2},
there is a further refinement, $\{p^{\prime\prime}\}$,
of $\{p\}$ such that $\{p^{\prime\prime}\}\sim \{g\circ f(p^\prime)\}$.
But then $\{p^{\prime\prime}\}$ is easily seen to be a set of base points 
also.
Now $\{p^{\prime\prime}\}$ is in 1{--}1 correspondence 
with $\{f(p^{\prime\prime})\}$, and $\{f(p^{\prime\prime})\times 0\}$
is a set of base points for $Y\times I$ by \fullRef{L.1.1.1}.
$\{f(p^{\prime\prime})\times 0\}$ is 
a subset of $\{f(p^{\prime})\times 0\}$ and is thus irreducible.
Hence $\{p^{\prime\prime}\}$ is easily seen to be irreducible, and 
therefore $X$ satisfies 1).

Let $\{p\}$ be an irreducible, locally finite set of base points for $X$.
We claim that there is an irreducible, locally finite set of base points $\{q\}$
for $Y\times I$ such that $\{p\}\sim\{g\circ\pi(q)\}$, where
$\pi\colon Y\times I\to Y$ is projection.

By the argument in \fullRef{L.1.1.2},
we see that we have a locally finite set
of paths $\{\lambda_p\}$ from $\{p\}$ to
$\{g\circ f(p)\}$.
However, $(g\circ f)^{-1}(g\circ f)(p)$ may contain more points of 
$\{p\}$ than just $p$.
But since $\{\lambda_p\}$ is locally finite, there are 
only finitely many such points, say $p_1$, \dots, $p_n$.
Let $q=f(p)\times 0$ and define $q_i=f(p)\times 1/i$ for $1\leq i\leq n$.
The resulting set of points, $\{q\}$ is easily seen to be locally finite, and 
by several applications of \fullRef{L.1.1.1},
$\{q\}$ is an irreducible set of
base points for $Y\times I$.

So suppose given $\{p\}$ and $\{p^\prime\}$, irreducible, locally
finite sets of base points for $X$.
Pick $\{q\}$ and $\{q^\prime\}$ as above to be irreducible, locally finite
sets of base points for $Y\times I$.
Since $Y$ is \prex{--}homogamous, $\{q\}\sim\{q^\prime\}$,
so $\{g\circ\pi(q)\}\sim\{g\circ\pi(q^\prime)\}$.
Thus $\{p\}\sim\{p^\prime\}$, so $X$ satisfies 2). 
\end{proof}

\medskip
\BEGIN{C.1.1.1.1.1}
A space which is the proper homotopy type of a \prex{--}homogamous space
is \prex{--}homogamous.
\end{Corollary}
\medskip
\BEGIN{C.1.1.1.2}
The mapping cylinder of a proper map whose range is \prex{--}homogamous  
is \prex{--}homogamous.
\end{Corollary}

\bigskip
\BEGIN{P.1.1.2}
Let $\{\mathcal O\}$ be a locally finite open cover of $X$. 
Further assume that each $\mathcal O$ is path connected and that
each $\Bar{\mathcal O}$ is compact.
Then $X$ is \prex{--}homogamous.
\end{Proposition}

\medskip
\BEGIN{C.1.1.2.1}
A locally compact, locally path connected, 
paracompact space is \prex{--}homogamous.
\end{Corollary}
\medskip
\BEGIN{C.1.1.2.2}
A locally compact CW complex is \prex{--}homogamous.
\end{Corollary}
\medskip
\BEGIN{C.1.1.2.3}
A paracompact topological manifold is \prex{--}homogamous.
\end{Corollary}
\bigskip
\begin{proof}
If $\{\mathcal O\}$ is the collection for $X$, $\{{\mathcal O}\times I\}$ is
a cover for $X\times I$ with the same properties, so, if we can show
1) and 2) hold for $X$, we are done.

Since each $\mathcal O$ is path connected, each path component of $X$ is open.
Also the complement of a path component is open, so each path component
is both open and closed.
Hence $X$ is \prex{--}homogamous \iff\ each path component is, so we assume
$X$ is path connected.

We claim $X$ is $\sigma${--}compact, i.e. the countable union 
of compact sets.
In fact, we will show $\{{\mathcal O}\}$ is at most countable.
As a first step, define a metric $d$ on $X$ as follows.
If $p\neq q$, look at a path $\lambda$ from $p$ to $q$.
$\lambda$ is compact, so it is contained in a finite union of ${\mathcal O}$'s.
Hence $\lambda$ is contained in a closed, compact set
so $\lambda$ intersects only finitely many ${\mathcal O}$'s.
Let $r(\lambda;p,q) =$ the number of ${\mathcal O}$'s that $\lambda$
intersects (non{--}empty).
Define $\displaystyle d(p,q)=\mathop{{\text{min}}}_{\lambda}
r(\lambda;p,q)$.
This is a natural number, so there is actually some path, $\lambda$,
such that $d(p,q)=r(\lambda;p,q)$.
If $p=q$, set $d(p,q)=0$.
$d$ is easily seen to be a metric.

Let us fix $p\in X$.
Then to each ${\mathcal O}$ we can associate a number 
$\displaystyle m({\mathcal O},p)=\mathop{{\text{min}}}_{q\in{\mathcal O}}d(p,q)$.
We claim that, for any $n$, $m({\mathcal O},p)\leq n$ for only finitely many
${\mathcal O}$.
For $n=0$ this is an easy consequence of the fact that 
$\{{\mathcal O}\}$ is locally finite.
Now induct on $n$.
Let ${\mathcal O}_1$, \dots,${\mathcal O}_k$ be all the ${\mathcal O}$'s such that
$m({\mathcal O},p)\leq n-1$.
Let $\displaystyle c=\mathop{\cup}_{i=1}^k \Bar{\mathcal O}_i$.
$c$ is compact.

Suppose $\Bar{\mathcal O}\cap c=\emptyset$. 
Then we claim $m({\mathcal O},p)\geq n+1$.
To see this, pick $q\in {\mathcal O}$, and any path $\lambda$ from $p$
to $q$.
If we can show $r(\lambda;p,q)\geq n+1$, we are done.
Let $[0,x]$ be the closed interval which is the first component of
$\lambda^{-1}(c)$, where $\lambda\colon I\to X$ is the path.
Since $c\cap \Bar{\mathcal O}=\emptyset$, 
$\lambda^{-1}(\Bar{\mathcal O})\geq s$, where $s>x$. 
Pick $x<t<s$.
Then $\lambda(t)\notin c$, so the path from $p$ to $\lambda(t)$
already intersects at least $n$ of the ${\mathcal O}$'s, so, from
$p$ to $q$ it must intersect at least $n+1$.

Therefore, if $m({\mathcal O},p)\leq n$, $\Bar{\mathcal O}\cap c\neq\emptyset$.
But since $\{{\mathcal O}\}$ is locally finite, there are only finitely many
${\mathcal O}$ for which this is true.
This completes the induction.

Hence the cover $\{{\mathcal O}\}$ is at most countable.
If $\{{\mathcal O}\}$ is finite, $X$ is compact and hence easily seen to satisfy
1) and 2).
Hence we assume $\{{\mathcal O}\}$ is infinite.

Enumerate $\{{\mathcal O}\}$, and set 
$\displaystyle C_k=\mathop{\cup}_{i=0}^k \Bar{\mathcal O}_i$.
Since $C_k$ is compact, there are but finitely many ${\mathcal O}$'s such that
$\Bar{\mathcal O}\cap C_k\neq\emptyset$.
Let $E$ be the union of $c$ and these ${\mathcal O}$'s.
Then $E$ is compact, as is $\partial E$, the frontier of $E$ in $X$.
Let $\{W_\alpha\}$ be the path components of $X-C_k$ not
contained entirely in $E$.

Look at $W_\alpha\cap \partial E$.
It might be empty, in which case  $W_\alpha$ is actually a component of $X$
since $\partial E$ separates the interior of $E$ and $X-E$.
But $X$ is connected, so $W_\alpha\cap\partial E\neq\emptyset$.
Now if $p\in\partial E$, $p\in{\mathcal O}$ with ${\mathcal O}\cap C_k=\emptyset$.
Now ${\mathcal O}$ is a path connected set missing $C_k$ with ${\mathcal O}$ not
contained entirely in $E$, so ${\mathcal O}\subseteq W_\alpha$ for
some $\alpha$.
Hence the $W_\alpha$ cover $\partial E$.

The $W_\alpha$ are disjoint, so, as $\partial E$ is compact, 
there are only finitely many of them.
Some $\Bar{W}_\alpha$ may be compact.
Set $D_k=E\cup ({\text{compact\ }}\Bar{W}_\alpha)$.
Then  $D_k$ is compact.

Since the $C_k$ are cofinal in the collection of all compact subsets of $X$,
we may assume, after refinement, that 
$C_0\subseteq D_0\subseteq C_1\subseteq D_1\subseteq\cdots
C_k\subseteq D_k\subseteq C_{k+1}\subseteq\cdots$

Now let $\{ p\}$ be a set of base points for $X$.
Let $\{W_{\alpha,k}\}$ be the set of unbounded path components 
of $X-C_k$, which we saw above was finite.
Since $\{p\}$ is a set of base points, in each $\{W_{\alpha,k}\}$
there are infinitely many $p\in \{p\}$ for which there exists an
${\mathcal O}\in \{{\mathcal O}\}$ such that $p\in{\mathcal O}\subseteq W_{\alpha,k}$.
We get a locally finite subset $\{p^\prime\}\subseteq\{p\}$ by
picking one element of $\{p\}\cap{\mathcal O}$ for each such non{--}empty
intersection as ${\mathcal O}$ runs over $\{{\mathcal O}\}$.
By the above remarks, this set is a set of base points.
It is clearly locally finite so $X$ satisfies 1).

Now let $\{p_k\}$ and $\{q_k\}$ be locally finite irreducible 
sets of base points (they are of necessity both countable).
Look at all the $p_k$'s in $D_0$. 
Join them by paths to some $q_\ell$ not in $D_0$.
Join the $q_k$'s in $D_0$ to some $p_\ell$'s not in $D_0$.
Note that the number of paths intersecting $C_0\leq
(\hbox{number of $p_k$ in $D_0$})+(\hbox{number of $q_k$ in $D_0$})$.

For the inductive step, assume we have joined all the $p_k$'s in $D_{n-1}$
to some $q_k$'s and vice versa.
Suppose moreover that the
\[\hbox{\small number of paths intersecting $C_{n-i}$}\leq
(\hbox{\small number of $p_k$ in $D_{n-i}$})+
(\hbox{\small number of $q_k$ in $D_{n-i}$})\]
for $1\leq i\leq n$.

Look at the $p_k$'s in $D_n-D_{n-1}$ which have not already been joined
to some $q_\ell$ in $D_{n-1}$.
Each of these lies in some $W_{\alpha,n-1}$; i.e. in an unbounded component
of $X-C_{n-1}$.
Join the $p_k$ in $W_{\alpha,n-1} \cap (D_n-D_{n-1})$ which
have not already been fixed up to some $q_\ell$ in $W_{\alpha,n-1}-D_{n}$
by a path in $W_{\alpha,n-1}$; i.e. outside of $C_{n-1}$.
(Recall there are an infinite number of $p_k$ [and $q_k$] in each
$W_{\alpha,\ell}$, so we can always do this.)
Do the same for the $q_k$ in $D_n-D_{n-1}$.

Now each of these new paths misses $C_{n-1}$, so the 
\[\hbox{\small number of paths intersecting $C_{n-i}$}\leq
(\hbox{\small number of $p_k$ in $D_{n-i}$})+
(\hbox{\small number of $q_k$ in $D_{n-i}$})\]
for $1\leq i\leq n$.
For $i=0$, the
\[\hbox{\small number of paths intersecting $C_{n}$}\leq
(\hbox{\small number of $p_k$ in $D_{n}$})+
(\hbox{\small number of $q_k$ in $D_{n}$})\ .\]
This completes the induction and shows $X$ satisfies 2).
\end{proof}

Local compactness and $\sigma${--}compactness are easily seen 
to be proper homotopy invariants so we define

\begin{xDefinition} A space is said to be \emph{homogamous} provided it is
locally compact, $\sigma${--}compact and \prex{--}homogamous.

Note now that any irreducible set of base points for an homogamous
space is countable.
\end{xDefinition}

\bigskip
\section{The \texorpdfstring{$\epsilon${--}$\Delta$}{e--D} construction}
\newHead{I.2}
In this section we describe our construction.
It will enable us to produce a proper homotopy functor on any
homogamous space from an ordinary homotopy functor (a homotopy
functor is a functor from the category of based topological spaces
and based homotopy classes of maps to some category).

Now our homotopy functor, say ${\bf H}$, takes values in some category
${\mathcal A}$.
Associated to any homogamous space, $X$, we have an irreducible set 
of locally finite base points, ${\bf I}$.
We also have a diagram scheme, ${\mathcal D}$, consisting of the closed,
compact subsets of $X$ (see the definition below for the definition
of a diagram scheme).
Our basic procedure is to associate an element in ${\mathcal A}$ to the
collection ${\bf H}(X-C,p)$, where $C$ is a closed compact subset of
$X$, and $p\in{\bf I}$.
In order to be able to do this, we must impose fairly strenuous conditions
on our category ${\mathcal A}$, but we prefer to do this in two stages.

\begin{xDefinition}[see \cite{btwentyfive} page 42]
A \emph{diagram scheme} is a triple ${\mathcal D}=(J,M,d)$, where $J$ is a
set whose elements are called vertices, $M$ is a set whose elements
are called arrows, and $d\colon M\to J\times J$ is a map.
Given a diagram scheme ${\mathcal D}$ and a category ${\mathcal A}$, a diagram
over ${\mathcal D}$ is a map from $J$ to the objects of ${\mathcal A}$
$( j\mapsto A_j )$ and a map from $M$ to the morphisms of ${\mathcal A}$
such that, if $d(m)=(i,j)$ $m$ goes to an element of
${\Homx}(A_i, A_j)$.
\end{xDefinition}

\insetitem{Notation}
$\bigl[{\mathcal D}, {\mathcal A}\bigr]$ denotes the category 
of all diagrams in ${\mathcal A}$ over ${\mathcal D}$.
(A map between diagrams over ${\mathcal D}$ is a collection of morphisms
$f_j\colon A_j\to B_j$ such that $f_j\circ m=\Bar m\circ f_i$, where
$m\in{\Homx}(A_i, A_j)$, and $\Bar m\in{\Homx}(B_i, B_j)$
correspond to the same element in $M$).
If ${\mathcal I}$ is an index set (i.e. a set) ${\mathcal A}^{\mathcal I}$ denotes the
category whose objects are sets of objects in ${\mathcal A}$ indexes by
${\mathcal I}$.
The morphisms are sets of morphisms in ${\mathcal A}$ indexed by ${\mathcal I}$.
Finally, if ${\mathcal A}$ and ${\mathcal B}$ are categories, 
$\bigl\{{\mathcal A}, {\mathcal B}\bigr\}$ is the category of covariant functors
from ${\mathcal A}$ to ${\mathcal B}$ (see \cite{btwentyfive} page 63).

\begin{xDefinition}
A category ${\mathcal A}$ is \emph{weakly regular} with respect to an
index set ${\mathcal I}$ provided:
\begin{enumerate}
\item[1)] ${\mathcal A}$ has products and zero objects.
\item[2)] Let ${\mathcal F}({\mathcal I})=\bigl\{ T\ \vert\ T\subseteq {\mathcal I}
\hbox{ and $T$ is finite}\bigr\}$.
If $\{ G_i \}$ is an object in ${\mathcal A}^{\mathcal I}$, 
each $T\in{\mathcal F}({\mathcal I})$ induces an endomorphism of $\{ G_i \}$ by
\[\begin{cases}G_i\ \to\ G_i& \text{if the identity if } i\notin T\\
G_i\ \to\ G_i& \text{if the zero map if }i\in T\ .\\\end{cases}
\hbox{This induces
a unique map}\]
$\displaystyle
X_T\colon \mathop{\bigtimes}_{i\in{\mathcal I}}G_i \longrightarrow
\mathop{\bigtimes}_{i\in{\mathcal I}}G_i$.
We require that there exist an object 
$\displaystyle\mathop{\mu}_{i\in{\mathcal I}}(G_i)$ and a map
$\displaystyle
\mathop{\bigtimes}_{i\in{\mathcal I}}G_i\longrightarrow
\mathop{\mu}_{i\in{\mathcal I}}(G_i)$ which is the coequalizer of the
family of morphisms $X_T$ for all $T\in{\mathcal F}({\mathcal I})$.
\end{enumerate}

\end{xDefinition}

We easily check
\BEGIN{L.1.2.1}
$\mu\colon {\mathcal A}^{\mathcal I}\to {\mathcal A}$ is a functor when ${\mathcal A}$
is a weakly regular category with respect to ${\mathcal I}$.\qed
\end{Lemma}

\begin{xExamples}
The categories of groups, abelian groups, rings and pointed sets are all
weakly regular with respect to any index set ${\mathcal I}$.
$\mu$ is each case is given as follows. 
We define an equivalence relation $R$ on 
$\displaystyle\mathop{\bigtimes}_{i\in{\mathcal I}}G_i$ by $x R y$ \iff\
(the $i^{\text{\thx}}$ component of $x$)=(the $i^{\text{\thx}}$ component of $y$)
for all but finitely many $i\in{\mathcal I}$.
Then $\displaystyle
\mathop{\mu}_{i\in{\mathcal I}}(G_i)=\left(\mathop{\bigtimes}_{i\in{\mathcal I}}G_i\right)/R$.
\end{xExamples}

\BEGIN{L.1.2.2}
If ${\mathcal D}$ is a diagram scheme, and if ${\mathcal A}$ is a weakly regular
category with respect to ${\mathcal I}$,
then $\bigl[{\mathcal D}, {\mathcal A}\bigr]$ is also weakly regular with respect
to ${\mathcal I}$.
\end{Lemma}
\begin{proof}
$\bigl[{\mathcal D}, {\mathcal A}\bigr]$ is easily seen to have a zero object.
$\bigl[{\mathcal D}, {\mathcal A}\bigr]$ has products, for to each object in
$\bigl[{\mathcal D}, {\mathcal A}\bigr]^{\mathcal I}$, $\bigl(\{G_{i j},\{m_i\}\bigr)$,
we associate the diagram 
$\displaystyle\Bigl(\ \mathop{\bigtimes}_{i\in{\mathcal I}} G_{i j},
\mathop{\bigtimes}_{i\in{\mathcal I}} m_i\Bigr)$.
It is not hard to check that this diagram has 
the requisite universal properties.

To see condition 2), to $\bigl\{G_{i j}\bigr\}$ associate
$\displaystyle\mathop{\mu}_{i\in{\mathcal I}}\bigl(G_{i j}\bigr)$.
Then $\displaystyle\mathop{\bigtimes}_{i\in{\mathcal I}} m_i$
induces $\displaystyle\mathop{\mu}_{i\in{\mathcal I}}(m_i)$, so we do
get a diagram.

To show it is a coequalizer, let $X_j$ be the objects of a diagram.
Set $H_j=\displaystyle\mathop{\bigtimes}_{i\in{\mathcal I}} G_{i j}$.
We are given $g_j\colon H_j\to X_j$ which commute with the diagram
maps.
If $T_1$, $T_2\in{\mathcal F}({\mathcal I})$ we also have 
$g_j\circ X_{T_1}=g_j\circ X_{T_2}$.
Hence by the universality for $\mu$ for ${\mathcal A}$, we get unique maps
$f_j\colon \displaystyle\mathop{\mu}_{i\in{\mathcal I}}(G_{i j})\to X_j$
such that

\hskip1.3in$\xymatrix@C70pt@R50pt{
H_j\ar[r]^-{g_j}\ar[d]&X_j\\
\displaystyle\mathop{\mu}_{i}(G_{i j})\ar[ru]_-{f_j}
}$
\vskip4pt\noindent
commutes.
If we have a map in ${\mathcal D}$ from $j$ to $k$, we get

\hskip1.3in$\xymatrix@C20pt@R30pt{
H_j\ar[rr]\ar[rd]\ar[dd]_-{g_j}&&H_k\ar[rd]\ar[dd]_<<<<<<<<<<<<<<{g_k}\\
&\displaystyle\mathop{\mu}_{i}(G_{i j})\ar@{.>}[rr]\ar[ld]^-{f_j}&&
\displaystyle\mathop{\mu}_{i}(G_{i k})\ar[ld]^-{f_k}\\
X_j\ar[rr]&&X_k\\
}$

\noindent with the front and back squares and both end triangles
commutative.
By the uniqueness of the map
$\displaystyle\mathop{\mu}_{i\in{\mathcal I}}\bigl(G_{i j}\bigr)\to X_k$,
the bottom square also commutes and we are done.
\end{proof}

Suppose given a functor $F\in\bigl\{{\mathcal A}, {\mathcal B}\bigr\}$.
If $F$ does not preserve products, it seems unreasonable to expect
it to behave well with respect to $\mu$, so assume $F$ preserves
products.
Then we get a natural map
$\mu\circ F^{\mathcal I}\to F\circ \mu$
($F^{\mathcal I}$ is the obvious element in 
$\bigl\{{\mathcal A}^{\mathcal I},{\mathcal B}^{\mathcal I}\bigr\}$).
$F$ preserves $\mu$ \iff\ this map is an isomorphism.

Now suppose ${\mathcal A}$ is complete with respect to a diagram scheme
${\mathcal D}$.
Then we have a functor
$\displaystyle\lim_{\mathcal D}\colon \bigl[{\mathcal D}, {\mathcal A}\bigr]\ \to\
{\mathcal A}$, the limit functor (see \cite{btwentyfive}, page 44).

\begin{xDefinition}
If ${\mathcal D}$ is a diagram scheme, and if ${\mathcal A}$ is a 
${\mathcal D}${--}complete, weakly regular category with respect to ${\mathcal I}$,
then we define 
$\epsilon\colon \bigl[{\mathcal D}, {\mathcal A}\bigr]^{\mathcal I}\to{\mathcal A}$
to be the composite $\displaystyle\lim_{\mathcal D}\circ \mu$.
\end{xDefinition}

\BEGIN{P.1.2.1}
Let ${\mathcal D}$ be a 
diagram scheme, and let ${\mathcal A}$ and ${\mathcal B}$ be  
two ${\mathcal D}${--}complete, weakly regular categories 
with respect to ${\mathcal I}$.
Let $F\in\bigl\{{\mathcal A}, {\mathcal B}\bigr\}$.
Then $\epsilon\colon \bigl[{\mathcal D}, {\mathcal A}\bigr]^{\mathcal I}\to{\mathcal A}$
is a functor.
If $F$ preserves products and $\displaystyle\lim_{\mathcal D}$, there is
a natural map
$\epsilon\circ F^{\mathcal I}_{\#}\to F\circ\epsilon$
(where $F_\#\colon\bigl[{\mathcal D}, {\mathcal A}\bigr]\
\to\ \bigl[{\mathcal D}, {\mathcal B}\bigr]$ is the induced functor).
If $F$ preserves $\mu$, this map is an isomorphism.
\end{Proposition}
\begin{proof}
Trivial. \end{proof}

\medskip
Unfortunately, the limit we are taking is an inverse limit, 
which is notorious for causing problems.
In some cases however (and in all cases in which we shall be interested)
it is possible to give a description of $\epsilon$ as a direct limit.
In fact, we will describe the $\Delta${--}construction as a direct limit 
and then investigate the relationship between the two.

\begin{xDefinition}
A lattice scheme ${\mathcal D}$ is a diagram scheme $(J, M, d)$ such that
$J$ is a partially ordered set with least upper and greatest lower
bounds for any finite subset of $J$.
We also require that $d\colon M\to J\times J$ be a monomorphism 
and that ${\Imx}\ d\ =\ \{ (j,k)\in J\times J\ \vert\ j>k\}$.
To the lattice scheme ${\mathcal D}$ and the index set ${\mathcal I}$, we
associate a diagram scheme ${\mathcal D}_{\mathcal I}$ as follows 
( ${\mathcal D}_{\mathcal I}$ is the diagram of 
``cofinal subsequences of ${\mathcal D}$'').
If $\displaystyle\alpha\in\mathop{\bigtimes}_{i\in{\mathcal I}} J$,
define $J_\alpha= \{j\in J\ \vert\ j=
p_i(\alpha)\hbox{ for some $i\in{\mathcal I}$}\}$.
$p_i$ is just the $i^{\text{\thx}}$ projection, so $J_\alpha$ is just the subset of
$J$ we used in making up $\alpha$.
Define $\rho_\alpha\colon {\mathcal I}\ \to\ J$ by
$\rho_\alpha(i)=p_i(\alpha)$.
Set 
\[J_{\mathcal I}=\bigl\{ \displaystyle\alpha\in\mathop{\bigtimes}_{i\in{\mathcal I}} J\ \vert\
\hbox{$J_\alpha$ is cofinal in $J$ and 
$\displaystyle\mathop{\cup}_{j\leq k} \rho^{-1}_\alpha(j)$ is finite for all $k\in J$}\}\ .\]
(A subset of $J$ is cofinal \iff\ given any $j\in J$, there is an element
$k$ of our subset so that $k\geq j$.
$J_{\mathcal I}$ may be thought of as the set of 
``locally finite, cofinal subsets'' of $J$).
\end{xDefinition}

We say $\alpha\geq\beta$ \iff\ $p_i(\alpha)\geq p_i(\beta)$ in $J$ for all
$i\in{\mathcal I}$.
Set $M_{\mathcal I}=\{ (\alpha, \beta)\in J_{\mathcal I}\times J_{\mathcal I}\ \vert\ 
\alpha\geq\beta\}$ and let $d_{\mathcal I}$ be the inclusion.
Given $\alpha$, $\beta\in J_{\mathcal I}$, define $\gamma\in J_{\mathcal I}$
by $p_i(\gamma)=\hbox{least upper bound of $p_i(\alpha)$
and $p_i(\beta)$.}$
(It is not hard to see $\gamma\in J_{\mathcal I}$.)
Greatest lower bounds can be constructed similarly.
Hence ${\mathcal D}_{\mathcal I}$ is also a lattice scheme.

Now if $J$ does not have any cofinal subsets of cardinality 
$\leq {\text{card}}({\mathcal I})$, $J_{\mathcal I}=\emptyset$.
Since $J$ has upper bounds for finite sets, if $J$ has finite cofinal subsets,
then $J$ has cofinal subsets of cardinality $\geq{\text{card}}({\mathcal I})-N$,
where $N$ is some natural number, then the condition that
$\displaystyle\mathop{\cup}_{j\leq k} \rho^{-1}_\alpha(j)$ be finite forces
$J_{\mathcal I}=\emptyset$.
Empty diagrams are a nuisance, so we define an ${\mathcal I}${--}lattice
scheme as a lattice scheme with cofinal subsets of cardinality 
$= {\text{card}}({\mathcal I})$.

We can now define 
$\delta\colon \bigl[{\mathcal D}, {\mathcal A}\bigr]^{\mathcal I}\ \to\ 
\bigl[{\mathcal D}_{\mathcal I}, {\mathcal A}\bigr]$ as follows.
If $\bigl\{ d_i\}\in\bigl[{\mathcal D}, {\mathcal A}\bigr]^{\mathcal I}$, 
$\delta(d)$ has for objects
$\displaystyle\delta_\alpha=
\mathop{\bigtimes}_{i\in{\mathcal I}}G_{i p_i(\alpha)}$,
where $G_{i j}$ is the $j^{\text{\thx}}$ object in the diagram for $d_i$
($\alpha\in J_{\mathcal I}$, $j\in J$, $i\in{\mathcal I}$).
If $\alpha\geq\beta$, we define $\delta_\alpha\ \to\ \delta_\beta$
by the maps $G_{i p_i(\alpha)}\to G_{i p_i(\beta)}$ 
which come from the diagram $d_i$.

We can also define maps $\displaystyle\delta_\alpha\to \mathop{\mu}_{j}$
as follows.
Map $G_{i p_i(\alpha)} \to G_{i j}$ by the unique map in 
$d_i$ if $p_i(\alpha)\geq j$, and by the zero map if $j>p_i(\alpha)$.
(Notice that there are at most finitely many $i$ 
such that $p_i(\alpha)<j$ by
the second defining condition on $J_{\mathcal I}$.)
These maps induce a unique map
$\delta_\alpha \to \mathop{\bigtimes}_{i\in{\mathcal I}}G_{i j}$.
Composing with the projection, we get a unique map
$\displaystyle\delta_\alpha\to \mathop{\mu}_{i\in{\mathcal I}}(G_{i j})=
\mathop{\mu}_{j}$.

If $k\geq j$,\hskip20pt \lower 12pt\vtop{\hsize=1in\noindent
$\displaystyle\begin{matrix}%
\delta_\alpha&\longrightarrow&\mathop{\mu}_{j}\\
&\searrow&\big\downarrow\\
&&\mathop{\mu}_{k}\\
\end{matrix}$}\hskip20pt commutes as one easily checks.
If $\alpha\geq\beta$,

\hskip20pt \lower 12pt\vtop{\hsize=1in\noindent$\displaystyle\begin{matrix}%
\delta_\alpha\\
&&\searrow\\
\big\downarrow&&&\mathop{\mu}_{j}\\
&&\nearrow\\
\delta_\beta\\
\end{matrix}$}\hskip20pt also commutes.

\bigskip
\BEGIN{L.1.2.3}
$\delta\colon\bigl[{\mathcal D}, {\mathcal A}\bigr]^{\mathcal I}\ \to\
\bigl[{\mathcal D}_{\mathcal I}, {\mathcal A}\bigr]$ is a functor.
\end{Lemma}
\smallskip\begin{proof}
The proof is easy and can be safely left to the reader. \end{proof}

Now suppose ${\mathcal A}$ is ${\mathcal D}_{\mathcal I}${--}cocomplete.
Then we have a colimit functor $\displaystyle\colim_{{\mathcal D}_{\mathcal I}}$.

\begin{xDefinition}
If ${\mathcal D}$ is an ${\mathcal I}${--}lattice scheme, and if ${\mathcal A}$ is a
${\mathcal D}_{\mathcal I}${--}cocomplete, weakly regular category with respect
to ${\mathcal I}$, then we define
$\Delta\colon \bigl[{\mathcal D}, {\mathcal A}\bigr]^{\mathcal I}\to{\mathcal A}$
to be the composition 
$\displaystyle\colim_{{\mathcal D}_{\mathcal I}}\ \circ\ \delta$.
\end{xDefinition}

\bigskip
\BEGIN{P.1.2.2}
Let ${\mathcal D}$ be a diagram scheme and let ${\mathcal A}$ and ${\mathcal B}$ be
two ${\mathcal D}_{\mathcal I}${--}cocomplete, weakly regular categories with
respect to ${\mathcal I}$.
Let $F\in\bigl\{{\mathcal A}, {\mathcal B}\bigr\}$.
Then  $\Delta\colon \bigl[{\mathcal D}, {\mathcal A}\bigr]^{\mathcal I}\to{\mathcal A}$
is a functor.
There is always a natural map 
$\Delta\circ F_{\#}^{\mathcal I}\to F\circ\Delta$.
If $F$ preserves products and $\displaystyle\colim_{{\mathcal D}_{\mathcal I}}$,
this map is an isomorphism.
\end{Proposition}
\smallskip\begin{proof}
Trivial. \end{proof}

The maps we constructed from 
$\displaystyle\delta_\alpha\to \mathop{\mu}_{j}$ combine to give us a
natural transformation from  $\Delta$ to $\epsilon$ whenever both
are defined.
We would like to study this natural transformation in order to get information
about both $\Delta$ and $\epsilon$.
A $({\mathcal D}, {\mathcal I})${--}regular category is about the most general category
in which we can do this successfully, and it includes all the examples we
have in mind.

\begin{xDefinition}
A category ${\mathcal A}$ is said to be $({\mathcal D}, {\mathcal I})${--}regular
provided 
\begin{enumerate}
\item[1)] ${\mathcal A}$ is weakly regular with respect to ${\mathcal I}$
\item[2)] ${\mathcal A}$ has images and inverse images
\item[3)] There is a covariant functor $F$ from ${\mathcal A}$ to 
the category of pointed sets and maps such that
\begin{enumerate}
\item[a)] $F$ preserves kernels, images, products, 
limits over ${\mathcal D}$, increasing unions, and $\mu$
\item[b)] $F$ reflects kernels, images, and isomorphisms
\end{enumerate}
\item[4)] ${\mathcal A}$ is ${\mathcal D}${--}complete and 
${\mathcal D}_{\mathcal I}${--}cocomplete
\item[5)] ${\mathcal I}$ is countable
\end{enumerate}
\end{xDefinition}
\bigskip
\begin{xExamples}
The categories of groups, abelian groups, rings and pointed sets are
all $({\mathcal D}, {\mathcal I})${--}regular for any ${\mathcal I}${--}lattice scheme
\footnote[1]{For which ${\mathcal I}$ is countable.}.
The functor $F$ is just the forgetful functor.
\end{xExamples}
\bigskip
\BEGIN{L.1.2.4}
Let ${\mathcal A}$ be a $({\mathcal D}, {\mathcal I})${--}regular category.
Then $\displaystyle\mathop{\bigtimes}_{i\in{\mathcal I}}$ and
$\mu$ preserve kernels and images.
\end{Lemma}

\begin{proof}
$\displaystyle\mathop{\bigtimes}_{i\in{\mathcal I}}$ is known to preserve
kernels (Mitchell \cite{btwentyfive}, page 67, Corollary 12.3).

Since $F$ preserves images, if ${\Imx}(f)$ is the image of
$A\ \RA{f}\ B$, then $F(A)$ is onto $F\bigl({\Imx}(f)\bigr)$
and $F\bigl({\Imx}(f)\bigr)$ injects into $F(B)$.
Let $K_i$ be the image of $A_i\ \RA{\ f_i\ }\ B_i$.
Then, since $F$ preserves products, 
$\displaystyle\mathop{\bigtimes}_{i\in{\mathcal I}} K_i\ \to\ 
\displaystyle\mathop{\bigtimes}_{i\in{\mathcal I}}B_i$ is a monomorphism
since $\displaystyle\mathop{\bigtimes}_{i\in{\mathcal I}}$ is a monofunctor.
Since $F$ preserves monomorphisms,
$\displaystyle\mathop{\bigtimes}_{i\in{\mathcal I}} F(K_i)\ \to\ 
\mathop{\bigtimes}_{i\in{\mathcal I}} F(B_i)$ is seen to be a monomorphism
as $F$ also preserves products.
Since $F$ reflects images, 
$\displaystyle\mathop{\bigtimes}_{i\in{\mathcal I}} K_i$ is the image of
$\displaystyle\mathop{\bigtimes}_{i\in{\mathcal I}}A_i\ \to\
\mathop{\bigtimes}_{i\in{\mathcal I}} B_i$.

Let $K_i\ \to\ A_i\ \to\ B_i$ be kernels.
Then
\[\begin{matrix}%
\bigtimes K_i&\RA{\ h\ }& \bigtimes A_i&\RA{\ g\ }&\bigtimes B_i\\
\big\downarrow&&\big\downarrow&&\big\downarrow\\
\mu(K_i)&\RA{\ \ \ }&\mu(A_i)&\RA{\ \ \ }&\mu(B_i)\\
\end{matrix}\leqno(\diagramOne)\]
commutes.
Since $F$ reflects kernels, we need only show that
$F\bigl(\mu(K_i)\bigr)$ injects into $F\bigl(\mu(A_i)\bigr)$ and is onto
$F(f)^{-1}(0)$.
Since $F$ preserves $\mu$, we may equally consider
$\mu\bigl(F(K_i)\bigr)$, etc.
Since $F$ preserves products, we may as well assume the diagram
$(\diagramOne)$ is in the category of pointed sets.

We show $\mu(K_i)$ is onto $f^{-1}(0)$.
Let $x\in f^{-1}(0)\subseteq \mu(A_i)$.
Lift $x$ to $y\in\bigtimes A_i$, which we may do since
$\displaystyle\mathop{\bigtimes}_{i\in{\mathcal I}}$ is onto
$\displaystyle\mathop{\mu}_{i\in{\mathcal I}}$ in the category of
pointed sets.
Now $g(y)\in\bigtimes B_i$ can have only finitely many 
non{--}zero components since it goes to $0$ in $\mu(B_i)$.

Define $\Bar y$ by
\[p_i(\Bar y)=\begin{cases}p_i(y)& \text{if } p_i\bigl(g(y)\bigr) = 0\\
0& \text{if }p_i\bigl(g(y)\bigr) \neq 0\\\end{cases}\ .\]
Then $\Bar y$ also lifts $x$ and $g(\Bar y)=0$.
There is a $z\in\bigtimes K_i$ such that $h(z)=\Bar y$, so
$\mu(K_i)$ maps onto $f^{-1}(0)$.
A similar argument shows $\mu(K_i)$ injects into $\mu(A_i)$.
Hence $\mu$ preserves kernels.

Now let $K_i$ be the image of $A_i\ \to\ B_i$.
Then
\[\begin{matrix}%
\bigtimes A_i&\RA{\ \ \ }& \bigtimes K_i&\RA{\ \ \ }&\bigtimes B_i\\
\big\downarrow&&\big\downarrow&&\big\downarrow\\
\mu(A_i)&\RA{\ \ \ }&\mu(K_i)&\RA{\ \ \ }&\mu(B_i)\\
\end{matrix}\]
commutes.
By general nonsense, it suffices to prove the result assuming  
we are working in the category of pointed sets.

Since $\displaystyle\mathop{\bigtimes}_{i\in{\mathcal I}}$ preserves
images, $\displaystyle\mathop{\bigtimes}_{i\in{\mathcal I}} A_i\ \to\ 
\mathop{\bigtimes}_{i\in{\mathcal I}} K_i$ is onto, so
$\mu(A_i)\ \to\ \mu(K_i)$ is easily seen to be onto.
Since $\mu$ preserves kernels, $\mu(K_i)$ injects into $\mu(B_i)$,
so $\mu(K_i)$ is the image of $\mu(A_i)\ \to\ \mu(B_i)$. \end{proof}

\bigskip
\BEGIN{T.1.2.1}
Let ${\mathcal A}$ be a $({\mathcal D}, {\mathcal I})${--}regular category.
Then $\epsilon$ preserves kernels and images.
\end{Theorem}
\begin{proof}
By Mitchell \cite{btwentyfive} (page 67, Corollary 12.2)
$\displaystyle\lim_{\mathcal D}$ preserves kernels so
$\epsilon = \displaystyle\lim_{\mathcal D}\ \circ\ \mu$ also does
using \fullRef{L.1.2.4} and general nonsense.

Now let $K_{i j}$ be the image of $A_{i j}\ \to\ B_{i j}$.
We claim that, if $x\in \epsilon(K_{i j})$, 
then there exists $\alpha\in J_{\mathcal I}$ such that
$x$ is in the image of $\delta_\alpha(K_{i j})$.
Assuming this for now we proceed as follows.

Since $\delta_\alpha(K_{i j})=
\displaystyle\mathop{\bigtimes}_{i\in{\mathcal I}} K_{i p_i(\alpha)}$, 
$\delta_\alpha$ preserves kernels and images by 
\fullRef{L.1.2.4}.
Hence

\[\begin{matrix}%
\delta_\alpha(A_{i j})&\RA{\hskip 10pt}&\delta_\alpha(K_{i j})
&\RA{\hskip 10pt}&\delta_\alpha(B_{i j})\\
\big\downarrow&&\big\downarrow&&\big\downarrow\\
\epsilon(A_{i j})&\RA{\hskip 10pt}&\epsilon(K_{i j})
&\RA{\hskip 10pt}&\epsilon(B_{i j})\\
\end{matrix}\]
commutes and $\delta_\alpha(K_{i j})$ is the image of
$\delta_\alpha(A_{i j})\ \to\ \delta_\alpha(B_{i j})$.
By the usual abstract nonsense, we may as well assume we are in the
category of pointed sets (note $F$ preserves $\epsilon$ by
\fullRef{P.1.2.1}.

Now using our claim we can easily get
$\epsilon(A_{i j})\ \to\ \epsilon(K_{i j})$ is onto.
$\epsilon(K_{i j})\ \to\ \epsilon(B_{i j})$ injects since $\epsilon$
preserves kernels.
Hence $\epsilon(K_{i j})$ is the image of 
$\epsilon(A_{i j})\ \to\ \epsilon(B_{i j})$. \end{proof}

We prove a stronger version of our claim than we have yet used.

\BEGIN{L.1.2.5}
Let ${\mathcal A}$ be the category of pointed sets.
Let $\{G_{i j}\}$ be an object in $\bigl[{\mathcal D}, {\mathcal A}\bigr]^{\mathcal I}$.
Then if $x\in\epsilon(G_{i j})$, there is an $\alpha\in J_{\mathcal I}$
such that $\delta_\alpha(G_{i j})$ contains $x$ in its image.
If $y$, $z\in\delta_\alpha(G_{i j})$ both hit $x$, then there is a 
$\beta\leq\alpha$ such that, in $\delta_\beta(G_{i j})$, the images of
$y$ and $z$ differ in only finitely many coordinates.
In fact, if there is a $j\in J$ such that $j\leq p_i(\alpha)$ for all
$i\in{\mathcal I}$ and such that $y$ and $z$ agree in 
$\displaystyle\mathop{\bigtimes}_{i\in{\mathcal I}}G_{i j}$, then $\beta$
can be chosen so that $y=z$ in $\delta_\beta(G_{i j})$.
\end{Lemma}

\begin{proof}
If $x\in \epsilon(G_{i j})$, there exist unique 
$a_j\in\displaystyle\mathop{\mu}_{i\in{\mathcal I}}(G_{i j})$ such that
$x$ hits $a_j$.
Since $\displaystyle\mathop{\bigtimes}_{i\in{\mathcal I}}$ is onto
$\displaystyle\mathop{\mu}_{i\in{\mathcal I}}$, we may lift $a_j$ to
$b_j\in\displaystyle\mathop{\bigtimes}_{i\in{\mathcal I}} G_{i j}$.
Since $J$ has countable cofinal subsets, let the natural numbers
$j=1$, $2$, \dots\ be one such.
Since ${\mathcal I}$ is countable (and infinite or our result is easy) we also
assume it to be the natural numbers.

Now look at $b_2$ and $b_1$. Since they agree in $\mu(G_{i1})$, 
$b_2$ projected into $\bigtimes G_{i j}$ differs from $b_1$ in only
finitely many coordinates.
Let ${\mathcal I}_1\subseteq{\mathcal I}$ be the finite subset which indexes these
unequal coordinates, together with the element $1\in{\mathcal I}$.

Next look at the pairs $(b_3,b_2)$ and $(b_3,b_1)$.
As before, projected into $\bigtimes G_{i2}$, $b_3$ and $b_2$ agree
in all but finitely many coordinates.
In $\bigtimes G_{i1}$, $b_3$ and $b_1$ differ in only 
finitely many coordinates.
Set ${\mathcal I}_2\subseteq{\mathcal I}$ to be the finite subset of ${\mathcal I}$
which indexes the unequal coordinates of $(b_3, b_2)$ or $(b_3,b_1)$
which lie in ${\mathcal I}-{\mathcal I}_1$, together with the smallest integer in
${\mathcal I}-{\mathcal I}_1$.

Define ${\mathcal I}_k$ to be the finite subset of ${\mathcal I}$ which indexes
the unequal coordinates of
$(b_k, b_{k-1})$, \dots, $(b_k, b_2)$, $(b_k,b_1)$ which lie in
${\mathcal I}-\bigl({\mathcal I}_{k-1}\ \cup\ \cdots\ \cup\ 
{\mathcal I}_2\ \cup\ {\mathcal I}_1\bigr)$,
together with the smallest integer in 
${\mathcal I}-\bigl({\mathcal I}_{k-1}\ \cup\ \cdots\ \cup\ 
{\mathcal I}_2\ \cup\ {\mathcal I}_1\bigr)$.

Then $\displaystyle{\mathcal I}=\mathop{\cup}_{k=1}^\infty {\mathcal I}_k$
as a disjoint union. 
Define $\alpha$ by $p_i(\alpha)=k$, where $i\in{\mathcal I}_k$.
Since ${\mathcal I}$ is countable, but not finite, and since each ${\mathcal I}_k$
is finite, $\alpha\in J_{\mathcal I}$.
Define $y\in\delta(G_{i j})$ by $p_i(y)=p_i(b_{p_i(\alpha)})$.
A chase through the definitions shows $y$ hits each $a_j$ through the
map $\delta_\alpha(G_{i j})\ \to\ 
\displaystyle\mathop{\mu}_{i\in{\mathcal I}}(G_{i j})$.
Thus $y$ hits $x$ in $\epsilon(G_{i j})$.

Now suppose $y$, $z\in \delta_\alpha(G_{i j})$ both map to $x$.
Then they map to the same element in each
$\displaystyle\mathop{\mu}_{i\in{\mathcal I}}(G_{i j})$.
Let $a_j$ be the image of $y$ in 
$\displaystyle\mathop{\bigtimes}_{i\in{\mathcal I}}G_{i j}$ under the map
$\delta_\alpha(G_{i j})\ \to\ 
\displaystyle\mathop{\bigtimes}_{i\in{\mathcal I}}G_{i j}$
which we defined just before \fullRef{L.1.2.3}.
Set $b_j$ to be the image of $z$ in 
$\displaystyle\mathop{\bigtimes}_{i\in{\mathcal I}} G_{i j}$.
Then $a_j$ and $b_j$ differ in only finitely many coordinates.

Let ${\mathcal I}_1$ be the finite subset of ${\mathcal I}$ which indexes the
unequal coordinates of $a_1$ and $b_1$.
If there is a $j\leq p_i(\alpha)$ for all $i\in{\mathcal I}$ such that $y$
and $z$ agree in $\displaystyle\mathop{\bigtimes}_{i\in{\mathcal I}} G_{i j}$,
we may assume $j=1$, so $a_1=b_1$, and ${\mathcal I}_1=\emptyset$.

Define ${\mathcal I}_k$ as the finite subset of ${\mathcal I}$ which indexes the
unequal coordinates of $(a_k,b_k)$ which lie in
${\mathcal I}-\bigl({\mathcal I}_{k-1}\ \cup\ \cdots\ \cup\ 
{\mathcal I}_2\ \cup\ {\mathcal I}_1\bigr)$.
Define $\beta$ by
\[p_i(\beta)=\begin{cases}k-1& \text{if }i\in{\mathcal I}_k \text{ for some } k\geq2\\
p_i(\alpha)& \text{if }i\notin{\mathcal I}_k \text{ for any }k\geq2\\\end{cases}\ .\]
Note $p_i(\beta)\leq p_i(\alpha)$, since $i\in{\mathcal I}_k$, this says
$p_i(a_k)\neq p_i(b_k)$.
But if $k>p_i(\alpha)$, $p_i(a_k)=0=p_i(b_k)$
by the definition of our map from $\delta_\alpha$ to 
$\displaystyle\mathop{\bigtimes}_{i\in{\mathcal I}}$.
Hence $k\leq p_i(\alpha)$, so $p_i(\beta)\leq p_i(\alpha)$.

Let $\Bar y$ be the projection of $y$ into $\delta_\beta(G_{i j})$,
and let $\Bar z$ be the projection of $z$ into $\delta_\beta(G_{i j})$.
$p_i(\Bar y)= p_i(a_{p_i(\beta)})$ and $p_i(\Bar z)= p_i(b_{p_i(\beta)})$.
If $p_i(a_{p_i(\beta)})\neq  p_i(b_{p_i(\beta)})$, then $i\notin{\mathcal I}_k$
for any $k\geq2$, since $i\in{\mathcal I}_k$ for $k\geq2$ says that
$p_i(a_k)\neq p_i(b_k)$ but $p_i(a_{k-1})=p_i(b_{k-1})$.
If $i\notin{\mathcal I}_k$ for any $k$, it says that $p_i(y)=p_i(z)$.
Thus $p_i(\Bar y)=p_i(\Bar z)$ if $i\notin {\mathcal I}_1$.
Hence they agree for all but finitely many $i\in{\mathcal I}$.
In fact, if ${\mathcal I}_1=\emptyset$, $\Bar y=\Bar z$. 
\end{proof}

We can now describe $\epsilon(G_{i j})$ as a colimit (direct limit).
Let $\displaystyle\mathop{\mu}_{\alpha}(G_{i j})$ be the $\mu$
functor applied to $\bigl\{ G_{i\> p_i(\alpha)}\bigr\}$.
Then the map $\delta_\alpha(G_{i j})\ \to\ \epsilon(G_{i j})$ factors
through $\displaystyle\mathop{\mu}_{\alpha}(G_{i j})$.

\bigskip
\BEGIN{T.1.2.2}
Let ${\mathcal A}$ be a $({\mathcal D}, {\mathcal I})${--}regular category.
Then the natural map \\
$\displaystyle\colim_{{\mathcal D}_{\mathcal I}}\ \mu_{\alpha}\ \to\ \epsilon$
is an isomorphism.
Hence $\epsilon$ is a cokernel, kernel preserving functor.
\end{Theorem}
\begin{proof}
Let us first show $F$ preserves $\displaystyle\colim_{{\mathcal D}_{\mathcal I}}$;
i.e. we must show that the natural map 
\[\displaystyle
\colim\ F(A_\alpha)\ \RA{f}\ F(\colim_{{\mathcal D}_{\mathcal I}} A_\alpha)\]
is an isomorphism.
To do this, we first compute ${\Imx}(f)$.
If ${\Imx}(f_\alpha)$ is the image of $F(A_\alpha)\ \to\ 
\colim\> F(A_\alpha)\ \to\ F(\colim\> A_\alpha)$, then by Mitchell
\cite{btwentyfive} (Proposition 2.8, page 46),
$\displaystyle{\Imx}(f)=\mathop{\cup}_\alpha\ {\Imx}(f_\alpha)$.
Let ${\Imx}(g_\alpha)$ be the image of 
$A_\alpha\ \to\ \colim A_\alpha$.
Then, since $F$ preserves images, $F\bigl({\Imx}(g_\alpha)\bigr)=
{\Imx}(f_\alpha)$, so
$\displaystyle\mathop{\cup}_\alpha {\Imx}(f_\alpha) = 
\mathop{\cup}_\alpha F\bigl({\Imx}(g_\alpha)\bigr)$.

Now $\{\alpha\}$ has a cofinal subsequence (which is countable and, if
${\mathcal I}$ is finite, it is also finite) $\{\alpha_i\}$ such that
$\alpha_o < \alpha_1< \cdots < \alpha_n < \cdots$\ .
Therefore
$\displaystyle\mathop{\cup}_\alpha {\Imx}(f_\alpha) =
\mathop{\cup}_{i=0}^\infty{\Imx}(f_{\alpha_i}) =
\mathop{\cup}_{i=0}^\infty F\bigl({\Imx}(g_{\alpha_i})\bigr)
$
since $\{\alpha_i\}$ is cofinal.

Again by Mitchell \cite{btwentyfive} (Proposition 2.8, page 46),
$\colim A_\alpha = 
\mathop{\cup}_\alpha {\Imx}(g_\alpha)=
\mathop{\cup}_{i=0}^\infty{\Imx}(g_{\alpha_i})
$.
Thus $F(\colim\ A_\alpha)=
F\bigl(\mathop{\cup}_{i=0}^\infty{\Imx}(g_{\alpha_i})\bigr)$.
Since ${\mathcal A}$ is a $({\mathcal D}, {\mathcal I})${--}regular category, 
the natural map
\[\displaystyle\mathop{\cup}_{i=0}^\infty F\bigl({\Imx}(g_{\alpha_i})
\bigr)\subseteq F\bigl(
\mathop{\cup}_{i=0}^\infty {\Imx}(g_{\alpha_i})\bigr)
\] is an isomorphism.
Thus the map 
$\displaystyle\mathop{\cup}_\alpha{\Imx}(f_\alpha)\subseteq
F\bigl(\mathop{\cup}_\alpha{\Imx}(g_\alpha)\bigr)$
is an isomorphism.
But this map is just the natural map
$\displaystyle\colim_{{\mathcal D}_{\mathcal I}}\ F(A_\alpha)\ \to\ 
F(\colim_{{\mathcal D}_{\mathcal I}}\ A_\alpha)$.

The natural map $\colim \mu_\alpha\ \to\ \epsilon$
is the map which comes from the maps $\mu_\alpha\ \to\ \epsilon$.
To show it is an isomorphism, it is enough to show it is for pointed sets
by the result above and the fact that $F$ reflects isomorphisms.
But this is exactly what \fullRef{L.1.2.5} says.

Now $\epsilon$ preserves kernels by \fullRef{T.1.2.1},
and it preserves cokernels since colimits preserve cokernels by
Mitchell \cite{btwentyfive} (page 67, Corollary 12.2 dualized). 
\end{proof}

\medskip
\BEGIN{C.1.2.2.1}
Let $\bigl\{ G_{i j}^n\bigr\}\in\bigl[{\mathcal D}, {\mathcal A}\bigr]^{\mathcal I}$
be a collection of exact sequences in a $({\mathcal D}, {\mathcal I})${--}regular
category ${\mathcal A}$ (i.e. there are maps
$f^n_{i j}\ \colon\ G^n_{i j}\ \to\ G^{n-1}_{i j}$ 
which are maps of diagrams
such that ${\Imx}(f^n_{i j})={\kerx}(f^{n-1}_{i j})$\ ).
Then the sequence
\[\cdots \to\ \epsilon(G^n_{i j})\ \RA{\ \epsilon(f^n_{i j})\ }\ 
\epsilon(G^{n-1}_{i j})\ \to\ \cdots\]
is also exact.
\end{Corollary}
\medskip
\BEGIN{C.1.2.2.2}
Let ${\mathcal A}$ be a $({\mathcal D}, {\mathcal I})${--}regular abelian category.
Let $\bigl\{ G^\ast_{i j}, f^\ast_{i j}\bigr\}$ 
be a collection of chain complexes in $\bigl[{\mathcal D}, {\mathcal I}\bigr]^{\mathcal I}$.
Then $\bigl\{\epsilon(G^\ast_{i j}), \epsilon(f^\ast_{i j})\bigr\}$ is
a chain complex, and
$H_\ast\bigl(\epsilon(G^\ast_{i j})\bigr) = 
\epsilon\bigl(H_\ast(G^\ast_{i j})\bigr)$, where $H_\ast$ is the homology
functor (see Mitchell \cite{btwentyfive}, page 152).
\end{Corollary}

\renewcommand{\proofname}{Proofs}
\begin{proof}
The first corollary is easily seen to be true.
(It is, in fact, a corollary of \fullRef{T.1.2.1}.)

The second corollary is almost as easy.
If $\bigl\{ Z^n_{i j}\bigr\}$ are the $n${--}cycles, and if 
$\bigl\{ B^{n+1}_{i j}\bigr\}$ are the $(n+1)${--}boundaries,
$0\to\ B^{n+1}_{i j}\ \to\ Z^n_{i j}\ \to\ H_n(G^\ast_{i j})\to0$
is exact.
Applying $\epsilon$, we get 
$0\to\ \epsilon(B^{n+1}_{i j})\ \to\ \epsilon(Z^n_{i j})\ \to\ 
\epsilon\bigl(H_n(G^\ast_{i j})\bigr)\to0$
is exact.
But as $\epsilon$ preserves kernels and images, 
$\epsilon\bigl(Z^n_{i j}\bigr)$ is the collection of $n${--}cycles for
$\epsilon(G^\ast_{i j})$ and
$\epsilon\bigl(B^{n+1}_{i j}\bigr)$ is 
the collection of $(n+1)${--}boundaries.
Hence 
$H_\ast\bigl(\epsilon(G^\ast_{i j})\bigr)\ \to\ 
\epsilon\bigl(H_n(G^\ast_{i j})\bigr)$ is an isomorphism.
\end{proof}
\renewcommand{\proofname}{Proof}

Now suppose $J$ has a unique minimal element $j_0$.
Then we get a square
\[\begin{matrix}%
\Delta(G_{i j})&\RA{\hskip1in}&\displaystyle \mathop{\bigtimes}_{i\in{\mathcal I}} G_{i j_0}\\
\Bigg\downarrow&&\Bigg\downarrow\\
\epsilon(G_{i j})&\RA{\hskip1in}&\displaystyle\mathop{\mu}_{i\in{\mathcal I}}(G_{i j_0})\\
\end{matrix}\]

\BEGIN{T.1.2.3}
In a $({\mathcal D}, {\mathcal I})${--}regular category, the above diagram is a 
pullback in the category of pointed sets, so if $F$ reflects pullbacks 
the above square is a pullback.
\end{Theorem}
\begin{proof}
As we showed in the proof of \fullRef{T.1.2.2}
that $F$ and $\displaystyle\colim_{{\mathcal D}_{\mathcal I}}$ commute, we
have $F\bigl(\Delta(G_{i j})\bigr) = \Delta\bigl(F(G_{i j})\bigr)$,
so we may work in the category of pointed sets.

The omnipresent \fullRef{L.1.2.5} can be used 
to show the above square is a pullback.
The pullback is the subset of 
$\epsilon(G_{i j})\times\bigl(\displaystyle\mathop{\bigtimes}_{i\in{\mathcal I}}
G_{i j_0}\bigr)$
consisting of pairs which project to the same element in
$\displaystyle\mathop{\mu}_{i\in{\mathcal I}}(G_{i j_0})$.

Given any element, $x$, in $\epsilon(G_{i j})$ 
we can find $\alpha\in J_{\mathcal I}$ such that 
the element is in the image of $\delta_\alpha(G_{i j})$.
Lift the image of $x$ in $\displaystyle\mathop{\mu}_{i\in{\mathcal I}}(G_{i j_0})$
to $y\in\displaystyle\mathop{\bigtimes}_{i\in{\mathcal I}}G_{i j_0}$.
Let $z\in\delta_\alpha(G_{i j_0})$ be an element which hits $x$. 
Then $y$ pushed into $\bigtimes G_{i j_0}$ and $z$ agree, except in
finitely many places.
It is then easy to find $\beta\in J_{\mathcal I}$ with $\beta\leq\alpha$
and an element $q\in\delta_\beta(G_{i j})$such that $q$ hits $x$
and $y$.
This says precisely that our square is a pullback.
\end{proof}

\begin{xRemarks}
In all our examples, $F$ reflects pullbacks. The analogues of
Corollaries \shortFullRef{C.1.2.2.1}
and \shortFullRef{C.1.2.2.2}
may be stated and proved by the reader for the $\Delta$ functor.
\end{xRemarks}

\BEGIN{T.1.2.4}
In a $({\mathcal D}, {\mathcal I})${--}regular category, $\epsilon(G_{i j})=0$
\iff\ given any $j\in J$ there exists a $k\geq j$ such that 
$G_{i k}\ \to\ G_{i j}$ is the zero map for all but finitely many $i$.
\end{Theorem}
\begin{proof}
Suppose given $j$ we can find such a $k$.
The we can produce a cofinal set of $j$'s, $j_0\leq j_1\leq \cdots$\ ,
such that the map $\mu(G_{i j_k})\ \to\ \mu(G_{i j_{k-1}})$ is
the zero map.
Hence $\epsilon=0$.

Conversely, suppose for some $j_0$ that no such $k$ exists.
This means that for every $k\geq j_0$ there are infinitely many
$i$ for which $G_{i k}\ \to\ G_{i j_0}$ is not the zero map.

As usual, it suffices to prove the result for pointed sets, 
so assume we have $Z_{i k}\in G_{i k}$ which goes non{--}zero
into $G_{i j_0}$.
Pick $j_0\leq j_1\leq j_2\leq \cdots$ a countable cofinal subsequence 
of $J$.
We define an element $\alpha$ of ${\mathcal D}_{\mathcal I}$ as follows.
Well order ${\mathcal I}$.
Then $\alpha(i)=j_0$ until we hit the first element of ${\mathcal I}$ for
which a $Z_{i j_k}$ is defined.
Set $\alpha(i)=j_k$ for this $i$ and continue defining $\alpha(i)=j_k$
until we hit the next element of ${\mathcal I}$ for which a $Z_{i j_{k_1}}$
is defined with $k_1\geq k$.
Set $\alpha(i)=j_{k_1}$ until we hit the next $Z_{i j_{k_2}}$ 
with $k_2\geq k_1$.
Continuing in this fashion is seen to give an element of ${\mathcal D}_{\mathcal I}$.
Define $Z_\alpha$ by $Z_{i\>\alpha(i)}=0$ unless $i$ is one of the
distinguished elements of ${\mathcal I}$, in which case set 
$Z_{i\>\alpha(i)}=Z_{i j_k}$, where $j_k=\alpha(i)$.

Then $Z_\alpha\in\Delta_\alpha(G_{i j})$. 
It is non{--}zero in $\displaystyle\mathop{\mu}_{i\in{\mathcal I}}(G_{i j_0})$
by construction, so $\epsilon(G_{i j})\neq0$.
\end{proof}

\bigskip
\section{Proper homotopy functors and their relations}
\newHead{I.3}
\medskip
We begin by clarifying the concept of an ordinary homotopy functor.
A homotopy functor is a functor, $h$, from 
the category of pointed topological spaces to some category, ${\mathcal A}$.
Given a space $X$ and two base points $p_1$ and $p_2$, and a path
$\lambda$ from $p_1$ to $p_2$, there is a natural transformation
$\alpha_\lambda\colon h(X,p_1)\ \to\ h(X,p_2)$ which is an
isomorphism and depends only on the homotopy class of $\lambda$
rel end points.
Furthermore, $h(X,p)\ \to\ h(X\times I, p\times t)$ given by 
$x\mapsto (x,t)$ is an isomorphism for $t=0$ and $1$.

For any homotopy functor we are going to associate a proper homotopy 
functor defined on the category of homogamous spaces and countable
sets of locally finite irreducible base points.

To be able to do this in the generality we need, we shall have to digress
momentarily to discuss the concept of a covering functor.

Let $X$ be a homogamous space, and let $\{ x_i\}$ be a countable,
locally finite, irreducible set of base points for $X$. (From now on we
write just ``set of base points'' for ``countable, locally finite, irreducible
set of base points.'')
Let ${\mathcal D}_X$ be some naturally defined collection of subsets of $X$
(by naturally defined we mean that if $f\colon X \to Y$ is a proper map,
$f^{-1}{\mathcal D}_Y\subseteq{\mathcal D}_X$).
${\mathcal D}_X$ is a diagram with arrows being inclusion maps.
Assume ${\mathcal D}_X$ is an $\{x_i\}${--}lattice, 
and assume $\emptyset\in{\mathcal D}_X$.
\medskip

\begin{xDefinition}
A \emph{covering functor} for ${\mathcal D}_X$ is a functor, $S$, which
assigns to each $\pi_1(X-C,x_i)$ a subgroup $S\pi_1(X-C,x_i)$
subject to
\[\begin{matrix}%
S\pi_1(X-C,x_i)&\subseteq& \pi_1(X-C,x_i)\\
\Big\downarrow&&\Big\downarrow\\
S\pi_1(X-D,x_i)&\subseteq& \pi_1(X-D,x_i)\\
\end{matrix}\]
commutes whenever $D\subseteq C$, $x_i\notin C$, and where the
vertical maps are induced by inclusion ($C$ and $D$ are any elements
of ${\mathcal D}_X$).
\end{xDefinition}
\medskip
\begin{xRemarks}
We have two examples for ${\mathcal D}_X$ in mind.
In this section we can use the set of all closed, compact subsets
of $X$ for ${\mathcal D}_X$.
For cohomology however, we will have to use 
the set of open subsets of $X$ with compact closure for ${\mathcal D}_X$.
\end{xRemarks}

\begin{xExamples}
There are three useful examples we shall define.
\begin{enumerate}
\item[1)] no covering functor (the subgroup is the whole group)
\item[2)] the universal covering functor (the subgroup is the
zero group)
\item[3)] the universal cover of $X$ but no more covering functor
(the subgroup is the kernel of
$\pi_1(X-C,x_i)\ \to\ \pi_1(X,x_i)$\ ).
\end{enumerate}
\end{xExamples}
\medskip
\begin{xDefinition}
A compatible covering functor for ${\mathcal D}_X$ is a covering functor
$S$ such that, for any $C\in{\mathcal D}_X$, the cover of the component
of $X-C$ containing $x_i$ corresponding to $S\pi_1(X-C,x_i)$ exists.
We write $(X,\coverFA{})$ for 
a compatible covering functor for ${\mathcal D}_X$ 
(which is inferred from context) to denote a collection of pointed spaces
$\bigl( (\coverFA{X-C})^i,\hat x_i\bigr)$, where $(\coverFA{X-C})^i$
is the covering space of the component of $X-C$ containing $x_i$, and
$\hat x_i$ is a lift of $x_i$ to this cover such that
$\pi_1\bigl((\coverFA{X-C})^i,\hat x_i\bigr)=S\pi_1(X-C,x_i)$.
Notice this notation is mildly ambiguous since if we change the $\hat x_i$
we get a different object.
As the two objects are homeomorphic this tends to cause no problems
so we use the more compact notation.

We say $(X,\coverFA{})\leq (X,\coverFB{})$
provided the subgroup of $\pi_1(X-C,x_i)$ corresponding to 
$ \coverFB{}$
contains the one corresponding to $\coverFA{}$.
Hence any $(X,\coverFA{})\leq (X, {{\text{no\ cover}}})$,
and if the universal covering functor is compatible with ${\mathcal D}_X$,
$(X, {\text{universal\ cover}}) \leq (X,\coverFA{})$.
\end{xDefinition}

Now the no covering functor is compatible with any ${\mathcal D}_X$.
If $X$ is semi{--}locally 1{--}connected, the universal cover of $X$
but no more is compatible with any ${\mathcal D}_X$.
If ${\mathcal D}_X$ is the collection of closed, compact subsets of $X$,
and if $X$ is locally 1{--}connected, the universal covering functor is
compatible with ${\mathcal D}_X$, as is any other covering functor.
Hence a CW complex is compatible with any covering functor
(see Lundell and Weingram \cite{btwentyone}  page 67, Theorem 6.6)
for ${\mathcal D}_X$.

We can now describe our construction.
Let $(X,\coverFA{})$ be a covering functor for $X$.
Assume from now on that our homotopy functor takes values in a
$({\mathcal D}_X, \{x_i\})${--}regular category for all homogamous $X$ with
base points $\{x_i\}$.
We apply the $\epsilon$ and $\Delta$ constructions to the collection
\[G_{i C}=\begin{cases}h\bigl((\coverFA{X-C})^i, \hat x_i\bigr)& \text{if } x_i\in X-C\\
0&\text{if }x_i\notin X-C\ .\end{cases}\]

If $D\subseteq C$ there is a unique map
$\bigl((\coverFA{X-D})^i, \hat x_i\bigr)\ \to\
\bigl((\coverFA{X-C})^i, \hat x_i\bigr)$ 
if $x_i\in X-D$ by taking the lift of the
inclusion which takes $\hat x_i$ in 
$(\coverFA{X-D})^i$ to $\hat x_i$ in
$(\coverFA{X-C})^i$.
Hence we get a map $G_{i D}\ \to\ G_{i C}$.
We denote these groups by
$\epsilon(X\Colon h,\{\hat x_i\},\coverFA{})$ and 
$\Delta(X\Colon h,\{\hat x_i\},\coverFA{}\ )$.

\BEGIN{T.1.3.1}
Let $\{ x_i\}$ and $\{ y_i\}$ be two sets of base points for $X$,
$X$ homogamous with countable base points.
Let $\Lambda$ be a locally finite collection of paths giving an equivalence
between $\{x_i\}$ and $\{y_i\}$.
Then there are natural transformations 
$\alpha_\Lambda\colon
\epsilon(X\Colon h,\{\hat x_i\},\coverFA{})\ \to\ 
\epsilon(X\Colon h,\{\hat y_i\},\coverFB{})$ and
$\alpha_\Lambda\colon
\Delta(X\Colon h,\{\hat x_i\},\coverFA{})\ \to\ 
\Delta(X\Colon h,\{\hat y_i\},\coverFB{})$
which are isomorphisms and depend only on the proper homotopy
class of $\Lambda$ rel end points. ($\coverFB{}$ is the covering
functor induced by the set of paths $\Lambda$.)
\end{Theorem}
\begin{proof}
Define $\alpha_\Lambda$ as follows.
By relabeling if necessary we may assume $x_i$ goes to $y_i$
by a path in $\Lambda$.
Map $h\bigl((\coverFA{X-C})^i, \hat x_i\bigr)$ to 
$h\bigl((\coverFB{X-C})^i, \hat y_i\bigr)$ by the zero map if the path from
$x_i$ to $y_i$ hits $C$.
If the path misses $C$, map $h\bigl((\coverFA{X-C})^i, \hat x_i\bigr)$ to 
$h\bigl((\coverFB{X-C})^i, \hat y_i\bigr)$ by lifting the path from
$x_i$ to $y_i$ into $(\coverFA{X-C})^i$ beginning at $\hat x_i$, say
it now ends at $z$, and then map $\bigl((\coverFA{X-C})^i,z\bigr)$
to $\bigl((\coverFB{X-C})^i,\hat y_i\bigr)$ by the unique homeomorphism
covering the identity which takes $z$ to $\hat y_i$.
This defines a homomorphism $\alpha_\Lambda$ on $\epsilon$
and $\Delta$.

If by $\Lambda^{-1}$ we mean the collection of paths from $y_i$ to $x_i$
given by the inverse of the path from $x_i$ to $y_i$, 
we can also define $\alpha_{\Lambda^{-1}}$.

$\alpha_{\Lambda}\ \circ\ \alpha_{\Lambda^{-1}}$ takes
$h\bigl((\coverFB{X-C})^i, \hat y_i\bigr)$ to itself by the zero
map if the path hits $C$ and by the identity otherwise.
Since all but finitely many paths miss $C$, this induces the identity on 
$\epsilon$.
Since the empty set is is the minimal element of ${\mathcal D}_X$,
$\alpha_{\Lambda}\ \circ\ \alpha_{\Lambda^{-1}}$ is the identity on
$\displaystyle\mathop{\bigtimes}_{i\in{\mathcal I}} h(X,y_i)$ and
$\mu\bigl(h(X,\hat y_i)\bigr)$.
Hence it is also the identity on $\Delta$.
A similar argument shows
$\alpha_{\Lambda^{-1}}\ \circ\ \alpha_{\Lambda}$
is the identity, so they are both isomorphisms.

The same sort of argument shows $\alpha_\Lambda$ depends only on
the proper homotopy type of $\Lambda$. 
It can be safely left to the reader. 
\end{proof}

If $h$ is actually a homotopy functor on the category of pairs 
(or $n${--}ads) we can define \\
$\gamma\bigl(X,A\Colon h,\{\hat x_i\},\coverFA{}\ \bigr)$
for the pair $(X,A)$ (where $\gamma$ denotes 
$\epsilon$ or $\Delta$) using
$G_{i C}= h\bigl((\coverFA{X-C})^i, (\coverFA{A}\ \cap\ (\coverFA{X-C})^i)\
\cup\ (\hat x_i), \hat x_i\ \bigr)$ where
$\coverFA{A}\ \cap\ (\coverFA{X-C})^i= \pi^{-1}\bigl(A\ \cap\ (X-C)^i\bigr)$,
$\pi\colon (\coverFA{X-C})^i\ \to\ (X-C)^i$, if $x_i\notin C$ and is $0$
otherwise ( for $n${--}ads use
$h\bigl((\coverFA{X-C})^i, (\coverFA{A}_1\ \cap\ (\coverFA{X-C})^i)\
\cup\ \{\hat x_i\},\cdots,(\coverFA{A}_{n-1}\ \cap\ (\coverFA{X-C})^i)\
\cup\ \{\hat x_i\},
 \hat x_i\ \bigr)$ or $0$).

Now suppose we have a connected sequence of homotopy functors 
$h_\ast$; i.e. each $h_n$ is a homotopy functor on some category
of pairs and we get long exact sequences.
By applying our construction to $(X,A)$, one would hope to get
a similar long exact sequence for the $\epsilon$ or $\Delta$ theories.

Several problems arise with this naive expectation.
To begin, we can certainly define groups which fit into a long exact
sequence.
Define $\gamma\bigl(A;X\Colon h_\ast,\{\hat x_i\},\coverFA{}\ \bigr)$
where $\gamma = \epsilon$ or $\Delta$ from 
$G_{i C} = h_\ast\Bigl(\bigl(\coverFA{A}\ \cap\ (\coverFA{X-C})^i\bigr)\ 
\cup\ \{\hat x_i\}, \hat x_i\Bigr)$ if $x_i\notin C$ and $0$ if $x_i\in C$.
Then \fullRef{C.1.2.2.1}, or its unstated analogue
for \fullRef{T.1.2.3}, shows we get a 
long exact sequence
\[\begin{aligned}%
\cdots\ \to\ \gamma\bigl(A;X\Colon h_n,&\{ \hat x_i\},\coverFA{}\ \bigr)
\ \to\ 
\gamma\bigl(X\Colon h_n,\{ \hat x_i\},\coverFA{}\ \bigr)\ \to\ 
\gamma\bigl(X, A\Colon h_n,\{ \hat x_i\},\coverFA{}\ \bigr)\\
&\to
\gamma\bigl(A;X\Colon h_{n-1},\{ \hat x_i\},\coverFA{}\ \bigr)\ \to\
\cdots\\
\end{aligned}\]
The problem of course is to describe $\gamma(A;X\Colon {\text{etc.}})$ 
in terms of $A$.

We clearly have little hope unless $A$ is homogamous, and 
for convenience we insist $A\subseteq X$ be a proper map.
Such a pair is said to be homogamous, and for such a pair we begin
to describe $\gamma(A;X\Colon {\text{etc.}})$.

Pick a set of base points for $A$, and then add enough new points to
get a set of base points for $X$.
Such a collection is a set of base points for $(X,A)$.
Two such are equivalent provided the points in $X-A$ can be made to
correspond via a locally finite collection of paths in $X$ all of which lie in
\setcounter{footnote}{0}
$X-A$\footnote{And two sets of base points for $A$ are equivalent.}.
A set of base points for $(X,A)$ is irreducible provided any subset 
which is also a set of base points for $(X-A)$ has the same cardinality
\footnote{And the set of base points for $A$ is irreducible.}.
(Note an irreducible set of base points for $(X,A)$ is not always an
irreducible set of base points for $X$.
$(S^1, S^0)$ is an example.)
We can construct $\epsilon$ and $\Delta$ groups for $X$ based on
an irreducible set of base points for $(X,A)$, and whenever we have
a pair, we assume the base points are an irreducible set of base points
for the pair.
If $X$ has no compact component, then any irreducible set of base
points for $(X,A)$ is one for $X$.
Over the compact components of $X$, the $\Delta$ group is just
the direct product of $h(\tilde X,p)$ for one $p$ in each component of $A$.
As in the absolute case, we drop irreducible and write 
``set of base points'' for ``irreducible set of base points''.

With a set of base points for $(X,A)$, there is a natural map
\[\gamma\bigl(A\Colon h, \{\hat x_i\}, \coverFA{}_F\ \bigr)\ \to\ 
\gamma\bigl(A; X\Colon h, \{\hat x_i\}, \coverFA{}\ \bigr)\ ,\]
where $\coverFA{}_F$ is the covering functor over 
$A$ induced as follows.
Let ${\mathcal D}(X)$ denote the following category.
The objects are closed compact subsets $C\subseteq X$.
The morphisms are the inclusions.
Given $A\subseteq X$ a closed subset, there is a natural map
${\mathcal D}(X)\ \to\ {\mathcal D}(A)$ given by $C\mapsto C\cap A$.
A lift functor $F\colon {\mathcal D}(A)\ \to\ {\mathcal D}(X)$ is a functor
such that ${\mathcal D}(A)\ \RA{F}\ {\mathcal D}(X)\ \to\ {\mathcal D}(A)$
is the identity and such that the image of $F$ is cofinal in ${\mathcal D}(X)$.
$\coverFA{}_F$ is the covering functor whose subgroups are the
pullbacks of
\[\begin{matrix}%
&&S\pi_1\Bigl(\big(X-F(C)\bigr)^i, x_i\Bigr)\\
&&\downlabeledarrow[\Bigg]{\cap}{}\\
\pi_1\bigl((A-C)^i, x_i\bigr)&\RA{\hskip 40pt}&
\pi_1\Bigl(\big(X-F(C)\bigr)^i, x_i\Bigr)\\
\end{matrix}\]
for $x_i\in A-C$, $C\in{\mathcal D}(A)$.
The existence of our natural map
$\gamma(A\Colon \cdots \coverFA{}_F)\ \to\ 
\gamma(A\Colon \cdots \coverFA{}\ )$
presupposes $\coverFA{}_F$ is compatible of $A$, but this is always
the case since the appropriate cover of $(A-C)^i$ is sitting in
$\bigl((\coverFA{X-F(C)}\bigr)^i$.
We denote this natural map by $\tau(A,X)$.

Notice first that $\tau(A,X)$ is a monomorphism since each map is.
Moreover, $\tau(A,X)$ is naturally split.
The splitting map is induced as follows.
We need only define it on some cofinal subset of ${\mathcal D}(X)$, so we 
define it on $\bigl\{ F(C)\ \vert\ C\in{\mathcal D}(A)\ \bigr\}$.
$h_\ast\Bigl(\bigl(\coverFA{A}\ \cap\ (\coverFA{X-F(C)})^i\bigr)\ \cup\ 
\{\hat x_i\},\ \hat x_i\ \Bigr)$ goes to $0$ if $x_i\notin A$, and it
goes to to
$h_\ast\bigl((\coverFA{A-C})^i,\hat x_i\bigr)$ if $x_i\in A$, where
$\coverFA{}$ in this last case is the cover given by the covering
functor $\coverFA{}_F$.
$\coverFA{A}\ \cap\ (\coverFA{X-F(C)})^i$ is just several disjoint copies
of $(\coverFA{A-C})^i$ union covers of other components of $A-C$.
The map collapses each of these covers of other components of $A-C$ 
to $\hat x_i$ and on the copies of $(\coverFA{A-C})^i$ it is just the
covering projection. 
At this point, this is all we can say about $\tau(A,X)$.
This map however has many more properties and we shall return to it again.

Now let $f\colon X\ \to\ Y$ be a proper map between homogamous spaces.
We have the mapping cylinder $M_f$.
$(M_f, X)$ is an homogamous pair (\fullRef{C.1.1.1.2}).
Let $\{x_i\}$ be a set of base points for $(M_f, X)$. 
We also have the homogamous pair $(M_f, Y)$.
By \fullRef{L.1.1.2} a set of base points for $Y$ is also a
set of base points for $(M_f, Y)$.
If $\{y_i\}$ is such a set, $\tau(Y, M_f)$ is an isomorphism.
This is seen by showing the splitting map is a monomorphism.
But if we use the lift functor
$F(C)=I\times f^{-1}(C)\ \cup\ C\subseteq M_f$ this is not hard to see.

Given a covering functor on $M_f$, it induces covering functors on $X$ and
$Y$, and these are the covering functors we shall use.
Given a covering functor on $Y$, we can get a covering functor on $M_f$
as follows.
The subgroups to assign to 
$\pi_1( M_f-f^{-1}(C)\times I\ \cup\ C)$ are the subgroups for
$\pi_1(Y-C)$.
One can then assign subgroups to all other required sets in 
such a way as to get a covering functor.
If we use the obvious lift functor for $Y$, the induced cover is the original.

By taking the cofinal collection $F(C)$, it is also not hard to see
$\gamma(M_f, Y\Colon h_n, \{ \hat y_i\}, \coverFA{}\ )=0$.
We define
$f\colon \gamma(X\Colon h_n, \{\hat x_i\},\coverFA{}\ )\ \to\ 
\gamma(Y\Colon h_n, \{\hat x_i\},\coverFA{}\ )$ if no component of
$Y$ is compact by
\[\begin{aligned}%
\gamma(X\Colon{\text{etc.}})\ \RA{\tau(X,M_f)}\ \gamma&(X;M_f\Colon{\text{etc.}})\ 
\RA{\hskip 10pt}\ \gamma(M_f\Colon h_n,\{\hat x_i\},\coverFA{}\ )\ 
\RA{\alpha_\Lambda}\\&\gamma(M_f\Colon h_n, \{\hat y_i\},\coverFA{}\ )
\ \LA{\ \cong}\ \gamma(Y\Colon {\text{etc.}})\\
\end{aligned}\]
Notice that this map may depend on the paths used to join 
$\{\hat x_i\}$ to $\{\hat y_i\}$.
If $f$ is properly $1/2${--}connected, 
(i.e. $f$ induces isomorphisms on $H^0$
and $H^0_{{\text{end}}}$: compare this definition and the one in 
\cite{beleven}) there is a natural choice of paths.

This choice is obtained as follows.
Take a set of base points $\{x_i\}$ for $X$.
By \fullRef{L.1.1.1}, $\bigl\{f(x_i)\bigr\}$
is a set of base points for $Y$.
Let $\{ x^\prime_i\}\subseteq \{x_i\}$ be any subset obtained by
picking precisely one element of $\{ x_i\}$ in each 
$f^{-1}f(x_i)$.
By \fullRef{L.1.3.1} below, $\{x_i^\prime\}$
is a set of base points for $X$.
Thus we can always find a set of base points for $X$ on which $f$ is 
$1${--}$1$ and whose image under $f$ is a set of base points for $Y$.
Take such a set of points as a set of base points for $(M_f,X)$.
Take their image in $Y$ as a set of base points for $(M_f,Y)$.
The paths joining these two sets are just the paths
\[\lambda_{x_i}(t) = \begin{cases}x_i\times t& 0\leq t\leq 1\\
f(x_i)& t=1\ .\\\end{cases}
\]

Given a properly $1/2${--}connected map $f$, 
we can get another definition of the induced map.
Pick a set of base points $\{x_i\}$ as in the last paragraph.
Then we have
\[f_\ast\colon \gamma(X\Colon h,\{\hat x_i\}, \coverFA{}\ )
\ \RA{\hskip 10pt}\ 
\gamma(Y\Colon h,\{\longHatFA{f(x_i)} \}, \coverFA{}\ )\]
defined by taking $h\bigl((\coverFA{X-C})^i, \hat x_i\bigr)\ \to\ 
h\bigl((\coverFA{Y-F(C)})^i, \longHatFA{f(x_i)}\bigr)$ 
by $f$, where $F$ is
a lift functor which splits ${\mathcal D}(Y)\ \to\ {\mathcal D}(X)$
and $F$ is the lift functor used to get the covering functor for $X$
from the one for $Y$.
One sees easily the two definitions of $f_\ast$ agree.

Now suppose we consider $i\colon A\subseteq X$ for an homogamous pair.
Then we can define $i_\ast$ as above. 
It is not hard to see 
\[\begin{matrix}%
\gamma(A;X\Colon {\text{etc.}})&\longrightarrow&\gamma(X\Colon{\text{etc.}})\\
\uplabeledarrow[\Bigg]{}{\tau(A,X)}&\hskip10pt\nearrow
 \lower 8pt\hbox{$i_\ast$}\\
\gamma(A\Colon {\text{etc.}})\\
\end{matrix}\]
commutes, where the paths we use in defining $i_\ast$ are $\lambda_{x_i}(t)= x_i\times t$ in
$A\times I\ \cup\ X\times 1 = M_i$.

\medskip
\BEGIN{L.1.3.1}
If $f\colon X\ \to\ Y$ is a proper map which induces epimorphisms on
$H^0$ and $H^0_{{\text{end}}}$, then, if $\bigl\{ f(p)\bigr\}$ is a set
of base points for $Y$, $\{ p\}$ is a set of base points for $X$.
\end{Lemma}
\begin{proof}
Since $f$ is an epimorphism on $H^0$, each path component of $X$ 
has a point of $\{p\}$ in it.

Now define a cochain in $S^0(X)$ for some closed compact set 
$D\subseteq X$,  $\varphi_D$ as follows.
$\varphi_D(q)=1$ if $q$ is in a path component of $X-D$ 
with no point of $\{ p\}$ in it and is $0$ otherwise.
$\delta\varphi_D = 0$ in $S^1_{{\text{end}}}$.

Since $f$ is an epimorphism on $H^0_{{\text{end}}}$, there must be a chain in
$S^0(Y)$, $\psi$, such that $f^\ast\psi=\varphi$ in $S^0_{{\text{end}}}(X)$.
But this means there is some closed compact set $C\subseteq X$ such
that $f^\ast\psi$ and $\varphi$ agree for any point in $X-C$.
Hence there is a closed, compact set $E\subseteq Y$ such that
$f^{-1}(E)\supseteq C\ \cup\ D$.
There is also a closed, compact $F\subseteq Y$ 
such that there is an $f(p)$
in each path component of $Y-E$ which is not contained in $F$.
$\psi$ restricted to $Y-E$ must be $0$ since some component of
$X-D$ which is not contained in $f^{-1}(E)$ has a point of $\{p\}$  in it.
Hence $\varphi$ restricted to $X-f^{-1}(E)$ is $0$, so we are done.
\end{proof}

\begin{xDefinition}
An homogamous pair $(X,A)$ is properly $0${--}connected if the inclusion
induces monomorphisms on $H^0$ and $H^0_{{\text{end}}}$.
We have already defined properly $1/2${--}connected.
If $(X,A)$ is properly $0${--}connected we can choose a set of base
points for the pair to be a set of base points for $A$.
We say $(X,A)$ is \emph{properly $n${--}connected}, $n\geq 1$
provided it is properly $1/2${--}connected, and, with base points chosen
as above,
$\Delta(X,A\Colon \pi_k, \{ x_i\}, {\text{no\ cover}}\ )=0$, $1\leq k\leq n$.
It is said to be \emph{properly $n${--}connected at $\infty$} provided it
is properly $1/2${--}connected and
$\epsilon(X,A\Colon \pi_k, \{ x_i\}, {\text{no\ cover}}\ )=0$, $1\leq k\leq n$.
\end{xDefinition}

\BEGIN{P.1.3.1}
If $(X,A)$ is properly $1/2${--}connected, and if 
\[i_\ast\colon \Delta(A\Colon \pi_1,\{ x_i\}, {\text{no\ cover}}\ )\ \to\ 
\Delta(X\Colon \pi_1,\{ x_i\}, {\text{no\ cover}}\ )\] is onto,
$(X,A)$ is properly $1${--}connected and conversely.
\end{Proposition}
\begin{proof}
If $(X,A)$ is properly $1/2${--}connected,
\[\Delta(A\Colon \pi_0,\{ x_i\}, {\text{no\ cover}}\ )\ \to\
\Delta(X\Colon \pi_0,\{ x_i\}, {\text{no\ cover}}\ )\]
is seen to be an isomorphism by applying 
\fullRef{T.1.2.4} to the kernel and cokernel
of this map, together with the definition of a set of base points.

Hence $\Delta(A;X \Colon \pi_1)\ \to\ 
\Delta(X\Colon \pi_1 )\ \to\ \Delta(X, A\Colon \pi_1)\to 0$ is exact.
\[\begin{matrix}%
\Delta(A;X\Colon \pi_1)&\longrightarrow&\Delta(X\Colon \pi_1)\\
\uplabeledarrow[\Bigg]{}{\tau(X,A)}&
\hskip14pt\nearrow\lower 4pt\hbox{$i_\ast$}\\
\Delta(A\Colon \pi_1)\\
\end{matrix}\]
commutes, and $i_\ast$ is an epimorphism.
Hence $\Delta(X,A\Colon \pi_1)=0$, so $(X,A)$ is properly $1${--}connected.
The converse follows trivially from 
\fullRef{P.1.3.2} below and the definitions.
\end{proof}

\BEGIN{P.1.3.2}
Let $(X,A)$ be a properly $1${--}connected pair.
Then $\tau(A,X)$ is an isomorphism if the base points for the pair
are a set of base points for $A$.
We may use any lift functor to induce the covering functor.
\end{Proposition}
\begin{proof}
If $\tau$ is an isomorphism on the $\epsilon$ objects, we need only show
$h(\coverFA{A},\hat x_i)= h(\coverFA{A}\ \cap\ \coverFA{X}, \hat x_i)$.
But $\coverFA{A}=\coverFA{A}\ \cap \coverFA{X}$ if $\pi_1(A)\to\pi_1(X)$
is onto, so if we can show the result for 
the $\epsilon$ objects we are done.

We need only show $\tau$ is onto. 
By \fullRef{T.1.2.4} applied to the cokernels of
the maps inducing $\tau$, we need only show 
that for each $C\in{\mathcal D}(X)$, there is a $D\supseteq C$ in ${\mathcal D}(X)$
such that
\[\begin{matrix}%
h\bigl((\coverFA{A-D})^i, \hat x_i)&\RA{\hskip10pt}&
h\biggl(\Bigl(\coverFA{A}\ \cap\ \bigl(\coverFA{X-F(D)}\bigr)^i\Bigr)\ 
\cup \{\hat x_i\},\hat x_i\biggr)\\
\downarrow[\Big]&&\downlabeledarrow[\Big]{i_\ast}{}\\
h\bigl((\coverFA{A-C})^i, \hat x_i)&\RA{\ \tau_\ast\ }&
h\biggl(\Bigl(\coverFA{A}\ \cap\ \bigl(\coverFA{X-F(C)}\bigr)^i\Bigr)\ 
\cup \{\hat x_i\},\hat x_i\biggr)\\
\end{matrix}\]
satisfies ${\Imx}\ i_\ast\ \subseteq {\Imx}\ \tau_\ast$
for all $x_i\notin D$.

We saw $\coverFA{A}\ \cap\ \bigl(\coverFA{X-F(C)}\bigr)^i$ 
was just some copies of $(\coverFA{A-C})^i$, together with covers
of components of $A-C\subseteq \bigl(X-F(C)\bigr)^i$.
Since $(X,A)$ is properly $1/2${--}connected, we can find $D$ so that
${\Imx}\ i_\ast\subseteq h\bigl({\text{copies\ of\ }}
(\coverFA{A-C})^i\bigr)$; i.e. we can find $D$ so that
$\bigl(X-F(D)\bigr)^i\ \cap\ (A-C) =
\bigl(X-F(D)\bigr)^i\ \cap\ (A-C)^i$.

Since $(X,A)$ is properly $1${--}connected, we can find 
$D_1\supseteq D$ so that
\[\begin{matrix}%
\pi_1\bigl( X-F(D_1), A-D_1, x_i\bigr)\\
\Big\downarrow\\
\pi_1\bigl( X-F(C), A-C, x_i\bigr)\\
\end{matrix}\]
is zero for all $x_i\notin D_1$.
But this says all the copies of
$(\coverFA{A-D_1})^i$ in 
$\coverFA{A}\ \cap\ \bigl(\coverFA{X-F(D_1)}\bigr)^i$ go to the same copy
of $(\coverFA{A-C})^i$ in 
$\coverFA{A}\ \cap\ \bigl(\coverFA{X-F(C)}\bigr)^i$, namely 
the one containing $\hat x_i$.
\end{proof}

\fullRef{T.1.2.4} can also be used to get

\begin{Named Theorem}[The subspace principle]
Let $(X,A)$ be an arbitrary homogamous pair.
Then \\
$\gamma(A;X \Colon h, \{\hat x_i\},\coverFA{}\ )=0$ \iff\ 
$\gamma(A\Colon h, \{\hat x_i\},\coverFA{}\ )=0$ provided, for the if part,
\begin{enumerate}
\item[1)]
if $A_\alpha$ is a collection of disjoint subsets of $A$,
$\displaystyle h(\mathop{\cup}_{\alpha}  \coverFA{A}_\alpha\ \cup\ p,p)
\cong \mathop{\oplus}_\alpha h(\coverFA{A}_\alpha\ \cup\ p,p)$.
\item[2)] if $E\subseteq B$ are subsets of $A$, and if there is a 
$q\in \coverFA{E}$ such that $h(\coverFA{E}, q)\ \to\ h(\coverFA{B},q)$
is the zero map, then
$h(\coverFA{E}\  \cup\ p, p)\ \to\ h(\coverFA{B}\  \cup\ p, p)$ is the zero
map for any $p$. 
$h$ need only be natural on subsets of $A$.
\end{enumerate}
\end{Named Theorem}

\begin{proof}
Only if is clear as $\tau(A,X)$ is naturally split, so we concentrate 
on the if part.

$\gamma(A\Colon h, \{\hat x_i\},\coverFA{}\ )=0$ implies by
\fullRef{T.1.2.4} that we can find a cofinal sequence
$C_0\subseteq C_1\subseteq\cdots$ of closed, compact subsets
of $A$ such that
$h\bigl((\coverFA{A-C_j})^i, \hat x_i\bigr)\ \to\ 
h\bigl((\coverFA{A-C_{j-1}})^i, \hat x_i\bigr)$
is the zero map for all $x_i\notin C_j$.
If $\gamma=\Delta$, $h(A)$ is also zero.

We then claim
$h\Bigl(\coverFA{A}\ \cap\ \bigl((\coverFA{X-F(C_j)}\bigr)^i
\ \cup\ \{\hat x_i\}, \hat x_i\Bigr)\ \to\ 
h\Bigl(\coverFA{A}\ \cap\ \bigl((\coverFA{X-F(C_{j-1})}\bigr)^i
\ \cup\ \{\hat x_i\}, \hat x_i\Bigr)$
is the zero map, and, if $\gamma=\Delta$,
$h(\coverFA{A}\ \cap\ \coverFA{X}, \hat x_i)=0$.

This last is easy since $\coverFA{A}\ \cap\ \coverFA{X}$ is the disjoint union
of copies of $\coverFA{A}$.
Now $\displaystyle\coverFA{A}\ \cap\ \bigl(\coverFA{X-F(C_j)}\bigr)^i=
\mathop{\cup}_\beta \mathop{\cup}_{\alpha_\beta} Z_{\alpha_\beta}$
where $\beta$ runs over the path components of $A-C_j$ in
$\bigl(X-F(C_j)\bigr)^i$, and $\alpha_\beta$ runs over 
the path components of $\pi^{-1}\bigl((A-C_j)^\beta\bigr)$ where
\[\pi\colon \bigl(\coverFA{X-F(C_j)}\bigr)^i \to\ 
\bigl(X-F(C_j)\bigr)^i\] 
is the covering projection and $(A-C_j)^\beta$
is the component of $A-C_j$ corresponding to $\beta$.
$Z_{\alpha_\beta}$ is the $\alpha_\beta^{\text{\thx}}$ component
of $\pi^{-1}\bigl((A-C_j)^i\bigr)$.

Similarly $\coverFA{A}\ \cap\ \bigl(\coverFA{X-F(C_{j-1})}\bigr)^i=
\displaystyle\mathop{\cup}_{\alpha}\mathop{\cup}_{\alpha_b}
Z_{\alpha_b}$.
The map we are looking at is just the map induced on the direct
sum by the maps
$h( Z_{\alpha_\beta}\ \cup\ \{\hat x_i\}, \hat x_i)\ \to\ 
h( Z_{a_b}\ \cup\ \{\hat x_i\}, \hat x_i)$
for the unique $a_b$ such that $Z_{a_b}$
is mapped into by $Z_{\alpha_\beta}$.
$Z_{\alpha_\beta}=(A-C_j)^\beta$, so if $\hat x_i\in Z_{\alpha_\beta}$,
the map is the zero map since it is then a map of the form
$h\bigl((\coverFA{A-C_j})^i, \hat x_i\bigr)\ \to\ 
h\bigl((\coverFA{A-C_{j-1}})^i, \hat x_i\bigr)$ which we know to be zero.

If $\hat x_i\notin Z_{\alpha_\beta}$, the map is now a map of the form
\par\noindent
$h\bigl((\coverFA{A-C_j})^\beta\  \cup \ \{\hat x_i\}, \hat x_i)
\ \to\ 
h\bigl((\coverFA{A-C_{j-1}})^\beta\  \cup \ \{\hat x_i\}, \hat x_i)$,
which is still zero by the properties of $h$.
\end{proof}

We now investigate the invariance of our construction.
\medskip
\BEGIN{T.1.3.2}
Let $f$, $g\colon X\ \to\ Y$ be properly homotopic maps 
between homogamous spaces.
Then there is a set of paths $\Lambda$ such that
\[\begin{matrix}%
\gamma(X\Colon h, \{\hat x_i\},\coverFA{}\ )&\RA{\ f_\ast\ }&
\gamma(Y\Colon h, \{\hat y_i\},\coverFA{}\ )\\
&\hbox to 0pt{\hss${g_\ast}$}\searrow&\downlabeledarrow[\Big]{\alpha_\Lambda}{}\\
&&\gamma(Y\Colon h, \{\hat y_i\},\coverFA{}\ )\\
\end{matrix}\]
commutes.
\end{Theorem}
\begin{proof}
Let $F \colon X\times I\to Y$ be the homotopy, and let $M_F$ be its
mapping cylinder.
Then it is possible to pick paths so that
\[\begin{matrix}%
\gamma(M_F\Colon h, \{\hat x_i\}\times 0,\coverFA{}\ )&\RA{\ \ \ }&
\gamma(M_F\Colon h, \{\hat y_i\},\coverFA{}\ )\\\noalign{\vskip4pt}
\Big\downarrow&&\downlabeledarrow[\Big]{\alpha_\Lambda}{}\\
\gamma(M_F\Colon h, \{\hat x_i\}\times 1,\coverFA{}\ )&\RA{\ \ \ }&
\gamma(M_F\Colon h, \{\hat y_i\},\coverFA{}\ )\\
\end{matrix}\]
commutes, where the horizontal maps are the maps induced by the
paths joining $\{\hat x_i\}\times t$ to $\{\hat y_i\}$ ($t=0$, $1$) and
the left hand vertical map is induced by the canonical path 
$\hat x_i\times 0$ to $\hat x_i\times 1$ in $X\times I\ \to\ M_F$.

It is now a chase of the definitions to show the desired diagram 
commutes.
\end{proof}

\BEGIN{C.1.3.2.1}
Let $f\colon X\to Y$ be a proper homotopy equivalence between
two homogamous spaces.
Then $f_\ast$ is an isomorphism. 
\end{Corollary}
\begin{proof}
There is a standard derivation of the corollary from the theorem.
\end{proof}

\BEGIN{C.1.3.2.2}
A proper homotopy equivalence between homogamous spaces 
is properly $n${--}connected for
all $n$ (i.e. its mapping cylinder modulo its domain is 
a properly $n${--}connected pair).
\end{Corollary}
\begin{proof}
$(M_f, X)$ is clearly properly $1/2${--}connected.
$i_\ast\colon \Delta(X\Colon \pi_1)\ \to\ \Delta(Y\Colon \pi_1)$ is onto, so it is
easy to show $(M_f,X)$ is properly $1${--}connected.
Then
$\gamma(X; M_f\Colon \pi_k)\ \cong\ \gamma(X\Colon \pi_k)\ \cong\ 
\gamma(M_f\Colon \pi_k)$, so
$\gamma(M_f, X\Colon \pi_k)=0$.
\end{proof}

\BEGIN{C.1.3.2.3}
If $f\colon X\to Y$ is a proper homotopy equivalence,
\[\gamma( M_f, X\Colon h, \{\hat x_i\}, \coverFA{}\ )=0\ .\]
\end{Corollary}
\begin{proof}
Since $f$ is properly $1${--}connected, $\gamma(X;M_f\Colon h, {\text{etc.}})\ \cong\ 
\gamma(X\Colon h, {\text{etc.}})$ by \fullRef{P.1.3.2}
$\gamma(X\Colon h, {\text{etc.}})\ \cong\ \gamma(M_f\Colon h, {\text{etc.}})$
by \fullRef{C.1.3.2.1}.
Hence 
$\gamma(M_f, X\Colon h, {\text{etc.}})=0$.
\end{proof}

In the other direction we have
\BEGIN{T.1.3.3}(Proper Whitehead)
Let $f\colon X\to Y$ be properly $n${--}connected.
Then for a locally finite CW complex, $K$, of dimension $\leq n$,
$f_\#\colon [K, X]\ \to\ [K,Y]$ is an epimorphism.
If $f$ is properly $(n+1)${--}connected, $f_\#$ is a bijection.
\end{Theorem}
\begin{xRemarks}
$[K,X]$ denotes the proper homotopy classes of proper maps of $K$
to $X$. 
For a proof of this result, see \cite{beleven} Theorem 3.4 and note
the proof is valid for $X$ and $Y$ homogamous.
\end{xRemarks}

\begin{xDefinition}
An homogamous space $Z$ is said to satisfy $D n$ provided the statement
of \fullRef{T.1.3.3} holds for $Z$ in place of $K$ and for 
each properly $n${--}connected map $f$ between homogamous spaces.
\end{xDefinition}

\BEGIN{P.1.3.3}
Let $Z$ be properly dominated by a space satisfying $D n$.
Then $Z$ satisfies $D n$.
\end{Proposition}
\begin{proof}
We leave it to the reader to modify the proof of 
\fullRef{P.1.1.1} to show $Z$ is homogamous 
\iff\ it is properly dominated by an homogamous space. 
Let $K$ be a space satisfying $D n$ and properly dominating $Z$.
Then $[Z,X]$ us a natural summand of $[K,X]$ 
for any homogamous $X$ and the result follows.
\end{proof}

We finish this section by proving a proper Hurewicz and a proper
Namioka theorem.

\begin{xDefinition}
A (as opposed to the) universal covering functor for $X$ is a
covering functor $\coverFA{}$ such that
$\epsilon(X\Colon \pi_1,\coverFA{}\ )=\Delta(X\Colon \pi_1,\coverFA{}\ )=0$.
Note that if the universal covering functor is
compatible with $X$, then it is a universal covering functor for $X$.
There are other examples however.
\end{xDefinition}

We start towards a proof of the Hurewicz theorem.
The proof mimics Spanier \cite{bthirtyfive} pages 391{--}393.
We first prove
\medskip
\BEGIN{L.1.3.2}
Suppose ${\mathcal G}=\{G_{i j}\}$ is a system of singular chain complexes on
spaces $X_{i j}$.
Suppose the projection maps $G_{i j}\ \to\ G_{i\>j-1}$ are induced
by continuous maps of the spaces $X_{i j}\ \to\ X_{i\>j-1}$.
Assume $i\geq 0$, $j\geq 0$.

Assume we are given a system $C=\{C_{i j}\}$, where each $C_{i j}$
is a subcomplex of $G_{i j}$ which is generated by the singular simplices 
of $G_{i j}$ which occur in $C_{i j}$.
Also assume that the projection $G_{i j}\ \to\ G_{i\>j-1}$
takes $C_{i j}\ \to\ C_{i\>j-1}$.

Lastly assume that to every singular simplex 
$\sigma\colon\Delta^q\to X_{i j}$ for $j\geq n$ ($n$ is given at the
start and held fixed throughout) there is assigned a map
$P_{i j}(\sigma)\colon
\Delta^q\times I\ \to\ X_{i\>j-n}$ which satisfies
\begin{enumerate}
\item[a)]
$P_{i j}(\sigma)(z,0)=\Bar\sigma(x)$, where
$\Bar\sigma\colon \Delta^q\ \RA{\sigma}\ X_{i j}\ \RA{\text{projection}}\ X_{i\>j-n}$.
\item[b)]
Define $\sigma_1\colon \Delta^q\ \to\ X_{i\>j-n}$ by
$\sigma_1(z)=P_{i j}(\sigma)(z,1)$.
Then we require that $\sigma_1\in C_{i\>j-n}$,
and, if $\sigma\in C_{i j}$, then $\sigma_1=\Bar\sigma$.
\item[c)] If $e^k_q\colon \Delta^{q-1}\ \to\ \Delta^q$
omits the $k^{\text{\thx}}$ vertex, then
$P_{i j}(\sigma)\ \circ\ (e^k_q\times 1) = P_{i j}(\sigma^{(k)})$.
\end{enumerate}
Then $\epsilon(C )\subseteq \epsilon({\mathcal G})$ is an homology equivalence.
(Compare Spanier \cite{bthirtyfive}, page 392, Lemma 7).
\end{Lemma}

\begin{proof}
Let $\alpha(i,k)\colon C_{i k}\subseteq G_{i k}$ be the inclusion, and
let $\tau(i,k)\colon G_{i k}\ \to\ C_{i\>k-n}$ be defined by
$\tau(i,k)(\sigma)=\Bar\sigma_1$ and extend linearly. 
(Here we must assume $k\geq n$).
Define $\rho_r\colon G_{i k}\ \to\ G_{i\>k-r}$ to be the
projection.

One easily checks that condition c) makes $\tau(i,k)$ into a chain map.
$\tau(i,k)\ \circ\ \alpha(i,k)\colon C_{i k}\penalty-2000 \to\ C_{i\>k-n}$ is
just the map induced on the $C_{i k}$ by $\rho_n$ on the $G_{i k}$.
This follows from condition b).

We claim $\alpha(i, k-n)\ \circ\ \tau(i,k)\colon G_{i k}\ \to\ G_{i\>k-n}$
is chain homotopic to $\rho_n$.
To show this, let $D_g\colon S(\Delta^q)\ \to\ S(\Delta^q\times I)$
be a natural chain homotopy between $\Delta(h_1)$ and
$\Delta(h_0)$, where $h_0$, $h_1\colon \Delta^q\ \to\ \Delta^q\times I$
are the obvious maps ($S$ is the singular chain functor).

Define a chain homotopy $D_{i k}\colon S(X_{i k})\ \to\ S(X_{i\>k-n})$
by $D_{i k}(\sigma)= S\bigl(P_{i k}(\sigma)\bigr)\bigl( D_q(\xi_q)\bigr)$
(where $\xi_q \colon \Delta^q\subset\Delta^q$ is the identity) where
$\sigma$ is a $q${--}simplex.
One checks, using c) and the naturality of $D_q$ that
$\partial D_{i k}+ D_{i k}\partial = \rho_n-
\alpha(i, k-n)\ \circ\ \tau(i,k)$.

By definition, $\displaystyle\epsilon({\mathcal G})=
\mathop{\lim}_{{\Atop{\leftarrow}{k}}}\ \mu(G_{i k})$ and
$\displaystyle\epsilon(C)=
\mathop{\lim}_{{\Atop{\leftarrow}{k}}}\ \mu(C_{i k})$.

Since
\[\begin{matrix}%
C_{i k}&\RA{\alpha(i,k)}&G_{i k}\\
\Big\downarrow&&\Big\downarrow\\
C_{i\>k-1}&\RA{\alpha(i,k-1)}&G_{i\>k-1}\\
\end{matrix}\]
commutes, we get a chain map 
$\alpha\colon\epsilon(C)\ \to\epsilon({\mathcal G})$,
which is just the inclusion.

Since $\tau(i,k)\ \circ\ \alpha(i,k)=\rho_n$,
\[\begin{matrix}%
C_{i k}&\RA{\alpha(i,k)} & G_{i k} &\RA{\tau(i,k)}&C_{i\>k-n}\\
\noalign{\vskip4pt}
\downlabeledarrow[\Big]{\rho_1}{}&&\downlabeledarrow[\Big]{\rho_1}{}
&&\downlabeledarrow[\Big]{\rho_1}{}\\
C_{i\>k-1}&\RA{\alpha(i,k-1)} & G_{i\>k-1} &\RA{\tau(i,k-1)}&C_{i\>k-1-n}\\
\end{matrix}\]
commutes along the outside square.
Unfortunately the right{--}hand square may not commute 
as we have made no stipulation as to the behavior of $P_{i j}$
with respect to $\rho_1$.
Similarly 
\[\begin{matrix}%
G_{i k}&\RA{\tau(i,k)}&C_{i\>k-n}&\RA{\alpha(i, k-n)}&G_{i\>k-n}\\
\downlabeledarrow[\Big]{\rho_1}{}&&&&\downlabeledarrow[\Big]{\rho_1}{}\\
G_{i\>k-1}&\RA{\tau(i,k-1)}&C_{i\>k-1-n}&\RA{\alpha(i, k-1-n)}&G_{i\>k-1-n}\\
\end{matrix}\]
may not commute.
However, since $\alpha(i, k-n)\ \circ\ \tau(i,k)$ 
is chain homotopic to $\rho_n$,
\[\begin{matrix}%
H_\ast(G_{i k})&\RA{H_\ast\bigl(\tau(i,k)\bigr)}&H_\ast(C_{i\>k-n})&
\RA{H_\ast\bigl(\alpha(i, k-n)\bigr)}&H_\ast(G_{i\>k-n})\\
\downlabeledarrow[\Big]{H_\ast(\rho_1)}{}&&&&\downlabeledarrow[\Big]{H_\ast(\rho_1)}{}\\
H_\ast(G_{i\>k-1})&\RA{H_\ast\bigl(\tau(i,k-1)\bigr)}&
H_\ast(C_{i\>k-1-n})&\RA{H_\ast\bigl(\alpha(i, k-1-n)\bigr)}&
H_\ast(G_{i\>k-1-n})\\
\end{matrix}\]
does commute.

Define $\beta(i,k)\colon G_{i k}\ \to\ C_{i\> k-2n}$ for $k\geq 2n$
by
$\beta(i,k)=\tau(i,k-n)\ \circ\ \alpha(i, k-n)\ \circ\ \tau(i,k)$.
We claim
\[\begin{matrix}%
H_\ast(G_{i k})& \RA{H_\ast\bigl(\beta(i,k)\bigr)}&H_\ast(C_{i\> k-2n})\\
\downlabeledarrow[\Big]{H_\ast(\rho_1)}{}&&
\downlabeledarrow[\Big]{H_\ast(\rho_1)}{}\\
H_\ast(G_{i\>k-1})& \RA{H_\ast\bigl(\beta(i,k-1)\bigr)}
&H_\ast(C_{i\> k-1-2n})\\
\end{matrix}\]
commutes.
To see this, look at
\[\begin{matrix}%
H_\ast(G_{i k})&\RA{H_\ast(\tau)}&
H_\ast(C_{i\>k-n})&\RA{H_\ast(\alpha)}&
H_\ast(G_{i\> k-n})&\RA{H_\ast(\tau)}&
H_\ast(G_{i\>k-2n})\\
\downlabeledarrow[\Big]{H_\ast(\rho_1)}{}&{\text{I}}&
\downlabeledarrow[\Big]{H_\ast(\rho_1)}{}&{\text{II}}&
\downlabeledarrow[\Big]{H_\ast(\rho_1)}{}&{\text{III}}&
\downlabeledarrow[\Big]{H_\ast(\rho_1)}{}\\
H_\ast(G_{i\>k-1})&\RA{H_\ast(\tau)}&
H_\ast(C_{i\>k-1-n})&\RA{H_\ast(\alpha)}&
H_\ast(G_{i\> k-1-n})&\RA{H_\ast(\tau)}&
H_\ast(G_{i\>k-1-2n})\\
\end{matrix}\]

The square II commutes since it already does on the chain level.
Similarly, the square II+III commutes.
The square I+II commutes on the homology level.
The desired commutativity is now a diagram chase.

Now define $\tau\colon \epsilon\bigl(H_\ast({\mathcal G})\bigr)\ \to\
\epsilon\bigl(H_\ast(C)\bigr)$ using the $H_\ast(\beta)$'s.
We also have $H_\ast(\alpha)\colon \hfill\penalty-10000
H_\ast\bigl(\epsilon(C)\bigr)\to\
H_\ast\bigl(\epsilon({\mathcal G})\bigr)$.
By \fullRef{C.1.1.2.2} we have
$H_\ast(\alpha)\colon \epsilon\bigl(H_\ast(C)\bigr)\ \to\
\epsilon\bigl(H_\ast({\mathcal G})\bigr)$.
$\tau\ \circ\ H_\ast(\alpha)$ and $H_\ast(\alpha)\ \circ\ \tau$
are both induced from the maps $H_\ast(\rho_{2n})$, and
hence are the identities on the inverse limits.
\end{proof}

\medskip
\BEGIN{L.1.3.3}
Let $X$ be an homogamous space.
Then we can find a countable, cofinal collection of closed,
compact sets $C_j\subseteq X$ with $C_j\subseteq C_{j+1}$.
Let $G_{i j}=S\bigl((\coverFA{X-C_j})^i, \hat x_i\bigr)$, the singular
chain groups on $(\coverFA{X-C_j})^i$.
Let $C_{i j}=S\bigl((\coverFA{X-C_j})^i,
\coverFA{A}\ \cap\ (\coverFA{X-C_j})^i, \hat x_i\bigr)^n$.
(see Spanier \cite{bthirtyfive}, page 391 for a definition).

Suppose $(X,A)$ is properly $1${--}connected and 
properly $n${--}connected at $\infty$ for $n\geq 0$.
Then the inclusion map
$\epsilon(C)\subseteq\epsilon({\mathcal G})$ is an homology equivalence.
(Notice that if we pick a set of base points $x_i$ for $A$, they are
a set for the pair, and $\epsilon({\mathcal G})=
\epsilon(X\Colon H_\ast,\{\hat x_i\},\coverFA{}\ )$.)
\end{Lemma}
\begin{proof}
Let $r=\min(q,n)$.
Then we produce for every $\sigma\in G_{i j}$ a map
\[
P_{i j}(\sigma) \colon \Delta^q\times I\ \to\ 
\bigl((\coverFA{X-C_{j-r}})^i,\hat x_i\bigr)
\]
which satisfies
\begin{enumerate}
\item[a)]
$P_{i j}(\sigma)(z,0)=\Bar\sigma\colon\Delta^q\ \RA{\sigma}\ 
\bigl((\coverFA{X-C_{j}})^i\bigr)\ \RA{\text{projection}}\ 
\bigl((\coverFA{X-C_{j-r}})^i\bigr)$.
\item[b)] If $\sigma_1(z)=P_{i j}(\sigma)(z,1)$, $\sigma_1\in C_{i\>j-r}$, 
and if $\sigma\in C_{i j}$, 
\[P_{i j}(\sigma)\colon
\Delta^q\times I\ \RA{{\text{proj}}}\ \Delta^q\ \RA{\sigma}\ 
\bigl((\coverFA{X-C_j})^i\bigr)\ \RA{{\text{projection}}}\ 
\bigl((\coverFA{X-C_{j-r}})^i\bigr)\]
\item[c)] $P_{i j}(\sigma)\ \circ\ (C^k_q\times 1)=
\begin{cases} ({\text{projection\ 1\ step}}\ )\ \circ\ P_{i j}(\sigma^{(k)})
& q\leq n\\
P_{i j}(\sigma^{(k)})&q>n\\\end{cases}$
\end{enumerate}

From such a $P$ it is easy to see how to get a $P$ as required by our
first lemma.
We remark that $C_{i j}$ and $G_{i j}$ satisfy all the other requirements 
to apply the lemma.
Hence \fullRef{L.1.3.2} will then give us the
desired conclusion.

We define $P_{i j}$ by induction on $q$.
Let $q=0$.
Then $\sigma\in G_{i j}$ is a map
$\sigma\colon \Delta^0\to \bigl((\coverFA{X-C_j})^i\bigr)$.
Since the point $\sigma(\Delta^0)$ lies in the same path component of
$(\coverFA{X-C_j})^i$ as $\hat x_i$, there is a path joining them.
Let $P_{i j}(\sigma)$ be such a path.
If $\sigma(\Delta^0)=\hat x_i$, $P_{i j}(\sigma)$ 
should be the constant path.
This defines $P_{i j}$ for $q=0$, and $P$ is easily seen to satisfy a){--}c).

Now suppose $P_{i j}$ is defined for all $\sigma$ of degree $<q$, 
$0<q\leq n$ so that it has properties a){--}c).

If $\sigma\in C_{i j}$, b) defines $P(\sigma)$, and $P$ then satisfies
a) and c).
So suppose $\sigma\notin C_{i j}$.
a) and c) define $P_{i j}$ on $\Delta^q\times 0\ \cup\ 
\dot\Delta^q\times I$; i.e. we get a map
$f\colon \Delta^q\times 0\ \cup\ \dot\Delta^q\times I\ \to\
(\coverFA{X-C_{j-q+1}})^i$.
There is a homeomorphism $h\colon E^q\times I\approx \Delta^q\times I$
such that $h(E^q\times 0)=\Delta^q\times 0\ \cup\ \dot\Delta^q\times I$;
$h(S^{q-1}\times 0)=\dot\Delta^q\times 1$; and
$h(S^{q-1}\times I\ \cup\ E^q\times 1)=\Delta^q\times 1$.
Let $g\colon (E^q,S^{q-1})\ \to\ 
\bigl((\coverFA{X-C_{j-q+1}})^i,
\coverFA{A}\ \cup\ (\coverFA{X-C_{j-q+1}})^i\bigr)$
be defined by $g=f\ \circ\ h$.

Because $q\leq n$ and $(X,A)$ is properly $n${--}connected at $\infty$,
we could have chosen (and did) the $C_j$ so that
\[\pi_q(X-C_k,A\cap(X-C_k),\ast)\ \to
\pi_q(X-C_{k-1},A\cap(X-C_{k-1}),\ast)\]
is the zero map for $q\leq n$.
Thus we get a homotopy
\[H\colon (E^q,S^{q-1})\times I\ \to\ 
\bigl((\coverFA{X-C_{j-q}})^i, \coverFA{A}\ \cap\ (\coverFA{X-C_{j-q}})^i\bigr)\]
between $\rho_1\ \circ\ g$ and an element of $C_{i\>j-q}$.

Define $P_{i j}(\sigma)$ to be the composite 
$\Delta^q\times I\ \RA{h^{-1}\times{\text{id}}}\ E^q\times I\ 
\RA{H}\ 
(\coverFA{X-C_{j-q}})^i$.
$P_{i j}$ clearly satisfies a) and b). Since $h$ was chosen carefully
c) is also satisfied.

In this way $P$ is defined for all simplices of degree $\leq n$.
Note that a singular simplex of degree $>n$ is in $C_{i j}$ \iff\ every
proper face is in $C_{i j}$.

Suppose that $P$ has been defined for all degrees $<q$, where $q>n$.
If $\sigma\in C_{i j}$, we define $P_{i j}(\sigma)$
by b) as usual.
It satisfies a) and c).
So suppose $\sigma\notin C_{i j}$.
Then a) and c) define a map
$f\colon \Delta^q\times 0\ \cup\ \dot\Delta^q\times I\ \to\
(\coverFA{X - C_{j-n}})^i$.
By the homotopy extension property we can extend $f$ to some map
$P(\sigma)\colon \Delta^q\times I\ \to\ (\coverFA{X - C_{j-n}})^i$.
It clearly satisfies a) and c).
It also satisfies b) since every proper face of $\sigma_1$ is in $C_{i\> j-n}$.
Hence we have defined our $P$.
\end{proof}

\bigskip
Now define $\epsilon^{(n)}(X,A\Colon H_q,\coverFA{}\ )$ to be
\[\begin{aligned}%
\epsilon\Biggl(H_q\biggl( &S\Bigl( (\coverFA{X-C_j})^i,
\coverFA{A}\ \cap\ (\coverFA{X-C_j})^i,\hat x_i\Bigr)^n /\\
&S\Bigl( (\coverFA{X-C_j})^i,
\coverFA{A}\ \cap\ (\coverFA{X-C_j})^i,\hat x_i\Bigr)^n\ \cap\ 
S\Bigl(\coverFA{A}\ \cap\ (\coverFA{X-C_j})^i, \hat x_i\Bigr)
\biggr)\Biggr)\ .\\\end{aligned}\]

Then there are natural maps
\[\epsilon^{(n)}(X,A\Colon H_q,\coverFA{}\ )\ \to\ 
\epsilon^{(n-1)}(X,A\Colon H_q,\coverFA{}\ )\ \to\ 
\cdots \ \to\ \epsilon(X,A\Colon H_q,\{\hat x_i\},\coverFA{}\ )\ .\]
\bigskip
\BEGIN{L.1.3.4}
Assume $(X,A)$ is a properly $1${--}connected pair which is 
properly $n${--}connected at $\infty$ for some $n\geq0$.
Then the natural map
$\epsilon^{(n)}(X,A\Colon H_q,\coverFA{}\ )\ \to\ 
\epsilon(X,A\Colon H_q,\{ \hat x_i\},\coverFA{}\ )$ is an isomorphism for all $q$.
\end{Lemma}
\begin{proof}
We have the following commutative diagram
\[\begin{aligned}%
&\begin{matrix}0\to&S\bigl((\coverFA{X-C_j})^i, {\text{etc.}}\bigr)^n
\ \cap\ S\bigl(\coverFA{A}\ \cap\ (\coverFA{X-C_j})^i\bigr)
&\to&S\bigl((\coverFA{X-C_j})^i, {\text{etc.}}\bigr)\\
&\downlabeledarrow[\Big]{\alpha}{}&&\downlabeledarrow[\Big]{\beta}{}\\
0\ \hbox to 0pt{$\RA{\hskip 2cm}$\hss}\hfill&
S\bigl( \coverFA{A}\ \cap\ (\coverFA{X-C_j})^i\bigr)&
\hbox to 1cm{\hss$\RA{\hskip 2.5cm}$}&
S\bigl((\coverFA{X-C_j})^i\bigr)\\
\end{matrix}\\\noalign{\vskip 14pt}
&\begin{matrix}\to&({\text{the\ quotient\ complex}})&\to 0\\%
&\downlabeledarrow[\Big]{\gamma}{}\\
\to&S\bigl((\coverFA{X-C_j})^i, \coverFA{A}\ \cap\ (\coverFA{X-C_j})^i\bigr)&
\to0\\
\end{matrix}\\
\end{aligned}\]
where (the quotient complex) was used in defining
$\epsilon^{(n)}(X,A\Colon H_q,\coverFA{}\ )$.

Now, since $(X,A)$ is properly $1${--}connected, the subspace groups
$\epsilon^{(n)}(A;X\Colon H_q,\coverFA{}\ )$ and \\
$\epsilon(A; X\Colon H_q,\coverFA{}\ )$ are the absolute groups.
Since $(X,A)$ is properly $1${--}connected and properly $n${--}connected
at infinity for all $n$, \fullRef{L.1.3.3} says
$\epsilon(\alpha)$ is an isomorphism on homology.
Similarly, \fullRef{L.1.3.3} says 
$\epsilon(\beta)$ is an isomorphism on homology.
Thus $\epsilon(\gamma)$ is an isomorphism on homology as asserted.
\end{proof}

\BEGIN{T.1.3.4}
Suppose $(X,A)$ is properly $1${--}connected and 
properly $(n-1)${--}connected at $\infty$ for some $n\geq 2$.
Then the Hurewicz map
\[\epsilon(X,A\Colon \pi_n^\prime,\{\hat x_i\},\coverFA{}\ )\ \to\ 
\epsilon(X,A\Colon H_n,\{\hat x_i\},\coverFA{}\ )\]
is an isomorphism, where $\pi^\prime_n\bigl(
(\coverFA{X-C_j})^i, \coverFA{A}\ \cap\ (\coverFA{X-C_j})^i \ \cup\ 
\hat x_i,\hat x_i\bigr)$ is $\pi_n$ quotiented out by the action of
$\pi_1\bigl( \coverFA{A}\ \cap\ (\coverFA{X-C_j})^i \ \cup\ 
\hat x_i,\hat x_i\bigr)$.
\end{Theorem}

\begin{proof}
The usual Hurewicz theorem contains the fact that
\[\pi^\prime_n\bigl(
(\coverFA{X-C_j})^i, \coverFA{A}\ \cap\ (\coverFA{X-C_j})^i ,\hat x_i\bigr)\to
H_n^{(n-1)}\bigl((\coverFA{X-C_j})^i, \coverFA{A}\ \cap\ 
(\coverFA{X-C_j})^i ,\hat x_i\bigr)\] is an isomorphism.
Thus $\epsilon(X,A\Colon \pi^\prime_n, \{\hat x_i\},\coverFA{}\ )\ \to\ 
\epsilon^{(n-1)}(X,A\Colon H_n, \coverFA{}\ )$ is an isomorphism.
But \fullRef{L.1.3.4} says 
$\epsilon^{(n-1)}(X,A\Colon H_n,\coverFA{}\ )\ \to\ 
\epsilon(X,A\Colon H_n, \{\hat x_i\},\coverFA{}\ )$ is an isomorphism.
\end{proof}

\BEGIN{T.1.3.5}
Suppose that $\epsilon(A\Colon \pi_1,\{ \hat x_i\},\coverFA{}_F\ )=0$
where $\coverFA{}_F$ is the cover over $A$ induced by the lift functor
from a cover $\coverFA{}$ over $X$.
Then the natural surjection\\
$\epsilon(X,A\Colon \pi_n, \{\hat x_i\},\coverFA{}\ )\ \to\ 
\epsilon(X,A\Colon \pi^\prime_n; \{\hat x_i\},\coverFA{}\ )$
is an isomorphism.
\end{Theorem}
\begin{proof}
Set $G_{i j}=\pi_n\bigl( (\coverFA{X-C_j})^i,
\coverFA{A}\ \cap\ (\coverFA{X-C_j})^i\ \cup\ \{\hat x_i\},\hat x_i\bigr)$
and $H_{i j}=\pi^\prime_n\bigl( (\coverFA{X-C_j})^i,
\coverFA{A}\ \cap\ (\coverFA{X-C_j})^i\ \cup\ \{\hat x_i\},\hat x_i\bigr)$.
Define $K_{i j}$ to be the kernel of $G_{i j}\ \to\ H_{i j}\ \to\ 0$.
$K_{i j}$ is generated by elements of the form $x-\alpha x$ where
$x\in \pi_n\bigl( (\coverFA{X-C_j})^i,
\coverFA{A}\ \cap\ (\coverFA{X-C_j})^i\ \cup\ \{\hat x_i\},\hat x_i\bigr)$
and $\alpha\in \pi_1\bigl(
\coverFA{A}\ \cap\ (\coverFA{X-C_j})^i\ \cup\ \{\hat x_i\},\hat x_i\bigr)$.

Since $\epsilon(A\Colon \pi_1,\{\hat x_i\},\coverFA{}_F\ )=0$,
the subspace principle says that we can assume the map
$\pi_1\bigl(
\coverFA{A}\ \cap\ (\coverFA{X-C_j})^i\ \cup\ \{\hat x_i\},\hat x_i\bigr)\ \to\ 
\pi_1\bigl(
\coverFA{A}\ \cap\ (\coverFA{X-C_{j-1}})^i\ \cup\ \{\hat x_i\},\hat x_i\bigr)$
is the zero map.
Then $K_{i j}\ \to\ K_{i\>j-1}$ takes $x - \alpha x$ to
$i_\ast(x)-i_\ast(\alpha x)=i_\ast(x)-i_\#(\alpha)\cdot i_\ast( x)=
i_\ast(x)-i_\ast(x)=0$, so this map is the zero map.
\end{proof}

\setcounter{footnote}{0}
\BEGIN{T.1.3.6}
Let $(X,A)$ be a properly $1${--}connected pair.
Then, for any covering functor $\coverFA{}$ on $X$, the natural map
$\epsilon(X,A\Colon \pi_n,\{\hat x_i\},\coverFA{}\ )\ \to\ 
\epsilon(X,A\Colon \pi_n,\{\hat x_i\},{\text{no\ cover}}\ )$ 
is an isomorphism\footnote{For $n\geq 2$.}.
\end{Theorem}
\begin{proof}
We have
\[\begin{matrix}%
\cdots\to\epsilon(A\Colon \pi_k,\{\hat x_i\},\coverFA{}_F\ )&\too&
\epsilon(X\Colon \pi_k,\{\hat x_i\},\coverFA{}\ )&\too&
\epsilon(X,A\Colon \pi_k,\{\hat x_i\},\coverFA{}\ )\to\cdots\\
\Big\downarrow&&\Big\downarrow&&\Big\downarrow\\
\cdots\to
\epsilon(A\Colon {\text{ditto}},{\text{no\ cover}})&\too&
\epsilon(X\Colon {\text{ditto}},{\text{no\ cover}})&\too&
\epsilon(X,A\Colon {\text{ditto}},{\text{no\ cover}})\to\cdots\\
\end{matrix}\]
commutes.
The first two maps are clearly isomorphisms for $k\geq2$, so the third
is for $k\geq3$.
Moreover
\[\begin{matrix}%
\epsilon(A\Colon \pi_1,\{\hat x_i\},\coverFA{}_F\ )&\to&
\epsilon(X\Colon \pi_1,\{\hat x_i\},\coverFA{}\ )\\
\Big\downarrow&&\Big\downarrow\\
\epsilon(A\Colon \pi_1,\{\hat x_i\},{\text{no\ cover}}\ )&\to&
\epsilon(X\Colon \pi_1,\{\hat x_i\},{\text{no\ cover}}\ )\\
\end{matrix}\]
is a pullback since it is obtained as the $\epsilon$ 
construction applied to pullbacks.
Hence the theorem remains true for $k=2$.
\end{proof}

\BEGIN{C.1.3.6.1}
Suppose $(X,A)$ is a properly $1${--}connected pair which is
$(n-1)${--}connected at $\infty$ for some $n\geq2$.
If $n=2$ assume $\epsilon(A\Colon \pi_1\{\hat x_i\},{\text{no\ cover}})\ \to\ \\
\epsilon(X\Colon \pi_1\{\hat x_i\},{\text{no\ cover}})$ is an isomorphism.
Then the Hurewicz map
\[\epsilon(X,A\Colon \pi_n,\{\hat x_i\},{\text{no\ cover}})\ \to\ 
\epsilon(X,A\Colon H_n,\{\hat x_i\},\coverFA{}\ )\] is an isomorphism
 where\ \ $\coverFA{}$\quad  is any universal covering functor for $X$.
\end{Corollary}

\BEGIN{T.1.3.7}
Theorems \shortFullRef{T.1.3.4}, \shortFullRef{T.1.3.5} 
and \shortFullRef{T.1.3.6} are true (after appropriate
changes) with $\Delta$ instead of $\epsilon$.
They are also true for the absolute groups.
\end{Theorem}
\begin{proof} Easy.\end{proof}
\bigskip
Now suppose $(X,A)$ is a locally compact CW pair.
Then we might hope to improve our Hurewicz theorems by 
getting information about the second non{--}zero map 
(see \cite{bfortytwo}).
We do this following Hilton \cite{bthirteen}.

\begin{xDefinition}
Two proper maps $f$, $g\colon X\ \to\ Y$ are said to be \emph{properly
$n${--}homotopic} if for every proper map $\phi\colon K\ \to\ X$,
where $K$ is a locally compact CW complex of dimension $\leq n$,
$f \circ\ \phi$ is properly homotopic to $g\ \circ\ \phi$.
$X$ and $Y$ are of the \emph{same proper $n${--}homotopy type}
provided there exist proper maps $f\colon X\ \to\ Y$ and 
$g\colon Y\ \to\ X$ such that $f\ \circ\ g$ and $g\ \circ\ f$ are properly
$n${--}homotopic to the identity.
Two locally compact CW complexes, $K$ and $L$, are said to be of
the same \emph{proper $n${--}type} \iff\ $K^n$ and $L^n$ have the same
proper $(n-1)${--}type.
A proper cellular map $f\colon K\ \to\ L$ is said to be a {\sl
proper $n${--}equivalence} provided there is a proper map
$g\colon L^{n+1}\ \to\ K^{n+1}$ with $f\vert_{K^{n+1}}\ \circ\ g$ and
$g\ \circ\ f\vert_{K^{n+1}}$ properly $n${--}homotopic to the identity.
\end{xDefinition}

A \emph{proper $J_m${--}pair}, $(X,A)$, is a properly $1${--}connected,
locally compact CW pair such that the maps
$\Delta(X^{n-1}\ \cup\ A, A\Colon \pi_n,\{\hat x_i\},{\text{no\ cover}})\to\ 
\Delta(X^{n}\ \cup\ A, A\Colon \pi_n,\{\hat x_i\},{\text{no\ cover}})$ are zero
for $2\leq n\leq m$.
A \emph{proper $J_m${--}pair at $\infty$} is the obvious thing.

\BEGIN{L.1.3.5}
The property of being a proper $J_m${--}pair is an invariant of 
proper $m${--}type.
\end{Lemma}
\begin{proof} See Hilton \cite{bthirteen}.
\end{proof}

\BEGIN{T.1.3.8}
Let $(X,A)$ be a proper $J_m${--}pair at $\infty$.
Then the Hurewicz  map
$h_n\colon \epsilon(X,A\Colon \pi_n,\{\hat x_i\},{\text{no\ cover}})\to
\epsilon(X,A\Colon H_n,\{\hat x_i\},\coverFA{}\ )$, 
where\ $\coverFA{}$\quad
is a universal covering functor for $X$, satisfies
$h_n$ is an isomorphism for $n\leq m$ and $h_{m+1}$ is an epimorphism.
\end{Theorem}
\begin{proof}
See Hilton \cite{bthirteen}.
\end{proof}
\BEGIN{C.1.3.8.1}
The same conclusions hold for a proper $J_m${--}pair with the $\Delta$
groups.
\end{Corollary}

\BEGIN{C.1.3.8.2}
Let $(X,A)$ be a properly $(n-1)${--}connected, locally compact CW pair,
for $n\geq 2$.
If $n=2$ let $\Delta(A\Colon \pi_1,\{\hat x_i\},{\text{no\ cover}})\ \to\ 
\Delta(X\Colon \pi_1,\{\hat x_i\},{\text{no\ cover}})$ be an isomorphism.
Then the Hurewicz map
$h_n\colon \Delta(X,A\Colon \pi_n,\{\hat x_i\},{\text{no\ cover}})\to
\Delta(X,A\Colon H_n,\{\hat x_i\},\coverFA{}\ )$, 
is an isomorphism, where\ $\coverFA{}$\quad
is a universal covering functor for $X$.
$h_{n+1}$ is an epimorphism.
\end{Corollary}
\begin{proof}
In section 5 we will see there is a locally finite $1${--}complex 
$T\subseteq A$ such that $(A,T)$ is a proper $1/2${--}equivalence
and $\Delta(T\Colon \pi_k.\{\hat x_i\}),{\text{no\ cover}})=0$ for $k\geq1$.
Then $(T,T)$ is certainly a proper $J_n${--}complex.
$(T,T)\subseteq (X,A)$ is a proper $(n-1)${--}equivalence, so
$(X,A)$ is a $J_n${--}complex by \fullRef{L.1.3.5}.
\end{proof}

\BEGIN{T.1.3.9}
(Namioka \cite{btwentyeight})
Let $\phi\colon (X,A)\ \to\ (Y,B)$ be a map of pairs of
locally compact CW complexes.
Let $\phi\vert_X$ and $\phi\vert_A$ be properly $n${--}connected,
$n\geq 1$. ($\phi\vert_X$ and $\phi\vert_A$ should induce 
isomorphisms on $\Delta(\quad\Colon \pi_1,\{\hat x_i\},{\text{no\ cover}})$
if $n=1$).
Then the Hurewicz map
$h_{n+1}\colon
\Delta\bigl((M_\phi\Colon M_{\phi\vert_A}, X)\Colon 
\pi_{n+1},\{\hat x_i\}, {\text{no\ cover}}\bigr)\ \to\ 
\Delta\bigl((M_\phi\Colon M_{\phi\vert_A}, X)\Colon
H_{n+1},\{\hat x_i\}, \coverFA{}\ \bigr)$, where \ $\coverFA{}$\quad
is a universal covering functor of $M_\phi$, is an epimorphism.
\end{Theorem}
\begin{proof}
$(M_\phi\Colon M_{\phi\vert_A},X)$ is a triad and the groups in question
are the proper triad groups.
The reader should have no trouble defining these groups.
We can pick a set of base points for $(X,A)$ and it will also be a set for
our triad.

The triad groups fit into a long exact sequence
\[\cdots\to\ 
\Delta\bigl((M_{\phi\vert_A}); M_\phi)\ \to\ 
\Delta\bigl((M_{\phi}), X)\ \to\ 
\Delta\bigl((M_\phi\Colon M_{\phi\vert_A},X)\bigr)\ \to\cdots \quad,\]
where again we get the subspace groups.
Since $\phi\vert_{A}$ and $\phi\vert_{X}$ are properly $n${--}connected,
$h_m$ for $(M_\phi,A)$ is an isomorphism for $m\leq n$ and 
an epimorphism for $m=n+1$.
By the subspace principle, $h_m$ for 
$\bigl((M_{\phi\vert_A}, A ); M_\phi\bigr)$ is an isomorphism for
$m\leq n$ and an epimorphism for $m=n+1$.
The strong version of the $5${--}lemma now shows the triad $h_n$ an
isomorphism and the triad $h_{n+1}$ an epimorphism.
\end{proof}

\insetitem{Notation}
$\Delta_\ast(X,A\Colon \coverFA{}\ )$ will hereafter denote 
$\Delta(X,A\Colon H_\ast,\{\hat x_i\},\coverFA{}\ )$
for some set of base points for the pair $(X,A)$.
Similar notation will be employed for homology $n${--}ad groups,
subspace groups, etc.
\medskip
We conclude this section with some definitions and computations.

\begin{xDefinition}
An homogamous space $X$ is said to have \emph{monomorphic ends},
provided \\
$\displaystyle\Delta(X\Colon \pi_1,\{ x_i\},{\text{no\ cover}})\to
\mathop{\bigtimes}_{i\in{\mathcal I}}\ \pi_1(X,x_i)$ is a monomorphism
(equivalently $\displaystyle\epsilon\ \to\ \mathop{\mu}_{i}$ 
is a monomorphism).
A space has \emph{epimorphic ends} provided the above map is onto,
and \emph{isomorphic ends} if the map is an isomorphism.
\end{xDefinition}

As examples, if $X$ is an homogamous space which is not compact,
$X\times\R$ has one, isolated end (see \cite{bthirtytwo}) which is
epimorphic.
$X\times\R^2$ has isomorphic ends.
These results use Mayer{--}Vietoris to compute the number of ends of
$X\times\R$ and van{--}Kampen to yield the $\pi_1$ information,
using the following pushout
\[\begin{matrix}%
(X-C)\times Y-D)&\to&X\times (Y-D)\\
\Big\downarrow&&\Big\downarrow\\
(X-C)\times Y&\to&X\times Y\ -\ C\times D\\
\end{matrix}\]

In fact, this diagram shows that if $X$ and $Y$ are not compact
(but are path connected),
$X\times Y$ has one end, which is seen to be epimorphic since
$\pi_1(X\times Y\ -\ C\times D)\ \to\ \pi_1(X\times Y)$ is
easily seen to be onto.
If $X$ has epimorphic ends, $\pi_1(X-C,p)\to\pi_1(X,p)$ 
must always be onto, so if $X$ and $Y$ have epimorphic ends,
$X\times Y$ has one isomorphic end.

Monomorphic ends are nice for then the third example of covering functor
that we gave (the universal cover of $X$ bit no more) becomes a 
universal covering functor.
Farrell and Wagoner (\cite{bnine} or \cite{beleven}) then showed \
that a proper map $f\colon X\ \to\ Y$, $X$, $Y$ locally compact CW,
with $X$ having monomorphic ends is a proper homotopy equivalence
provided it is a properly $1${--}connected map; a homotopy equivalence;
and $f^\ast\colon H^\ast_\cmpsup(\tilde Y)\ \to\ 
H^\ast_\cmpsup(\tilde X)$ is an isomorphism 
where $\tilde{\ }$ denotes the universal cover (coefficients
are the integers).

\section{Proper cohomology, coefficients and products}
\newHead{I.4}
In attempting to understand ordinary homotopy theory, cohomology theory 
is an indispensable tool.
In ordinary compact surgery, the relationship between homology and
cohomology in Poincar\'e duality spaces forms the basis of many of the results.
To extend surgery to paracompact objects, we are going to need 
a cohomology theory.

If one grants that the homology theory that we 
constructed in section 3 is the right one, 
then the correct cohomology theory is not hard to intuit.
To be loose momentarily, in homology we associate to each compact set
$C$ the group $H_\ast(\coverFA{X-C})$.
If $M-C$ is a manifold with boundary, Lefschetz duality tells us this is
dual to 
$H^\ast_\cmpsup(\coverFA{\Bar{X-C}},\coverFA{\partial C})$, 
where $\Bar{M-C}$ is the closure of $M-C$.
If $C\subseteq D$, we have a map
$H_\ast(\coverFA{M-D})\ \to\ H_\ast(\coverFA{M-C})$, 
so we need a map
$H^\ast_c(\coverFA{\Bar{M-D}},\coverFA{\partial D})\ \to\
H^\ast_c(\coverFA{\Bar{M-C}},\coverFA{\partial C})$.
A candidate for this map is
\[
H^\ast_\cmpsup(\coverFA{\Bar{M-D}},\coverFA{\partial D})\ \RA{\text{tr}}\ 
H^\ast_\cmpsup(\dottedBar{\Bar{M-C}},\dottedBar{\partial D})\ \LA{ex}\ 
H^\ast_\cmpsup(\coverFA{\Bar{M-C}},\coverFA{\Bar{D - C}})\ \RA{\text{inc}}\ 
H^\ast_\cmpsup(\coverFA{\Bar{M-C}},\coverFA{\partial C})
\]
where $\dottedBar{\Bar{M - D}}=\pi^{-1}(M - D)$ 
($\pi\colon \coverFA{M-C}\ \to\ M-C$), 
${\text{inc}}$ is the map induced by inclusion, ${\text{tr}}$ is the trace and
${\text{ex}}$ is an excision map.

The first problem that arises is that ${\text{ex}}$ need not be an isomorphism.
This problem is easily overcome.
We define ${\mathcal O}(X)$ to be the category whose objects are 
open subsets of $X$ whose closure (in $X$) is compact.
If $U$, $V\in{\mathcal O}(X)$, there is a morphism $U\ \to\ V$ \iff\ 
$\Bar U\subseteq V$ or $U=V$.
${\mathcal O}(X)$ will be our diagram scheme.
Note we have a functor ${\mathcal O}(X)\ \to {\mathcal D}(X)$ which sends 
$U\mapsto \Bar U$.
Since $X$ is locally compact, this functor has a cofinal image ($X$ is
homogamous, hence locally compact).

The second problem which arises concerns covering functors.
Since $X - U$, $U\in{\mathcal O}(X)$ is closed, it is hard to get conditions
on $X$ so that $X-U$ has arbitrary covers.
There are two solutions to this problem. 
We can restrict ${\mathcal O}(X)$ (e.g. if $X$ is an homogamous CW complex,
and if we pick sets $U$ so that $X-U$ is a subcomplex, then we  
always have covers),  or we can ignore the problem.
We choose the latter alternative, and when we write $\coverFA{}$ is
a covering functor for $X$, we mean \ $\coverFA{}$\quad is
compatible with $X-U$ for each $U\in{\mathcal O}(X)$.
It is not hard to see that if $X$ is locally $1${--}connected, then
universal covering functors exist despite the fact that the universal
covering functor need not.

Now we could have defined homology and homotopy groups 
using ${\mathcal O}(X)$ instead of ${\mathcal D}(X)$.
Given a covering functor for ${\mathcal O}(X)$ there is an obvious one for
${\mathcal D}(X)$.
It is not hard to show that the homology and homotopy groups for $X$
are the same whether one uses ${\mathcal O}(X)$ or ${\mathcal D}(X)$.

\begin{xDefinition}
$\Delta_\ast(X;A_1,\cdots, A_n\Colon \coverFA{}\ , \Gamma)$,
where $\Gamma$ is a local coefficient system on $X$, denotes the
$\Delta${--}construction applied to $G_{i U}=
H_\ast\bigl((\coverFA{X-U})^i;
\coverFA{A}_1\ \cap\ (\coverFA{X-U})^i, \cdots, 
\coverFA{A}_n\ \cap\ (\coverFA{X-U})^i\Colon i^\ast\Gamma\bigr)$,
where the homology group is the ordinary (singular) $n${--}ad
homology group with coefficients $i^\ast\Gamma$, where
$i^\ast\Gamma$ is the local system induced from $\Gamma$ by
the composite $(\coverFA{X-U})^i\ \RA{\pi}\ X-U\subseteq X$.
\end{xDefinition}

\begin{xDefinition}
$\Delta^\ast(X\Colon \coverFA{}\ ,\Gamma)$ is the $\Delta${--}construction
applied to \[G_{i U}= H_\cmpsup^\ast\bigl((\coverFA{X-U})^i,
\coverFA{\partial U}\ \cap\ (\coverFA{X-U})^i; i^\ast\Gamma\bigr)\]
($\partial U = {\text{frontier\ of\ }} U {\text{\ in\ }}X$).

$\Delta^\ast(X,A\Colon \coverFA{}\ ,\Gamma)$ is the $\Delta${--}construction
applied to
\[G_{i U} = H^\ast_\cmpsup\bigl( (\coverFA{X-U})^i;
\coverFA{\partial U}\ \cap\ (\coverFA{X-U})^i,
\coverFA{A}\ \cap\ (\coverFA{X-U})^i; i^\ast\Gamma\bigr)\ .\]
\end{xDefinition}

\insetitem{{\bf Caution}} $(X,A)$ must be a proper pair (i.e. $A\subseteq X$ is proper)
before $H^\ast_\cmpsup(X,A)$ makes sense. 
A similar remark applies to $n${--}ads.

$\Delta^\ast(X; A_1, \cdots, A_n\Colon \coverFA{}\ ,\Gamma)$ is defined
similarly. 

In our definition we have not defined our maps
$G_{i V}\ \to\ G_{i U}$ if $\Bar U\subseteq V$.
If $\Bar{X - V} = \pi^{-1}(X-V)$, where $\pi\colon (\coverFA{X-U})^i\ \to\
X-U$, then the map is the composite 
\[\begin{aligned}%
H^\ast_\cmpsup&\bigr((\coverFA{X-V})^i, \coverFA{\partial V}\ \cap\ 
(\coverFA{X-V})^i; i^\ast\Gamma\bigr)\ \RA{{\text{tr}}}\ 
H^\ast_\cmpsup(\Bar{X-V},\Bar{\partial V}; \Gamma_1)
\ 
\rightlabeledarrow{\ ex \ }{\cong} 
\ \\
&H^\ast_c\bigr((\coverFA{X-U})^i, (\coverFA{V-U})\ \cap\ (\coverFA{X-U})^i;
\Gamma_2\bigr)\ \RA{\ {\text{inc}}\ }\ 
H^\ast_c\bigr((\coverFA{X-U})^i, \coverFA{\partial U}\ 
\cap\ (\coverFA{X-U})^i; i^\ast\Gamma\bigr)\\\end{aligned}\]
where $\Gamma_1$ and $\Gamma_2$ are the obvious local systems.
A similar definition gives the map in the pair and $n${--}ad cases.

Once again we get long exact sequences modulo 
the usual subspace difficulties.\\
We let $\Delta^\ast(A; X\Colon \coverFA{}\ ,\Gamma)$ denote the subspace
group with a similar notation for the sub{--}$n${--}ad groups.
Again we get a subspace principle. 
Lastly, the cohomology groups are ``independent'' of base points
(compare \fullRef{T.1.3.2}) and are invariant under
proper homotopy equivalence.
The proofs of these results should be easy after section 3, and hence
they are omitted.

One reason for the great power of cohomology is that we have various
products.
The first product we investigate is the cup product.

\BEGIN{T.1.4.1}
There is a natural bilinear pairing, the cup product
\[H^m(X,A;\Gamma_1)\times
\Delta^n(X,B\Colon \coverFA{}\ ,\Gamma_2)\ \to\
\Delta^{m+n}(X;A,B\Colon \coverFA{}\ ,\Gamma_1\otimes\Gamma_2)\ .\]

If $\{ A, B\}$ us a properly excisive pair, the natural map
\[\Delta^n(X,A\ \cup\ B\Colon \coverFA{}\ ,\Gamma_1\otimes\Gamma_2)\ \to\
\Delta^n(X;A,B\Colon \coverFA{}\ ,\Gamma_1\otimes\Gamma_2)\]
is an isomorphism, so we get the ``usual'' cup product.
\end{Theorem}
\begin{proof}
Given $\varphi\in H^m(X,A;\Gamma_1)$, define, for any $U\in{\mathcal O}(X)$,
$\varphi_U\in H^m\bigl((\coverFA{X-U})^i, 
\coverFA{A}\ \cap\ (\coverFA{X-U})^i; i^\ast\Gamma_1\bigr)$ as the image of
$\varphi$ under the composite 
\[H^m(X,A;\Gamma_1)\ \to\ H^m\bigl((X-U), (A-U);\Gamma_1)
\ \RA{\ \pi^\ast\ }\ 
H^m\bigl((\coverFA{X-U})^i, 
\coverFA{A}\ \cap\  (\coverFA{X-U})^i,i^\ast\Gamma_1\bigr)\ .\]
One then checks that if $G_{i U}=
H^n_\cmpsup\bigl((\coverFA{X-U})^i, 
\coverFA{B}\ \cap\  (\coverFA{X-U})^i, \coverFA{\partial U}\ \cap\ 
(\coverFA{X-U})^i, i^\ast\Gamma_2\bigr)$ and if
$H_{i U} = $ the corresponding group for 
$\Delta^{m+n}(X;A,B\Colon \coverFA{}\ ,
\Gamma_1\otimes\Gamma_2)$, then
\[\begin{matrix}%
G_{i U}&\RA{\ \cup\>\varphi_U\ }&H_{i U}\\
\Big\downarrow&&\Big\downarrow\\
G_{i V}&\RA{\ \cup\>\varphi_V\ }&H_{i V}\\
\end{matrix}\]
commutes.
Hence the maps $\cup\>\varphi_U$ give us a map
$\Delta^n(X,B\Colon \coverFA{}\ ,\Gamma_2)\ \to\
\Delta^{m+n}(X; A, B\Colon \coverFA{}\ ,\Gamma_1\otimes\Gamma_2)$.
One easily checks this map gives us a natural bilinear pairing.

Now we have a natural map
$\Delta^\ast(X, A\ \cup\ B)\ \to\ \Delta^\ast(X; A,B )$.
We get a commutative diagram
\[\begin{matrix}%
\cdots\to&\Delta^\ast(X, A\ \cup\ B)&\to& \Delta^\ast(X,A)&\to&
\Delta^\ast(A\ \cup\ B, A;X)&\to\cdots\\
&\Big\downarrow&&\Big\downarrow&&\Big\downarrow\\
\cdots\to&\Delta^\ast(X; A, B)&\to& \Delta^\ast(X,A)&\to&
\Delta^\ast(B, A\ \cap\ B;X)&\to\cdots\\
\end{matrix}\]
where the rows are exact.
$\{A, B\}$ a properly excisive pair implies
$\Delta^\ast(A\ \cup\ B, A)\ \to\ \Delta^\ast(B, A\ \cap\ B)$
is an isomorphism for a set of base points in $A\ \cap\ B$ which is
a set for $A$, $B$, and $A\ \cup\ B$.
The subspace principle now shows the right hand map is an isomorphism.
The middle map is the identity, so the left hand map is an isomorphism.
This establishes the last part of our claim.
\end{proof}

For completeness we give the definition of a properly excisive pair.

\begin{xDefinition}
A pair $\{A, B\}$ of homogamous spaces is said to be properly excisive
with respect to some covering functor \ $\coverFA{}$\ ,
provided
\[\Delta^\ast(A\ \cup\ B\Colon A, B\Colon \coverFA{}\ ,\Gamma)=
\Delta_\ast(A\ \cup\ B\Colon A, B\Colon \coverFA{}\ ,\Gamma)=
H^0_{{\text{end}}}(\coverFA{A\ \cup\ B}; \coverFA{A}, \coverFA{B}; \Gamma)=0\]
for any local system $\Gamma$.
\medskip
The pair is properly excisive if it is properly excisive with respect to all
covering functors compatible with $A\ \cup\ B$.
\end{xDefinition}
\medskip

The other product of great importance is the cap product.
We get two versions of this ( Theorems \shortFullRef{T.1.4.2} 
and \shortFullRef{T.1.4.3}).

\BEGIN{T.1.4.2}
There is a natural bilinear pairing, the cap product
\[\Delta^m(X,A\Colon \coverFA{}\ ,\Gamma_1)\times
H^\locf_{n+m}(X; A, B;\Gamma_2)\ \to\ 
\Delta_{n}(X,B\Colon \coverFA{}\ ,\Gamma_1\otimes\Gamma_2)\]

If $\{A,B\}$ is a properly excisive pair, we can define the ``usual'' 
cap product.
\end{Theorem}
\begin{proof}
Let $C\in H^\locf_{n+m}(X;A,B;\Gamma_2)$.
Define \[C_U\in H^\locf_{n+m}\bigl((\coverFA{X-U})^i;
\coverFA{A}\ \cap\ (\coverFA{X-U})^i, \coverFA{B}\ \cap\ (\coverFA{X-U})^i,
\coverFA{\partial U}\ \cap\ (\coverFA{X-U})^i; i^\ast\Gamma_2\bigr)\]
as the image of $C$ under the composite
\[\begin{matrix}%
H^\locf_{n+m}(X; A, B;\Gamma_2)\\
\downarrow\\
H_{n+m}^\locf(X;A,B,U;\Gamma_2)& 
\leftlabeledarrow{\text{\ ex\ }}{\cong}&
H_{n+m}^\locf(X-U;A-U,B-U, \partial U;\Gamma_2)\\
&&\downlabeledarrow{{\text{tr}}}{}\\
\hbox to 0pt{$H^\locf_{n+m}\bigl((\coverFA{X-U})^i;
\coverFA{A}\ \cap\ (\coverFA{X-U})^i, \coverFA{B}\ \cap\ (\coverFA{X-U})^i,
\coverFA{\partial U}\ \cap\ (\coverFA{X-U})^i; i^\ast\Gamma_2\bigr)$ .\hss}
\\\end{matrix}\]

One can check that $\cap\ C_U$ satisfies the necessary commutativity 
relations to define a map
\[\Delta^m(X,A\Colon \coverFA{}\ ,\Gamma_1)\ \to\ 
\Delta_n(X,B\Colon \coverFA{}\ ,\Gamma_1\otimes\Gamma_2)\ .\]

If $\{A, B\}$  is properly excisive, $H^\ast_\cmpsup(A\ \cup\ B;A, B)=0$
follows from $\Delta^\ast=0$.
Universal coefficients shows $H^\locf_\ast(A\ \cup\ B;A,B)=0$, so the
standard exact sequence argument shows
$H^\locf_\ast(X;A,B;\Gamma_2)\cong 
H^\locf_\ast(X,A\ \cup\ B;\Gamma_2)$.
\end{proof}

\BEGIN{T.1.4.3}
There is a natural bilinear pairing, the cap product
\[H^m(X,A;\Gamma_1)\times
\Delta_{n+m}(X;A,B\Colon \coverFA{}\ ,\Gamma_2)\ \to\ 
\Delta_n(X,B\Colon \coverFA{}\ ,\Gamma_1\otimes\Gamma_2)\ .\]

If $\{A,B\}$ is a properly excisive pair, we can define the
``usual'' cap product.
\end{Theorem}
\begin{proof}
Given $\varphi\in H^m(X,A;\Gamma_1)$, define, for any $U\in{\mathcal O}(X)$,
$\varphi_U\in H^m\bigl((\coverFA{X-U})^i, 
\coverFA{A}\ \cap\ (\coverFA{X-U})^i; i^\ast\Gamma_1\bigr)$ as the image of
$\varphi$ under the composite 
\[H^m(X,A;\Gamma_1)\ \to\ H^m\bigl((X-U), (A-U);\Gamma_1)
\ \RA{\ \pi^\ast\ }\ 
H^m\bigl((\coverFA{X-U})^i, 
\coverFA{A}\ \cap\  (\coverFA{X-U})^i,i^\ast\Gamma_1\bigr)\ .\]
One checks again that the necessary diagrams commute.
The statement about $\{A,B\}$ follows from the $5${--}lemma and
the subspace principle.
\end{proof}

We will also need a version of the slant product for our theory.
To get this we need to define a group for the product of two ``ads''.
As usual we apply the $\Delta${--}construction to a particular situation.
Pick a set of base points for $X$ and a set for $Y$.
Our indexing set is the Cartesian product of these two sets.
Our diagram is ${\mathcal O}(X)\times{\mathcal O}(Y) = 
\{ U\times V\subseteq X\times Y\ \vert\ U\in{\mathcal O}(X),
\ V\in{\mathcal O}(Y)\}$.

\[\begin{aligned}%
G^{i\times j}_{U\times V}=
H_\ast\bigl(&(\coverFA{X-U})^i\times(\coverFC{Y-V\ })^j;
\bigl(\coverFA{A}_1\ \cap\ (\coverFA{X-U})^i\bigr)\times (\coverFC{Y-V\ })^j,
\cdots, \\&
\bigl(\coverFA{A}_n\ \cap\ (\coverFA{X-U})^i\bigr)\times (\coverFC{Y-V\ })^j,
(\coverFA{X-U})^i\times\bigl(\coverFA{B_1}\ \cap\ 
\coverFC{Y-V})^j,\cdots,\\&
(\coverFA{X-U})^i\times\bigl(\coverFA{B_m}\ \cap\ \coverFC{Y-V})^j;
i^\ast\Gamma\bigr)\quad .\\
\\\end{aligned}\]
The resulting groups will be denoted
$\Delta_\ast\bigl( (X;A_1,\cdots, A_n)\times (Y;B_1,\cdots, B_m)\Colon 
\coverFA{}\ , \coverFC{\hskip10pt }\ ,\Gamma)$ ($\Gamma$ is 
some local system on $X\times Y$).

\BEGIN{T.1.4.4}
There is a natural bilinear pairing, the slant product
\[\begin{matrix}%
\slantp\colon H^m(\coverFC{Y\ }; \coverFC{B_1},\cdots, \coverFC{B_m}; 
\Gamma_1)\times
\Delta_{n+m}\bigl( (X;A_1,\cdots, A_n)\times (Y;B_1,\cdots, B_m)\Colon 
\coverFA{}\ , \coverFC{\hskip10pt }\ ,\Gamma_2)\\
\downarrow\\
\Delta_n(X;A_1,\cdots, A_n,\coverFA{}\ ,\Gamma_1\otimes\Gamma_2)
\\
\end{matrix}\]
\end{Theorem}
\begin{proof}
For $\varphi\in H^m(\coverFC{Y\ }; \coverFC{B_1},\cdots, \coverFC{B_m}; 
\Gamma_1)$,
define $\varphi_V$ by analogy with the definition in 
\fullRef{T.1.4.1}.
These give us the necessary maps.
\end{proof}

\BEGIN{C.1.4.4.1}
If $d\colon X\ \to\ X\times X$ is the diagonal, and \\if 
$C\in \Delta_{n+m}(X; A, B\Colon \coverFA{}\ ,\Gamma_1)$, and if
$\varphi\in H^m(X,A;\Gamma_2)$, then
\[\varphi\ \cap\ C\ =\ \varphi\slantp d_\ast C\ .\quad\mathqed\]
\end{Corollary}

Using our slant product, we can define the cap product of 
\fullRef{T.1.4.3} on the chain level.
There are two basic chain groups we would like to use,
For an homogamous CW complex we would like to use the cellular chains,
and when $X$ is a paracompact manifold with a locally finite handlebody
decomposition, we want to use the chains based on the handles.
We do the former case and leave the reader to check the theory still
holds in the latter.

If $X$ is an homogamous CW complex,
 we define
\[P_\ast(X;A,B\Colon \coverFA{}\ ,\Gamma)=
\Delta_\ast(X^\ast; X^{\ast-1}, A^\ast, B^\ast\Colon \coverFA{}\ , \Gamma)\]
(where $A^\ast=A\ \cap X^\ast$, etc.) for $\ast\geq 2$.
If $\ast=0$ or $1$, we must use subspace groups
\[\Delta_\ast\bigl((X^\ast;X^{\ast-1},A^\ast,B^\ast);X\Colon 
\coverFA{}\ ,\Gamma)\ .\]
$A$ and $B$ are subcomplexes.
Similarly define
\[\begin{aligned}%
P^\ast(A;A,B\Colon \coverFA{}\ ,\Gamma)&= 
\Delta^\ast(X^\ast; X^{\ast-1}, A^\ast, B^\ast\Colon \coverFA{},\Gamma)\\
C^\locf_\ast(X\Colon A,B)&=
H^\locf_\ast(X^\ast; X^{\ast-1}, A^\ast, B^\ast)\\
C^\ast(X\Colon A,B)&= H^\ast(X^\ast; X^{\ast-1}, A^\ast, B^\ast)\quad .\\
\end{aligned}\]
The triple $(X^\ast, X^{\ast-1}, X^{\ast-2})$ gives us a boundary map
$P_\ast\ \to\ P_{\ast-1}$, $P^\ast\ \to\ P^{\ast+1}$, etc.
This boundary map makes the above objects into chain complexes 
($\partial\partial=0$), and by \fullRef{C.1.2.2.2},
the homology of these complexes is just what one expects.

A diagonal approximation 
\[h_\ast\colon P_\ast(X; A, B\Colon \coverFA{}\ ,\Gamma)\ \to\ 
\Delta_\ast\Bigl(\bigl( (X,A)\times (X,B)\bigr)^\ast\Colon 
\coverFA{}\ , \Gamma\Bigr)\]
is a cellular approximation to $d\colon X\ \to\ X\times X$ with a homotopy
between $d$ and the cellular map,
$H\colon X\times I\ \to\ X\times X$, such that $\pi_1\ \circ H$ and
$\pi_2\ \circ H$ are proper.
$\bigl( (X,A)\times (X,B)\bigr)^\ast$ is just
$\displaystyle\mathop{\cup}_k (X,A)^k\times (X,B)^{\ast-k}$.
Any two such diagonal approximations are cellularly homotopic so that
the homotopy composed with the projections is proper.

\BEGIN{T.1.4.5}
Given any diagonal approximation $h$, there is a bilinear pairing
\[B_h\colon
C^m(X,A;\Gamma_1)\times 
P_{n+m}(X;A,B\Colon \coverFA{}\ ,\Gamma_2)\ \to\ 
P_n(X,B\Colon \coverFA{}\ ,\Gamma_1\otimes\Gamma_2)\ .\]

If $f\in C^m(X,A;\Gamma_1)$ and 
$c\in P_{n+m}(X;A,B\Colon \coverFA{}\ ,\Gamma_2)$, then
\[\partial B_h(f,c)= (-1)^{n} B_h(\delta f, c) + B_h(f, \partial c)\ .\]
Hence we get an induced pairing on the homology level.
Any two $B_h(f,\quad)$ are chain homotopic, so the pairing on homology
does not depend on the diagonal approximation.
This pairing on homology is the cap product of 
\fullRef{T.1.4.3}.
\end{Theorem}
\begin{proof}
Consider the element $h_\ast(c)\in
\Delta_{n+m}\Bigl(\bigl((X,A)\times(X,B)\bigr)^{n+m},{\text{etc.}}\Bigr)$.
The group 
\[\Delta_{n+m}\bigl((X^n;X^{n-1},A^n)\times (X^m; X^{m-1}, B^m)\Colon 
\coverFA{}\ ,\Gamma_2)\] 
lies as a natural summand of this first group.
Let $p^n_m$ be this projection.
Then $B_h(f,c) = f\Big\vert_{p^n_m\bigl(h_\ast(c)\bigr)}$.
The rest of the proof involves checking this definition has all the asserted 
properties.
\end{proof}

We also want to define the cap product of 
\fullRef{T.1.4.2} on the chain level.
Unfortunately, there is no slant product of the needed type, 
so we must use brute force.

\BEGIN{T.1.4.6}
Given any diagonal approximation $h$, there is a bilinear pairing
\[B_h\colon
P^m(X,A\Colon \coverFC{}\ ,\Gamma_1)\times
C^\locf_{n+m}(X;A,B;\Gamma_2)\ \to\ 
P_n(X,B\Colon \coverFC{}\ ,\Gamma_1\otimes\Gamma_2)\ .\]

If $f\in P^m(X,A\Colon \coverFC{}\ ,\Gamma_1)$ and 
$c\in C^\locf_{n+m}(X;A,B;\Gamma_2)$, then
\[\partial B_h(f,c)=(-1)^{n}B_h(\delta f,c) + B_h(f, \partial c)\quad .\]
Hence we get an induced pairing (independent of $h$) on the homology
level.
This pairing is the cap product of \fullRef{T.1.4.2}.
\end{Theorem}
\begin{proof}
Let $c\in C^\locf_{n+m}(X;A,B;\Gamma_2)$.
Define \[c_U\in C^\locf_{n+m}\bigl( (\coverFC{X-U})^i;
\coverFC{A\ }\ \cap\ (\coverFC{X-U})^i, \coverFC{B\ }\ \cap\ (\coverFC{X-U})^i,
\coverFC{\partial U\ }\ \cap\ (\coverFC{X-U})^i; i^\ast\Gamma_2\bigr)\]
by excision and trace as in the proof of \fullRef{T.1.4.2}.
We define $B_h(\quad, c)$ from the maps\hfill\\\noindent
$H^m_\cmpsup\bigl( \coverFC{(X-U)^i\ \cap\ X^m};
\coverFC{(X-U)^i\ \cap\ X^{m-1}}, {\text{etc.}}\bigr)\ \RA{\ \slantp b_U\ }\ 
$\hfil\par\noindent$H_n\bigl( \coverFC{(X-U)^i\ \cap\ X^n};
\coverFC{(X-U)^i\ \cap\ X^{n-1}}, {\text{etc.}} \bigr)$
where $\slantp$ is the slant product and $b_U$ is the homology class
which is the image of $c_U$ under the following composite.
\[\begin{aligned}%
H^\locf_{n+m}\bigl(& \coverFC{(X-U)^i\ \cap\ X^m\ }; {\text{etc.}} \bigr)
\RA{\ h_\ast\ }
H^\locf_{n+m}\Bigl( \bigl( (\coverFC{X-U})^i\times(\coverFC{X-U})^i\bigr)
\cap
\coverFC{(X\times X)^{n+m}}; {\text{etc.}} \Bigr)\\
&\RA{\ p^n_m\ }\ 
H^\locf_{n+m}\Bigl( \bigl( \coverFC{(X-U)^i\ \cap\ X^m\ }\bigr) \times 
\bigl( \coverFC{(X-U)^i\ \cap\ X^n\ }\bigr); {\text{etc.}}\Bigr)
\\\end{aligned}\]
(superscript $i$ denotes a component containing $\hat x_i$,
and superscripts $n$, $m$, and $n+m$ denote skeletons.)

Note in passing that $h_\ast({\text{tr}}\ b_U) \neq {\text{tr}}(h_\ast b_U)$,
which is why we are unable to define a general slant product like
\fullRef{T.1.4.4} to cover this case.

The rest of the proof involves verifying diagrams commute and
verifying our equation.
\end{proof}

\medskip
Lastly we prove the Browder lemma, which will be essential in our
study of Poincar\'e duality.
\BEGIN{T.1.4.7}
Let $(X,A)$ be a proper pair, and let $c\in H^\locf_{n}(X,A;\Gamma_2)$.
Then
\[\begin{aligned}\begin{matrix}%
\Delta^{\ast-1}(A;X\Colon \coverFA{}\ ,\Gamma_1)&\to&
\Delta^\ast(X,A\Colon \coverFA{}\ ,\Gamma_1)&\to\\
\downlabeledarrow[\Big]{\ \cap\ (-1)^{n-1}\partial c}{}
&&\downlabeledarrow[\Big]{\cap\ c}{}\\\noalign{\vskip4pt}
\Delta_{n-\ast}(A;X\Colon \coverFA{}\ ,\Gamma_1\otimes\Gamma_2)&\to&
\Delta_{n-\ast}(X\Colon \coverFA{}\ ,\Gamma_1\otimes\Gamma_2)&\to&
\\\end{matrix}\\\noalign{\vskip20pt}
\begin{matrix}%
\Delta^\ast(X\Colon \coverFA{}\ ,\Gamma_1)&\to&
\Delta^\ast(A;X\Colon \coverFA{}\ ,\Gamma_1)&\\
\downlabeledarrow[\Big]{\cap\ c}{}&&
\downlabeledarrow[\Big]{\cap\ \partial c}{}\\\noalign{\vskip4pt}
\Delta_{n-\ast}(X,A\Colon \coverFA{}\ ,\Gamma_1\otimes\Gamma_2)&\to&
\Delta_{n-1-\ast}(A;X\Colon \coverFA{}\ ,\Gamma_1\otimes\Gamma_2)
&
\\\end{matrix}\\\end{aligned}
\]
commutes.
\end{Theorem}

\begin{proof}
The usual Browder lemma (see section 1) says that the corresponding
diagram commutes for ordinary homology and cohomology with compact
supports.
Commutativity is then trivial for the above diagram. 
(While we have not defined a cap product for subspace groups,
the reader should have no difficulty writing down the necessary maps.)
\end{proof}

\bigskip
\section{Chain complexes and simple homotopy type}
\newHead{I.5}
In our $\Delta${--}construction as applied to the homology or homotopy
functors, we still have some structure that we have not utilized. 

As an example of this extra structure, let us consider $\epsilon(X\Colon \pi_1)$.
This is an inverse limit\\
$\displaystyle\lim_{\to} \mu\bigl( \pi_1(X-C, x_i)\Bigr)$.
Now many of the $\pi_1(X-C, x_i)$ are isomorphic.
(Unfortunately this isomorphism are not unique but depends on paths
joining $x_i$ to $x_j$.)
Our $\epsilon${--}construction makes no use of this fact.
In order to be able to make effective use of this extra structure, we
need a way to choose the above isomorphisms.

We will do this through the concept of a tree.
A \emph{tree for an homogamous space} $X$ will be a $1${--}dimensional, 
locally finite, simplicial complex, $T$, such that
\begin{enumerate}
\item[1)] $\Delta(T\Colon \pi_k)=0$ for $k>0$
\item[2)] If $T^\prime\subseteq T$ is a subcomplex of $T$, $T^\prime$
has the proper homotopy type of $T$ \iff\ $T=T^\prime$.
\end{enumerate}

\noindent
(This last condition is to insure that
$\displaystyle\begin{matrix}&&\circ&&\circ\\
&&\vert&&\vert\\
\circ&\vrule width 10pt height 4.4pt depth -3.9pt&
\circ&\vrule width 10pt height 4.4pt depth -3.9pt&
\circ&\vrule width 10pt height 4.4pt depth -3.9pt&
\circ&\cdots
\\
\end{matrix}$
is not a tree for $\R^2$ but rather
$\circ\ \vrule width 10pt height 4.4pt depth -3.9pt\ 
\circ\ \vrule width 10pt height 4.4pt depth -3.9pt\ 
\circ\ \cdots$
is.)
We also require a map $f\colon T\ \to\ X$ 
which is properly $1/2${--}connected.

Two trees $(T,f)$ and $(S,g)$ are equivalent provided there is a
proper homotopy equivalence $h\colon T\ \to\ S$
with $g\ \circ h$ properly homotopic to $f$.

A space $X$ is said to have a tree provided $X$ is homogamous and
there is a tree for $X$.
Any locally path connected homogamous space has a tree.
To see this,  let $\{ p\}$ be a set of base points for our space $X$.
Let $\{C_i\}$ be a cofinal collection of compact subsets of $X$.
We can assume $X$ is path connected since we can do each path
component separately.
We may assume $\{ p\} \ \cap\ C_0\neq\emptyset$.
Pick a point $p_0\in\{ p\}\ \cap\ C_0$.
Look at the components of $X-C_0$ with a point of $\{p\}$ in them.
As we showed in the proof of \fullRef{P.1.1.2},
there are only finitely many such components of $X- C_0$.
The components whose closure is not compact are called essential 
components.
We may assume $\{p\}\ \cap\ ({\text{each\ essential\ component\ of\ }}
X-C)\ \cap\ C_1 \notin C_1\neq \emptyset$
since this is true for some compact set.
Let $p_1^{\alpha_1}$, $p_1^{\alpha_2}$, \dots , $p_1^{\alpha_n}$
be a subset of $\{p\}\ \cap\ C_1$, one for each essential component
of $X-C_0$.
Join $p_1^{\alpha_i}$ to $p_0$ by a path $\lambda_{1,i}$.
Now look at the essential components of $X-C_1$.
Pick $p_2^{\alpha_1}$, \dots , $p_2^{\alpha_m}$
(which we may assume are in $C_2$), one for each essential component
of $X-C_1$.
Each $p_2^{\alpha_i}$ lies in an essential component of $X-C_0$, so
pick paths $\lambda_{2,i}$ which join $p_2^{\alpha_i}$ to the
appropriate element in $\{p_1^\alpha\}$.
These paths should lie in $X-C_0$.
(Since $X$ is locally path connected, the components of $X-C_0$ are
path connected.)
Continue in this fashion to get $\{ p^\alpha_j\}$, one for each 
essential component of $X-C_{j-1}$.
$\{p^\alpha_j\}$ may be assumed to lie in $X-C_j$.
We can also get paths $\lambda_{j,\alpha_i}$ which join
$p^{\alpha_i}_j$ to the appropriate $p^\alpha_{j-1}$ and which
lie in $X-C_{j-1}$.

Now $T$ has $\{ p^\alpha_j \}$ for vertices and
$( p^{\alpha_i}_j, p^{\alpha_\ell}_k)$ is a $1${--}simplex \iff\ $k=j-1$
and $\lambda_{j,\alpha_i}$ joins $p^{\alpha_i}_j$ to $p^{\alpha_\ell}_k$.
The map $f\colon T\ \to\ X$ is the obvious one.

We claim $H_1(T;\Z)=0$, and in fact, if $H_1$ is computed from the
simplicial chains then there are no $1${--}cycles.
This is fairly clear , so it will be left to the reader.
Now any locally finite $1${--}complex with $H_1(T;\Z)=0$ satisfies
$\Delta(T\Colon \pi_k)=0$ for $k>0$.
One shows $f$ is properly $1/2${--}connected by showing that
$Z^0_{\text{end}}(X)\ \to\ Z^0_{\text{end}}(T)$ is an isomorphism
($Z^0_{\text{end}}$ are the $0${--}cycles in $S^0_{\text{end}}$).
But this follows from our construction.
Lastly suppose $T^\prime\subseteq T$ is a connected subcomplex,
and suppose $p\in T- T^\prime$.
Now by definition $p$ is in an essential component of $X-C_i$ for
all $i\leq n$ for some $n$.
Since each essential component of $X-C_i$ has infinitely many base points
in it, let $\{q\}$ be the set of base points in the component of $X-C_n$
containing $p$.
Then $\{q\}\subseteq T-T^\prime$, as is easily seen.
Hence $H^0_{\text{end}}(T)\ \to\ H^0_{\text{end}}(T^\prime)$ has a kernel,
and so $T^\prime\subseteq T$ is not a proper homotopy equivalence.
Hence $X$ has a tree.

From now on in this section we restrict ourselves to the category of
homogamous CW complexes.
We will use \hCWx\ complex to denote objects in this category.

Given $X$, an \hCWx\  complex, we have the category ${\mathcal C}(X)$ whose
objects are all sets $A\subseteq X$ such that 
\begin{enumerate}
\item[1)] $A$ is a subcomplex
\item[2)] $A$ is connected
\item[3)] There exists an element of ${\mathcal O}(X)$, $U$, such that
$A$ is an essential component of $X-U$.
\end{enumerate}

\noindent The morphisms are inclusions.

Now given a tree $(T,f)$ for $X$, we get a functor
${\mathcal C}(X)\ \RA{{\mathcal C}(f)}\ {\mathcal C}(T)$ ($f$ is always assumed to be
cellular).

\begin{xDefinition}
A \emph{lift of} ${\mathcal C}(f)$ is a covariant functor 
$F\colon {\mathcal C}(T)\ \to\ {\mathcal C}(X)$ such that
${\mathcal C}(f)\ \circ F$ is the identity and such that the image of $F$ 
is cofinal.
The set of all such lifts is a diagram scheme by defining $F\leq G$
\iff\ $F(A)\subseteq G(A)$ for all $A\in{\mathcal C}(T)$.
We denote this diagram scheme by ${\mathcal L}(f)$.
\end{xDefinition}

\begin{xDefinition}
A \emph{tree of rings} is a covariant functor 
$R\colon{\mathcal C}(T)\ \to\ {\mathcal R}$, where ${\mathcal R}$ is the category of
all rings (rings have units and all ring homomorphisms preserve units).
A \emph{tree of modules} over $R$ is a collection of modules $M_A$,
$A\in{\mathcal C}(T)$, where $M_A$ is a unitary $R_A${--}module.
A tree of right (left) $R${--}modules requires each $M_A$ to be a right
(left) $R_A${--}module.
If $A\subseteq B$ in ${\mathcal C}(T)$, there is a unique map
$p_{AB}\colon M_A\ \to\ M_B$, which is an
$R(A\subseteq B)${--}linear map; i.e. if $f\colon R_A\ \to\ R_B$ is the
ring homomorphism associated to $A\subseteq B$ by $R$,
\[p_{AB}(a\cdot\alpha + b\cdot\beta)=
p_{AB}(a)\cdot f(\alpha) + p_{AB}(b)\cdot f(\beta)\]
for $\alpha$, $\beta\in R_A$; $a$, $b\in M_A$.
\end{xDefinition}

An $R${--}module homomorphism $f\colon M\ \to\ M^\prime$ 
is a set of maps $f_A\colon M_A\ \to\ M^\prime_A$
for each $A\in{\mathcal C}(T)$ such that
\begin{enumerate}
\item[1)] $f_A$ is an $R_A${--}module homomorphism
\vskip 4pt
\item[2)] For $A\subseteq B$,
\lower 12pt\vtop {\hsize=1.75in$\displaystyle\begin{matrix}%
M_A&\RA{\ f_A\ }&M^\prime_A\\
\downlabeledarrow[\big]{p_{AB}}{}&&\downlabeledarrow[\big]{p^\prime_{AB}}{}\\
M_B&\RA{\ f_B\ }&M^\prime_B\\\noalign{\vskip6pt}\end{matrix}$\vss}
commutes, where the vertical maps come from the tree structures
on $M$ and $M^\prime$.
\end{enumerate}
\medskip

\begin{xExample}
Given an \hCWx\  complex $X$ with a tree $(T,f)$ and given $F\in{\mathcal L}(f)$,
we get a tree of rings from $R_A=\Z\bigl[\pi_1\bigl( F(A),f(p)\bigr)\bigr]$
where if $A\neq T$, $p$ is the vertex of $\partial A$, the set
theoretic frontier of $A$.
If $A=T$ pick a vertex for a base point and use it.
This will be the tree of rings we will consider for our geometry,
and we will denote it by $\Z\pi_1$.

The tree of $\Z\pi_1${--}modules we will consider will be 
various chain modules. 
The basic idea is given by
$M_A= H_i\bigl( \coverFA{F(A)^i\ },  \coverFA{F(A)^{i-1}\  }, f(p)\bigr)$,
where \ $\coverFA{}$\quad denotes the universal cover of $F(A)$,
and $\coverFA{F(A)^i}$ is $\pi^{-1}$ of the $i${--}skeleton of $F(A)$
in $\coverFA{F(A)}$. ($\pi\colon \coverFA{F(A)}\ \to\ F(A)$).

Now given an $R${--}module $M$, we can form $\Delta(M)$ by applying
the $\Delta${--}construction with index set the vertices of $T$, and
with diagram scheme ${\mathcal O}(X)$.
Given $U\in{\mathcal O}(X)$, there are finitely many $A\in{\mathcal C}(T)$ for
which $A\ \cap\ \Bar U = {\text{a\ vertex}}$.
Set \[G_{p U}=\begin{cases}M_A& \text{if }p\in A\\0&\text{otherwise}\end{cases}\]
for some $A$ such that $A\ \cap\ \Bar U = {\text{a\ vertex}}$.
An $R${--}module homomorphism $f\colon M\ \to\ M^\prime$
clearly induces a map $\Delta(f)\colon \Delta(M)\ \to\ \Delta(M^\prime)$.
An $R${--}module homomorphism, $f$, which induces an isomorphism 
$\Delta(f)$ is said to be a \emph{strong equivalence} and 
the two modules
are said to be \emph{strongly equivalent}.
Note that this relation on $R${--}modules seems not to be symmetric.
Nevertheless we can define two $R${--}modules to be \emph{equivalent}
\iff\ there is a (finite) sequence of $R${--}modules $M=M_0$, $M_1$, \dots,
$M_n=M^\prime$ such that either $M_{i}$ is strongly equivalent to 
$M_{i+1}$ or $M_{i+1}$ is strongly equivalent to $M_{i}$.

We tend only to be really interested in the equivalence class of $M$ 
(indeed, we are often interested merely in $\Delta(M)$ ).
The relation of equivalence is not however very nice.
We would like $M$ to be equivalent to $M^\prime$ \iff\ there were
``maps'' $f\colon M\ \to\ M^\prime$ and $g\colon M^\prime\ \to\ M$
whose composites were the identity.
To do this properly we need a short digression.
\end{xExample}

\begin{xDefinition}
A functor $F$ which assigns to each $A\in{\mathcal C}(T)$ 
a cofinal subcomplex of $A$, $F(A)$, such that $F(A)\subseteq F(B)$ 
whenever $A\subseteq B$ and such that $F(T)=T$ will be called a
\emph{shift functor}.
${\mathcal S}(T)$ will denote the set of all shift functors on $T$.
${\mathcal S}(T)$ is partially ordered via $F\geq G$ \iff\ $F(A)\subseteq G(A)$
for all $A\in{\mathcal C}(T)$.
Define $(F\cap G)(A)=F(A)\ \cap\ G(A)$, and one checks it is 
a shift functor.
$F\cap G\geq F$ and $F\cap G\geq G$.
\end{xDefinition}

Given a tree of $R${--}modules and a shift functor $F$,
we get a tree of $R${--}modules, $M_F$, in a natural way; i.e.
$F$ is going to induce a functor from the category of $R${--}modules
to itself.
$M_F$ is defined as follows.
Let $A\in{\mathcal C}(T)$.
Then $F(A)=\displaystyle\mathop{\cup}_{i=1}^n A_i$, 
with $A_i\in{\mathcal C}(T)$.
$\displaystyle(M_F)_A=\mathop{\oplus}_{i=1}^n\ M_{A_i}\otimes R_A$,
where the tensor product is formed using the homomorphisms 
$R_{A_i}\ \to\ R_A$.
Note that there is an $R_A${--}module map
$(M_F)_A\ \to\ M_A$.
$(p_F)_{A B}\colon \displaystyle\mathop{\oplus}_{i=1}^n\ M_{A_i}\otimes R_A
\ \to\ \mathop{\oplus}_{i=1}^n\ M_{B_i}\otimes R_B$
is defined as follows.
Since $A\subseteq B$, $F(A)\subseteq F(B)$, so each $A_i$ is
contained in a unique $B_j$.
Let $p_{i j}$ be $p_{A_i B_j}$ if $A_i\subseteq B_j$ and $0$ otherwise.
$f_{i j}$ is the map $R_{A_i}\ \to\ R_{B_j}$ if $A_i\subseteq B_j$ 
and $0$ otherwise.
$g$ is the map $R_A\ \to\ R_B$.
Then $(p_F)_{A B}=\displaystyle\mathop{\oplus}_{i=1}^n\ \mathop{\oplus}_{j=1}^m\ 
p_{i j}\otimes f_{i j}\otimes g$.
Notice that
\lower 16pt\hbox{$\displaystyle\begin{matrix}%
(M_F)_A&\RA{\hskip 40pt}& M_A\\
\downlabeledarrow[\Big]{(p_F)_{A B}}{}&&\downlabeledarrow[\Big]{p_{A B}}{}\\
(M_F)_B&\RA{\hskip 40pt}& M_B\\\end{matrix}$\hskip 20pt}
commutes.
If $f\colon M\ \to\ M^\prime$ is a map, 
$\displaystyle(f_F)_A=\mathop{\oplus}_{i=1}^n f_{A_i}\otimes g_{A_i A}$, 
where $g_{A_i A}\colon R_{A_i}\ \to\ R_A$, defines a map
$f_F\colon M_F\ \to\ M^\prime_F$ so that
\lower 12pt\hbox{%
$\displaystyle\begin{matrix}%
M_F&\RA{f_F}&M^\prime_F\\
\big\downarrow&&\big\downarrow\\
M&\RA{\ f\ }&M^\prime\\\end{matrix}$\hskip10pt} commutes.
For the natural map of $M_F$ into $M$ we write $M_F\subseteq M$.
If $G\geq F$ there is a natural map $M_G\ \to\ M_F$ induced by
the inclusion of each component of $G(A)$ in $F(A)$.
\medskip

\BEGIN{L.1.5.1}
$M_F\subseteq M$ is a strong equivalence.
\end{Lemma}
\begin{proof}
We must show $\Delta(M_F)\ \to\ \Delta(M)$ is an isomorphism.
Suppose $B\in{\mathcal C}(T)$ and $B\subseteq F(A)$.
Then
\lower 16pt \hbox{$\begin{matrix}%
(M_F)_B&\RA{\hskip10pt}& M_B \\
\big\downarrow&&\big\downarrow\\
(M_F)_A&\RA{\hskip10pt}& M_A\\\end{matrix}$}
commutes and there is a map
$h\colon M_B\ \to\ (M_F)_A$
so that the resulting triangles commute.
But then clearly $\Delta(M_F)\cong\Delta(M)$.
\end{proof}
\medskip
As motivation for our next definition we prove

\BEGIN{L.1.5.2}
Let $f\colon M\ \to\ N$ be a strong equivalence.
Then there is a shift functor $F$ and a map $N_F\ \to\ M$ such that
\lower 12pt\hbox{$\begin{matrix}&&M\\
&\nearrow&\downarrow\\
N_F&\subseteq&N\\\end{matrix}$}
commutes.
\end{Lemma}
\begin{proof}
By \fullRef{T.1.2.4} applied to kernel and cokernel,
$f$ is a strong equivalence \iff\ for any $A\in{\mathcal C}(T)$ there is a
$U\in{\mathcal O}(T)$ such that for any $B\in{\mathcal C}(T)$ with
$B\subseteq A-U$
\[\begin{matrix}%
M_B&\RA{\ f_B\ }&N_B\\
\downlabeledarrow[\Big]{p^M_{A B}}{}&&\downlabeledarrow[\Big]{p^N_{A B}}{}\\
M_A&\RA{\ f_A\ }&N_A\\\end{matrix}\]
satisfies
\begin{enumerate}
\item[1)] ${\kerx}\ f_B\ \subseteq\ {\kerx}\ p^M_{A B}$ and
\item[2)] ${\Imx}\ p^N_{A B}\ \subseteq\ {\Imx}\ f_A$.
\end{enumerate}
\medskip
For each $A\in{\mathcal C}(T)$, pick such an element in ${\mathcal O}(T)$, $U_A$.
Now let $F(A)= A-\displaystyle\mathop{\cup}_{A\subseteq D} U_D$.
$F$ is easily seen to be a shift functor, and for any $B\in{\mathcal C}(T)$
with $B\subseteq F(A)$, 1) and 2) hold.

Now look at
\lower 25pt\hbox{$\begin{matrix}%
M_{A_2}&\RA{\hskip10pt}& N_{A_2}\\
\big\downarrow&&\downlabeledarrow[\big]{p}{}\\
M_{A_1}&\RA{f_{A_1}}& N_{A_1}\\
\downlabeledarrow[\big]{q}{}&&\big\downarrow\\
M_{A}&\RA{\hskip10pt}& N_{A}\\\end{matrix}$}
where $A_1\subset F(A)$, $A_2\subset F(A_1)$.
Then there exists a map $h\colon N_{A_2}\ \to\ M_A$
defined by $h(x)= q(f_{A_1})^{-1}p(x)$ for all $x\in N_{A_2}$.
By properties 1) and 2), $h$ is well{--}defined, and if 
$g\colon R_{A_2}\ \to\ R_A$ is the homomorphism given by the tree,
$h$ is easily seen to be $g${--}linear.

Define a shift functor $F\ \circ\ G$ by $F\ \circ\ G(A)=
\displaystyle\mathop{\cup}_{i=1}^n F(A_i)$,
where $G(A)=\displaystyle\mathop{\cup}_{i=1}^n A_i$.
Then one checks that the $h$ defined above yields a map 
$N_{F\circ G}\ \to\ M$.
\end{proof}

\bigskip
\begin{xDefinition}
A $T${--}map $f\colon M\ \to\ N$ is a map $M_F\ \to\ N$,
where $F\in{\mathcal S}(T)$.
$M_F\ \to\ N$ induces a natural map $M_G\ \to\ N$ for all $G\geq F$.
We say $f$ is defined on $M_G$ for all $G\geq F$.
Two $T${--}maps $f$, $g\colon M\ \to\ N$ are equal provided that,
for some $F\in{\mathcal S}(T)$ such that $f$ and $g$ are defined on
$M_F$, the two maps $M_F\ \to\ N$ are equal.
\end{xDefinition}

\begin{xRemarks}
If $f$ is defined on $M_F$, and if $g$ is defined on $M_G$, $f$ and
$g$ are both defined on $M_{F\cap G}$.
With this remark it is easy to see equality of $T${--}maps is an
equivalence relation.
It is also easy to see how to add or subtract two $T${--}maps,
and it is easy to check that if $f_1=f_2$ and $g_1=g_2$, then
$f_1\pm g_1 = f_2\pm g_2$.

Hence, if ${\Homx}_T(M,N)$ is the set of equivalence classes of
$T${--}maps from $M$ to $N$, ${\Homx}_T(M,N)$ has the
structure of an abelian group.
An equivalence class of $T${--}maps is called a \emph{map{--}germ}.

We can compose two $T${--}maps $f\colon M\ \to\ N$ and
$g\colon N\ \to\ P$ as follows.
$g$ is defined on $N_G$ and $f$ is defined on $M_F$.
Hence $f\colon N_F\ \to\ P$ is an actual map, and we define the
$T${--}map $g\ \circ\ f$ to be the map
$g\ \circ\ f_G\colon (M_F)_G\ \to\ N_G\ \to\ P$.
Note $(M_F)_G=M_{F\circ G}$.
One can check that the map{--}germ $g\ \circ\ f$ is well{--}defined.
\end{xRemarks}
\medskip
Hence \fullRef{L.1.5.2} becomes

\BEGIN{L.1.5.3}
$M$ and $N$ are equivalent \iff\ they are $T${--}equivalent.
\end{Lemma}
\begin{proof}
If $M$ and $N$ are equivalent, \fullRef{L.1.5.2} 
shows how to get $T${--}maps $M\ \to\ N$ and $N\ \to\ M$ using
the sequence of strong equivalences.

If $M$ and $N$ are $T${--}equivalent, we have map $T${--}maps
$f\colon M\ \to\ N$ and $g\colon N\ \to\ M$ such that 
$f\ \circ\ g ={\text{id}}_N$ and $g\ \circ\ f = {\text{id}}_M$.
Now a $T${--}map $f\colon M\ \to\ N$ induces a unique map
$\Delta(f)\colon \Delta(M)\ \to\ \Delta(N)$ via
$\Delta(f)=\Delta(f)\ \circ\ \Delta({\text{inc}})^{-1}$
where $f$ is defined on $M_F$ and ${\text{inc}}\colon M_F\subseteq M$.
It is clear that $\Delta(f)$ depends only on the map{--}germ of $f$.
Hence in our case, $g$ induces an equivalence of $M$ and $N$ by
$N\supseteq N_G\ \RA{\ g\ }\ M$.
\end{proof}
\medskip
Also useful is
\medskip
\BEGIN{L.1.5.4}
Let $f$ and $g$ be $T${--}maps.
Then $f=g$ \iff\ $\Delta(f)=\Delta(g)$.
\end{Lemma}
\begin{proof}
$f=g$ \iff\ $f-g=0$.
$\Delta(f-g)=\Delta(f)-\Delta(g)$
Thus we need only show $h=0$ \iff\ $\Delta(h)=0$.
Since $\Delta(h)$ depends only on the map{--}germ, and since
$\Delta(0)=0$, one way is easy.

So assume we are given a $T${--}map $h\colon M\ \to\ N$ with
$\Delta(h)=0$. 
We may as well assume that $h$ is an actual map, since otherwise 
set $M=M_H$ and proceed.
We have a submodule ${\kerx}\ h\subseteq M$ defined in 
the obvious way.
Since ${\kerx}\ h\subseteq M$ is a strong
equivalence, \fullRef{L.1.5.2} says we
can find $F$ such that $M_F\ \to\ {\kerx}\ h\subseteq M$.
But then $M_F\ \to\ N$ is the zero map.
\end{proof}
\bigskip
\begin{xDefinition}
If $R$ is a tree of rings, let ${\mathcal M}_R$ be the category of trees of
$R${--}modules and germs of maps.
Let ${\mathcal M}_{\Delta(R)}$ be the category of $\Delta(R)${--}modules.
\end{xDefinition}
\bigskip
\BEGIN{P.1.5.1}
${\mathcal M}_R$ is an abelian category.
The natural functor \[\Delta\colon {\mathcal M}_R\ \to\ {\mathcal M}_{\Delta(R)}\]
is an exact, additive, faithful functor.
\end{Proposition}
\begin{proof}
The functor just takes $M$ to $\Delta(M)$ and $[f]$ to $\Delta(f)$
($[f]$ denotes the map{--}germ of $f$).
$\Delta$ is additive more or less by definition, and faithful by 
\fullRef{L.1.5.4}.

$\Delta$ preserves kernels: Let $M\ \RA{\ [g]\ }\ N$ be a map{--}germ
in ${\mathcal M}_R$.
We can find $G$ such that $M_g\ \RA{\ g\ }\ N$ is a representative.
Clearly any  kernel for $[g]$ is equivalent to 
${\kerx}\ g\subseteq M_G$, where ${\kerx}\ g$ 
is the obvious submodule.
But $\Delta({\kerx}\ g)$ is clearly a kernel for $\Delta(g)$.

An entirely similar argument shows $\Delta$ preserves cokernels,
so $\Delta$ is exact.

To see ${\mathcal M}_R$ is normal and conormal, take representatives for
the germs and construct the quotient or the kernel module.

${\mathcal M}_R$ has pullback and pushouts, again by finding representatives
for the germs and constructing the desired modules.
Now by \cite{btwentyfive}, Theorem 20.1 (c), page 33, ${\mathcal M}_R$
is abelian.
\end{proof}

We want to do stable algebra, and for this we need an analogue of
finitely{--}generated projective.
Projective is easy, we just insist that a projective $R${--}module 
is projective in the category ${\mathcal M}_R$ (see \cite{btwentyfive}, 
pages 69{--}71 for definitions and elementary properties).

For the analogue of finitely{--}generated, we first produce the analogue
of a finitely{--}generated, free module.

\begin{xDefinition}
Let $T$ be a tree and let $S$ be a set.
A \emph{partition of} $S$ is a functor $\pi\colon {\mathcal C}(T)\ \to\ 2^S$ 
(where $2^S$ is the category of subsets of $S$ and inclusion maps)
satisfying
\begin{enumerate}
\item[1)] $\pi(T)=S$.
\item[2)] If $A\ \cap\ B=\emptyset$, $\pi(A)\ \cap\ \pi(B)=\emptyset$
( $A$, $B\in{\mathcal C}(T)$).
\item3)] Let $A_i\in{\mathcal C}(T)$, $i=1$, \dots, $n$.
If $\displaystyle T-\mathop{\cup}_{i=1}^n A_i$ is compact,
$\displaystyle \pi(T)-\mathop{\cup}_{i=1}^n \pi(A_i)$ is finite.
\item[4)] Let $s\in S$.
Then there exist $A_i\in{\mathcal C}(T)$, $i=1$, \dots, $n$ such that
$\displaystyle T - \mathop{\cup}_{i=1}^n A_i$ is compact and
$s\notin \pi(A_i)$ for any $i=1$, \dots, $n$.
\end{enumerate}
\end{xDefinition}

\bigskip
\begin{xDefinition}
Let $R$ be a tree of rings over $T$.
Let $\pi$ be a partition of $S$.
The \emph{free $R${--}module based on $\pi$}, $F_\pi$, is the tree
of $R${--}modules defined by $(F_\pi)_A$ is the free $R_A${--}module
based on $\pi(A)$, and if $A\subseteq B$,
$p_{A B}\colon (F_\pi)_A\ \to\ (F_\pi)_B$ is induced by the inclusion
$\pi(A)\subseteq\pi(B)$.
\end{xDefinition}

\begin{xDefinition}
A tree of $R${--}modules, $M$, is said to be locally{--}finitely generated
\iff\ there is a set of generators, $S$, and a partition $\pi$, of $S$,
such that there is an epimorphism $F_\pi\ \to\ M$.
\end{xDefinition}

Let us briefly discuss partitions.
If $\pi$ and $\rho$ are two partitions of a set $S$, we say 
$\pi\subseteq\rho$ \iff\ $\pi(A)\subseteq\rho(A)$ for all $A\in{\mathcal C}(T)$.
(Hence we could talk about the category of partitions, but we shall
largely refrain.)
Two partitions are \emph{equivalent} \iff\ there exist a finite sequence
$\pi=\pi_0$, $\pi_1$, \dots, $\pi_n=\rho$ of partitions with
$\pi_{i}\subseteq\pi_{i+1}$, or $\pi_{i+1}\subseteq\pi_{i}$.
(This is clearly an equivalence relation.)
Given two sets $X$ and $Y$, and partitions $\pi$ and $\rho$,
$\pi\cup\rho$ is the partition $X\cup Y$ given by $(\pi\cup\rho)(A)=
\pi(A)\cup\rho(A)$.
\bigskip
\BEGIN{L.1.5.5}
Let $R$ be a tree of rings over $T$, and let $X$ and $Y$ be sets.
Then if $\pi$ and $\pi^\prime$ are equivalent partitions of $X$,
$F_\pi$ is isomorphic to $F_{\pi^\prime}$ in ${\mathcal M}_R$.
If $\rho$ is a partition of $Y$,
$F_{\pi\cup\rho}=F_\pi\oplus F_\rho$ ($X$ and $Y$ disjoint).
\end{Lemma}
\begin{proof}
To show the first statement we need only show it for 
$\pi\subseteq \pi^\prime$.
In this case there is a natural map $f\colon F_\pi\ \to\ F_{\pi^\prime}$.
For each $A\in{\mathcal C}(T)$, $(F_\pi)_A\ \to\ (F_{\pi^\prime})_A$ is
injective, so $f$ is a monomorphism. 
If $\pi\subseteq \pi^\prime$, then $\pi^\prime(A)-\pi(A)$ has only
finitely many elements.
To see this, observe we can find $A_i\in{\mathcal C}(T)$, $i=1$, \dots, $n$
such that $A\ \cap\ A_i=\emptyset$, and
$\displaystyle T - \mathop{\cup}_{i=1}^n A_i - A$ is compact.
Then by 2) $\displaystyle\pi^\prime(A)\subseteq \pi^\prime(T)-
 \mathop{\cup}_{i=1}^n \pi^\prime(A_i)$, so
$\displaystyle\pi^\prime(A)-\pi(A)\subseteq \pi^\prime(T)-
 \mathop{\cup}_{i=1}^n \pi^\prime(A_i)-\pi(A)\subseteq
\pi(T)- \mathop{\cup}_{i=1}^n \pi(A_i)-\pi(A)$, which is finite.
Since $\pi^\prime(A)-\pi(A)$ is finite, 
$f_A\colon (F_\pi)_A\ \to (F_{\pi^\prime})_A$ has finitely generated 
cokernel, so when the $\Delta${--}construction is applied to it,
4) guarantees that $\Delta(f)$ is onto, so $f$ is an equivalence.
The second statement is the definition of $\pi\cup\rho$ and
$F_\pi\oplus F_\rho$.
\end{proof}

It is not hard to see that if we have a partition of $S$ for the tree $T$,
then $S$ has at most countably many elements if $T$ is infinite, and
at most finitely many if $T$ is a point.
In the case $S$ is infinite, we have a very handy countable infinite
set lying around, namely the vertices of $T$.
There is an obvious partition, $\pi$, where $\pi(A)=\{ p\ \vert\ 
p{\text{\ is\ a\ vertex\ of\ }}A\ \}$.
Denote $F_\pi$ by $F^{(1)}$.
If $T$ is a point, let $F^{(1)}$ denote the free module on one 
generator; i.e. still $F_\pi$ for the above partition $\pi$.
$F^{(n)}= F^{(n-1)}\oplus F^{(1)}$ for $n\geq 2$.
\bigskip
\BEGIN{L.1.5.6}
Let $\pi$ be any partition of a set $S$ for the tree $T$, and let $R$
be a tree of rings.
Then $F_\pi\oplus F^{(1)}$ is equivalent to $F^{(n)}$ for some
$n\geq 1$.
If $T$ is infinite, $n$ can be chosen to be $1$.
\end{Lemma}
\begin{proof}
If $T$ is a point, this is obvious, so assume $T$ is infinite.
$F_\pi\oplus F^{(1)}$ is just $F_{\pi\cup\rho}$, where $\rho$
is the standard partition on $V$, the vertices of $T$.
Since $V\ \cup\ S$ is infinite (and countable), there is a $1${--}$1$
correspondence $\alpha\colon V\ \cup\ S\ \to\ V$.
Any such $\alpha$ induces an equivalence of categories
$\alpha\colon 2^{V\cup S}\ \to\ 2^V$.
We show that we can pick $\alpha$ so that 
$\alpha\ \circ\ (\pi\ \cup\ \rho)$ is equivalent to $\rho$.
(We will show in \fullRef{L.1.5.7} that 
$\alpha\ \circ\ (\pi\ \cup\ \rho)$ is a partition for any $\alpha$.)
Our $\alpha$ is defined by picking a strictly increasing sequence
of finite subcomplexes, $C_0\subseteq C_1\subseteq\cdots$, so
that $\displaystyle\mathop{\cup}_{i=0}^\infty C_i=T$.
Let $A_k(i)$ be the essential components of $T-C_i$.
Set $A_1(-1)=T$, and let $\displaystyle K_{k i}=
(\pi\cup\rho)\bigl(A_k(i)\bigr) -\mathop{\cup}_\ell\ (\pi\cup\rho)\bigl(
A_\ell(i+1)\bigr)$.
Note $K_{k i}\ \cap\ K_{k^\prime i}=\emptyset$ and 
$K_{k i}\ \cap\ K_{k^\prime i+1}=\emptyset$ by 2), so
$K_{k i}\ \cap\ K_{\ell j}\neq \emptyset$ \iff\ $k=\ell$ and $i=j$.

Now $K_{k i}$ is finite.
We define $\alpha$ on $K_{k i}$ by induction on $i$.
Let $\displaystyle L_{k i}=\rho\bigl(A_k(i)\bigr) -\mathop{\cup}_\ell
\rho\bigl(A_\ell(i+1)\bigr)$, and  note that the cardinality of $K_{k i}$
is greater than or equal to the cardinality of $L_{k i}$.
Define $\alpha$ on $K_{1\ -1}$ by mapping some subset of it to
$L_{1\ -1}$ and mapping any left over elements to any elements of $V$
($\alpha$ should be injective).

Suppose $\alpha$ defined on $k_{k\>i-1}$ so that $\alpha(K_{k j})
\subseteq\rho\bigl(A_k(j)\bigr)$ for $j\leq i-1$.
We need only define $\alpha$ on $K_{k i}$ so that $\alpha(K_{k i})
\subseteq\rho\bigl(A_k(i)\bigr)$ to be done.
Look at 
$\displaystyle M = L_{k i} - 
\mathop{\cup}_{\Atop{\scriptscriptstyle{\text{all\ }}\ell}{j\leq i-1}}
\ {\Imx}\ \alpha(K_{\ell j})$.
Map some subset of $K_{k i}$ to $M$.
Map the rest of $K_{k i}$ to any elements of $\rho\bigl(A_k(i)\bigr)$ at all.

By 4), $\displaystyle V\cup S=
\mathop{\cup}_{\Atop{\scriptscriptstyle{\text{all\ }}k}{{\text{all\ }}i}} K_{k i}$
and 
$\displaystyle S= 
\mathop{\cup}_{\Atop{\scriptscriptstyle{\text{all\ }}k}{{\text{all\ }}i}} L_{k i}$
(as disjoint unions as we saw).
Since $\alpha$ is onto each $L_{k i}$, and since it injects when restricted to
each $K_{k i}$, $\alpha$ is $1${--}$1$.
Furthermore, $\tau=\alpha\ \circ\ (\pi\cup\rho)$ satisfies
$\tau\bigl(A_k(i)\bigr)\subseteq \rho\bigl(A_k(i)\bigr)$ by construction. 

Set $\lambda(A)=\tau(A)\cap\rho(A)$.
We claim $\lambda$ is a partition. 
Clearly $\lambda$ is a functor ${\mathcal C}(T)\ \to\ 2^V$.
1) and 2) are trivial and 4) is not much harder (\ 1), 2) and 4) hold for
the intersection of any two partitions, it is only 3) which might fail).
To show 3), note $\lambda\bigl(A_k(i)\bigr)=\tau\bigl(A_k(i)\bigr)$.
If $\displaystyle T-\mathop{\cup}_{j=1}^n B_j$ is compact, there is
a minimal $i$ such that $B_j$ contains $A_k(i)$ for some $k$
(perhaps several, say $k=1$, \dots, $m$).
Then
$\displaystyle \mathop{\cup}_{k=1}^m \lambda\bigl(A_k(i)\bigr)\subseteq
\lambda(B_j)$.
$\displaystyle\lambda(T)-\mathop{\cup}_{j=1}^n\lambda(B_j)\subseteq
\lambda(T)-\cup\ \lambda\bigl(A_k(i)\bigr)=
\tau(T)-\cup\ \tau\bigl(A_k(i)\bigr)$.
The last two unions are over all the $A_k(i)\subseteq B_j$ for 
$j=1$, \dots, $n$. 
This last set is finite, so 3) holds.
Hence $\lambda$ is a partition and thus $\tau$ is equivalent to $\rho$.

The map from $F_{\pi\cup\rho}\ \to F_\tau$ induced by $\alpha$
is the obvious map: $(F_{\pi\cup\rho})_A\ \to\ (F_\tau)_A$ is the
isomorphism induced by the equivalence of bases $\alpha\colon
(\pi\cup\rho)(A)\leftrightarrow\tau(A)$.
\fullRef{L.1.5.5} completes the proof modulo
the proof of \fullRef{L.1.5.7}.
\end{proof}

\bigskip
\BEGIN{L.1.5.7}
Let $X$ and $Y$ be two (disjoint) sets, and let $\pi$ be a partition
of $X$ for the tree $T$.
Any $1${--}$1$ correspondence $\alpha\colon X\ \to\ Y$ induces a
partition $\alpha\ \circ\ \pi$ of $Y$ for the tree $T$.
\end{Lemma}
\begin{proof} The easy proof is omitted.\end{proof}
\bigskip
\BEGIN{L.1.5.8}
$F_\pi$ is projective.
\end{Lemma}
\begin{proof}
By \fullRef{L.1.5.6} and standard nonsense, it is
enough to prove the result for $F^{(1)}$.
By Mitchell \cite{btwentyfive} Proposition 14.2, page 70, we need only show
$M\ \RA{\ [f]\ }\ F^{(1)}$ splits whenever $[f]$ is an epimorphism
(note ${\mathcal M}_R$ is abelian by \fullRef{P.1.5.1} 
so we may apply Mitchell).

By taking a representative for $[f]$, we may as well assume that we
have a map $f\colon M\ \to\ F=F^{(1)}$ which is an epimorphism.
Now there is a partition $\pi$ with $\pi\subseteq \rho$
($\rho$ the standard partition for $F^{(1)}$), such that the
inclusion of $(F_\pi)_A$ in $F_A$ lies in the image of $M_A$
under $f_A$; i.e. define $\pi(A)=\{ x\in\rho(A)\ \vert\ 
x\in{\Imx}\ f_A\ \}$.
Since $f$ is an epimorphism, one can easily check $\rho(A)-\pi(A)$
is finite, and from this result one easily deduces $\pi$ is a partition.

Now pick a base point $\ast\in T$.
This choice immediately partially orders all the vertices of $T$ by saying
$p\geq q$ provided the minimal path from $p$ to $\ast$ hits $q$.
$A_p\in{\mathcal C}(T)$ for each $p$ a vertex of $T$, $p\neq\ast$, is defined
as the unique $A\in{\mathcal C}(T)$ such that $q\in A_p$ implies $q\geq p$.

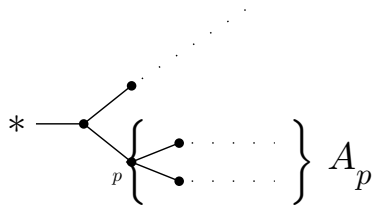
\IpdfLRT{diagram151}{2.5in}{2in}

\vskip -60pt\noindent Given a partition $\pi$, define a new partition $\tau$ by
$\displaystyle\tau(A) = \mathop{\cup}_{A_p\subseteq A}\ \pi(A_p)$
 (again $\pi(A)-\tau(A)$ is finite, $\tau(A)\subseteq \pi(A)$, so one
can check $\tau$ is a partition).
Since $\tau\subseteq\pi$, $(F_\tau)_A\subseteq F_A$ lies in 
${\Imx\ }(f_A)$.

Now given any vertex $v$ of $T$, there is a unique $p$ such that 
$v\in\tau(A_p)$ and $v\in\tau(A)$ \iff\ $A_p\subseteq A$, unless 
$v\notin \tau(A_p)$ for any $A_p$ (there are only finitely many of
the latter).
To see this, set $A=\displaystyle\mathop{\cap}_{v\in\tau(A_p)} A_p$.
Now $A_p\ \cap\ A_q\neq\emptyset$ implies $A_p\subseteq A_q$
(or $A_q\subseteq A_p$).
By 4) the intersection runs over finitely many objects, so $A=A_p$ for
some $p$.
This $A_p$ has the properties we claimed.

Define $x_v\in M_{A_p}$ to be any element such that $f_{A_p}(x_v)$
hits the image of the generator in $(F_\pi)_{A_p}$ corresponding to $v$.
Define $h\colon F_\tau\ \to\ M$ by $h_A\colon (F_\tau)_A\ \to\ M_A$
takes the generator corresponding to $v$ to $p_{A_p\>A}(x_v)$.
We extend linearly.
Notice that if the generator corresponding to $v$ lies in $(F_\tau)_A$,
$A_p\subseteq A$, so $p_{A_p\>A}$ makes sense.

It is not hard to check that the $h_A$ induce a map 
$h\colon F_\tau\ \to\ M$, and $f\ \circ\ h\colon F_\tau\ \to\ F$
is just the inclusion.
\end{proof}

\medskip
If ${\mathcal P}_R$ is the category of locally{--}finitely generated trees of
projective $R${--}modules, we have
\medskip
\BEGIN{L.1.5.9}
Let $0\to\ M\ \to\ N\ \to\ Q\ \to0$ be a short exact sequence
of $R${--}modules.
Then, if $N$, $Q\in{\mathcal P}_R$, $M\in{\mathcal P}_R$.
If $M$, $Q\in{\mathcal P}_R$, $N\in{\mathcal P}_R$.
Lastly, any $P\in{\mathcal P}_R$ is a summand of a locally{--}finitely generated
free module.
\end{Lemma}
\begin{proof} The proof is easy.\end{proof}

\medskip\begin{xRemarks}
${\mathcal P}_R$ is a suitable category in which to do stable algebra (see 
Bass \cite{bone}).
${\mathcal P}_R$ has a product, the direct sum. ${\mathcal P}_R$ is also a full
subcategory of ${\mathcal M}_R$, which is abelian by 
\fullRef{P.1.5.1}.
Hence we may use either of Bass's definitions of the $K${--}groups.
Note ${\mathcal P}_R$ is semi{--}simple (Bass \cite{bone}) so
the two definitions agree.
\end{xRemarks}

\insetitem{Notation}
$K_0(R)=K_0({\mathcal P}_R)$ and $K_1(R)=K_1({\mathcal P}_R)$ for
$R$ a tree of rings.
\bigskip
Given a map of trees of rings $R\ \to\ S$ ($R_A\ \to\ S_A$
takes units to units) we can define $M\otimes_R S$ for
$M$ a right $R${--}module by taking $(M\otimes_R S)_A=
M_A\otimes_{R_A} S_A$.
$\otimes$ induces a functor ${\mathcal M}_R\ \to\ {\mathcal M}_S$.
The only non{--}trivial part of this is to show $\otimes$ is well{--}defined
on map{--}germs.
But since
\[\begin{matrix}%
{\mathcal M}_R&\rightlabeledarrow{\hskip40pt\otimes_R S \hskip40pt}{}&{\mathcal M}_S\\
\Big\downarrow&&\Big\downarrow\\
{\mathcal M}_{\Delta(R)}&
\rightlabeledarrow{\hskip 40pt\otimes_{\ \Delta(R)} \Delta(S)\hskip40pt}{}&{\mathcal M}_{\Delta(S)}\\
\end{matrix}\]
commutes, this is easy.
$\otimes$ is, as usual, an additive, right exact functor.

Now given a partition $\pi$, $F^R_\pi\otimes_R S= F^S_\pi$, where
$F^R_\pi$ is the free $R${--}module based on $\pi$ ( $F^S_\pi$similarly).
Hence it is easy to see that $\otimes$ takes ${\mathcal P}_R$ to ${\mathcal P}_S$.
$\otimes$ is cofinal in the sense of Bass \cite{bone}, so we get a 
relative group $K_0(f)$, where $f\colon R\ \to\ S$ is the map of trees
of rings.
There is an exact sequence
\[K_1(R)\ \to\ K_1(S)\ \to\ K_0(f)\ \to\  K_0(R)\ \to\ K_0(S)\ .\]

We denote by $K_i(T)$, $i=0$, $1$, the result of applying the 
$K${--}groups to ${\mathcal P}_T$, where ${\mathcal P}_T$ is the category of
locally{--}finitely generated projective modules 
over the tree of rings $``T''$,
where $(``T'')_A=\Z$ for all $A$ and $p_{A B}={\text{id}}$.
There is always a functor ${\mathcal P}_T\ \to\ {\mathcal P}_R$
induced by the unit map $``T''\ \to\ R$.
The relative $K_0$ of this map will be called the reduced $K_1$
of $R$, written $\Bar{K}(R)$.

\medskip\begin{xRemarks}
If the tree of rings is a point the functor 
${\mathcal M}_R\ \to\ \Delta({\mathcal M}_R)$ induces a functor
${\mathcal P}_R\ \to\ {\mathcal P}_{\Delta(R)}$, where ${\mathcal P}_{\Delta(R)}$
is the category of finitely{--}generated projective $\Delta(R)${--}modules.
This functor induces an isomorphism on $K_0$ and $K_1$.
For the compact case ($T={\text{pt}}$.),
torsions lie in quotients of $K_1({\mathcal P}_{\Delta(R)})$.
This, together with \fullRef{P.1.5.2} below is
supposed to motivate our choice of ${\mathcal P}_R$ as the category
in which to do stable algebra.
\end{xRemarks}

\medskip\begin{xDefinition}
Let $W$ be an \hCWx\  complex of finite dimension.
Let $X$ and $Y$ be subcomplexes.
Let $(T,f)$ be a tree for $W$.
Lastly let $F\in {\mathcal L}(f)$.
Then $\Z \pi_1(W,F,f)$ is the tree of rings we had earlier as an example.
Pick a locally{--}finite set of paths, $\Lambda$, from the cells of $W$
to the vertices of $f(T)$ (the paths all begin at the barycenter 
of each cell).
\end{xDefinition}

\medskip
$C_\ast(W;X,Y\Colon \Lambda,F)$ is the tree of $\Z\pi_1(W,F,f)${--}modules
given at $A$ by
\[H_\ast\left(
\coverFC{F(A)\ }^{\lower5pt\hbox{$\scriptstyle\ast$}}; \coverFC{F(A)\ }^{\lower5pt\hbox{$\scriptstyle\ast-1$}},
\coverFC{F(A)^\ast\cap X\hskip 6pt },
\coverFC{F(A)^\ast\cap Y\hskip 6pt }
\right)
\]
where \ $\coverFC{\hskip10pt}$\quad is the universal cover of $F(A)$, so,
for example, $\coverFC{F(A)^\ast\cap Y\hskip 6pt}$ 
is the part of the universal
cover of $F(A)$ lying over $Y\ \cap\ ({\text{the\ }}\ast{\text{{-}skeleton\ of\ }}F(A)$.
In each $\coverFC{F(A)}$ pick a base point covering the vertex
$\partial A$.
These choices give us maps 
$\coverFC{F(A)}\ \to\ \coverFC{F(B)}$ whenever $A\subseteq B$.
\medskip
$C^\ast(W;X,Y\Colon \Lambda,F)$ is defined from the cohomology groups
\[H^\ast_\cmpsup\left(
\coverFC{F(A)\ }^{\lower5pt\hbox{$\scriptstyle\ast$}}; \coverFC{F(A)\ }^{\lower5pt\hbox{$\scriptstyle\ast-1$}},
\coverFC{\partial F(A)\ }^{\lower5pt\hbox{$\scriptstyle\ast$}},
\coverFC{F(A)^\ast\cap X\hskip 6pt}, 
\coverFC{F(A)^\ast\cap Y\hskip 6pt}
\right) .
\]
The maps are the ones we defined in section 4.
\bigskip
\BEGIN{P.1.5.2}
$C_\ast(W;X,Y\Colon\,\Lambda,F)$ ($C^\ast(W;X,Y\Colon\,\Lambda,F)$) is a 
locally{--}finitely generated, free, right (left) $\Z\pi_1(W,F,f)${--}module.
If $G\in{\mathcal L}(f)$ satisfies $G\geq F$, there is an induced map
$\Z\pi_1(W,F,f)\ \to\ \Z\pi_1(W,G,f)$. 
$C_\ast(W;X,Y\Colon\,\Lambda,F)\otimes_{\Z\pi_1(W,F,f)}\Z\pi_1(W,G,f)$
is equivalent to $C_\ast(W;X,Y\Colon\,\Lambda,G)$.
$\Z\pi_1(W,G,f)\otimes_{\Z\pi_1(W,F,f)}C^\ast(W;X,Y\Colon\,\Lambda,F)$
is equivalent to \\$C^\ast(W;X,Y\Colon\,\Lambda,G)$.
The $\Delta${-}functor applied to $C_\ast(W;X,Y\Colon\,\Lambda,F)$
is \\$P_\ast(W;X,Y\Colon\,\coverFC{\hskip10pt})$;
$\Delta\bigl(C^\ast(W;X,Y\Colon\,\Lambda,F)\bigr) = 
P^\ast(W;X,Y\Colon\,\coverFC{\hskip10pt})$ (the $P$ were defined in section 4) ).
\end{Proposition}
\begin{proof}
The assertions are all fairly obvious.
Note in passing  that the set $S$ for $C_\ast$ ($C^\ast$) is the set of
all $\ast${--}cells in $W-(X\ \cup\ Y)$.
\end{proof}

\medskip
\BEGIN{P.1.5.3}
The choice of paths $\Lambda$ determines a basis for $C_\ast$ ($C^\ast$).
\end{Proposition}
\begin{proof}
Let $S$ be the set of all $\ast${--}cells in $W-(X\ \cup\ Y)$.
Partition $S$ by $\pi(A)=$ the set of all $\ast${--}cells in $W-(X\ \cup\ Y)$
such that the cell and its associated path both lie in $F(A)$.
$\pi$ is seen to be a partition, and $F_\pi$ is equivalent to $C_\ast$.
The path also determines a lift of the cell into $\coverFC{F(A)}$,
so each $(F\pi)_A$ is based.
\end{proof}
\medskip
Apparently our tree of rings and modules is going to depend on the
lift functor we choose.
This is not the case and we proceed to prove this.
Given a shift functor $F$ and a tree of rings $R$, $R_F$ is the tree 
of rings given by 
$\displaystyle (R_F)_A=\mathop{\oplus}_{i=1}^n\ R_{A_i}$
where the $A_i$ are the essential components of $F(A)$.
$p_{A B}$ is just $\oplus p_{i j}$, where $p_{i j}$ is the projection
$p_{A_i B_j}$ where $A_i\subseteq B_j$.

We now redefine $M_F$. 
$M_F$ is going to be an $R_F${--}module.
$\displaystyle (M_F)_A=\mathop{\oplus}_{i=1}^n\ M_{A_i}$
with the obvious $R_F${--}module structure.
Note $M_F\otimes_{R_F} R$ is just our old $M_F$.

Now a $T${--}map of rings is just a map $R_F\ \to\ S$.
As in the case of modules, we can define a map{--}germ between
two rings.

\medskip
\BEGIN{L.1.5.10}
The maps $K_i(R_F)\ \to\ K_i(R)$, $i=0$, $1$, are isomorphisms.
\end{Lemma}
\begin{proof} 
$M\mapsto M_F$, $f\mapsto f_F$ defines a functor 
${\mathcal P}_R\ \to\ {\mathcal P}_{R_F}$.
Using this functor, one checks ${\mathcal P}_{R_F}\ \to\ {\mathcal P}_R$
is an equivalence of categories. 
The result is now easy.
\end{proof}
\medskip
Hence given a map{--}germ $f\colon R\to S$, we get well{--}defined
induced maps $K_i(R)\ \to\ K_i(S)$, $i=0$, $1$, and 
$\Bar K_1(R)\ \to\ \Bar K_1(S)$.
\medskip
\BEGIN{L.1.5.11}
Let $f\colon R\ \to\ S$ be a map such that $\Delta(f)$ is an
isomorphism.
Then there is a shift functor $F$ and a map $g\colon S_F\ \to\ R$
such that 
\[\begin{matrix}%
&&R\\
&\hbox to 0pt{\hss$\scriptstyle g$}\nearrow&\downlabeledarrow{}{f}\\
S_F&\to&S
\end{matrix}\]
commutes.
\end{Lemma}

\begin{proof}
The proof is just like that of \fullRef{L.1.5.2}.
\end{proof}

\medskip
\BEGIN{L.1.5.12}
Let $[f]\colon R\ \to\ S$ be a map{--}germ such that $\Delta(f)$
is an isomorphism.
Then the maps 
$K_0(R)\ \to\ K_0(S)$; $K_1(R)\ \to\ K_1(S)$; and 
$\Bar K_1(R)\ \to\ \Bar K_1(S)$ are isomorphisms.
\end{Lemma}
\begin{proof}
This proof is easy and will be left to the reader.
\end{proof}
\medskip\begin{xRemarks}
By \fullRef{L.1.5.12}, the $K${--}groups we get
will not depend on which lift functor we use.
Let $\displaystyle K_i(X\Colon f)=
\lim_{\Atop{\longrightarrow}{F\in{\mathcal L}(f)}}\ 
K_i\bigl(\Z\pi_1(X,F,f)\bigr)$.
Since all the maps in our direct limit are isomorphisms, $K_i(X\Colon f)$
is computable in terms of $K_i\bigl(\Z\pi_1(X,F,f)\bigr)$ for any $F$.
$\Bar K_1(X\Colon f)$ is defined similarly.
\end{xRemarks}
\bigskip

\begin{xDefinition}
A stably free (s{--}free) tree of $R${--}modules is an element, $P$,
of ${\mathcal P}_R$ such that $[P]$ is in the image of $K_0(T)$.
Let $P$ be an s{--}free\ $R${--}module.
An s{--}basis for $P$ is an element $F\in{\mathcal R}_T$ and an isomorphism
$b\colon F\otimes_T R\ \to\ P\oplus F_1\otimes_T R$, where 
$F_1\in{\mathcal P}_T$.
\end{xDefinition}

Two s{--}bases $b\colon F\otimes_T R\ \to\ P\oplus F_1\otimes_T R$
and $c\colon F_2\otimes_T R\ \to\ P\oplus F_3\otimes_T R$ are
\emph{equivalent} ($B\sim c$) \iff\ 
$0=(F\oplus F_3, (b\oplus {\text{id}}_{F_3})\ \circ\ {\tw\ }\circ\
(c\oplus {\text{id}}_{F_1})^{-1}, F_2\oplus F_1)$ in $\Bar K_1(R)$,
where 
${\tw\ }\colon (P\oplus F_1\otimes_T R)\oplus F_3\otimes_T R\ \to\ 
 (P\oplus F_3\otimes_T R)\oplus F_1\otimes_T R$ is the obvious map.

\medskip
We can now give an exposition of torsion following 
Milnor \cite{btwentythree}.
Given a short exact sequence 
$0\to\ E\ \RA{\ i\ }\ F\ \RA{\ p\ }\ G\ \to0$ and s{--}bases $b$ for $E$
and $c$ for $G$, define an s{--}basis $b c$ for $F$ by picking a
splitting $r\colon G\to F$ for $p$ and then taking the composition
$F_1\oplus F_2\ \RA{\ (b,c)\ }\ (E\oplus F_3)\oplus (G\oplus F_4)\ 
\RA{\ h\ }\ F\oplus(F_3\oplus F_4)$, where
$h(e,x,g,z)$ goes to $\bigl( i(e)+r(g), x, z\bigr)$.
It is not hard to check that this s{--}basis does not depend on the choice 
of splitting map.

We use Milnor's formulation.
Let $F_0\subseteq F_1\subseteq$ \dots $\subseteq F_k$ and suppose
each $F_i/F_{i-1}$ has an s{--}basis $b_i$. 
Then $b_1b_2\cdots b_k$ is seen to be well{--}defined; 
i.e. our construction is associative.

Let $E$ and $F$ be submodules of $G$.
Then $E+F$ is the submodule of $G$ generated by $E$ and $F$.
$E\ \cap\ F$ is the pullback of
\lower \baselineskip\hbox{$\begin{matrix}%
&&E\\
&&\downarrow\\
F&\to&G
\end{matrix}$}

\smallskip
\BEGIN{L.1.5.13}
(Noether) The natural map
$E/(E\cap F)\ \to\ (E+F)/F$ is an isomorphism.
\end{Lemma}
\begin{proof}
Apply the ordinary Noether isomorphism to each term.
\end{proof}
\medskip
Now let $E/(E\cap F)$ have an s{--}basis $b$ and let $F/(E\cap F)$
have an s{--}basis $c$. 
Base $(E+F)/F$ by $b$ composed with the Noether map (we will
continue to denote it by $b$). 
Similarly base $(E+F)/E$ by $c$.
Then $b c\sim c b$ as s{--}bases for $(E+F)/(E\cap F)$.
\medskip\begin{xDefinition}
Let $b$ and $c$ be two s{--}bases for $P$.
Then $[b/c]\in\Bar K_1(R)$ is defined as follows:
if $F\ \RA{\ b\ }\ P\oplus F_1$; $G\ \RA{\ c\ }\ P\oplus F_2$, then
$[b/c]= (F\oplus F_2, h, G\oplus F_1)$ where 
$\vrule height 12pt depth 0pt width 0pt
h\colon F\oplus F_2\ \RA{\ b\oplus{\text{id}}\ }\ 
(P\oplus F_1)\oplus F_2\ \to\ (P\oplus F_2)\oplus F_1\ 
\RA{\ c^{-1}\oplus {\text{id}}\ }\ G\oplus F_1
\vrule height 14pt depth 0pt width 0pt$.
Two s{--}bases are equivalent \iff\ $[b/c]=0$.
The formulas
$[b/c]+[c/d]=[b/d]$ and $[b/c]+[d/e]=[b d/c e]$ are easy to derive
from the relations in the relevant $K_1$.
\end{xDefinition}

\setcounter{footnote}{0}
We next define a torsion for chain complexes.
A \emph{free} chain complex is a set of s{--}free modules, $P_n$,
together with map{--}germs $\partial_n\colon P_n\ \to\ P_{n-1}$
such that $\partial_n\ \circ\ \partial_{n-1}=0$\footnote{Probably should have been called s{--}free.}.
A \emph{finite} free chain complex is one with only finitely many non{--}zero
$P_n$.
A \emph{positive} free chain complex has $P_n=0$ for $n<0$.
\medskip\begin{xDefinition}
Let $\{P_n, \partial_n\}$ be a finite free chain complex.
Let $P_n$ be s{--}based by $c_n$, and suppose each homology group
$H_i$ is s{--}free and s{--}based by $h_i$.
\end{xDefinition}

The sequences $0\to\ B_{n+1}\ \to\ Z_n\ \to\ H_n\ \to0$ and
$0\to\ Z_n\ \to\ P_n\ \to\ B_n\ \to0$, where 
$B_n={\Imx}\ (P_n\ \to\ P_{n-1})$ and 
$Z_n={\kerx} (\partial_n)$, are short exact.
Let $b_n$ be an s{--}basis for $B_n$, which exists by an inductive
argument.
\[\tau(P_\ast)=\sum_{n}(-1)^n\bigl[ b_n h_n b_{n-1}/c_n\bigr]
\in\Bar K_1(R)\ .\]

It is easy to show $\tau(P_\ast)$ does not depend on the choice of $b_n$.
Let $0\to\ P^\prime_\ast\ \to\ P_\ast\ \to\ 
P^{\prime\prime}_\ast\ \to0$ be a short exact sequence of finite free 
chain complexes.
There is a long exact sequence 
\[\begin{matrix}H_\ast(P^\prime)&\RA{\hskip20pt}&H_\ast(P)\\
\hskip10pt\hbox to 0pt{$\hbox to 0pt{\hss$\scriptstyle\partial$}\nwarrow$\hss}
&&\hskip-20pt\hbox to0pt{\hss$\swarrow$}\\
&H_\ast(P^{\prime\prime})\\\end{matrix}\]
Suppose each homology group is s{--}based.
Then we have a torsion associated to ${\mathcal H}$, where
\[{\mathcal H}_{3n} = H_n(P^{\prime})\ ,\quad
{\mathcal H}_{3n-1} = H_n(P)\ ,\quad
{\mathcal H}_{3n-2} = H_n(P^{\prime\prime})\ ,\]
since ${\mathcal H}$ is acyclic.

\bigskip
\BEGIN{T.1.5.1}
$\tau(P_\ast)=\tau(P^\prime_\ast)+\tau(P^{\prime\prime}_\ast)+
\tau({\mathcal H})$.
\end{Theorem}
\begin{proof} See Milnor \cite{btwentythree}, Theorems 3.1 and 3.2.
\end{proof}
\medskip
We next describe the algebraic Subdivision Theorem of 
Milnor \cite{btwentythree} (Theorem 5.2).
Given a chain complex $C_\ast$, suppose it is filtered by
$C_\ast^{(0)}\subseteq C_\ast^{(1)}\subseteq \cdots
\subseteq C_\ast^{(n)}= C_\ast$ such that the homology groups
$H_i\bigl( C^{(\lambda)}/C^{(\lambda-1)}\bigr)=0$ for
$i\neq \lambda$. ($C_\ast^{(-1)}=0$). 

Then we have a chain complex $(\Bar C_\ast, \Bar\partial_\ast)$ 
given by $\Bar C_\lambda= H_\lambda
\bigl( C^{(\lambda)}/C^{(\lambda-1)}\bigr)$ and $\Bar\partial$ is given
by the boundary in the homology exact sequence of the triple
$\bigl( C^{(\lambda)}, C^{(\lambda-1)}, C^{(\lambda-2)}\bigr)$.
There is a well{--}known canonical isomorphism
$H_i(\Bar C)\ \RA{\ \cong\ }\ H_i(C)$ (see Milnor, Lemma 5.1).

Now suppose each $C^{(\lambda)}_i/C^{(\lambda-1)}_i$ has an s{--}basis
$c^\lambda_i$: each $\Bar C_\lambda$ has an s{--}basis $\Bar c_\lambda$:
each $H_i(\Bar C)$ has an s{--}basis $h_i$.
Assume $C_\ast$ is a finite complex.
Then so is $\Bar C_\ast$.

Each $C^{(\lambda)}/C^{(\lambda-1)}$ has a torsion.
If $C_i$ is s{--}based by $c^0_i c^1_i$ \dots $c^n_i$, and $H_i(C)$ is
based by $h_i$ composed with the canonical isomorphism, then the
torsion of $C$ is defined.
Lastly, the torsion of $\Bar C$ is also defined.
\bigskip
\BEGIN{T.1.5.2}
(Algebraic Subdivision Theorem)
\[\tau(C)= \tau(\Bar C\> )+\sum_{\lambda=0}^n
\tau\bigl(C^{(\lambda)}/C^{(\lambda-1)}\bigr)\quad .\]
\end{Theorem}
\begin{proof}
The proof is the same as Milnor's \cite{btwentythree}, Theorem 5.2.
One does the same induction, but one just shows
$\displaystyle \tau(C^{(k)})= \tau(\Bar C^{(k)})+\sum_{\lambda=0}^k
\tau\bigl(C^{(\lambda)}/C^{(\lambda-1)}\bigr)$ 
(notation is the same as Milnor's).
\end{proof}
\medskip
Now let $(K,L)$ be a pair of finite dimensional \hCWx\  complexes with $L$
a proper deformation retract of $K$.
We have the modules $C_\ast(K, L\Colon \Lambda, F)$.
The exact sequence of a triple makes $C_\ast$ into a chain complex,
whose homology is zero since $L$ is a proper deformation retract of
$K$.
The paths $\Lambda$ give us a basis for $C_\ast$ up to sign; i.e.
we must orient each cell, which we can do arbitrarily. 
$\tau(K,L\Colon \Lambda,f)\in \Bar K_1\bigl(\Z\pi_1(K,F,f)\bigr)$
is the torsion of this complex with the basis given by $\Lambda$.
We proceed to show that it does not depend on the choice of signs.

Let $\tau^\prime$ be the torsion with a different choice of signs.
Then, by \fullRef{L.1.5.14} below,
$\tau^\prime - \tau=
\displaystyle\sum_\ast (-1)^\ast[c_\ast/c^\prime_\ast]$
where $c_\ast$ and $c^\prime_\ast$ are maps $F_\pi\ \to\ C_\ast$,
one with the signs for $\tau$ and the other 
with the signs for $\tau^\prime$.
But $c^{-1}_\ast\ \circ\ c^\prime_\ast\colon F_\pi\ \to\ F_\pi$
lies in the image of ${\mathcal P}_T\ \to\ {\mathcal P}_R$, and so
$[c_\ast/c^\prime_\ast]=0$ in $\Bar K_1(R)$.
\medskip
\BEGIN{L.1.5.14}
Let $C_\ast$ be a finite chain complex. 
Let $c_\ast$ and $c^\prime_\ast\colon F_\pi\ \to\ C_\ast$ be two free
bases for $C_\ast$. 
Suppose $H_\ast(C)$ is s{--}based.
Let $\tau$ and $\tau^\prime$ be the torsions from the bases
$c_\ast$ and $c^\prime_\ast$ respectively.
Then $\tau-\tau^\prime = 
\displaystyle\sum_{\ast}(-1)^\ast \bigl[c_\ast/c^\prime_\ast\bigr]$.
\end{Lemma}
\begin{proof}
This is a fairly dull computation.
\end{proof}
\medskip
Now suppose $G$ is a different lift functor with $F\leq G$.
Then by \fullRef{P.1.5.2}, the basis
$c_\ast\colon F_\pi\ \to\ C_\ast(F)$ goes to $c_\ast\colon
F_\pi\ \to\ C_\ast(G)$ under $\otimes_{\Z\pi_1(F)}\Z\pi_1(G)$.
Let $c^\prime_\ast\colon F_\rho\ \to\  C_\ast(G)$ be the usual
basis.
Then $\pi\subseteq\rho$, and $F_\pi\ \to\ F_\rho
\ \RA{\ c^\prime_\ast\ }\ 
C_\ast(G)$ is just $c_\ast$.
The inclusion $F_\pi\ \to\ F_\rho$ lies in the image of ${\mathcal P}_T$
in ${\mathcal P}_R$, so $[c_\ast/c^\prime_\ast]=0\in \Bar K_1\bigl(
\Z\pi_1(K,G,f)\bigr)$.
Hence
$i_\ast \tau(K,L\Colon F,\Lambda) - \tau(K,L\Colon G,\Lambda)=0$ where
$i_\ast\colon \Bar K_1\bigl(\Z\pi_1(K,F,f)\bigr)\ \to\ 
\Bar K_1\bigl(\Z\pi_1(K,G,f)\bigr)$.
Therefore we can define $\tau(K,L\Colon \Lambda)\in \Bar K_1(K\Colon f)$.

$\tau(K,L\Colon \Lambda)$ depends strongly on $\Lambda$.
We would like this not to be the case, so we pass to a quotient of 
$\Bar K_1$.

\medskip\begin{xDefinition}
Let $G$ be a tree of groups with associated tree of rings $\Z G$.
The \emph{Whitehead group} of $G$, 
$\wh(G)=\Bar K_1(\Z G)/\bigl(\Delta(G)\bigr)$, 
where $\bigl(\Delta(G)\bigr)$ is the subgroup generated by all
objects of the form $( F^{(1)},[g], F^{(1)})$ where $[g]$ is the
map{--}germ of $F^{(1)}$ to itself induced by any element 
$g\in\Delta(G)$ as follows:
$g$ can be represented by a collection $\{g_p\}$, where 
$g_p\in G_{A(p)}$, with $p\in A(p)$ and $\{ A(p) \}$ cofinal and
locally finite.
Define a partition, $\pi$, of the vertices of $T$ by 
$\pi(A)=\{\ p\in T\ \vert\ A(p)\subseteq A\ \}$.
$\pi$ is seen to be a partition and $\pi\subseteq \rho$, 
the standard partition.
Define a map $g \colon F_\pi\ \to\ F_\pi$ by
$g_A\colon (F_\pi)_A\ \to\ (F_\pi)_A$ takes $e_p$ to
$e_p\cdot f_{A_p\>A}(g_p)$ where 
$f_{A B}\colon (\Z G)_A\ \to\ (\Z G)_B$.
It is not hard to show this is a well{--}defined map{--}germ.
What we have actually done is to construct a homomorphism
$\Delta(G)\ \to\ \Bar K_1(\Z G)$ defined by 
$g\mapsto \bigl(F^{(1)},[g], F^{(1)}\bigr)$.
Be definition
$\Delta(G)\ \to\ \Bar K_1(\Z G)\ \to\ \wh(G)\ \to0$ is exact.
\end{xDefinition}

Given a homomorphism $f\colon G\ \to\ H$ between two trees of groups,
we clearly get a commutative square
\[\begin{matrix}%
\Delta(G)&\to&\Bar K_1(\Z G)\\
\downarrow&&\downarrow\\
\Delta(H)&\to&\Bar K_1(\Z H)\hbox to 0pt{$\quad ,$\hss}\\
\end{matrix}\]
so we get a homomorphism $\wh(G)\ \to\ \wh(H)$.

\medskip
\BEGIN{L.1.5.15}
Let $f\colon G\ \to\ H$ be a map between two trees of groups for
which $\Delta(f)$ is an isomorphism.
Then $\wh(G)\ \to\ \wh(H)$ is an isomorphism.
\end{Lemma}
\begin{proof}
$\Delta(f)\colon \Delta(G)\ \to\ \Delta(H)$ is also an isomorphism,
so apply \fullRef{L.1.5.12} 
and the $5${--}lemma.
\end{proof}
\medskip
We can now define $\wh(X\Colon f)$ as
$\displaystyle\lim_{
\Atop{\longrightarrow}{F\in{\mathcal L}(f)}
}\ 
\wh\bigl(\Z\pi_1(X,F,f)\bigr)$.
\bigskip
\BEGIN{P.1.5.4}
Let $(K,L)$ be a pair of finite dimensional \hCWx\  complexes with $L$
a proper deformation retract of $K$.
Then, if $\Lambda$ and $\Lambda^\prime$ are two choices of
paths, $\tau(K,L\Colon \Lambda)=\tau(K,L\Colon \Lambda^\prime)$ in $\wh(X\Colon f)$,
Hence we can then define $\tau(K,L)\in\wh(X\Colon f)$.
\end{Proposition}
\begin{proof}
We can pick any lift functor we like, say $F$.
$C_\ast(K,L\Colon \Lambda, F) = C_\ast(K,L\Colon \Lambda^\prime, F)$,
and each is naturally based.
Let $\pi_\ast$ be the partition associated to $\Lambda$ (see
\fullRef{P.1.5.3} ) and let $\pi^\prime_\ast$ be the
partition associated to $\Lambda^\prime$.
Let $\rho_\ast$ be the partition $\rho_\ast(A)=\{ e\ \vert\ e$
is a $\ast${--}cell in $F(A)$ and the path for $e$ in $\Lambda$ lies in
$F(A)$ and the path for $e$ in $\Lambda^\prime$ also lies in $F(A)\ \}$.
$\rho_\ast=\pi_\ast\ \cap\ \pi_\ast^\prime$.

The basis $F_{\rho_\ast}\ \to\ C_\ast$ is equivalent to the basis
$F_{\pi_\ast}\ \to\ C_\ast$.
Similarly $F_{\rho_\ast}\ \to\ C^\prime_\ast$ is equivalent to the basis
$F_{\pi^\prime_\ast}\ \to\ C^\prime_\ast$.
($C_\ast=C_\ast(\cdots, \Lambda)$; 
$C^\prime_\ast=C_\ast(\cdots, \Lambda^\prime)$.)

$\tau^\prime- \tau = \tau(K,L\Colon \Lambda^\prime) -\tau(K,L\Colon \Lambda)=
\displaystyle\sum_\ast (-1)^\ast [\pi_\ast/\pi_\ast^\prime]$, by
\fullRef{L.1.5.14}.
If we can show $[\pi_\ast/\pi_\ast^\prime]$ is in the image of
$\Delta(\pi_1)$ we are done.
But this is not hard to see ( $\wh(\ )$ was defined by factoring out
such things).
\end{proof}

\medskip
Having defined a torsion, we prove it invariant under subdivision.
We follow Milnor \cite{btwentythree}.
\bigskip
\BEGIN{T.1.5.3}
The torsion $(K,L)$ is invariant under subdivision of the pair $(K,L)$;
$(K,L)$ a finite dimensional \hCWx\  pair.
\end{Theorem}
\begin{proof}
Following Milnor \cite{btwentythree} we prove two lemmas.
\medskip
\BEGIN{L.1.5.16}
Suppose that each component of $K-L$ has compact closure and is simply 
connected.
If $L$ is a proper deformation retract of $K$, then $\tau(K,L)=0$.
\end{Lemma}
\begin{proof}
(Compare Milnor \cite{btwentythree} Lemma 7.2).
Let $f\colon T\ \to\ K$ be a tree.
We wish to find a set of paths $\Lambda$ so that the boundary maps
in $C_\ast(K,L\Colon F,\Lambda)$ come from ${\mathcal P}_T$.
Let $\{ M_i\}$ be the components of $K-L$.
Pick a point $q\in M_i$ and join $\{ q_i\}$ to $T$ by a locally
finite set of paths $\lambda_i$.
Now join each cell in $M_i$ to $q_i$ by a path lying in $M_i$.
Let $\Lambda$ be the set of paths gotten by following the path from the
cell to a $q_i$ and then following the path $\lambda_i$.
Clearly $\Lambda$ is a locally finite set of paths joining the cells of $K-L$
to $T$.

Let $e$ be a cell of $K-L$.
Then if $f$ is a cell of $\partial e$, to compute the coefficient of $f$ 
in $\partial e$ we join the barycenter of $f$ to the barycenter of $e$
by a path in $e$ and look at the resulting loop.
The path from $e$ and the path from $f$ hit the same $q_i$, and since
$\pi_1(M_i,q_i)=0$, the coefficient is $\pm1$, so the boundary maps come
from ${\mathcal P}_T$.
\end{proof}

\medskip
\BEGIN{L.1.5.17}
Suppose that $H_\ast\bigl(C_\ast(K,L\Colon \Lambda)\bigr)$ is not $0$,
but is a free $\Z\pi_1(K)${--}module with a preferred basis.
Suppose each basis element can be represented by a cycle lying over
a single component of $K-L$.
Assume as before that each component of $K-L$ is compact and 
simply connected.
Then $\tau(K,L)=0$.
\end{Lemma}
\begin{proof}
Pick a set of paths as in \fullRef{L.1.5.16} 
so that the boundary maps come from ${\mathcal P}_T$.
Look at a cycle $z$, representing a basis element of $H_\ast$.
What this means is the following. 
Let $C_0\subseteq C_1\subseteq$ \dots be an increasing sequence
of compact subcomplexes with $\cup\ C_i=K$ and $M_i\subseteq C_i$.
Then $z\in H_\ast(\coverFA{K-C_i}, \coverFA{L-C_i})$ for a maximal
$C_i$.
Then $z$ is represented by a cycle lying in some component of 
$\pi^{-1}(M_{i+1})$, where $\pi\colon \coverFA{K-C_i}\ \to\ K-C_i$.
All the lifted cells of $M_{i+1}$ lie in a single component of 
$\pi^{-1}(M_{i+1})$, so let $g\in \pi_1(K-C_i)$ be such that $g z$ also lies
in this distinguished component.

Then the torsion computed with this altered basis is zero since it again
comes from $\wh(T)=0$.
But the new basis for $H_\ast$ is clearly equivalent to the old one 
in $\wh(K)$.
\end{proof}
\medskip
The proof of \fullRef{T.1.5.3} now follows Milnor's proof of
Theorem 7.1 word for word except for a renumbering of 
the requisite lemmas.
\end{proof}

\medskip
\BEGIN{L.1.5.18}
If $M\subseteq L\subseteq K$, where both $L$ and $M$ are 
proper deformation retracts of $K$, then 
$\tau(K,L)=\tau(K,M)+ i_\ast\tau(L,M)$, where 
$i_\ast\colon\wh(L\Colon f)\ \to\ \wh(K\Colon i \circ f)$ is the map induced by
$i\colon L\subseteq K$. (Note the tree must be in $L$.)
\end{Lemma}
\begin{proof}
This is a simple application of \fullRef{T.1.5.1}.
\end{proof}
\medskip
Let $f\colon X\ \to\ Y$ be a proper, cellular map between two 
finite dimensional \hCWx\  complexes. 
Let $M_f$ be the mapping cylinder. 
$Y$ is a proper deformation retract of $M_f$ and we have
\medskip
\BEGIN{L.1.5.19}
$\tau(M_f,Y)=0$ in $\wh(M_f, t)$ where $t\colon T\ \to\ Y$
is a tree for $Y\subseteq M_f$.
\end{Lemma}
\begin{proof}
Word for word Milnor \cite{btwentythree} Lemma 7.5.
\end{proof}
\bigskip\begin{xDefinition}
For any cellular proper homotopy equivalence $f\colon X\ \to\ Y$,
$X$ and $Y$ as above, there is a torsion, $\tau(f)$, defined as follows.
Let $t\colon T\ \to\ Y$ be a tree for $Y$.
Then, as in \fullRef{L.1.5.19}, $t$ is also
a tree for $M_f$ under $T\ \to\ Y\subseteq M_f$.
$\tau(f)=r_\ast\tau(M_f,X)\in\wh(Y\Colon t)$ where 
$r_\ast\colon \wh(M_f\Colon t)\ \to\ \wh(Y\Colon t)$, where $r$ is the retraction.
\end{xDefinition}

Just as in Milnor we have
\BEGIN{L.1.5.20}
If $i\colon L\ \to\ K$ is an inclusion map
$\tau(i)=\tau(K,L)$ if either is defined.
\end{Lemma}
\medskip
\BEGIN{L.1.5.21}
If $f_0$ and $f_1$ are properly homotopic, $\tau(f_0)=\tau(f_1)$.
\end{Lemma}
\medskip
\BEGIN{L.1.5.22}
If $f\colon X\ \to\ Y$ and $g\colon Y\ \to\ Z$ are cellular proper
homotopy equivalences, then
\[\tau(g\ \circ\ f)=\tau(g)+g_\ast\tau(f)\ ,\] where $t\colon T\ \to\ Y$
is a tree for $Y$ and $g_\ast\colon \wh(Y\Colon t)\ \RA{\ \cong\ }\ 
\wh(Z\Colon g \circ t)$.
\end{Lemma}

\medskip\begin{xRemarks}
It follows from \fullRef{L.1.5.21} that we may define
the torsion of any proper homotopy equivalence between finite dimensional
\hCWx\  complexes, since we have a proper cellular approximation
theorem \cite{beleven}.
\end{xRemarks}
\medskip
Now in \cite{bthirtythree}, Siebenmann defined the notion of simple
homotopy type geometrically.
In particular, he got groups $\sieb(X)$ associated to any locally compact
CW complex.
If $X$ is finite dimensional, we can define a map 
$\tau\colon\sieb(X)\ \to\ \wh(X\Colon f)$ by choosing a tree 
$f\colon T\ \to\ X$.
If $g\colon X\ \to\ Y$ is an element of $\sieb(X)$, $g$
goes to $\tau(M_{g^{-1}},Y)$ where $g^{-1}\colon Y\ \to\ X$
is a proper homotopy inverse for $g$.

$\tau$ is additive by \fullRef{L.1.5.22} 
and depends only on the proper homotopy class of $g$ by
\fullRef{L.1.5.21}.
That $\tau$ is well{--}defined reduces therefore to showing that
$g$ a simple homotopy equivalence implies $\tau(g)=0$.
We defer for the proof to Farrell{--}Wagoner \cite{bten},
where it is also proved $\tau$ is an isomorphism.
The inverse for $\tau$ is easy to describe. 
Let $\alpha\in\wh(X\Colon f)$ be an automorphism of $F^{(n)}$ for some $n$.
Wedge $n$ $2${--}spheres to each vertex of the tree.
Attach $3${--}cells by $\alpha$ to get an \hCWx\  complex Y with $Y-X$
$3${--}dimensional. 
Then $i\colon X\subseteq Y$ is an element of $\sieb(X)$ and
$\tau(i)=\alpha$.
Again we defer to \cite{bten} for the proof that this map is well{--}defined.

In \cite{bthirtythree} Siebenmann also constructs an exact sequence
\[0\to\ \wh^\prime \pi_1(X)\ \to\ \sieb(X)\ \to\ 
K_0\pi_1E(X)\ \to\ K_0 \pi_1(X)\ .\]
We have
\[\begin{matrix}%
&&\sieb(X)\\
&\hskip 60pt\hbox to 0pt{$\nearrow$\hss}&&\searrow\\
0\to&\wh^\prime \pi_1(X)&\uplabeledarrow[\Big]{\tau^{-1}}{}&&
K_0\pi_1E(X)&\to&K_0 \pi_1(X)\\
&\hskip 60pt\hbox to 0pt{$\hbox to 0pt{\hss$\scriptstyle\alpha$}\searrow$\hss}&&
\nearrow\hbox to 0pt{$\scriptstyle\beta$\hss}\\
&&\wh(X\Colon f)\\
\end{matrix}\]
commutes.
Farrell and Wagoner describe $\alpha$ and $\beta$ and prove this diagram
commutes.
They show that the bottom row is exact, so $\tau^{-1}$ is an
isomorphism.

Note now that if $g\colon T\ \to\ X$ is another tree for $X$, we
have natural maps \\$\wh(X\Colon f)\ 
{\crB}\ \wh(X\Colon g)$
which take $\tau(X,Y)$ computed with $f$ to $\tau(X,Y)$ computed
with $g$ and vice{--}versa.
This shows $\wh(X\Colon f)$ does not really depend on the choice of tree.
We content ourselves with remarking that the map
$\wh(X\Colon f)\ \to\ \wh(X\Colon g)$ is not easy to describe algebraically.

In \cite{bthirtythree} Siebenmann derives some useful formulas which
we name
\begin{enumerate}
\item[1)] Sum formula
\item[2)] Product formula
\item[3)] Transfer formula
\end{enumerate}

Note if $\pi\colon \coverFA{Y}\ \to\ Y$ is a cover, $\pi$ induces
$\pi^\ast\colon \sieb(Y)\ \to\ \sieb(\coverFA{Y})$.
We are unable to say much about this map algebraically.
The product formula is algebraically describable however.
\medskip
\BEGIN{L.1.5.23}
Let $C_\ast$ be an s{--}based, finite chain complex over the tree of
rings $R$.
Let $D_\ast$ be an s{--}based, finite chain complex on the ring $S$
(the tree of rings over a point).
Then $(C\otimes D)_\ast$ is defined.
If $C_\ast$ is acyclic with torsion $\tau$, $(C\otimes D)_\ast$
is acyclic with torsion $\chi(D)\cdot i_\ast \tau(C)\in \wh(R\times S)$
where $(R\times S)_A= R_A\times S$, and 
$i_\ast\colon \wh(R)\ \to\ \wh(R\times S)$ is the obvious split 
monomorphism.
If $D_\ast$ is acyclic, then so is $(C\otimes D)_\ast$, and if
$\tau(D)=0$, then $\tau(C\otimes D)=0$.
\end{Lemma}
\begin{proof}
The first formula is Siebenmann's product formula and is proved
by induction on the number of cells in $D_\ast$.
The second formula is new, but it is fairly easy.
It basically requires the analysis of maps
$\wh(S)\ \to\ \wh(R\times S)$ of the form 
$D_\ast\ \to\ P\otimes D_\ast$ for $P$ and s{--}based $R${--}module.
These maps are homomorphisms, and so, if $\tau(D_\ast)=0$,
$\tau(P\otimes D_\ast)=0$.
But $\tau\bigl((C\otimes D)_\ast\bigr)=\displaystyle
\sum_k(-1)^k\tau(C_k\otimes D_\ast)$.
(There is evidence for conjecturing that the map
$\wh(S)\ \to\ \wh(R\times S)$ is always $0$\setcounter{footnote}{0}\footnote{If the tree is infinite.}.)
\end{proof}

\medskip
We conclude this section by discussing the notion of duality.
In particularly, we would like a functor 
$\du\colon {\mathcal M}_R\ \to\ {\mathcal M}^\ell_{\mathcal R}$
which generalizes the usual duality $P\ \to\ {\Homx}(P, R)$ in
the compact case.
Up until now, ${\mathcal M}_R$ has denoted without prejudice either the
category of right or left $R${--}modules.
We now fix it to be the category of right $R${--}modules.
${\mathcal M}^\ell_R$ then denotes the category of left $R${--}modules.

Actually, we are really only interested in
$\du\colon {\mathcal P}_R\ \to\ {\mathcal P}^\ell_R$.
Hence we begin by discussing a functor
$\du\colon {\mathcal F}_R\ \to\ {\mathcal F}^\ell_R$, where ${\mathcal F}_R$
is the category of locally{--}finitely generated free right $R${--}modules.
$\du$ will satisfy

\begin{enumerate}
\item[1)] $\du$ is a contravariant, additive, full faithful functor
\item[2)] $\du\du$ is naturally equivalent to the identity.
\end{enumerate}

By this last statement we mean the following.
Given $\du\colon{\mathcal F}_R\ \to\ {\mathcal F}^\ell_R$ there will be another
obvious duality $\du\colon {\mathcal F}^\ell_R\ \to\ {\mathcal F}_R$.
The composition of these two is naturally equivalent to the identity.

We proceed to define $\du$.
If $F_A$ is a free right $R_A${--}module based on the set $A$, there is
also a free left $R_A${--}module based on the same set,
$F_A^\ast$.
$F^\ast_A$ can be described as ${\Homx}^\cmpsup_{R_A}(F_A,R_A)$,
where ${\Homx}^\cmpsup_{R_A}$ is the set of all $R_A${--}linear 
homomorphisms which vanish on all but finitely many generators.
${\Homx}^\cmpsup_{R_A}(F_A,R_A)$, is easily seen to have 
the structure of a left $R_A${--}module.

Let $A\subseteq B$, and let $f\colon R\ \to\ S$ be a ring homomorphism.
Then we have
\[\begin{aligned}%
{\Homx}^\cmpsup_{R_A}(F_A,R_A)\ \to\ 
{\Homx}^\cmpsup_{R_B}(F_A\otimes R_B,R_B)\ \LA{\ {\text{ex}}\ }\\
{\Homx}^\cmpsup_{R_B}(F_B/F_{B-A},R_B)\ \to\ 
{\Homx}^\cmpsup_{R_B}(F_B,R_B)\\
\end{aligned}\]

The map ${\text{ex}}$ is an isomorphism since $0\to\ F_A\otimes R_B\ \to\ 
F_B\ \to\ F_{B-A}\ \to0$ is split exact.
Thus we get a well{--}defined homomorphism 
\[{\Homx}^\cmpsup_{R_A}(F_A,R_A)\ \to\
{\Homx}^\cmpsup_{R_B}(F_B,R_B)\ .\]

Now given $F_\pi$, let $F^\ast_\pi$ be the tree of left modules over
the tree of rings $R$ defined by
$(F^\ast_\pi)_A={\Homx}^\cmpsup_{R_A}(F_{\pi(A)},R_A)$, and use
the map discussed above to define $p_{A B}$.

Given a map $f\colon F_\pi\ \to\ F_\rho$, define
$f^\ast\colon F_\rho^\ast\ \to\ F_\pi^\ast$ by 
\[(f^\ast)_A={\Homx}(f_A) \colon
{\Homx}^\cmpsup_{R_A}(F_{\rho(A)},R_A)\ \to\ 
{\Homx}^\cmpsup_{R_A}(F_{\pi(A)}, R_A)\ .\]
We must check that $(f^\ast)_A$ is defined and that the requisite
diagrams commute.
This last is trivial, so we concentrate on the first objective.
To this end, let $\alpha\in {\Homx}^\cmpsup_{R_A}(F_{\rho(A)},R_A)$.
We must show ${\Homx}(f_A)(\alpha)$ lies in
${\Homx}^\cmpsup_{R_A}(F_{\pi(A)},R_A)\subseteq 
{\Homx}_{R_A}(F_{\pi(A)},R_A)$.
Since $\alpha$ has compact support, $\alpha$ vanishes on the
generators corresponding to a subset $S\subseteq \rho(A)$ with
$\rho(A)-S$ finite.
Hence there is a $B\in{\mathcal C}(T)$ so that $\rho(B)\subseteq S$;
i.e. $\alpha$ vanishes on generators corresponding to $\rho(B)$.
Let $\Bar F_{\pi(B)} = F_{\pi(B)}\otimes_{R_B} R_A$; let 
$\Bar F_{\rho(B)} = F_{\rho(B)}\otimes_{R_B} R_A$; and let
$\Bar f_B=f_b\otimes {\text{id}}$. 
Then
\[\begin{matrix}%
{\Homx}(F_{\rho(A)}, R_A)&\RA{\ {\Homx}(f_A)\ }&
{\Homx}(F_{\pi(A)}, R_A)\\
\downlabeledarrow[\Big]{i}{}&&\downlabeledarrow[\Big]{j}{}\\
{\Homx}(\Bar F_{\rho(B)}, R_A)&\RA{\ {\Homx}(\Bar f_B)\ }&
{\Homx}(\Bar F_{\pi(B)}, R_A)\\
\end{matrix}\]
commutes.
$\alpha$ is in the kernel of $i$, so ${\Homx}(f_A)(\alpha)\in{\kerx}\ j$.
But this means ${\Homx}(f_A)(\alpha)$ has compact support.

There is a natural map $F\ \to\ F^{\ast\ast}$ induced by the natural
inclusion of a module into its double dual.
This map is an isomorphism\setcounter{footnote}{0}\footnote{Since $F$ is finitely-generated and free.} and
\[\begin{matrix}%
F&\RA{\hskip 10pt}&F^{\ast\ast}\\
\downlabeledarrow{f}{}&&\downlabeledarrow{}{f^{\ast\ast}}{}\\
G&\RA{\hskip 10pt}&G^{\ast\ast}\\
\end{matrix}\]
commutes.

$\du$ is clearly contravariant and a functor.
If $\pi\subseteq\rho$, one sees $F_\rho^\ast\ \to\ F^\ast_{\pi}$
is an equivalence.
Hence we can define $\du$ for map{--}germs.
$(f+g)^\ast= f^\ast+g^\ast$ is easy to see, so $\du$ is additive.
Since $\du\du$ is naturally isomorphic to the identity, $\du$ must be both
faithful and full, so 1) is satisfied.

We next define the subcategory on which we wish to define $\du$.
Let $\Bar{\mathcal M}_R$ be the full subcategory of ${\mathcal M}_R$
such that $M\in \Bar{\mathcal M}_R$ \iff\ there exists 
$f\colon F_{\rho}\ \to\ F_{\pi}$ with ${\cokerx\ }f\cong M$.
Note ${\mathcal P}_R\subseteq \Bar{\mathcal M}_R$.
We define $\du \colon \Bar{\mathcal M}_R\to {\mathcal M}^\ell_R$ by
$M^\ast={\kerx}(f^\ast)$.

Given $M$, $N\in \Bar{\mathcal M}_R$, a map $g\colon M\ \to\ N$, and
resolutions $F_\rho\ \to\ F_\pi\ \to\ M\ \to0$ and
$F_\alpha\ \to\ F_\beta\ \to\ M\ \to0$, note that we can compare
resolutions.
That is, we can find $h$ and $f$ so that
\[
\begin{matrix}%
&F_\rho&\RA{\ f\ }&F_{\alpha}\\
&\downarrow&&\downarrow\\
\hbox to 0pt{\hss A)\hskip 20pt}&F_\pi&\RA{\ h\ }&F_{\beta}\\
&\downarrow&&\downarrow\\
&M&\RA{\ g\ }&N\\
&\downarrow&&\downarrow\\
&0&&0\\
\end{matrix}\] commutes.
\par\noindent Define $g^\ast\colon N^\ast\ \to\ M^\ast$ by
\[
\begin{matrix}%
&0&&0\\
&\downarrow&&\downarrow\\
&N^\ast&\RA{\ g^\ast\ }&M^\ast\\
&\downarrow&&\downarrow\\
\hbox to 0pt{\hss B)\hskip 20pt}&
F^\ast_\beta&\RA{\ h^\ast\ }&F^\ast_{\pi}\\
&\downarrow&&\downarrow\\
&F^\ast_\alpha&\RA{\ f^\ast\ }&F^\ast_{\rho}\ .\\
\end{matrix}\]

We first note that the definition of $g^\ast$ does not depend on $h$
and $f$, for if we pick $h_1$ and $f_1$ such that A) commutes, there 
is a commutative triangle
\[\begin{matrix}%
&&F_\alpha\\
&\hbox to 0pt{\hss $\scriptstyle p$}\nearrow&\downarrow\\\noalign{\vskip 6pt}
F_\pi&\hbox to 0pt{\hss$\RA{h-h_1}$\hss}& F_\beta
\hbox to 0pt{\quad .\hss}\\
\end{matrix}\]
Dualizing, we get
\[\begin{matrix}%
F^\ast_\beta&\hbox to 0pt{\hss$\RA{h^\ast-h^\ast_1}$\hss}
& F^\ast_\pi\\\noalign{\vskip6pt}
\downarrow&\nearrow\hbox to 0pt{$\scriptstyle p^\ast$\hss}\\
F^\ast_\alpha&&\hbox to 0pt{\quad .\hss}\\
\end{matrix}\]
Now this triangle shows that the map we get from $f_1$ and $h_1$
is the same as we got from $f$, $h$.

To show $M^\ast$ does not depend on the resolution is now done
by comparing two resolutions and noting $({\text{id}})^\ast={\text{id}}$.

Unfortunately, $(M^\ast)^\ast$ may not even be defined, so we have
little hope of proving a result like 2).
One useful result that we can get however is

\medskip
\BEGIN{L.1.5.24}
Let $f\colon P\ \to\ M$ be an epimorphism with $M\in\Bar{\mathcal M}_R$
and $P\in{\mathcal P}_R$.
Then $f^\ast\colon M^\ast\ \to\ P^\ast$ is a monomorphism.
\end{Lemma}
\begin{proof}
The proof is easy.
\end{proof}
\medskip
If we restrict ourselves to ${\mathcal P}_R$, we can get 1) and 2) to hold.
It is easy to see $P^\ast\in{\mathcal P}^\ell_R$ for $P\in{\mathcal P}_R$.
Now the equation $(P\oplus Q)^\ast=P^\ast\oplus Q^\ast$ is easily
seen since direct sum preserves kernels.
Thus $(P\oplus Q)^{\ast\ast}=P^{\ast\ast}\oplus Q^{\ast\ast}$, so it is
not hard to see $P\ \to\ P^{\ast\ast}$ must be an isomorphism since
if $P$ is free the result is known.
Lastly, $\du$ is natural, i.e. if $f\colon R\ \to\ S$ is a map,
\lower \baselineskip\hbox{$\begin{matrix}%
\Bar{\mathcal M}_R&\RA{\hskip10pt}&\Bar{\mathcal M}_S\\
\downlabeledarrow{\du}{}&&\downlabeledarrow{\du}{}\\
{\mathcal M}^\ell_R&\RA{\hskip10pt}&{\mathcal M}^\ell_S\\
\end{matrix}$} commutes.
That $\Bar{\mathcal M}_R$ hits $\Bar{\mathcal M}_S$ follows since $\otimes$ 
is right exact.
\medskip\begin{xDefinition}
Let $\{ M_i, \partial_i\}$ be a chain complex with
$M_i\in\Bar{\mathcal M}_R$.
Then $\{ M^\ast_i,\partial^\ast_i\}$ is also a chain complex.
The \emph{cohomology} of $\{ M_i, \partial_i\}$ is defined
as the homology of $\{M^\ast_i, \partial^\ast_i\}$.
\end{xDefinition}
\medskip
\BEGIN{P.1.5.5}
Let $(X,Y)$ be an \hCWx\  pair; let $F$ be a lift functor; and let $\Lambda$
be a set of paths.
Then $\{C_\ast(X,Y\Colon F,\Lambda), \partial_\ast\}$ is a chain complex 
as we saw.
Its dual is 
$\{C^\ast(X,Y\Colon F,\Lambda), \delta^\ast\}$.
Hence the cohomology of a pair is the same as the cohomology of
its chain complex.
\end{Proposition}
\begin{proof} Easy.
\end{proof}

\medskip
Notice that our geometric chain complexes lie in ${\mathcal P}_R$.
For such complexes we can prove
\bigskip
\BEGIN{T.1.5.4}
Let $\{ P_r, \partial_r\}$ be a finite chain complex in ${\mathcal P}_R$.
$H_k(P)=0$ for $k\leq n$ \iff\ there exist maps
$D_r\colon P_r\ \to\ P_{r+1}$ for $r\leq n$ with 
$D_{r-1}\partial_r + \partial_{r+1}D_r={\text{id}}_{P_r}$.
\end{Theorem}
\begin{proof} Standard.
\end{proof}

\medskip
\BEGIN{C.1.5.4.1}
(Universal Coefficients).
With $\{ P_r,\partial_r\}$ as above, $H_k(P)=0$ for $k\leq n$
implies $H^k(P)=0$ for $k\leq n$.
$H^k(P)=0$ for $k\geq n$
implies $H_k(P)=0$ for $k\geq n$.
\end{Corollary}
\begin{proof} Standard.
\end{proof}
\bigskip
Now suppose $\{P_r,\partial_r\}$ is a chain complex in ${\mathcal P}_R$.
Then ${\cokerx}\ \partial_{r+1}\in\Bar{\mathcal M}_R$.
By \fullRef{L.1.5.24},
${\kerx\ }\delta^{r}=({\cokerx\ }\partial_{r+1})^\ast$.
Now

\topD{12}{$\begin{matrix}%
&&P_{r-1}\\
&\hbox to 0pt{\hss$\scriptstyle\partial _r$}\nearrow&\uparrow\\\noalign{\vskip 6pt}
P_r&\to&{\cokerx\ }\partial_{r+1}\\
&&\uparrow\\
&&H_r(P)\\
&&\uparrow\\
&&0\\\end{matrix}
$}{5}
commutes and is exact.
If $H_r(P)\in\Bar{\mathcal M}_R$, applying duality to this 
diagram yields
\\
\topD{12}{$\begin{matrix}%
P^\ast_{r-1}&\RA{\ \alpha\ }&{\kerx}(\delta^r)&\RA{\ \beta\ }&
\bigl(H_r(P)\bigr)^\ast\\
&\hbox to 0pt{\hss$\scriptstyle \delta^{r-1}$}\searrow&\downarrow\\\noalign{\vskip 3pt}
&&P^\ast_r\\
\end{matrix}$\ .}{5}
By definition,
${\cokerx\ }\alpha = H^r(P)$. $\beta\ \circ\ \alpha=0$, so there
is a unique, natural map
$H^r(P)\ \to\ \bigl(H_r(P)\bigr)^\ast$.
\medskip
\BEGIN{C.1.5.4.2}
With $\{P_r,\partial_r\}$ as above, if $H_k(P)=0$ for $k<n$,
$H_n(P)\in\Bar{\mathcal M}_R$.
If $H_n(P)\in{\mathcal P}_R$, the natural map
$H^n(P)\ \to\ \bigl(H_n(P)\bigr)^\ast$ is an isomorphism.
\end{Corollary}
\begin{proof}
By induction, one shows $Z_n\in{\mathcal P}_R$, and since
$P_{n+1}\ \RA{\ \partial_{n+1}\ }\ Z_n\ \to\ H_n(P)\ \to0$ is exact,
it is not hard to see $H_n(P)\in\Bar{\mathcal M}_R$.
If $H_n(P)\in{\mathcal P}_R$,
$0\to\ \bigl(H_n(P)\bigr)^\ast\ \to\ Z^\ast_n\ \to\ P^\ast_{n+1}$
is exact, so $H^n(P)\cong \bigl(H_n(P)\bigr)^\ast$.
\end{proof}
\bigskip
\BEGIN{T.1.5.5}
With $\{ P_r,\partial_r\}$ as above, suppose $H_k(P)=0$ for $k<n$ 
and $H^k(P)=0$ for $k>n$.
Then $H_n(P)\in{\mathcal P}_R$ and the natural map
$H^n(P)\ \to\ \bigl(H_n(P)\bigr)^\ast$ is an isomorphism.
In $K_0(R)$,
$\bigl[H_n(P)\bigr]=(-1)^n\chi(P)$, where
$\chi(P)\in K_0(R)$ is $\displaystyle\sum_r (-1)^r[P_r]$.
\end{Theorem}

\begin{proof}
Since $H_k(P)=0$ for $k<n$, the sequence
$\cdots\to\ P_{n+1}\ \to\ P_n\ \to\ P_{n-1}\ \to\cdots$
splits up as
\[\begin{aligned}
\cdots\to\ P_{n+1}\ \to\ Z_n&\ \to 0\\
0\to\ Z_n&\ \to\ P_n\ \to\ P_{n-1}\ \to\cdots\\\end{aligned}\]
The second sequence is exact, and
$\cdots\to\ P_{n+1}\ \to\ Z_n\ \to H_n(P)\ \to 0$
is exact since $H_k(P)=0$ for $k>n$ by 
\fullRef{C.1.5.4.1}.

By \fullRef{C.1.5.4.2}, $H_n\in\Bar{\mathcal M}_R$.
Dualizing, we get
$\cdots\leftarrow\ P^\ast_{n+1}\ \leftarrow\ Z^\ast_n\ \leftarrow\ 
(H_n)^\ast\ \leftarrow 0$ is exact by 
\fullRef{C.1.5.4.1} and 
\fullRef{L.1.5.24}.
As in the proof of
\fullRef{T.1.5.4}, we get a chain retraction up to
$D\colon P^\ast_{n+1}\ \to\ Z^\ast_n$.
This shows $(H_n)^\ast\in{\mathcal P}_R$.
But
\[\begin{aligned}%
\cdots\leftarrow\ P^\ast_{n+1}\ \leftarrow\ &Z^\ast_n\ \leftarrow\ 
(H_n)^\ast\ \leftarrow 0\\
0\leftarrow\ &Z^\ast_{n}\ \leftarrow\ P^\ast_n\ \leftarrow\ 
P^\ast_{n-1}\ \leftarrow\cdots\\
\end{aligned}\]
splice together to give the cochain complex.
$H^n\ \to (H_n)^\ast$ is now easily seen to be an isomorphism.

Now $\displaystyle
\sum_{r\geq n+1} (-1)^r[P_r] +(-1)^n[Z_n] + (-1)^{n-1}[H_n]=0$
and $\displaystyle
\sum_{r\leq n} (-1)^r[P_r] +(-1)^{n+1}[Z_n] =0$
in $K_0(R)$ by Bass \cite{bone}, Proposition 4.1, Chapter VIII.
Summing these two equations shows
$\chi(P)+(-1)^{n-1}[H_n]=0$.
\end{proof}

Now let us return and discuss the products we defined in section 4.
We defined two versions of the cap product on the chain level
(see Theorems \shortFullRef{T.1.4.5} and \shortFullRef{T.1.4.6}).
Notice that the maps we defined on $P_\ast(X;A,B)$ and $P^\ast(X;A,B)$
actually come from maps on the tree modules
$C_\ast(X;A,B\Colon \Lambda,F)$ and $C^\ast(X;A,B\Colon \Lambda,F)$.
Thus if $f$ is a cocycle in $C^m(X;A;\Gamma)$, and if $h$ is a
diagonal approximation, \fullRef{T.1.4.5} yields a
chain map
$C_{\ast+m}(X;A,B\Colon \Lambda,F)\ \RA{\ \capf f h\ }\ 
C_{\ast}(X,B\Colon \Lambda,F)$.
Note that in order for this to land in the asserted place, $\Gamma$
pulled up to the universal cover of $X$ must just be ordinary
integer coefficients.

$\capf f h$ dualizes to 
$\cupf f h\colon C^\ast(X,B\Colon \Lambda,F)\ \to\ 
C^{m+\ast}(X\Colon A,B\Colon \Lambda,F)$.
Since we did not define cup products on the chain level, we may
take this as a definition.
Nevertheless we assert that on homology $f\cup_h$ induces 
the cup product of \fullRef{T.1.4.1}.
This follows from the duality relations we wrote down between
ordinary cohomology and homology (see the discussion around
the universal coefficient theorems in section 1)
\setcounter{footnote}{0}\footnote{See page \pageref{Universalcoefficientx}.}.

Now one easily sees $\du$ induces a map $\wh(\du)\colon
\wh(R)\ \to\ \wh^\ell(R)$, where $\wh^\ell(R)$ is the group formed
from left modules.
If $\capf f{}$ (or $\cupf f{}$) is a chain equivalence, we can compare 
$\tau(\cupf f{})$ and $\tau(\capf f{})$.
We get $\wh(\du)\bigl(\tau(\capf f{})=(-1)^m\tau(\cupf f{})$
by definition.

Next we study the cap product of \fullRef{T.1.4.6}.\footnote
{We did not get this quite right in the original so there are some changes
here.}\\
A cycle $c\in C^\locf_m(X;A,B;\Gamma)$ yields maps
$C^\ast(X,A\Colon \Lambda,F)\ \to\ C_{m-\ast}(X,B\Colon \Lambda,F)$.
$C^\ast$ is a left module while $C_\ast$ is a right module, so
$\capf c h$ is not a map of tree modules.
If $\Gamma$ has all its groups isomorphic to $\Z$, which it must to
yield the asserted product, we get a homomorphism
$\wo\colon \Gamma(X)\ \to\ \cy2={\Autx}(\Z)$ 
given by the local system.
We can make $C^\ast$ into a right module (or $C_\ast$ into a left module)
by defining $m_A\cdot a=\Bar{a}\cdot m_A$, where $m_A\in(C^\ast)_A$,
$a\in(\Z\pi_1)_A$ and $\Bar{\barnone}$ \ 
is the involution on $(\Z\pi_1)_A$ induced by
$g\in (\pi_1)_A$ goes to $\wo(g)\cdot g^{-1}$, where $\wo(g)\in
\cy2=\{1,-1\}$ is the image of $g$ under the composition
$(\pi_1)_A\ \to\ \pi_1(X)\ \RA{\ \wo\ }\ \cy2$.
With this right module structure, $\capf{c}{h}$ is a right module map.
It is not however the case that $\capf{c}{h}$ is a chain map.
The requisite diagrams commute up to sign, but they only commute in
half the dimensions.
To overcome this annoyance, alter the boundary maps in 
$C^\ast(X,A\Colon \Lambda,F)$ to be 
$\delta^\ast_{(m)}=(-1)^{\ast+m}\delta^\ast$ where $\delta^\ast$ are
the duals of the $\partial$ boundary maps.
Let $\twistedC{\ast}{\wo}{m}$ be 
$C^\ast$ with our right module structure
and boundary maps $\delta^\ast_{(m)}$.

Given any finite, projective chain or cochain complex, 
$\{ P_\ast, \partial_\ast\}$ ( or $\{ P^\ast, \delta^\ast\}$ ), 
we can get a new complex $\{ P_\ast, (-1)^{\ast+m}\partial_\ast\}$.
There are evident chain isomorphisms among the three complexes
and these isomorphisms are simple (even measured in $\Bar{K}_1$).
Given a complex, its $(\wo,m)${--}dual is formed by taking the
dual modules, converting them to modules of the same sideness as
the original using $\wo$ 
and altering the dual boundary maps by $(-1)^{\ast+m}$.
The double $(\wo,m)${--}dual of a complex is chain isomorphic to the
original complex. 
The map on a particular module is just the isomorphism
$P\ \to\ P^{\ast\ast}$.

Now $\capf ch\colon \twistedC{\ast}{\wo}{m}(X,A)\ \to\ 
C_{m-\ast}(X,B)$
is a chain map.
If we $(\wo,m)${--}dualize, we get a map
$(\capf c h)^\ast\colon \twistedC{m-\ast}{\wo}{m}(X,B)\ \to\ 
\bigl(\twistedC{\ast}{\wo}{m}(X,A)\bigr)^{(\wo,m){-}{\text{dualized}}}$.
$\twistedC{\ast}{\wo}{m}$ $(\wo,m)${--}dualized is just 
$C_{\ast}$, and $(\capf c h)^\ast=\capf c h$.

The involution is seen to induce an isomorphism $\wh^\ell(G)\ \to\ \wh(G)$,
and the composition
$\wh(G)\ \RA{\ \wh(\du)\ }\ \wh^\ell(G)\ \to\ \wh(G)$ is the map
induced by $\Z G\ \to\ \Z G$ 
via $\Bar{\barnone}$ (it is not hard to see that 
this map induces a map on the Whitehead group level). 
We will denote the map on $\wh(G)$ also by 
$\Bar{\barnone}$.

If $\capf h c$ is a chain isomorphism, either from
$\twistedC{\ast}{\wo}{m}(X,A)\ \to\ C_{m-\ast}(X,B)$ or from 
$\twistedC{m-\ast}{\wo}{m}(X,B)\ \to\ C_{\ast}(X,A)$, we can
compare the two torsions.
We get the confusing equation $\tau(\capf c h)=
(-1)^m\Bar{\tau(\capf c h) }$ where despite their similar appearance,
the two $\capf c h$'s are not the same (which is which is irrelevant).

We conclude by recording a notational convention.
We will sometimes have a map on homology such as $\cap c{}\colon
\Delta^\ast(M)\ \to\ \Delta_\ast(M)$.
If this map is an homology isomorphism we will often speak of the torsion
of $\capf c{}$ (or $\cupf f{}$, etc.).
By this we mean that there is a chain map, possibly after
twisting the cochain complex, (the chain map and twists
will be clear from context)
and these maps on the chain level are equivalences.
Note that by the usual nonsense, the torsions of these product maps
do not depend on a choice of cycle (or cocycle) within 
the homology (cohomology) class.
Nor do they depend on lift functor or choice of paths.
They are dependent on the tree at this stage of our discussion, but
this too is largely fictitious. 
A better proof of independence is given at the end of section 6.
Especially relevant for this last discussion are 
\fullRef{T.2.1.2} and 
the discussion of the Thom isomorphism theorem
in the appendix to Chapter 2.

\section{The realization of chain complexes}
\newHead{I.6}

In \cite{bthirtyseven} and \cite{bthirtyeight}, Wall discussed 
the problem of constructing a CW complex whose chain complex
corresponds to a given chain complex.
We discuss this same problem for locally compact CW complexes.
Throughout this section, complex will mean a finite dimensional, 
locally compact CW complex.

If we have a chain complex $A_\ast$, there are many conditions it
must satisfy if it is to be the chain complex of a complex.
Like Wall \cite{bthirtyeight} we are unable to find an algebraic
description of these conditions in low dimensions.
We escape the dilemma in much the same way.
\medskip

\begin{xDefinition}
A \emph{geometric chain complex} is a positive, finite, 
chain complex $A_\ast$ together with a $2${--}complex $K$,
a tree $f\colon T\ \to\ K$, and a lift functor $F\in{\mathcal L}(f)$
such that
\begin{enumerate}
\item[1)] each $A_k$ is a locally{--}finitely generated free
$\Z\pi_1(K,F,f)${--}module
\item[2)] each $\partial_k\colon A_k\ \to\ A_{k-1}$ is a map
(not a map{--}germ)
\item[3)] in dimensions $\leq 2$, $C_\ast(K\Colon F)=A_\ast$.
\end{enumerate}
\end{xDefinition}

For 3) to make sense, we must define equality for two free tree modules.
If $A$ is free and based on $(S,\pi)$ and if $B$ is free and based on 
$(R,\rho)$, $A=B$ \iff\ there exists a $1${--}$1$ map $\alpha\colon
S\leftrightarrow R$ such that $\alpha\ \circ\ \pi$ is equivalent to $\rho$.
One easily checks that this is an equivalence relation.

Notice that if $A_\ast$ is going to be the chain complex of some complex,
then all the above conditions are necessary.

Given two geometric chain complexes $A_\ast$ and $B_\ast$,
a map $f_\ast\colon A_\ast\ \to\ B_\ast$ is a map (not a germ) on
each $A_k$ and $\partial_k f_k = f_{k-1}\partial_k$ as maps.
\begin{xDefinition}
A map $f_\ast\colon A_\ast\ \to\ B_\ast$ between two geometric
chain complexes is \emph{admissible} provided
\begin{enumerate}
\item[1)] if $L$ is the $2${--}complex for $B_\ast$, $L=K$ wedged
with some $2${--}spheres in a locally finite fashion
\item[2)] $f_0$ and $f_1$ are the identity
\item[3)] $f_2$ is the identity on the $2${--}cells of $K$ and takes
any $2${--}sphere to its wedge point. 
(The tree for $L$ is just the tree for $K$.
The lift functor for $L$ is just $g^{-1}\bigl($ lift functor for $K \bigr)$,
where $g\colon K\ \to\ L$ is the collapse map.)
\end{enumerate}
\end{xDefinition}

\begin{xRemarks}
It seems unlikely that we really need such strong conditions on a map 
before we could handle it, but in our own constructions we usually get this, 
and these assumptions save us much trouble.
\end{xRemarks}
\medskip
The chief geometric construction is the following.
\medskip
\BEGIN{T.1.6.1}
Let $X$ be a connected complex. 
Let $A_\ast$ be a geometric chain complex with an admissible map
$f_\ast\colon A_\ast\to C_\ast(X)$ which is an equivalence.
Then we can construct a complex $Z$ and a proper, cellular map
$g\colon Z\ \to\ X$ so that $C_\ast(Z)=A_\ast$ and

\hskip2.5in$\xymatrix@C4pt{
A_\ast\ar[rr]^-{f_\ast}\ar@2@{-}[rd]&&C_\ast(X)\\
&C_\ast(Z)\ar[ru]_-{g_\ast}\\
}$

\noindent
commutes.
$g$ is a proper homotopy equivalence.
\end{Theorem}

\begin{proof}
We construct $Z$ skeleton by skeleton. 
Since $f_\ast$ is admissible, $Z^2=X^2$ wedge $2${--}spheres.
$g_2\colon Z^2\ \to\ X$ is just the collapse map onto $X^2$.
To induct, assume we have an $r${--}dimensional complex
$Z^r$ and $g_r\colon Z^r\ \to\ X$ so that
$C_\ast(Z_r)=A_\ast$ in dimensions $\leq r$ and
$(g_r)_\ast=f_\ast$ in these dimensions.
If we can show how to get $Z^{r+1}$ and $g_{r+1}$ we are
done since $A_\ast$ is finite.

Now $A_{r+1}$ is free, so pick generators $\{e_i\}$.
We have a map $\partial\colon A_{r+1}\ \to\ A_r$
and $C_r(Z^r)=A_r$.
Hence each $\partial e_i$ is an $r${--}chain in $Z^r$.
We will show that these $r${--}chains are locally finite and spherical 
(i.e. there is a locally finite collection of $r${--}spheres, and, after
subdivision, cellular maps $\cup S^r_i\ \to\ Z^r$ such that 
$\partial e_i$ is homologous to $S^r_i$, and, if $h_i$ is an
$(r+1)$-chain giving the homology, the $\{h_i\}$ 
may be picked to be locally finite.)
We will then attach cells by these spheres and extend the map.

Let us now proceed more carefully.
For each $i$, $\partial e_i\in A_r$ and $\partial e_i\in (A_r)_{W_i}$
for some $W_i\in{\mathcal C}(T)$ with $\{ W_i\}$ cofinal in the subcategory 
of ${\mathcal C}(T)$ consisting of all $A$ such that $\partial e_i\in (A_r)_A$.
Since $C_r(Z^r)=A_r$, $\partial e_i=c_i\in \bigl(C_r(Z^r)\bigr)_{B_i}$
for some $B_i\in {\mathcal C}(T)$ with $B_i\leq W_i$ ( we write
$B_i\leq W_i$ provided $B_i\subseteq W_i$  and $\{ B_i\}$ is cofinal in the
subcategory of all $A\in{\mathcal C}(T)$ for which $e_i\in (A_r)_A$.
$c_i$ is now a real geometric chain.
$\partial c_i=0$ since $\partial$ is actually a map.
Let $[c_i]$ be the homology class of $c_i$ in 
$H_r\bigl(\coverFA{F_r(B_i)}\bigr)$, where $F_r$ 
is the lift functor for $Z^r$.
Now $g_\ast[c_i]=0$ in $H_r\bigl(\coverFA{F(U_i)}\bigr)$, where $F$ is
the lift functor for $X$ and $U_i\leq B_i$.

Hence there is an $f_i\in H_{r+1}\bigl(\coverFA{g_r}\colon
 \coverFA{F_r(U_i)}\ \to\ 
\coverFA{F(U_i)}\bigr)$ with $f_i\mapsto [c_i]$.
But $g_r\colon Z^r\ \to\ X$ is properly $r${--}connected (it induces an
isomorphism of $\Delta(\quad\Colon \pi_1)$'s and $H^0_{{\text{end}}}$'s 
by assumption, so it is always 1-1/2{--}connected.
Hence the universal covering functor for $X$ 
is a universal covering functor for $Z^r$, so $\Delta(\quad\Colon \pi_k)=0$ \iff\ 
$\Delta(M_{g_r},X\Colon H_k,\coverFA{}\ )=0$ for $k\leq r$ by the
Hurewicz theorem.
But $\Delta(M_{g_r},X\Colon H_k,\coverFA{}\ )=0$ for $k\leq r$
\iff\ $\Delta\bigl(H_k(g_r)\bigr)$ is an isomorphism for $k<r$ and
an epimorphism for $k=r$, which it is.)
Hence the Hurewicz theorem gives us elements
$s_i\in \pi_{r+1}\bigl( g_r\colon F_r(V_i)\ \to\ F(V_i)\bigr)$
where $V_i\leq U_i$ and $s_i$ hits the image of $f_i$ in
$H_{r+1}\bigl( \coverFA{g}_r\colon \coverFA{F_r(V_i)}\ \to\ 
\coverFA{F(V_i)}\bigr)$ under the Hurewicz map.

Let $Z^{r+1} = Z^r\ \cup$ a collection of $(r+1)${--}cells, $\{e_i\}$
attached by $s_i$.
$g_{r+1}\colon Z^{r+1}\ \to\ X$ is $g_r$ on $Z^r$.
Since $g_r\ \circ\ s_i\colon S^r\ \to\ Z^r\ \to\ X$ are properly
null homotopic, choose a locally finite collection $\{Q_i\}$
of null homotopies of $g_r\ \circ\ s_i$ to zero in $F(V_i)$.
$g_{r+1}\colon Z^{r+1}\ \to\ X$ is defined by $Q_i$ on each $e_i$.
$g_{r+1}$ is obviously still proper.
$C_\ast(Z^r)\ \to\ C_\ast(Z^{r+1})$ induces 
an isomorphism for $\ast\leq r$.
$C_{r+1}(Z^{r+1})= A_{r+1}$ by taking the cell $e_i$ to the generator
$e_i$.
$F_{r+1}(B)= F_r(B)\ \cup\ ($ all cells $e_i$ for which the generator
$e_i$ lies in $B$ less those for which $g_{r+1}(e_i)\not\subseteq
F(B)$ ).

Then $g^{-1}_{r+1}\bigl(F(B)\bigr)\supseteq F_{r+1}(B)$.
Notice that if a cell $e$ does not attach totally in $F_r(B)$,
$g_{r+1}(e)\not\subseteq F(B)$, so $F_{r+1}(B)$ is a subcomplex.
$F_{r+1}(B)$ is cofinal in $B$, so $F_{r+1}$ is a lift functor.

Look at the chain map
$(g_{r+1})_\ast\colon C_{r+1}(Z^{r+1})\ \to\ C_{r+1}(X)$.
$e_i$ as a cell goes under $(g_{r+1})_\ast$ to the same element in
$\bigl(C_{r+1}(X)\bigr)_B$ as the generator $e_i$ does under $f_\ast$
for all $B\in{\mathcal C}(T)$ such that $e_i$ is a cell in $F_{r+1}(B)$.
Hence
\topD{12}{$\begin{matrix}%
A_{r+1}&\hbox to 0pt{\hss$\RA{\ f_{r+1}\ }$\hss}&\ C_{r+1}(X)\\
\parallel&\hskip 10pt\hbox to 0pt{$\nearrow\hbox to 0pt{$\scriptstyle g_{r+1}$\hss}$\hss}\\
\noalign{\vskip 4pt}
\hbox to 30pt{$C_{r+1}(Z^{r+1})$\hss}\\
\end{matrix}$}{4}
commutes.
\end{proof}

\begin{xDefinition}
A \emph{relative geometric chain complex} is a triple 
$(A_\ast, K, L)$ consisting
of a finite, positive chain complex $A_\ast$ 
and a pair of complexes $(K,L)$.
Understood is a tree and a lift functor.
Then each $A_k$ is a locally{--}finitely generated free 
$\Z\pi_1(K)$ tree module; each $\partial_k$ is a map; and in
dimensions $\leq2$, $A_\ast = C_\ast(K,L)$.

An \emph{admissible map} $f_\ast$ from
$(A_\ast, K, L)$ to $(B_\ast, K^\prime, L)$
is a chain map, not a germ, 
$f_\ast\colon A_\ast\ \to\ B_\ast$ and $K=K^\prime$ wedge a locally
finite collection of $2${--}spheres.
$f_0$ and $f_1$ are the identity, and $f_2$ is the map induced by the
collapse $K\ \to\ K^\prime$.
\end{xDefinition}
\medskip
\BEGIN{C.1.6.1.1}
Let $(X,Y)$ be a pair of complexes, $X$ connected.
Let $A_\ast$ be a relative geometric chain complex with an admissible 
map $f_\ast\colon A_\ast\ \to\ C_\ast(X,Y)$
which is an equivalence.
Then we can construct a complex $Z$ with $Y$ as a subcomplex and a 
proper cellular map $g\colon Z\ \to\ X$ which is the identity on $Y$
such that $C_\ast(Z,Y)=A_\ast$ and
\topD{12}{$\begin{matrix}%
A_{\ast}&\hbox to 0pt{\hss$\RA{\ f_{\ast}\ }$\hss}&\ C_{\ast}(X,Y)\\
\parallel&\hskip 10pt\hbox to 0pt{$\nearrow\hbox{$\scriptstyle g_{\ast}$\hss}$\hss}\\
\noalign{\vskip 4pt}
\hbox to 30pt{$C_{\ast}(Z,Y)$\hss}\\
\end{matrix}$}{4}
commutes. 
$g$ is a proper homotopy equivalence of pairs.
\end{Corollary}
\begin{proof}
The proof parallels the proof of \fullRef{T.1.6.1},
except that we must now use Namioka to show our elements
are spherical.
\end{proof}

Now let $f_\ast\colon A_\ast\ \to\ C_\ast(X)$ be an arbitrary
chain equivalence.
As in Wall \cite{bthirtyeight}, we would like to replace $A_\ast$
by an admissible complex with $f_\ast$ admissible while changing
$A_\ast$ as little as possible. 
Look at
\[\begin{matrix}%
\cdots\to&A_3&\RA{\partial}&A_2&\to&A_1&\to& A_0&\to&A_{-1}
\to\cdots\to0\\
&\downlabeledarrow{f_3}{}&&\downlabeledarrow{f_2}{}&&\downarrow
&&\downarrow\\
\cdots\to&C_3&\to&C_2&\to&C_1&\to& A_0&\to0\\
\end{matrix}\]
One might like to try the complex
\[\begin{matrix}%
\cdots\to&A_3&\RA{f_2\ \circ\ \partial}&C_2&\to&C_1&\to& C_0&\to0\\
&\downlabeledarrow{f_3}{}&&\downlabeledarrow{\text{id}}{}&&
\downlabeledarrow{\text{id}}{}&&\downlabeledarrow{\text{id}}{}\\
\cdots\to&C_3&\to&C_2&\to&C_1&\to& C_0&\to0\quad .\\
\end{matrix}\]
The top complex is clearly admissible, but unfortunately the map
is no longer an equivalence.
The cycles in $A_3$ are now bigger with no new boundaries, and
the boundaries in $C_2$ are smaller with no fewer cycles.

Note first that $X$ is not of great importance.
If we replace $X$ by something in its proper homotopy class,
we will not be greatly concerned.
Let $X^\prime$ be $X$ with $2${--}spheres wedged on 
to give a basis for $A_2$ and $3${--}cells attached to kill them.
Then $X^\prime$ has the same simple homotopy type as $X$, 
$C_k(X^\prime)=C_k(X)$ except for $k=2$, $3$, and
$C_k(X^\prime)=C_k(X)\oplus A_2$ for $k=2$, $3$.
Let $f^\prime_k=f_k$, $k\neq 2$, $3$, and let $f^\prime_3=(f_3,\partial)$
and $f^\prime_2=(f_2, {\text{id}})$.
Then $A_\ast\ \RA{\ f^{\prime}_\ast\ }\ C_\ast(X^\prime)$
is still an equivalence and now $f^\prime_2$ is a monomorphism.
Let $A^\prime_\ast$ be the complex
\[\cdots\to\ A_3\ \RA{\ f^\prime_2\ \circ\ \partial\ }\ 
C_2^\prime=C_2(X^\prime)\ \to\ C^\prime_1\ \to\ 
C^\prime_0\ \to0\quad.\]
Then $h_\ast\colon A_\ast\ \to\ C^\prime_\ast$ has homology in only one 
dimension:
$0\to\ H_2(h)\ \to\ H_2(A^\prime)\ \to\ H_2(C^\prime)\ \to0$.
Since $A_\ast\ \to\ A^\prime_\ast$ and since the composition
$A_\ast\ \to\ A^\prime_\ast\ \to\ C^\prime_\ast$ is an
equivalence, $H_2(A^\prime)=H_2(h)\oplus H_2(C^\prime)$.

Now by \fullRef{T.1.5.5}, $h_2(h)$ is s{--}free, provided
we can show $H^k(h)=0$ for $k\geq 3$.
But since we have a chain equivalence $A_\ast\ \to\ C^\prime_\ast$
we have a chain homotopy inverse in each dimension.
We then clearly get a chain homotopy inverse for 
$A^\prime_k\ \to\ C^\prime_k$, $k\geq 4$, and
$h_3\ \circ\ g^\prime_3$ chain homotopic to ${\text{id}}_{{C_3}^\prime}$.
But this implies $H^k(h)=0$, $k\geq 3$.

Since $H_2(h)$ is projective, we get a map 
$\rho=\partial\ \circ\ \rho^\prime$, where 
$\rho^\prime\colon H_2(h)\ \to\ C^\prime_3$ is given as follows.
Both $A_3$ and $C^\prime_3$ map into $C_2=C^\prime_2$, and
$0\to\ {\Imx} A_3\ \to\ {\Imx} C_3\ \to\ H_2(h)\ \to0$
is exact.
Split this map by $\sigma\colon H_2(h)\ \to\ {\Imx} C_3^\prime$
and note $ {\Imx}\ A_3\ \cap\ {\Imx}\ \sigma = \{1\}$.
Now $C^\prime_3\ \to\ {\Imx}\ C^\prime_3\ \to0$ is exact,
so we can lift $\sigma$ to $\rho^\prime\colon H_2(h)\ \to\ C^\prime_3$.
Since $\sigma$ is a monomorphism, note 
$ {\Imx}\ \rho^\prime\ \cap\ {\Imx}\ f_3 = \{1\}$.

Form $A^{\prime\prime}_\ast$ and $h^\prime_\ast$ by
\[\begin{matrix}%
\cdots\to&A_4&\to&A_3\oplus H_2(h)&
\rightlabeledarrow{\hskip 15pt \partial + \rho \hskip 15pt}{}&C^\prime_2&\to&C^\prime_1&\to
&C^\prime_0\to0\\
&\downlabeledarrow[\bigg]{f_4}{}&
&\downlabeledarrow[\bigg]{f_3 + \rho^\prime}{}
&&\downlabeledarrow[\bigg]{\text{id}}{}&&\downlabeledarrow[\bigg]{\text{id}}{}
&&\downlabeledarrow[\bigg]{\text{id}}{}\hfill
\\\noalign{\vskip4pt}
\cdots\to&C^\prime_4&\to&C^\prime_3&
\RA{\hskip 30pt}&
C^\prime_2&\to&C^\prime_1&\to
&C^\prime_0\to0\\
\end{matrix}\]

Note ${\kerx}(\partial+\rho)= ({\kerx}\ \partial, 0)$ since $\rho$
is a monomorphism and if $\rho(x)\in{\Imx}\ \partial$,
$\rho(x)=\{1\}$ as ${\Imx}\ A_3\ \cap\ {\Imx}\ \sigma=
\{1\}$.
Likewise note ${\Imx}(\partial+\rho)={\Imx}\ C^\prime_3$
since ${\Imx}\ A_3\ \oplus\ {\Imx}\ \rho=
{\Imx}\ C^\prime_3$.
Hence $h^\prime_\ast$ is an equivalence.

Note $\rho\colon H_2(h)\ \to\ C_2^\prime$ is a direct summand.
We split $\rho$ as follows.
\[\begin{matrix}%
\cdots\to&A_4&\to& A_3\oplus H_2(h)&
\rightlabeledarrow{\hskip 15pt \partial + \rho \hskip 15pt}{}& C^\prime_2&\to\cdots\\
&\uplabeledarrow[\big]{}{{\text{id}}}&&\uplabeledarrow[\big]{}{\alpha}
&&\uplabeledarrow[\big]{}{f_2}
\\\noalign{\vskip 4pt}
\cdots\to&A_4&\to& A_3&
\hbox to 32pt{\hss$\RA{\hskip 48pt}$}& A_2&\to\cdots\quad,\\
\end{matrix}\]
where $\alpha$ is inclusion on the first factor, commutes.
These maps must define a chain equivalence, so the dual situation
is also an equivalence.

\[\begin{matrix}%
\cdots\leftarrow&A^\ast_4&\leftlabeledarrow{\hskip 10pt \delta^3 \hskip 10pt}{}&
A^\ast_3\oplus \bigl(H_2(h)\bigr)^\ast&
\leftlabeledarrow{\hskip 20pt \partial^\ast+\rho^\ast\hskip 20pt}{}&(C^\prime_2)^\ast&
\leftarrow\cdots\\
&\bigg\downarrow&&\downlabeledarrow[\bigg]{\alpha^\ast}{}&&
\downlabeledarrow[\bigg]{f^\ast_2}{}
\\\noalign{\vskip 4pt}
\cdots\leftarrow&A^\ast_4&
\hbox to 22pt{$\leftlabeledarrow{\hskip 20pt \delta^3_A\hskip 20pt}{}$\hss}&
A^\ast_3&
\hbox to 42pt{\hss$\leftlabeledarrow{\hskip 30pt \delta^2_A\hskip 30pt}{}$}&A_2^\ast&
\leftarrow\cdots\\
\end{matrix}\]

${\kerx}\ \delta^3 = {\kerx}\ (\delta^3)_A \oplus
\bigl(H_2(h)\bigr)^\ast$, and ${\Imx}\ (\partial^\ast+\rho^\ast)=
{\Imx}\ \partial^\ast \oplus {\Imx}\ \rho^\ast$.
\[\begin{aligned}H_3({\text{Top\ complex}}) =&
{\kerx}\ \delta^3 / {\Imx}(\partial^\ast+\rho^\ast)\\
=&\Bigl({\kerx}\ (\delta_3)_A/{\Imx}\ \partial^\ast\Bigr)
\oplus
\Bigl(\bigl(H_2(h)\bigr)^\ast/{\Imx}\ \rho^\ast\Bigr)
\\\noalign{\medskip}
H_3({\text{Bottom\ complex}})=&
{\kerx}\ (\delta^3_A)/{\Imx}\ (\delta^2_A)\\
\end{aligned}\]
$H_3(\alpha)\colon H_3({\text{Top\ complex}})\ \to\ 
H_3({\text{Bottom\ complex}})$ is
\[\begin{matrix}%
{\kerx}\ (\delta^3_A)/{\Imx}\ \partial^\ast&&
{\kerx}\ (\delta^3_A)/{\Imx}\ (\delta^2_A)\\
\noalign{\vskip 6pt}
\oplus&\hbox to 20pt{\hss$\RA{\hskip 60pt}$\hss}&\oplus\\
\noalign{\vskip 6pt}
\bigl(H_2(h)\bigr)^\ast/{\Imx}\ \rho^\ast&&0\\
\end{matrix}\]
Hence, since $H_3(\alpha)$ is an isomorphism,  
${\Imx}\ \partial^\ast = {\Imx}\ (\delta^2_A)$,
and $\rho^\ast\colon (C^\prime_2)^\ast\ \to\ \bigl(H_2(h)\bigr)^\ast$
is onto.
$\bigl(H_2(h)\bigr)^\ast$ is projective so split $\rho^\ast$.
Dualizing splits $\rho\colon H_2(h)\ \to\ C^\prime_2$.

$H_2(h)$ may not be free (it is only s{--}free).
$A_3\oplus H_2(h)$ is often free, but we prefer to keep $A_3$.
Hence form $A^s_\ast$ and $f^s_\ast$ by
\[\begin{matrix}%
\cdots\to&A_4&\to&
A_3\oplus \bigl(H_2(h)_S\oplus F^{(n)}\bigr)&\to&
C^\prime_2\oplus F^{(n)}&\to\cdots\\
&\downlabeledarrow[\Big]{f_4}{}&&
\downlabeledarrow[\Big]{f_3+(\rho^\prime_S+0)}{}&&
\downlabeledarrow[\Big]{{\text{id}}+0}{}\\
\cdots\to&C^\prime_4&\hbox to 0pt{$\RA{\hskip 40pt}$\hss}&
C^\prime_3&
\hbox to 10pt{\hss\hbox to 50pt{\hss$\RA{\hskip 70pt}$}\hss}&
C^\prime_2&\to\cdots\\
\end{matrix}\]
where $S$ is a shift functor so that the map germs $\rho$ 
and $\rho^\prime$ are actual maps.

By wedging on $n$ $2${--}spheres at each vertex of the tree, we
see $A^s_\ast$ and $f^s_\ast$ are admissible.
Notice that exactly the same procedure makes a map 
$f\ast\colon A_\ast\ \to\ C_\ast(X,Y)$ admissible.

In section 3, \fullRef{P.1.3.3} we defined what it
meant by $X$ satisfies $D n$.
We briefly digress to prove
\medskip
\BEGIN{T.1.6.2}
The following are equivalent for $n\geq 2$, $X$ a complex
\begin{enumerate}
\item[1)] $X$ satisfies $D n$
\item[2)] $X$ is properly dominated by an $n${--}complex
\item[3)] $\Delta^k(X\Colon {\text{universal\ covering\ functor}})=0$
for $k>n$.
\end{enumerate}
\end{Theorem}
\begin{proof}
1) implies 2) as $X^n\subseteq X$ is properly $n${--}connected and
hence dominates $X$ if $X$ satisfies $D n$.
2) implies 3) by computing $\Delta^k$ from the cellular chain complex
of the dominating complex.

3) implies 2): Since $\Delta^k(X\Colon \quad)=0$ for $k>n$, by 
\fullRef{T.1.5.4}\ dualized to cohomology, we
get chain retracts

\[0\to\ C_r\ 
\crA\ 
C_{r-1}\ 
\crA\ 
\cdots\cdots\cdots\ \crA\  
C_{n+1}\ 
\crA\ 
 C_n\ 
\crA\ 
\cdots
\]
where $r={\text{dim}}\ X<\infty$.
By an induction argument, ${\Imx}\ \partial_{n+1}$ is s{--}free,
and $C_n={\Imx}\ \partial_{n+1}\oplus A_n$
(dualize everything to get these results in the cochain complex and then 
dualize back).
$A_n$ is s{--}free and
\[\begin{matrix}%
&0&\to&A_n&\to&C_{n-1}&\to\ \cdots\ \to&C_0&\to 0\\
&&&\hbox to 0pt{\hss$\scriptstyle r$}\downarrow&&\downlabeledarrow{{\text{id}}}{}&&
\downlabeledarrow{{\text{id}}}{}\\
\cdots \to&C_{n+1}&\to&C_n&\to&C_{n-1}&\to\ \cdots\ \to&
C_0&\to 0\\
\end{matrix}\] 
gives us an $n${--}dimensional chain complex and a chain equivalence.
$A_n$ is only s{--}free, so form
$0\to\ A_n\oplus F^{(m)}\ \to\ C_{n-1}\oplus F^{(m)}\ \to\ \cdots$
which is now a free complex.
If $n\geq 3$, the complex and the map are clearly admissible, so by
\fullRef{T.1.6.1} we get an $n${--}complex $Y$
and a proper homotopy equivalence $g\colon Y\ \to\ X$, so $X$
satisfies 2).

If $n=2$, $X$ has the proper homotopy type of a $3${--}complex by
the above, so we assume $X$ is a $3${--}complex.
Its chain complex is then 
$0\to\ C_3\ \to\ C_2\ \to\ C_1\ \to\ C_0\ \to0$
with $H^3(C)=0$.
Wedge $2${--}spheres to $X$ at the vertices of the tree to get a
chain complex
$0\to\ C_3\ \to\ C_2\oplus C_3\ \to\ C_1\ \to\ C_0\ \to0$.
Since $H^3(C)=0$, $C_2=C_3\oplus M$.
Let $j\colon C_3 \to\ C_2$ be the inclusion.
Then we have
\[\begin{matrix}%
\hbox to 0pt{\hss{\text{A:}}\hskip 40pt}
&&0\to&C_3&\RA{j}&C_2\oplus C_3&\to&C_1&\to&C_0\to0\\
\noalign{\vskip 4pt}
&&&\Big\uparrow&&\uplabeledarrow[\Big]{}{r}&&\uplabeledarrow[\Big]{}{{\text{id}}}
&&\uplabeledarrow[\Big]{}{{\text{id}}}\\
\noalign{\vskip 6pt}
\hbox to 0pt{\hss{\text{B:}}\hskip 40pt}
&&&0&\to&C_2&\to&C_1&\to&C_0\to0\\
\end{matrix}\]
where $r\colon ( C_3\oplus M)\ \to\ 
 ( C_3\oplus M)\oplus C_3$
is given by $r(x,y)=(0,y,x)$.
This is a chain equivalence between B and A.

Both A and B are the chain complexes for a space (A for $X\vee S^2$'s
and B for the $2${--}skeleton of $X$).
The chain  map is easily realized on the $1${--}skeleton as a map, 
and we show we can find a map
$g\colon X^2\ \to\ X\vee_j S^2_j$ realizing the whole chain map.

Let $\{ e_i\}$ be the two cells of $X^2$.
Their attaching maps determine an element in 
$\Delta( X^1\Colon \pi_1, \coverFA{}\ )$, where this group denotes the 
$\Delta${--}construction applied to the groups
$\pi_1\bigl(p^{-1}(X^1-C),\hat x_i\bigr)$, where 
$p\colon \coverFA{X}\ \to X$ is the projection for the universal cover
of $X$. (i.e. $\coverFA{}$\quad denotes the covering functor
over $X^1$ induced in the above manner form the universal covering
functor on $X$.)
Let $g_1\colon X^1\ \to\ X\vee_j S^2_j$ be the natural inclusion.
As in the proof of \fullRef{T.1.6.1}, the $\{e_i\}$
determine an element of $\Delta(g_1\Colon H_2, \coverFA{}\ )$.
Our two elements agree in $\Delta( X^1\Colon H_1, \coverFA{}\ )$.
The following diagram commutes and the rows are exact

\[\begin{matrix}%
&&\Delta(X^2\Colon \pi_2,\coverFA{}\ )&\zto&\Delta(g_1\Colon \pi_2,\coverFA{}\ )
&\zto&\Delta(X^1\Colon \pi_1,\coverFA{}\ )&\zto&
{\scriptstyle\Delta(X^2\Colon \pi_1,\coverFA{}\ )=}0\\
&&\downlabeledarrow[\Big]{h}{}&&\Big\downarrow&&\Big\downarrow\\
0={\scriptstyle\Delta(X^1\Colon H_2,\coverFA{}\ )}&\zto&
\Delta(X^2\Colon H_2,\coverFA{}\ )&\zto&\Delta(g_1\Colon H_2,\coverFA{}\ )
&\zto&\Delta(X^1\Colon H_1,\coverFA{}\ )\\
\end{matrix}\]
where $\Delta(X^1\Colon  H_1,\coverFA{}\ )$ and 
$\Delta(X^1\Colon H_2,\coverFA{}\ )$ are defined similarly to 
$\Delta(X^1\Colon \pi_1,\coverFA{}\ )$.
$X^1\subseteq X^2$ is properly $1${--}connected, so the subspace
groups are the groups asserted.
$h$ is an isomorphism by the Hurewicz theorem, so a diagram chase 
yields a unique element in $\Delta(g_1\Colon \pi_2,\coverFA{}\ )$
which hits our elements in both $\Delta(X^1\Colon \pi_1,\coverFA{}\ )$
and $\Delta(g_1\Colon H_2,\coverFA{}\ )$.
Use this element to extend the map to 
$g_2\colon X^2\ \to\  X\vee_j S^2_j$.
By our choices, $g_2$ induces an isomorphism of $\Delta(\pi_1)$'s.
Hence $g_2$ is a proper homotopy equivalence.
This 3) implies 2) for $n\geq 2$.

2) implies 1) is trivial.
\end{proof}

\medskip
\BEGIN{C.1.6.2.1}
If $X$ satisfies $D n$ for $n\geq 3$, $X$ has the proper homotopy
type of an $n${-}complex. \mathqed\end{Corollary}
\medskip
Combining our admissibility construction with 
\fullRef{T.1.6.1} gives
\bigskip
\BEGIN{T.1.6.3}
Let $f_\ast\colon A_\ast\ \to\ C_\ast(X)$ be a chain equivalence
for a complex $X$ ($A_\ast$ free, finite and positive).
Then there exists a complex $Y_0$ satisfying $D2$; a complex
$Y\supseteq Y_0$ such that $C_\ast(Y,Y_0)=A_\ast$
in dimensions greater than or equal to $3$; and a proper,
cellular homotopy equivalence $g\colon Y\ \to\ X$ such that
$g_\ast=f_\ast$ in dimensions greater than or equal to $4$.
The torsion of $g$ may have any preassigned value.
\end{Theorem}
\begin{proof}
Make $A_\ast$, $f_\ast$ admissible.
The new complex is
$\cdots\to\ A_4\ \to\ A_3\oplus (?)\ \to\ 
C_2\oplus(?)\ \to\ X_1\ \to\cdots$.
Construct a $Y$ from this complex as in \fullRef{T.1.6.1}.
When we pick a basis for $A_3\oplus (?)$, pick a basis for $A_3$
and one for $(?)$ and use their union.
Then there is a subcomplex $Y_0\subseteq Y$ whose chain complex is
$0\to\ (?)\ \to\ C_2\oplus (?)\ \to\ C_1\ \to\ C_0\ \to0$
The first $(?)$ is $H_2(h)\oplus F^{(m)}$.
It is not hard to show $\Delta^3(Y_0\Colon \coverFA{}\ )=0$, 
so $Y_0$ satisfies $D2$.
The remainder of the theorem is trivial except for the remark about torsion.
But for some $m\geq 0$, we can realize a given torsion by an
automorphism $\alpha\colon F^{(m)}\to\ F^{(m)}$.
Hence by altering the basis in $F^{(m)}$ we can cause our map to
have any desired torsion (we may have to take $m$ bigger,
although in the infinite case $m=1$ will realize all torsions).
\end{proof}

\bigskip
\BEGIN{T.1.6.4}
Let $f_\ast\colon A_\ast\ \to\ C_\ast(X,Z)$ be a chain map for a pair
$(X,Z)$ ($A_\ast$ free, finite and positive).
Then there exist complexes $Y_0$ and $Y$ such that
$Y\supseteq Y_0\supseteq Z$; $C_\ast(Y,Y_0)=A_\ast$
in dimensions greater than or equal to 3; a proper cellular homotopy
equivalence $g\colon Y\ \to\ X$ which is the identity on $Z$ such
that $g_\ast=f_\ast$ in dimensions greater than or equal to $4$.
The torsion of $g$ may have any preassigned value.
\end{Theorem}
\begin{proof} Use \fullRef{C.1.6.1.1}.
\end{proof}
\bigskip
We conclude this chapter by returning briefly to the question of 
the invariance of torsion for chain maps under change of trees.
The natural map $\wh(X\Colon f)\ \to\ \wh(X\Colon g)$ is a homomorphism,
so the property of being a simple chain equivalence is 
independent of the tree.
But now use \fullRef{T.1.6.3} to get a proper
homotopy equivalence $X\ \to\ X_\tau$ with torsion $\tau$.
Suppose given a chain map, say for example, $\capf c{}\colon
\Delta^\ast(X)\ \to\ \Delta_{m-\ast}(X)$ with torsion $-\tau$.
Then the composition $\Delta^\ast(X)\ \RA{\ \capf c{}\ }\ 
\Delta_{m-\ast}(X)\ \to\ \Delta_{m-\ast}(X_\tau)$ is simple, so a change
of trees leaves it simple.
But the second map is a proper homotopy equivalence of spaces,
and so is preserved by our change of tree map.
Hence so must be the torsion of $\capf c{}$.
(Note we are using our convention of writing chain equivalences on
the homology level.)

\chapter{Poincar\'e Duality Spaces}
\section{Introduction, definitions and elementary properties}
\newHead{II.1}
\bigskip
In this chapter we discuss the analogue of manifold in the proper homotopy
category.
We seek objects, to be called Poincar\'e duality spaces, which have the
proper homotopy attributes of paracompact manifolds.
To this end, we begin by discussing these attributes.

There is a well known Lefschetz duality between $H_\ast$ and
$H^\ast_\cmpsup$ or between $H^\ast$ and $H_\ast^\locf$
which is valid for any paracompact manifold with boundary (see for
instance Wilder \cite{bfortyfour}).
This duality is given via the cap product with a generator of 
$H^\locf_N$,
perhaps with twisted coefficients. 
This generator is called the {\sl fundamental class}.

Given any paracompact handlebody $M$, $M$ can be covered by an
increasing sequence of compact submanifolds with boundary.
Let $\{C_i\}$ be such a collection.
If $[M]\in H^\locf_{\dim M}(M;\Z^t)$ is the fundamental class, its
image in $H^\locf_{\dim M}(\Bar{M-C_i},\partial C_i; \Z^t)$ via the
inclusion and excision is the fundamental class for the pair
$(\Bar{M-C_i},\partial C_i)$.
A word about notation: $\Z^t$ occurring as a coefficient group will
always denote coefficients twisted by the first Stiefel{--}Whitney 
class of the manifold.

\bigskip
\BEGIN{T.2.1.1}
The fundamental class $[M]$ in $H^\locf_N(M;\Z^t)$ induces via
the cap product an isomorphism
\[\capf {[M]}{}\colon \Delta^{N-\ast}(M\Colon \coverFA{})\ \to\
\Delta_\ast(M\Colon \coverFA{})\]
where \ $\coverFA{}$\quad is any covering functor.

If $M$ has a boundary, we get a fundamental class
$[M]\in H^\locf_N(M, \partial M;\Z^t)$ and isomorphisms
\[\begin{matrix}%
\capf {[M]}{}\colon& \Delta^{N-\ast}(M,\partial M\Colon \coverFA{})&\to&
\Delta_\ast(M\Colon \coverFA{})\hfill\cr
\capf {[M]}{}\colon& \Delta^{N-\ast}(M\Colon \coverFA{})\hfill&
\hbox to10pt{\hss$\RA{\hskip 25pt}$}&
\Delta_\ast(M,\partial M\Colon \coverFA{})\cr
\end{matrix}\]
A similar result holds for a manifold $n${--}ad.
\end{Theorem}

\begin{proof}
The proof is easy.
On the cofinal subset of compact submanifolds with boundary of $M$,
$[M]$ induces, via inclusion and excision, the fundamental class 
for the pair ( $(n+1)${--}ad in general) $(\Bar{M-C_i}, \partial C_i)$,
where $C_i$ is a compact submanifold with boundary of $M$,
and $(\Bar{M-C_i}$ is the closure in $M-C_i$ in $M$.
$\partial C_i$ is equally the boundary of $C_i$ as a manifold or the
frontier of $C_i$ as a set.
By the definition of $\capf{[M]}{}$, it induces an isomorphism for each
base point and set $C_i$.
Hence it must on the inverse limit.
\end{proof}
\medskip
If one computes the homology and cohomology from chain complexes
based on a PL triangulation, on a handlebody decomposition, or on
a triangulation of the normal disc bundle, $\capf{[M]}{}$ induces a 
chain isomorphism.
We can ask for the torsion of this map.
We have
\BEGIN{T.2.1.2}
If $(M,\partial M)$ is a manifold with (possibly empty) boundary, and if
$[M]\in H^\locf_N(M,\partial M;\Z^t)$ is the fundamental class,
$\capf{[M]}{}\colon \Delta^{N-\ast}(M,\partial M\Colon \coverFA{} )\ \to\
\Delta_\ast(M\Colon \coverFA{})$ and 
$\capf{[M]}{}\colon \Delta^{N-\ast}(M\Colon \coverFA{})\ \to\
\Delta_\ast(M,\partial M\Colon \coverFA{})$
are simple equivalences, where \ $\coverFA{}$\quad
is the universal covering functor.
\end{Theorem}
\begin{proof}
Given a handlebody decomposition, the proof is easy.
The cap product with the fundamental class takes the cochain which
is $1$ on a given handle and zero on all the other handles to the
dual of the given handle.
Hence $\capf{[M]}{}$ takes generators in cohomology to generators
in homology (up to translation by the fundamental group).
The fact that the simple homotopy type as defined by a PL triangulation 
or by a triangulation of the normal disc bundle is the same as that defined
by a handlebody has been shown by Siebenmann \cite{bthirtyfour}.
\end{proof}

\medskip
We are still left with manifolds which have no handlebody decomposition.
Let $N=\CP^4\ \#\ S^3\times S^5\ \#\ S^3\times S^5$.
Then $\chi(N)=1$. $N\times M$ has $[N]\times [M]$ 
as a fundamental class.
For $M$ we use the simple homotopy type defined by a triangulation 
of the normal disc bundle.
Then $\capf{[N]\times[M]}{}$ is a simple homotopy equivalence \iff\ 
$\capf{[M]}{}$ is by \fullRef{L.1.5.23},
since $\capf{[N]}{}$ is known to induce a simple equivalence.
But $\capf{[N]\times[M]}{}$ is a simple equivalence since $N\times M$
has a handlebody structure (Kirby{--}Siebenmann \cite{beighteen}).
Note Theorems \shortFullRef{T.2.1.1} 
and \shortFullRef{T.2.1.2} now hold for arbitrary
paracompact manifolds.

With these two theorems in mind, we make the following definition.
\medskip\begin{xDefinition}
A locally finite, finite dimensional CW pair $(X,\partial X)$ with
orientation class $\wone\in H^1(X;\cy2)$
is said to satisfy
Poincar\'e duality with respect to $[X]$ and the covering functor
\ $\coverFA{}$\quad provided there is a class
$[X]\in H^\locf_N(X,\partial X; \Z^\wone)$ such that the maps
\[\begin{matrix}%
\capf{[X]}{}\colon &\Delta^{N-\ast}(X\Colon \coverFA{})\hfill&\to&
\Delta_\ast(X,\partial X\Colon \coverFA{})\cr\noalign{\vskip6pt}
\capf{[X]}{}\colon &
\Delta^{N-\ast}(X,\partial X\Colon \coverFA{} )&\to&
\Delta_\ast(X\Colon \coverFA{})\hfill\cr
\end{matrix}\]
are isomorphisms.
$\Z^\wone$ denotes integer coefficients twisted by the class $\wone$.
\end{xDefinition}

If $X$ is an $n${--}ad we require that all the duality maps be isomorphisms.

\bigskip\begin{xRemarks}
The two maps above are dual to one another, so if one is an isomorphism
the other is also.
\end{xRemarks}

Suppose \ $\coverFA{}$\quad is a regular covering functor for $X$,
and suppose \ $\coverFC{}$\quad is another regular
covering functor with $\coverFC{} > \coverFA{}$.
Then the chain and cochain groups have the structure of 
$\Z\pi^\prime_1(X\Colon F,f\Colon \coverFA{})${--}modules, when
$f\colon T\ \to\ X$ is a tree and $F\in{\mathcal L}(f)$.
The tree of groups $\pi^\prime_1(X\Colon F,f\Colon \coverFA{})$
is the tree given by
$(\pi^\prime_1)_A=\pi_1\bigl(F(A),p\bigr)/\pi_1\Bigl(\coverFA{F(A)},p\Bigr)$
where $p$ is the minimal vertex for $A$.
There is a map of rings $\Z\pi^\prime_1(X\Colon F,f\Colon \coverFA{} )\ \to\ 
\Z\pi^\prime_1(X\Colon F,f\Colon \coverFC{} )$, and the tensor product takes
$\Delta(X\Colon \coverFA{})$ to $\Delta(X\Colon \coverFC{})$.
Since $\capf{[X]}{}$ is an isomorphism for $\coverFA{}$,
we can get chain homotopy inverses, so under tensor product,
$\capf{[X]}{}$ still induces isomorphisms for $\coverFC{}$.
\medskip
As we have the Browder Lemma (\fullRef{T.1.4.7}),
with patience we can prove a variety of cutting and gluing theorems.
The following are typical.

\bigskip
\BEGIN{T.2.1.3}
Let $(X\Colon \partial_0X,\partial_1X)$ be a triad.
Then the following are equivalent.
\begin{enumerate}
\item[1)] $(X\Colon \partial_0X,\partial_1X)$ satisfies Poincar\'e duality with respect
to \hfil\penalty-10000
\hbox{$V\in H^\locf_N(X\Colon \partial_0X,\partial_1X\Colon \Z^\wone)$ and
$\coverFA{}$}. 
\item[2)] $(\partial_0X,\partial_{\{0,1\}}X)$ satisfies Poincar\'e duality 
with respect to\hfil\penalty-10000
\hbox{$\partial V\in H^\locf_{N-1}(\partial_0X\Colon 
\partial_{\{0,1\}}X\Colon \Z^\wone)$ and
$\coverFA{}$ where $\coverFA{}$ is induced from}
$\coverFA{}$ over $X$\hfil\penalty-1000 and $\wone$ is the orientation
class induced from $\wone$ over $X$.
Moreover, one of the maps
\[\begin{matrix}%
\capf{V}{}\colon&\Delta^\ast(X,\partial_1X\Colon \coverFA{})&
\to&\Delta_{N-\ast}(X,\partial X_0\Colon \coverFA{})\cr
\noalign{\vskip6pt}
\capf{V}{}\colon&\Delta^\ast(X,\partial_0X\Colon \coverFA{})&
\to&\Delta_{N-\ast}(X,\partial X_1\Colon \coverFA{})\cr
\end{matrix}\]
is an isomorphism. (Hence they are both isomorphisms.)
\item[3)] The same conditions as 2) but considering
$(\partial_1X,\partial_{\{0,1\}}X)$.
\end{enumerate}
\end{Theorem}

\begin{proof}
The proof is fairly standard.
We look at one of the sequences associated to a triple, say
\[\begin{matrix}%
\Delta^\ast(X\Colon \partial_0X,\partial_1X\Colon \coverFA{})&\to&
\Delta^\ast(X,\partial_1X\Colon \coverFA{})&\to&
\Delta^\ast(\partial_0X\Colon \partial_{\{0,1\}}X\Colon X\Colon \coverFA{})\cr
\downlabeledarrow[\Big]{}{\capf{V}{}}&&\downlabeledarrow[\Big]{}{\capf{V}{}}&&
\downlabeledarrow[\Big]{}{\capf{\partial V}{}}\cr\noalign{\vskip4pt}
\Delta_{N-\ast}(X\Colon \coverFA{})&
\hbox to 10pt{\hss$\RA{\hskip 26pt}$}&
\Delta_{N-\ast}(X,\partial_0X\Colon \coverFA{})&\to&
\Delta_{N-1-\ast}(\partial_0X\Colon X\Colon \coverFA{})\cr
\end{matrix}\]
1) implies both $\capf{V}{}$'s are isomorphisms, so the $5${--}lemma
shows $\capf{\partial V}{}$ is an isomorphism.
2) implies one of the $\capf{V}{}$'s is an isomorphism and that
\[\capf{\partial V}{}\colon \Delta^\ast(\partial_0X,\partial_{\{0,1\}}X\Colon 
\coverFA{})\ \to\ \Delta_{N-1-\ast}(\partial_0X\Colon \coverFA{})\]
is an isomorphism.
Hence we must investigate how the subspace groups depend 
on the absolute groups.
Make sure that the set of base points for $X$ contains 
a set for $\partial_0X$.
Then we have a diagram
\[\begin{matrix}%
\Delta^\ast(\partial_0X,\partial_{\{0,1\}}X\Colon \coverFA{})&
\to&\Delta^\ast(\partial_0X,\partial_{\{0,1\}}X\Colon X\Colon 
\coverFA{})\cr\noalign{\vskip6pt}
\downlabeledarrow[\big]{}{\cap\partial V}&&\downlabeledarrow[\big]{}{\cap\partial V}\cr
\noalign{\vskip6pt}
\Delta_{N-1-\ast}(\partial_0X \Colon \coverFA{})&
\to&\Delta_{N-1-\ast}(\partial_0X\Colon X \Colon \coverFA{})\cr
\end{matrix}\]
which commutes.
The horizontal maps are naturally split, so if
$\capf{\partial V}{}$ on the subspace groups is an isomorphism,
then it is also an isomorphism on the absolute groups.
Hence 1) implies 2) and 3).

Now if $\capf{\partial V}{}$ on the absolute groups is an isomorphism,
then it is also an isomorphism on the subspace groups by
\fullRef{T.2.1.3}.
Hence 2)  or 3) implies 1).
\end{proof}

\bigskip
\BEGIN{T.2.1.4}
Let $Z= Y\ \cup\ Y^\prime$ and set $X = Y\ \cap\ Y^\prime$.
Then any two of the following imply the third.
\begin{enumerate}
\item[1)] $Z$ satisfies Poincar\'e duality with respect to $[Z]$ and
$\coverFA{}$.
\item[2)] $(Y,X)$ satisfies Poincar\'e duality with respect to $\partial[Z]$
and $\coverFA{}$.
\item[3)] $(Y^\prime,X)$ satisfies Poincar\'e duality with respect to 
$\partial[Z]$ and $\coverFA{}$.
\end{enumerate}

where $\coverFA{}$ is a covering functor over $Z$, which
then induces $\coverFA{}$ over $Y$ and $Y^\prime$.
An orientation class over $Z$ which induces one over $Y$ and
$Y^\prime$ has been assumed in our statements.
\end{Theorem}
\begin{proof}
The reader should have no trouble proving this.
\end{proof}

\medskip
A map $\varphi\colon M\ \to\ X$, where $M$ and $X$ are locally compact
CW $n${--}ads which satisfy Poincar\'e duality with respect to $[M]$ and
$\coverFA{}$, and $[X]$ and $\coverFC{}$
respectively, is said to be {\sl degree $1$\/} provided it is a map
of $n${--}ads and
\begin{enumerate}
\item[1)] $\varphi^\ast(\coverFC{})=\coverFA{}$,
where $\varphi^\ast(\coverFC{})$ is the covering functor 
over $M$ induced by $\varphi$ from $\coverFC{}$
over $X$.
\item[2)] If $\wone\in H^1(X;\cy2)$ is the orientation class for $X$,
$\varphi^\ast\wone$ is the orientation class for $M$.
\item[3)] $\varphi_\ast[M]=[X]$.
\end{enumerate}

\bigskip
\BEGIN{T.2.1.5}
Let $\varphi\colon M\ \to\ X$ be a map of degree $1$ 
of Poincar\'e duality spaces.
Then the diagram
\[\begin{matrix}%
\Delta^r(M\Colon \coverFA{})&\LA{\ \varphi^\ast\ }&
\Delta^r(X\Colon \coverFA{})\cr
\downlabeledarrow[\Big]{}{\capf{[M]}{}}&&\downlabeledarrow[\Big]{}{\capf{[X]}{}}\cr
\noalign{\vskip6pt}
\Delta_{n-r}(M\Colon \coverFA{})&\RA{\ \varphi_\ast\ }&
\Delta_{n-r}(X\Colon \coverFA{})\cr
\end{matrix}\]
commutes. ($\coverFA{}$ over $M$ is the covering functor
induced from $\coverFA{}$ over $X$.)
Hence $\capf{[M]}{}$ induces an isomorphism on the cokernel of
$\varphi^\ast$, $K^r(M\Colon \coverFA{})$ onto the kernel of
$\varphi_\ast$, $K_{n-r}(M\Colon \coverFA{})$.
Thus if $\varphi$ is $k${--}connected, $\varphi_\ast$ and $\varphi^\ast$
are isomorphism for $r<k$ and $r>n-k$.

Similarly let $\varphi\colon (N,M)\ \to\ (Y,X)$
be a degree $1$ map of pairs.
Then $\varphi_\ast$ gives split surjections of homology groups with
kernels $K_\ast$, and split injections of cohomology with cokernels
$K^\ast$.
The duality map $\capf{[N]}{}$ induces isomorphisms
$K^\ast(N\Colon \coverFA{})\ \to\ 
K_{n-\ast}(N,M\Colon\coverFA{})$ and
$K^\ast(N, M\Colon \coverFA{})\ \to\ 
K_{n-\ast}(N\Colon \coverFA{})$.

Analogous results hold for $n${--}ads.
\end{Theorem}

\begin{proof}
The results follow easily from definitions and the naturality 
of the cap product.
\end{proof}

\bigskip
\section{The Spivak normal fibration}
\newHead{II.2}

One important attribute of paracompact manifolds is 
the existence of normal bundles.
In \cite{bthirtysix} Spivak constructed an analogue for these bundles
in the homotopy category.
Although he was interested in compact spaces, he was often forced 
to consider paracompact ones.
It is then not too surprising that his definition is perfectly adequate for
our problem. 
This is an example of a general principle in the theory of paracompact
surgery, namely that all bundle problems encountered are exactly the
same as in the compact case.
One does not need a ``proper'' normal bundle or a ``proper''
Spivak fibration.

\bigskip\begin{xDefinition}
Let $(X,\partial X)$ be a locally compact, finite dimensional CW pair.
Embed $(X,\partial X)$ in $(\HS^n,\R^{n-1})$, where $\HS^n$
is the upper half plane and $\R^{n-1}=\partial\HS^n$.
Let $(N;N_1,N_2)$ be a regular neighborhood of $X$ as a subcomplex
of $\HS^n$; i.e. $X\subseteq N$, $\partial X\subseteq N_2$ and
$N$ (resp. $N_2$) collapses to $X$ (resp. $\partial X$).
Let ${\mathcal P}(N_1,N,X)$ be the space of paths starting in $N_1$,
lying in $N$, and ending in $X$ endowed with the compact-open topology.
(If $A$, $B$, $C$ are spaces with $A$, $C\subseteq B$, a
similar definition holds for ${\mathcal P}(A,B,C)$.)
There is the endpoint map $w\colon {\mathcal P}(N_1,N,X)\ \to\ X$.
$w$ is a fibration and is called the {\sl Spivak normal fibration}.
Its fibre is called the {\sl Spivak normal fibre}.
\end{xDefinition}

Spivak showed that a necessary and sufficient condition for a finite complex
to satisfy Poincar\'e duality with respect to the universal covering functor
was that the Spivak normal fibre of the complex should have the
homotopy type of a sphere.
He also showed that if one started with a compact manifold, then the
normal sphere bundle had the same homotopy type as the Spivak
normal fibration, at least stably.
Before we can do this for paracompact manifolds, we will need to do some
work.

In practice, the fact that the Spivak normal fibration is constructed from
a regular neighborhood is inconvenient.
More convenient for our purposes is a {\sl semi{--}regular neighborhood}

\begin{xDefinition}
Let $(X,\partial X)$ be a pair of finite dimensional, locally compact 
CW complexes.
A {\sl semi{--}regular neighborhood\/} ( s-r neighborhood) is a manifold
triad $(M\Colon M_1, M_2)$ and proper maps
$i\colon X\ \to\ M$ and $j\colon \partial X\ \to\ M_2$
such that
\topD{10}{$\begin{matrix}%
X&\RA{\ i\ }&M\cr
\cup\vrule width .1pt depth 0pt height 6pt&&
\cup\vrule width .1pt depth 0pt height 6pt\cr
\partial X&\RA{\ j\ }&M_2\cr
\end{matrix}$}{4}
commutes, and such that
$i$ and $j$ are simple homotopy equivalences.
Lastly we require that $M$ be parallelizable (equivalent to being
stably parallelizable).
The definition for an $n${--}ad is similar: we have a manifold $(n+1)${--}ad
$(M\Colon M_1,\cdots M_{n})$ with a simple homotopy equivalence of $n${--}ads
$i\colon X\to \delta_1 M$ .
\end{xDefinition}

\bigskip
\BEGIN{T.2.2.1}
The fibration $w\colon{\mathcal P}(M_1,M,X)\ \to\ X$ is stably fibre
homotopy equivalent to the Spivak normal fibration.
\end{Theorem}
\begin{proof}
The proof needs 
\medskip
\BEGIN{L.2.2.1}
If $\bigl((M\Colon M_1,M_2), i, j \bigr)$ is an s-r neighborhood of $X$, then so is
\[\bigl((M\Colon M_1,M_2)\times(D^m,S^{m-1}), i\times c, j_1\bigr)\]
where $c$ denotes the constant map.
The $(n+1)$ structure on the product is $(M\times D^m\Colon 
M_1\times D^m\ \cup\ M\times S^{m-1},
M_2\times D^m,\cdots)$.
Let $\xi$ be ${\mathcal P}(M_1, M,X)\ \to\ X$ and let $\eta$ be
${\mathcal P}(M_1\times D^m\ \cup\ M\times S^{m-1},M\times D^m)\ \to\ X$.
Then $\xi\ast(m)$ is fibre homotopy equivalent to $\eta$, where $(m)$ is 
the trivial spherical fibration of dimension $m-1$ and $\ast$ denotes the 
fibre join.
\end{Lemma}
\begin{proof}
The first statement is trivial and the second is Spivak \cite{bthirtysix},
Lemma 4.3.
\end{proof}
\smallskip
Now if $(M\Colon M_1,M_2)$ is an s-r neighborhood of $X$, then for some $n$,
$(M\Colon M_1,M_2)\times (D^n,S^{n-1})$ is homeomorphic to a regular
neighborhood of $X$ in $\R^{n+m}$, where $m=\dim M$.
If we can show this then the lemma easily implies that $\xi$ is stably
equivalent to the Spivak normal fibration formed from this regular
neighborhood.

By crossing with $D^n$ if necessary, we may assume 
$\dim M\geq 2\dim X+1$, so we may assume $i$ and $j$ are embeddings.
Since $M$ is parallelizable, $(M,\partial M)$ immerses in $(\HS^m,\R^{m-1})$.
Since $m\geq 2\dim X+1$ we can subject $i$ and $j$ to a proper homotopy
so that $i\colon X\ \to\ M\subset \HS^m$ and
$j\colon \partial X\ \to\ M_2\subset \R^{m-1}$ become embeddings
on open neighborhoods $U$ and $U_2$, where $U$ is a neighborhood
of $X$ in $M$ and $U_2=M_2\ \cap\ U$ is a neighborhood of
$\partial X$.
In $U$ sits a regular neighborhood of $X$, $(N\Colon N_1, N_2)$.
Hence $(N\Colon N_1, N_2)\subseteq (M\Colon M_1, M_2)$ and excision gives a
simple homotopy equivalence $\partial N\subseteq \Bar{M-N}$ and
$\partial N_2\subseteq \Bar{M_2-N_2}$ (this uses the fact that $i$ and
$j$ are simple homotopy equivalences).
By the $s${--}cobordism theorem (see \cite{bthirtythree} or
\cite{bten}) these are products (assume $m\geq 6$) so
$(M\Colon M_1,M_2)$ is homeomorphic to a 
regular neighborhood of $X$ in $\R^m$.
\end{proof}

\medskip
\BEGIN{C.2.2.1.1}
The Spivak normal fibration is stably well defined.
\end{Corollary}
\medskip\begin{xRemarks}
By definition we have a Spivak normal fibration for any regular 
neighborhood, so we can not properly speak of ``the'' Spivak normal 
fibration.
By the corollary however they are all stably equivalent, so we will
continue to speak of the Spivak normal fibration when we really
mean any fibration in this stable class.
This includes fibrations formed from s-r neighborhoods.
\end{xRemarks}

Now, for finite complexes we know the complex satisfies Poincar\'e duality 
\iff\ the Spivak normal fibre has the homology of a sphere.
Unfortunately, this is not true for our case.
In fact, Spivak has already shown what is need to get the normal fibre a
sphere.
This information is contained in Theorems 
\shortFullRef{T.2.2.2} and \shortFullRef{T.2.2.3}

\begin{xDefinition}\footnote[1]{Expanded.}
A locally compact, finite dimensional  CW complex is a {\sl Spivak space\/}
provided the fibre of any Spivak normal fibration has the homology of
a sphere.
A {\sl Spivak pair\/} is a pair, $(X,\partial X)$, of locally compact,
finite dimensional CW complexes such that the fibre of any Spivak
normal fibration has the homology of a sphere, and such that the
Spivak normal fibration for $X$ restricted to $\partial X$ is the Spivak
normal fibration for $\partial X$. 
To be slightly more precise, given an s-r neighborhood for 
the pair $(X,\partial X)$, there is a natural fibre map from the
Spivak normal fibration for $\partial X$ to the Spivak normal fibration
for $X$ restricted to $\partial X$: it is this map we are requiring to
be an equivalence.
A {\sl Spivak $n${--}ad\/} is defined analogously.
Any Spivak space (pair, $n${--}ad) has a first Stiefel{--}Whitney class
and a fundamental class.
The first Stiefel{--}Whitney class of the Spivak space, $\wone$, is the
first Stiefel{--}Whitney class of the Spivak normal spherical fibration.%
There is a Thom isomorphism $H^\locf_{m+k}(M,\partial M;\Z)\ \to\ 
H^\locf_{m}(M,M_2;\Z^\wone)$.
Since a parallelizable manifold is oriented, we get a fundamental class
$[M]\in H^\locf_{m+k}(M,\partial M;\Z)$ and the fundamental class of
the Spivak space is the image of this class in 
$H^\locf_{m}(X,\partial X;\Z^\wone)\cong H^\locf_{m}(M,M_2;\Z^\wone)$.
We denote it by $[X]$ and note that it is defined up to sign.
A choice of sign will be called an {\sl orientation}.
\end{xDefinition}

\bigskip
\BEGIN{T.2.2.2}
The following are equivalent
\begin{enumerate}
\item[1)] $X$ is a Spivak space
\item[2)] 
$\capf{[X]}{}\colon 
H^\ast_\cmpsup(\tilde X)\ \to\ H_{N-\ast}(\tilde X)$ is an isomorphism
\item[3)] $\capf{[X]}{}\colon 
H^\ast(\tilde X)\ \to\ H^\locf_{N-\ast}(\tilde X)$ is an isomorphism
\end{enumerate}
\end{Theorem}

\begin{proof}
2) implies 3) thanks to the following commutative diagram.
\[\begin{matrix}%
0\to&{\text{ Ext}}\bigl(H^{\ast+1}_\cmpsup(\tilde X;\Z),\Z\bigr)&\to&
H^\locf_\ast(\tilde X;\Z) &\to&
{\Homx}\bigl(H^\ast_\cmpsup(\tilde X;\Z),\Z)&\to0\cr
&\uplabeledarrow[\bigg]{}{{\text{ Ext}}\bigl(\capf{[X]}{}\bigr)}&&
\uplabeledarrow[\bigg]{}{\capf{[X]}{}}&&
\uplabeledarrow[\bigg]{}{{\Homx}\bigl(\capf{[X]}{}\bigr)}\cr
0\to&{\text{ Ext}}\bigl(H_{N-1-\ast}(\tilde X;\Z),\Z\bigr)&\to&
H^{N-\ast}(\tilde X;\Z) &\to&
{\Homx}\bigl(H_{N-\ast}(\tilde X;\Z),\Z)&\to0\cr
\end{matrix}\]

3) implies 1) thanks to Spivak, Proposition 4.4, and the observation that
the Spivak normal fibration for $X$ pulled back over $\tilde X$ is
the Spivak normal fibration for $\tilde X$.
This observation is an easy consequence of 
\fullRef{T.2.2.1}, the definition of an s-r neighborhood,
and the fact that the transfer map $\sieb(X)\ \to\ \sieb(\tilde X)$
is a homomorphism.

1) implies 2) as follows.
Look at
\[\begin{matrix}%
H_{N-\ast}\bigl(D(\tilde X);\Z\bigr)&
\ZRA{0}{\hskip70pt}{10}&H_{N-\ast}(\tilde X;\Z)\cr
\uplabeledarrow[\Big]{}{\capf{U}{}}\cr
H_{N+k-\ast}\bigl(D(\tilde X), S(\tilde X);\Z)&\LA{\hskip30pt}&
H_{N+k-\ast}(\tilde N,\partial N;\Z)\cr
&&\uplabeledarrow[\Big]{}{\capf{[N]}{}}\cr
H^\ast_\cmpsup(\tilde X;\Z)&\ZLA{-4}{\hskip108pt}{0}
&H^\ast_\cmpsup(N)\cr
\end{matrix}\]
where $U$ is the Thom class for the normal disc fibration $D(\tilde X)$
with spherical fibration $S(\tilde X)$.
$(N,\partial N)$ is an s-r neighborhood for $\tilde X$.
The horizontal maps are induced by the inclusion $\tilde X\subseteq N$
and the proper homotopy equivalence $(N,\partial N)\ \to\ 
\bigl(D(\tilde X), S(\tilde X)\bigr)$.
All horizontal maps are isomorphisms.
The composite map $H^\ast_\cmpsup(\tilde X)$ to $H_{N-\ast}(\tilde X)$
is essentially the cap product with $U\cap [N]$, where $U\cap[N]$
should be actually be written $i_\ast\bigl(U\cap j_\ast[N]\bigr)$,
where $i_\ast\colon H^\locf_\ast\bigl(D(\tilde X)\bigr)\ \to\ 
H^\locf_\ast(\tilde X)$ and $j_\ast\colon
H^\locf(N, \partial N)\ \to\ H^\locf_\ast\bigl(D(\tilde X), S(\tilde X)\bigr)$.
$H^\locf_\ast\bigl(D(\tilde X)\bigr)$ is the homology group of the infinite
singular chains on $D(\tilde X)$ which project to give locally finite chains
on $\tilde X$.
( $H^\locf_\ast\bigl(D(\tilde X), S(\tilde X)\bigr)$ is similar.)
Now $[X]=U\cap[N]$ shows 2) is satisfied. 
[ 1) was used to get the Thom class $U$.]
\end{proof}

\medskip
\BEGIN{C.2.2.2.1}
\footnote[1]{New Corollary.}
For a Spivak space $X$, \\
$\capf{[X]}{}\colon 
H^\ast_\cmpsup(\tilde X;\Gamma)\ \to\ 
H_{N-\ast}(\tilde X;\Gamma^\wone)$ and
$\capf{[X]}{}\colon 
H^\ast(\tilde X;\Gamma)\ \to\ 
H^\locf_{N-\ast}(\tilde X;\Gamma^\wone)$ are isomorphisms for
any local coefficients $\Gamma$.
\end{Corollary}

\smallskip\begin{proof} The results follows from the universal coefficient
formulas.
\end{proof}

\bigskip
\BEGIN{T.2.2.3}
\footnote[2]{Restated Theorem and reworked proof.}
Fix an $[X]\in H^\locf_N(X,\partial X;\Z^\wone)$ and
consider the following four families of maps
\begin{enumerate}\item[]
\begin{enumerate}
\item[A)] $\capf{[X]}{}\colon 
H^\ast_\cmpsup(\tilde X,\widetilde{\partial X};\Z)\ \to\ 
H_{N-\ast}(\tilde X;\Z)$
\item[B)] $\capf{[X]}{}\colon 
H^\ast(\tilde X,\widetilde{\partial X};\Z)\ \to\ 
H^\locf_{N-\ast}(\tilde X;\Z)$
\item[C)] $\capf{[X]}{}\colon 
H^\ast_\cmpsup(\tilde X;\Z)\ \RA{\hskip16pt}\ 
H_{N-\ast}(\tilde X,\widetilde{\partial X};\Z)$
\item[D)] $\capf{[X]}{}\colon 
H^\ast(\tilde X;\Z)\ \RA{\hskip16pt}\ 
H^\locf_{N-\ast}(\tilde X,\widetilde{\partial X};\Z)$
\end{enumerate}
\end{enumerate}
Suppose $\partial X$ is a Spivak space with fundamental
class $\partial [X]$.
Then the following are equivalent.
\begin{enumerate}
\item[1)] $(X,\partial X)$ is a Spivak pair with fundamental class $[X]$
\item[2)] any one of the above four maps is an isomorphism
\end{enumerate}
\end{Theorem}

\begin{proof}
1) implies by definition that all four maps are isomorphisms.

By a diagram similar to the one in the proof of 
\fullRef{T.2.2.2}, A) implies D) and C) implies B).
Similar diagrams show D) implies A) and B) implies C).
In fact, 2) implies all four maps are isomorphisms.
If A) is an isomorphism, then the Browder lemma shows C) is too.
Conversely, if C) is an isomorphism the Browder lemma shows A) is.
This shows the claim.

D) implies that the Spivak normal fibre is 
a sphere by Spivak, Proposition 4.4.

Either 1) or 2) implies that the Spivak normal fibration 
has a spherical fibre, $S^{k-1}$,  so we have a Thom isomorphism
\[H^\locf_{N+k}(M,\partial M;\Z)\ \RA{\ \capf{U}{}\ }\ 
H^\locf_N(X,\partial X;\Z^\wone)\]
and that we can pick a fundamental class 
$[M]\in H^\locf_{N+k}(M,\partial M;\Z)$
with $[X]=[M]\cap U$. 

This shows 1) implies 2) and the only remaining point in 2) implies 1)
is to check that the Spivak normal fibration for $X$ 
restricted to $\partial X$, ${\mathcal P}(M_1, M, \partial X)$, 
is the Spivak normal fibration for $\partial X$, 
${\mathcal P}(M_1\cap M_2, M_2, \partial X)$.
There is an evident inclusion \[{\mathcal P}(M_1\cap M_2, M_2, \partial X)
\subseteq {\mathcal P}(M_1, M, \partial X)\ ,\leqno(\ast)\]
where both fibrations are spherical with fibre $S^{k-1}$.
Let $\wone$ and $U$ be the orientation class and the Thom class
for ${\mathcal P}(M_1, M, \partial X)$.
The following 
diagram of restrictions and Thom isomorphisms commutes
\[\begin{matrix}%
H^\locf_{N+k}(M,\partial M;\Z)&\RA{\ \capf{U}{}\ }&
H^\locf_N(X,\partial X;\Z^\wone)\cr
\big\downarrow&&\big\downarrow\cr
H^\locf_{N-1+k}(M_2,\partial M_2;\Z)&\RA{\ \capf{U}{}\ }& 
H^\locf_{N-1}(\partial X;\Z^\wone)
\end{matrix}\]
where the vertical maps are Poincar\'e isomorphisms and the top map is
a Thom isomorphism.
It follows that the bottom map is an isomorphism, so the Thom class of
one of our fibrations restricts to be the Thom class of the other,
which shows $(\ast)$ is a fibre homotopy equivalence.
\end{proof}

\medskip
\BEGIN{C.2.2.3.1}
\footnote[1]{New Corollary.}
For a Spivak pair $(X,\partial X)$, the analogues of all four maps
are isomorphisms with
any local coefficients $\Gamma$.
\end{Corollary}
\smallskip\begin{proof} The results follows from the universal coefficient
formulas.
\end{proof}

\bigskip
Now suppose $\xi$ is an arbitrary spherical fibration over a locally compact,
finite dimensional CW complex $X$.
The total space in general is not such a complex.
Our techniques apply best to such spaces however, and we want to
study these total spaces.
Hence we wish to replace any such space by a space with the proper
homotopy type of a locally compact, finite dimensional CW complex.

\medskip\begin{xDefinition}
Let $S(\xi)$ be the total space of a spherical fibration $\xi$ over
a locally compact, finite dimensional $n${--}ad $X$.
A \emph{\cwation}\ of $\xi$ is an $n${--}ad $Y$ and a proper map 
$Y\ \to\ X$ such that the following conditions are satisfied.
$Y$ has the proper homotopy type of a locally finite, finite dimensional 
$n${--}ad.
There are maps 
$\vrule width 0pt depth 10pt height 12pt
\displaystyle S(\xi)\ {\scriptscriptstyle\Atop%
{\Limitsarrow{10}{10}{\scriptscriptstyle g}{}{\rightarrowfill}}
{\Limitsarrow{10}{10}{}{\scriptscriptstyle h}{\leftarrowfill}}
}\ Y$ 
such that $h\ \circ\ g$ is a fibre map, fibre homotopic to the
identity and such that $g\ \circ\ h$ is properly homotopic to
the identity.
Lastly
\topD{12}{$\begin{matrix}%
 S(\xi)&{\scriptscriptstyle\Atop
{\Limitsarrow{10}{10}{\scriptscriptstyle g}{}{\rightarrowfill}}
{\Limitsarrow{10}{10}{}{\scriptscriptstyle h}{\leftarrowfill}}
}& Y\cr
&\hbox to 10pt{\hss$\searrow\hskip 20pt\swarrow\hbox to 0pt{$\scriptstyle{f}$\hss}$\hss}\cr
&X\cr
\end{matrix}$}{4}
should commute.
The pair $(M_f,Y)$ is seen to satisfy the Thom isomorphism for 
$\Delta^\ast$ and $\Delta_\ast$ theories. (See the appendix page \pageref{PoincareAppendix} for
a discussion of the Thom isomorphism in these theories.)
The simple homotopy type of $Y$ is defined by any locally compact,
finite dimensional CW complex having the same proper homotopy 
type as $Y$ and for which the Thom isomorphisms are simple
homotopy equivalences.
For a fibration $\xi$, $\bigl(D(\xi),C(\xi)\bigr)$ will denote the pair
$(M_f,Y)$ with this simple homotopy type.
Such a pair is said to be a {\sl simple \cwation}.
It has a {\sl Thom class\/} $U_\xi\in H^k\bigl(D(\xi),C(\xi);
\Z^{w_1(\xi)}\bigr)$, where $w_1(\xi)$ is the first Stiefel{--}Whitney
class of the spherical fibration $\xi$.
There is a Thom isomorphism with twisted coefficients.
\end{xDefinition}

\medskip\begin{xRemarks}
Any spherical fibration of dimension two or more has a simple \cwation.
The proof of this fact is long and is the appendix page \pageref{PoincareAppendix} to this chapter.
\end{xRemarks}

\bigskip
\BEGIN{T.2.2.4}
Let $\xi$ be any spherical fibration of dimension $>1$ over a Spivak
space $X$.
Then $\bigl(D(\xi), C(\xi)\bigr)$ is a Spivak pair 
with fundamental class $[\xi\>]\in H^\locf_{N+k}\bigl(
D(\xi), C(\xi);\Z^\omega\bigr)$, where 
$\omega(g)=\wone(g)\cdot \bigl(w_1(\xi)\bigr)(g)$, where
$\wone$ is the first Stiefel{--}Whitney class for $X$. 
If $U_{\xi}$ is a Thom class for $\xi$, $U_\xi \cap\ [\xi] = [X]$, 
a fundamental class for $X$. 

If $(X,\partial X)$ is a Spivak pair, then 
$\bigl(D(\xi)\Colon  D(\xi\vert_{\partial X}),C(\xi)\bigr)$ is a Spivak triad with
fundamental class $[\xi\>]$.
In general, a \cwation\ of a Spivak $n${--}ad has 
an $(n+1)${--}ad structure with fundamental class $[\xi\>]$.
Fundamental classes and Thom classes are related as in the absolute case. 
\end{Theorem}

\medskip\begin{proof}\footnote[1]{New proof.}
To show $\bigl(D(\xi), C(\xi)\cup D(\xi\>\vert_{\partial X})\bigr)$ 
is a Spivak pair, look at
\[\begin{matrix}%
H^\ast\big(\coverFC{D(\xi)}\bigr)&
\RA{\hskip 10pt\psi\hskip10pt}&
H^\locf_{N+k-\ast}\bigl(\coverFC{D(\xi)}, 
\coverFC{C(\xi)\cup D(\xi\>\vert_{\partial X})\quad}\hskip-5pt\bigr)\cr
\downlabeledarrow[\Big]{\cong}{\cap[X]}&&
\downlabeledarrow[\Big]{\cong}{U_\xi\cap}\cr
H^\locf_{N-\ast}(\coverFC{X\ },\coverFC{\partial X})
&\hbox to 40pt{$\RA{\hskip 60pt}$\hss}&
H^\locf_{N-\ast}\bigl(\coverFC{D(\xi)}, 
\coverFC{D(\xi\>\vert_{\partial X})}\bigr)\cr
\end{matrix}\]
The bottom horizontal map is an isomorphism, hence so is $\psi$.
We claim $\psi(x)=x\ \cap\ \psi(1)$, where 
$1\in H^0\bigl(\coverFC{D(\xi)}\bigr)$ is a choice of generator
which we pick so $U_\xi\cap \psi(1)=[X]$.
But from the diagram,
$U_\xi\ \cup\ \psi(x)=x\cap [X]=x\cap U_\xi\cap \psi(1)=
U_\xi\cap\bigl(x\cap \psi(1)\bigr)$ which proves the result.

Suppose we can show there is a class 
$[\xi\>]\in H^\locf_{N+k}\bigl(D(\xi), 
C(\xi)\cup D(\xi\vert_{\partial X});\Z^\omega\bigr)$
with $tr\bigl([\xi\>]\bigr)=\psi(1)$.
Note
\[\begin{matrix}%
H^\ast_\cmpsup\big(\coverFC{D(\xi)}\bigr)&
\RA{\hskip10pt \cap [\xi\>]\hskip 10pt}&
H_{N+k-\ast}\bigl(\coverFC{D(\xi)}, 
\coverFC{C(\xi)\cup D(\xi\>\vert_{\partial X}})\bigr)\cr
\downlabeledarrow[\Big]{\cong}{\cap[X]}&&
\downlabeledarrow[\Big]{\cong}{U_\xi\cap}\cr
H_{N-\ast}(X,\partial X)
&\hbox to 40pt{$\RA{\hskip 60pt}$\hss}&
H_{N-\ast}\bigl(\coverFC{D(\xi)}, 
\coverFC{D(\xi\>\vert_{\partial X})}\bigr)\cr
\end{matrix}\]
commutes and since three of the maps are isomorphisms, 
so is $\cap\ [\xi\>]$.

Since $\xi$ has dimension $>1$,  
$\coverFC{D(\xi\>\vert_{\partial X}) \cup C(\xi)\quad}$ 
is simply connected.
By the Browder lemma,
\[\cap \partial [\xi\>]\colon 
H^\ast\bigl(\coverFC{D(\xi\>\vert_{\partial X}) 
\cup C(\xi)\quad};\Z\bigr)\ \RA{\hskip10pt}\ 
H_{N+k-1-\ast}\bigl(\coverFC{D(\xi\>\vert_{\partial X}) 
\cup C(\xi)\quad};\Z\bigr)\]
is an isomorphism, so by \fullRef{T.2.2.2},
$D(\xi\>\vert_{\partial X}) \cup C(\xi)$ is a Spivak space with
fundamental class $\partial [\xi\>]$.
It follows from \fullRef{T.2.2.3} that
$\bigl(D(\xi), D(\xi\>\vert_{\partial X}) \cup C(\xi)\bigr)$ is a
Spivak pair with fundamental class $[\xi\>]$. 

Hence we are done if we can construct the class $[\xi\>]$
with the right trace.
To do this, look at the Thom isomorphism
\[H^\locf_{N+k-\ast}\bigl(D(\xi), 
C(\xi)\cup D(\xi\>\vert_{\partial X});\Z^\omega\bigr)\ 
\RA{\ U_\xi \cap\ }\ 
H^\locf_{N-\ast}\bigl(D(\xi), D(\xi\>\vert_{\partial X});\Z^\wone\bigr)
\]
There is a unique element $[\xi\>]\in H^\locf_{N+k}\bigl(D(\xi), 
C(\xi)\cup D(\xi\>\vert_{\partial X});\Z^\omega\bigr)$
such that $U_\xi\cap [\xi\>]=[X]$.
Since cap product commutes with trace,
\[\begin{matrix}%
H^\locf_{N+k-\ast}\bigl(D(\xi), 
C(\xi)\cup D(\xi\>\vert_{\partial X});\Z^\omega\bigr)&
\RA{\ U_\xi \cap\ }& 
H^\locf_{N-\ast}\bigl(D(\xi), D(\xi\>\vert_{\partial X});\Z^\wone\bigr)\cr
\downlabeledarrow[\big]{}{{\text{ tr}}}&&\downlabeledarrow[\big]{}{{\text{ tr}}}\cr
H^\locf_{N+k-\ast}\bigl(\coverFC{D(\xi)}, 
\coverFC{C(\xi)\cup D(\xi\>\vert_{\partial X})};\Z\bigr)&
\RA{\ U_\xi \cap\ }& 
H^\locf_{N-\ast}\bigl(\coverFC{D(\xi)}, 
\coverFC{D(\xi\>\vert_{\partial X})};\Z\bigr)\cr
\end{matrix}\]
commutes and it is easy to see ${\text{ tr}}\bigl([\xi\>]\bigr)=\psi(1)$
since 
\[\begin{aligned}%
H^\locf_{N}\bigl(D(\xi), D(\xi\>\vert_{\partial X});\Z^\wone\bigr)\cong
H^\locf_{N-\ast}\bigl(X, \partial X);\Z^\wone\bigr)\cong\Z\cr
H^\locf_{N}\bigl(\coverFC{D(\xi)}, \coverFC{D(\xi\>\vert_{\partial X})};
\Z^\wone\bigr)\cong
H^\locf_{N-\ast}\bigl(\coverFC{X}, \coverFC{\partial X});
\Z^\wone\bigr)\cong\Z\cr
\end{aligned}\]
and ${\text{ tr}}$ is an isomorphism by \fullRef{T.2.2.5} below.
\end{proof}

\medskip
\BEGIN{C.2.2.4.1}
Let $X$ be a locally compact, finite dimensional CW $n${--}ad.
Let $\xi$ be a spherical fibration of dimension $\geq 2$ over $X$.
Then if $D(\xi)$ is a Spivak $(n+1)${--}ad, $X$ is a Spivak $n${--}ad.
\end{Corollary}

\smallskip\begin{proof}
Let $[X]=U_\xi\ \cap\ [\xi\>]$. 
Then $\cap [X]$ induces isomorphisms
\[H^\ast_\cmpsup\bigl(\coverFA{D(\xi)};\Z\bigr)\ \to\
H_{N-\ast}\bigl(\coverFA{D(\xi)},\coverFA{D(\xi\>\vert_{\partial X})};\Z)\ ,\]
or equivalently
$H^\ast_\cmpsup\bigl(\coverFA{X};\Z\bigr)\ \to\
H_{N-\ast}\bigl(\coverFA{X},\coverFA{\partial X};\Z)$.
Inducting over the $n${--}ad structure of $X$ and applying Theorems
\shortFullRef{T.2.2.2} and \shortFullRef{T.2.2.3}, 
we get $X$ is a Spivak $n${--}ad.
\end{proof}

\bigskip
\BEGIN{T.2.2.5}
\hskip-5pt\setcounter{footnote}{0}\footnote{Expanded statement of theorem and added to proof.}
$X$ is a Spivak $n${--}ad \iff\ $\coverFA{X}$ is for any cover of $X$:
moreover ${\text{ tr}}\bigl([X]\bigr)=[\coverFA{X}]$.
If $X$ is an $n${--}ad and $Y$ is an $m${--}ad,
$X\times Y$ is a Spivak $(n+m-1)${--}ad \iff\ 
$X$ is a Spivak $n${--}ad and $Y$ is a Spivak $m${--}ad:
$[X\times Y]=[X]\times[Y]$.
\end{Theorem}

\medskip\begin{proof}
Our first statement if immediate from \fullRef{T.2.2.1},
since if $N$ is an s-r neighborhood for $X$, $\coverFA{N}$ is one
for $\coverFA{X}$. 
Since ${\text{ tr}}$ preserves the fundamental class of manifolds and
since $[X]=U_\xi\ \cap\ [N]$, 
$[\coverFC{X}]=U_\xi\ \cap\ [\coverFC{N}]$, we see
${\text{ tr}}\bigl([X]\bigr)=[\coverFA{X}]$.

The ``adic'' part of the product result is easy once we see
$(X\times Y, \partial X\times Y\ \cup X\times\partial Y)=
(X,\partial X)\times (Y,\partial Y)$ is a Spivak pair.
To see this we first need

\medskip
\BEGIN{L.2.2.2}
If $\nu_Z$ is the Spivak normal fibration 
for an finite dimensional $n${--}ad $Z$, and 
if $X$ and $Y$ are such complexes,
$\nu_X\ \ast\ \nu_Y=\nu_{X\times Y}$.
\end{Lemma}

\smallskip\begin{proof}
Let $D_Z$ be the ``disc'' fibration of $\nu_Z$.
Then
$\nu_X\ \ast\ \nu_Y=\nu_X\times D_Y\ \cup\ D_X\times\nu_Y
\subseteq D_X\times D_Y$ denotes the fibrewise join.
Let $(N\Colon N_1,N_2)$ be the s-r neighborhood from which $\nu_X$
was formed.
$(M\Colon M_1, M_2)$ is the corresponding object for $\nu_Y$.
Then $\nu_{X\times Y}$ has
${\mathcal P}(N\times M_1\ \cup\ N_1\times N, N\times M,X\times Y)$ for
total space.
$\nu_X\times D_Y$ consists of triples $e\in\nu_X$, $f\in\nu_Y$ and
$t\in[0,1]$ with $(f,0)=(g,0)$ if $f(1)=g(1)\in Y$ ($D_Y$ is the fibrewise
cone on $\nu_Y$).
There is a similar map for $D_X\times\nu_Y$ which agrees with the first
on $\nu_X\times\nu_Y$. 
Hence we get a fibre map
$\nu_X\ \ast\ \nu_Y\ \to\ \nu_{X\times Y}$.
Now $\nu_X$ restricts from a fibration $\nu^\prime_X$ over all of $N$.
$\nu^\prime_Y$ and $\nu^\prime_{X\times Y}$ are defined similarly,
we have a fibre map 
$\nu^\prime_X\ \ast\ \nu^\prime_Y\ \to\ \nu^\prime_{X\times Y}$, and
a homotopy equivalence $\nu_Z\subseteq \nu^\prime_Z$.
There is an initial point map $\nu^\prime_{X\times Y}\ \to\ 
N\times M_1\ \cup\ N_1\times M$, which is a homotopy equivalence.
$\nu^\prime\ \ast\ \nu^\prime_Y\ \to\ 
N\times M_1\ \cup\ N_1\times M$ via the composition is likewise a
homotopy equivalence.
Hence by Dold \cite{bseven}, 
$\nu^\prime\ \ast\ \nu^\prime_Y\ \to\ \nu^\prime_{X\times Y}$
is a fibre homotopy equivalence.
Hence so is $\nu\ \ast\ \nu_Y\ \to\ \nu_{X\times Y}$.
\end{proof}

Now $X\times Y$ is a Spivak ad \iff\ the fibre of $\nu_{X\times Y}$ 
has the homology of a sphere,
and, if $Z\subseteq X\times Y$ is a piece of the ``adic'' structure,
$\nu_{X\times Y}\vert_Z\cong \nu_Z$.
Since the fibre of the fibrewise join is the join of the fibres, the
fibre of $\nu_{X\times Y}$ has the homology of a sphere \iff\ 
the fibres of $\nu_X$ and $\nu_Y$ do.
The equation $[X\times Y]=[X]\times[Y]$ also follows.
If $\partial X=\partial Y=\emptyset$, we are done.

For any ``adic'' piece $Z\subseteq X\times Y$, there is a fibre map
$\nu_Z\subseteq \nu_{X\times Y}\vert_Z$.
Next suppose $\partial Y=\emptyset$ and $X$ and $Y$ are Spivak.
It follows easily from the lemma that $\partial X\times Y$ is a Spivak space,
and, since $\nu_X\vert_{\partial X}\cong \nu_{\partial X}$ is
an equivalence, that $X\times Y$ is a Spivak pair.
Conversely, if $X\times Y$ is a Spivak pair, $\partial X\times Y$
is a Spivak space so $\partial X$ and $Y$ are.
$\nu_X$ has fibres that have the homology of spheres and since
$\nu_X\ \ast\ \nu_Y\vert_{\partial X\times Y}\cong 
\nu_{\partial X}\ \ast\ \nu_Y$ it follows that $\nu_X\vert_{\partial X}\cong
\nu_{\partial X}$ so $(X,\partial X)$ is a Spivak pair.

\medskip
\begin{Lemma}
\footnote[2]{This was a remark tossed off in the original which
now seems a bit harder.}
If $(Y,X)$ and $(Y^\prime, X)$ are Spivak pairs, $Y\ \cup_X\ Y^\prime$
is a Spivak space.
\end{Lemma}

\smallskip\begin{proof}
We choose our neighborhoods with care. 
Embed $X$ in $\R^{K-1}$ for $K$ much bigger than the dimension of $X$. 
Then we can extend this embedding to embeddings $Y$ in $\HS^K$ and
$Y^\prime$ in $\HS^K$.
Let $(M\Colon M_1,M_2)$ be a regular neighborhood for $(Y,X)$ and
let $(M^\prime\Colon M^\prime_1,M^\prime_2)$ be a regular neighborhood for $(Y^\prime,X)$. 
Write $Z=Y\ \cup\ Y^\prime$ and note $M\ \cup\ M^\prime \subset \R^K$
is a regular neighborhood for $Z$.
Let $\nu_X={\mathcal P}(M_1\cap M_2,M_2,X)$; 
$\nu^\prime_X={\mathcal P}(M^\prime_1\cap M^\prime_2,M^\prime_2,X)$; 
$\nu_Y={\mathcal P}(M_1,M,Y)$; 
$\nu_{Y^\prime}={\mathcal P}(M^\prime_1,M^\prime,Y^\prime)$;
and $\nu_Z={\mathcal P}(M_1\ \cup_{h\vert_{M_1\cap M_2}}\ M^\prime_1,
M\ \cup_h\ M^\prime,Z)$.
$h$ gives a fibre equivalence between $\nu_X$ and $\nu^\prime_X$
and the natural maps $\nu_X\subseteq\nu_Y\vert_X$, 
$\nu^\prime_X\subseteq\nu_{Y^\prime}\vert_X$ are fibre
homotopy equivalences.
Let $\nu^\prime_Z$ be the pushout
\topD{10}{$\begin{matrix}%
\nu_x&\to&\nu_Y\cr
\downarrow&&\downarrow\cr
\nu_{Y^\prime}&\to&\nu^\prime_Z\cr
\end{matrix}$}{4}
and note we have a projection $\nu^\prime_Z\to Z$
and a fibre map $\nu^\prime_Z\subseteq \nu_Z$.
$\nu^\prime_Z\ \to\ Z$ is a Dold fibration \cite{bseven}
and the map $\nu^\prime_Z\subseteq \nu_Z$ is a homotopy
equivalence, both being homotopy equivalent to $Z$.
Hence the map is a fibre homotopy equivalence and the
fibres of $\nu_Z$ have the homology of spheres.
\end{proof}

Now assume $X$ and $Y$ are pairs.
Then $\partial X\times (Y,\partial Y)$ and $(X,\partial X)\times \partial Y$
are Spivak pairs.
Then $\partial (X\times Y)$ is a Spivak space by the last lemma.
Since the fibre of $\nu_{X\times Y}$ has the homology of a sphere,
$\bigl(X\times Y, \partial (X\times Y)\bigr)$ is a Spivak pair,
and so $X\times Y$ is a Spivak triad.

Conversely, if $X\times Y$ is a Spivak triad, $\partial X\times\partial Y$
is a Spivak space, so $\partial X$ and $\partial Y$ are.
$X\times\partial Y$ and $\partial X\times Y$ are Spivak pairs,
so $X$ and $Y$ must be Spivak pairs as well.

The general case follows by induction.
\end{proof}

\bigskip
\BEGIN{T.2.2.6}
Let $X$ be a Spivak $n${--}ad, and let $N$ be a regular neighborhood
for $X$.
If $\bigl(D(X),S(X)\bigr)$ is a simple \cwation\ for this normal fibration,
there is a proper map of $(n+1)${--}ads $g\colon N\ \to\ D(X)$
such that the composition $N\ \to\ D(X)\ \to\ N$ is a proper
homotopy inverse for $X\ \to\ N$.
If $[N]$ and $\bigl[D(X)\bigr]$ are the fundamental classes for $N$ and
$D(X)$ respectively, $g_\ast[N]=\bigl[D(X)\bigr]$.

$g$ is a homotopy equivalence of $(n+1)${--}ads (not necessarily a proper
homotopy equivalence)
$g$ is however properly $(\dim N - \dim X -1)$ connected.
\end{Theorem}

\medskip\begin{xRemarks}
If $[X]$ lives in $k${--}dimensional homology, then the normal fibration
has a simple \cwation\ if $\dim N -k \geq 3$.
\end{xRemarks}

\medskip\begin{proof}
To be momentarily sloppy, let $D(X)$ denote the total space of the 
normal disc fibration for $X$.
Since $X\ \to\ N$ is a proper homotopy equivalence, pick
an inverse $N\ \to\ X$. 
Pull $D(X)$ back over $N$.
It is also a disc fibration and so has a section (see Dold \cite{bseven}
Corollary 6.2).
Map $N\ \to\ D(X)$ by the section followed by the map into $D(X)$.
Under the composition $N\ \to\ D(X)\ \to\ X$, we just get our original
map.
But now we can take the map from the total space of the fibration to
the \cwation. 
Letting $D(X)$ be the disc \cwation\ again, we get a map $N\ \to\ D(X)$
so that the map $N\ \to\ D(X)\ \to\ X$ is a proper homotopy inverse
to our original map $X\ \to\ N$.
The map $N\ \to\ D(X)$ is easily seen to be proper and is $g$.

$g$ is a homotopy equivalence of $(n+1)${--}ads by construction.
$g_0\colon N\ \to\ D(X)$ is also a proper homotopy equivalence
($g\colon N\ \to\ D(X)$ is not necessarily a proper homotopy
equivalence of $n$-ads).
The following diagram commutes
\[\begin{matrix}%
H^\ast(N,\partial N)&\LA{\ g^\ast\ }&H^\ast\bigl(D(X),C(X)\bigr)\cr
\downlabeledarrow[\big]{}{\capf{[N]}{}}&&
\downlabeledarrow[\big]{}{g_\ast\capf{[N]}{}}\cr\noalign{\vskip4pt}
H^\locf_{N-\ast}(N)&\RA{(g_0)_\ast}&
H^\locf_{N-\ast}\bigl(D(X)\bigr)\cr
\end{matrix}\]
$\capf{[N]}{}$, $g^\ast$ and $(g_0)_\ast$ are all isomorphisms, so
$g_\ast[N]$ is also an isomorphism.
Therefore $g_\ast[N]=\pm\bigl[D(X)\bigr]$ and we may orient $N$ so that
$g_\ast[N]=\bigl[D(X)\bigr]$. 

The map $C(X)\ \to\ X$is properly $q${--}connected, where the normal
spherical fibration has fibre $S^{q}$.
This is seen from the fibration sequence $S^q\ \to\ S(\xi)\ \to\ X$,
where $S(\xi)$ is the total space of the normal spherical fibration, 
by noticing that
\[\begin{matrix}%
\Delta(S(\xi)\Colon \pi_k)&\to&\Delta(C(X)\Colon \pi_x)\cr
&\hbox to 0pt{\hss$\searrow\hskip 60pt\swarrow$\hss}\cr
&\Delta(X\Colon \pi_k)\cr
\end{matrix}\]
commutes, where $\Delta(S(\xi)\Colon \pi_k)$ is formed from the groups
$\pi_k(S(\xi)\vert_{X-C},\hat p)$, where $\hat p\in S(\xi)$ covers
$p\in X$ (i.e. just pick one $\hat p$ for each base point in $X$).
The horizontal map is an isomorphism since $C(X)$ is a \cwation.
The first vertical map is an isomorphism for $k<q$ and an epimorphism
for $k=q$, so we are done.

The map $N_1\subseteq N\ \to\ X$ is properly $r${--}connected, where
$r=(\dim N - \dim X - 1)$.
This is seen by showing that the map $N_1\subseteq N$ is properly
$r${--}connected.
But this is easy.
If $K$ is a locally compact complex with $\dim K\leq r$, any map of
$K\ \to\ N$ deforms properly by general position to a map whose
image lies in $N-X$, and so can be properly deformed into $N_1$.
Hence $\Delta(N_1\Colon \pi_k)\ \to\ \Delta(N\Colon \pi_k)$ is onto for $k\leq r$
and $1${--}$1$ for $k\leq r-1$.

Now $g_0\colon N\ \to\ D(X)$ is a proper equivalence so the map
is properly $r${--}connected.
Since $r\leq q$, $g\colon N_1\to\ C(X)$ is properly $r${--}connected.
If $X$ is a space, we are done.
If $(X,\partial X)$ is a pair, the regular neighborhood is $(N\Colon N_1,N_2)$
and the \cwation\ is $\bigl(D(X)\Colon C(X), D(\partial X)\bigr)$.
$C(X)\ \cap\ D(\partial X)= C(\partial X)$.
$g\colon N_1\ \cap\ N_2\ \to\ C(\partial X)$ 
is properly $r${--}connected as it is an example of the absolute case.
$N_2\ \to D(\partial X)$ is a proper homotopy equivalence, hence
properly $r${--}connected.
$g\colon N_1\ \to\ C(X)$ and $g\colon N\ \to\ D(X)$ we saw were 
properly $r${--}connected, so the case for pairs is done.
For the $n${--}ad case, just induct. 
\end{proof}

\medskip
We are now ready to define Poincar\'e duality spaces.

\medskip\begin{xDefinition}
A Spivak $n${--}ad is a {\sl Poincar\'e duality $n${--}ad\/} \iff\ the $g$ of
\fullRef{T.2.2.6} is a proper homotopy equivalence
of $(n+1)${--}ads for some regular neighborhood.
\end{xDefinition}

\medskip\begin{xRemarks}
A priori our definition depends on which regular neighborhood we have
used in \fullRef{T.2.2.6}.
In fact this is not the case as our next theorem demonstrates.
\end{xRemarks}

\bigskip
\BEGIN{T.2.2.7}
Let $X$ be a locally finite, finite dimensional CW $n${--}ad.
Then $X$ is a Poincar\'e space \iff\ $X$ satisfies Poincar\'e duality with respect
to $[\tilde X]\in H^\locf_{N}(\coverFA{X},\coverFA{\partial X};\Z)$
and with respect to a universal covering functor.

A pair $(X,\partial X)$ is a Poincar\'e pair \iff\ $\partial X$ is a Poincar\'e space
and $X$ satisfies Poincar\'e duality with respect to a universal covering
functor and a class $[\tilde X]\in H^\locf_N(\coverFA{X},
\coverFA{\partial X};\Z)$ such that 
$\partial[\tilde X]=[\coverFA{\partial X}]$.

A similar result holds for $n${--}ads.
\end{Theorem}

\medskip\begin{proof}
Since $\capf{[\tilde X]}{}\colon 
\Delta^\ast(X\Colon \coverFA{})\ \to\ 
\Delta_{N-\ast}(X\Colon \coverFA{})$ an isomorphism implies
$\capf{[\tilde X]}{}\colon 
H^\ast_\cmpsup(\coverFA{X})\penalty-1000 \to\ 
H^\locf_{N-\ast}(\coverFA{X})$ is an isomorphism, if $X$ satisfies Poincar\'e 
duality then, by \fullRef{T.2.2.2}, $X$ is a Spivak space.
Similarly, by \fullRef{T.2.2.3}, we may show
$(X,\partial X)$ is a Spivak space if $\partial X$ is a Poincar\'e 
and if $(X,\partial X)$ satisfies Poincar\'e duality.
In both cases, the fundamental class ,$[X]$, 
transfers up to give $\pm[\tilde X]$. 
Now look at
\[\begin{matrix}%
\Delta^\ast(X,\partial X\Colon \coverFA{})&\RA{\ r^\ast\ }&
\Delta^\ast(N,N_2\Colon \coverFA{})\cr
&&\downlabeledarrow[\big]{}{\capf{[N]}{}}\cr\noalign{\vskip4pt}
&&\Delta_{n+k-\ast}(N,N_1\Colon \coverFA{})&
\RA{\ g_\ast\ }&
\Delta_{n+k-\ast}(D(X),C(X)\Colon \coverFA{})\cr
&&&&\downlabeledarrow[\big]{}{U_{\nu}\cap}\cr\noalign{\vskip4pt}
&&&&\Delta_{n-\ast}(X\Colon \coverFA{})\cr
\end{matrix}\]
where $r\colon (N,N_2)\ \to\ (X,\partial X)$ is a proper homotopy inverse
for $(X,\partial X)\subseteq (N,N_2)$, and $U_\nu$ is the Thom class
for the normal fibration $\nu$.
By \fullRef{T.2.2.4}, the composition is just
$\capf{[X]}{}$, and $r^\ast$, $\capf{[N]}{}$, and $U_\nu\cap$ are 
all isomorphisms.
Hence $(X,\partial X)$ satisfies Poincar\'e duality \iff\ $g_\ast$ is an isomorphism.

If $g$ is a proper homotopy equivalence, $g_\ast$ is clearly an isomorphism.

If $(X,\partial X)$ satisfies Poincar\'e duality, and if $\dim N - \dim X\geq 3$,
$g$ is a proper homotopy equivalence by the Whitehead theorem.
To see this, first note \ $\coverFA{}$ is a universal covering
functor for both $N$ and $D(X)$.
Since $\dim N - \dim X  \geq 3$, $N_1\subseteq N$ and 
$C(X)\subseteq D(X)$ are at least properly $2${--}connected.
Since $\partial X$ is by hypothesis a Poincar\'e duality space,
$g\ast\colon\Delta_\ast(N_1\Colon \coverFA{})\ \to\ 
\Delta_\ast(C(X)\Colon \coverFA{})$ is an isomorphism.
By the connectivity of $N_1\subseteq N$ and $C(X)\subseteq D(X)$.
these groups are already the subspace groups for a wise choice
of base points.
By the Browder lemma
$g_\ast\Delta(N\Colon \coverFA{})\ \to\ 
\Delta_\ast(D(X)\Colon \coverFA{})$ is an isomorphism, and $g$ is
at least properly $2${--}connected, so the Whitehead theorem applies
to show that $g$ is a proper homotopy equivalence.
\end{proof}

\medskip\begin{xRemarks}Note that the proof shows that if $X$
is Poincar\'e, $g$ must be a proper homotopy equivalence whenever 
$\dim N - \dim X\geq 3$.
\end{xRemarks}

\medskip
We have seen that manifolds satisfy Poincar\'e duality 
with respect to any covering functor.
The Thom isomorphism theorem also holds for any covering functor.
Hence it is easy to see

\medskip
\BEGIN{C.2.2.7.1}
A Poincar\'e duality $n${--}ad satisfies Poincar\'e duality with respect to
any covering functor.
\end{Corollary}

\medskip\begin{xDefinition}
The torsion of the equivalence $\capf{[X]}{}\colon 
\Delta^\ast(X,\partial X\Colon \coverFA{})\ \to\ 
\Delta_{N-\ast}(X\Colon \coverFA{})$ is defined to be the
torsion of the Poincar\'e duality space $X$ ($\coverFA{}$
is the universal covering functor).
Since $\bigl(D(X), C(X)\bigr)$ is a simple \cwation,
and since $\capf{[N]}{}$ is a simple equivalence 
(\fullRef{T.2.1.2}),
$\tau(X)=(-1)^{N+k}\tau(g)$, where $\tau(X)$ is the torsion of $X$
and everything else comes from the diagram in the proof of 
\fullRef{T.2.2.7}.
A simple Poincar\'e $n${--}ad is one for which all the duality maps are simple.
\end{xDefinition}

\bigskip\begin{xExamples}
By Theorems \shortFullRef{T.2.1.1} and \shortFullRef{T.2.1.2}, any paracompact manifold
$n${--}ad is a simple Poincar\'e $n${--}ad.
There are also examples of Spivak spaces which are not Poincar\'e duality
spaces.
One such is the following.
Let $X$ be a finite complex whose reduced homology 
with integer coefficients is zero, but which is not contractible.
(The dodecahedral manifold minus an open disc is such an example.)
Look at $\open{C}(X\vee S^2)$, the open cone on $X\vee S^2$.
The obvious map
$\R^3=\open{C}(S^2)\ \to \ \open{C}(X\vee S^2)$ is seen to induce
isomorphisms on $\pi_1$, $H_\ast$ and $H^\ast_\cmpsup$.
Since $\R^3$ is a Spivak space, so is $\open{C}(X\vee S^2)$.
$\open{C}(X\vee S^2)$ is not a Poincar\'e duality space as $X\vee S^2$
is not a Poincar\'e duality space.
\end{xExamples}

\medskip
In the other direction, we have as an application of 
a theorem of Farrell{--}Wagoner \cite{bnine}

\medskip
\BEGIN{T.2.2.8}
Let $X$ be a locally compact complex with monomorphic ends.
Then $X$ is a Poincar\'e duality space \iff\ $X$ is a Spivak space.
\end{Theorem}

An analogous result is true for $n${--}ads.

\medskip
\BEGIN{C.2.2.8.1}
Let $X$ be a Spivak $n${--}ad.
Then $X\times \R^2$ is a Poincar\'e duality $n${--}ad.
\end{Corollary}

\medskip\begin{proof}
We only prove $X$ Spivak implies $X$ Poincar\'e.
If $X$ has monomorphic ends,  and if $N$ is an s-r neighborhood with
$\dim N - \dim X\geq 3$, $\partial N$ has monomorphic ends.
$C(X)$ also has monomorphic ends.
The $g$ of \fullRef{T.2.2.6} is at least properly
$2${--}connected.
Hence by \cite{bnine} we need only prove $g$ induces isomorphisms
on $H_\ast$ and $H^\ast_\cmpsup$.
But $g\ast[N]=\bigl[D(X)\bigr]$, and $g$ on homology is an
isomorphism since it is a homotopy equivalence.
Since $N$ and $C(X)$ are both Spivak spaces, 
\fullRef{T.2.2.1} shows $g$ induces isomorphisms
on $H^\ast_\cmpsup$.

To show the corollary, observe that if $X$ is not compact, $\times\R^2$
has monomorphic ends.
It is a Spivak space by \fullRef{T.2.2.5}, so, in
this case, we are done.
If $X$ is compact, $X$ is already a Poincar\'e duality space, so the result
will follow from the next theorem.
\end{proof}

\bigskip
\BEGIN{T.2.2.9}
Let $X$ be a Poincar\'e duality $n${--}ad,  and 
let $Y$ be a Poincar\'e duality $m${--}ad. 
Then $X\times Y$  is a Poincar\'e duality $(n+m-1)${--}ad.
If $X$ or $Y$ is compact, the converse is true.
\end{Theorem}

\medskip\begin{proof}
From \fullRef{L.2.2.2} we have
$C(X\times Y)=D(X)\times C(Y)\ \cup C(X)\times D(Y)\subseteq
D(X)\times D(Y)=D(X\times Y)$.
If $N$ is an s-r neighborhood for $X$ and if $M$ is one for $Y$,
$N\times Y$ is one for $X\times Y$.
Hence we have $g\times f\colon N\times M\ \to\ D(X)\times D(Y)$
is a map of $(n+m+1)${--}ads.
It is a proper homotopy equivalence if $f$ and $g$ are.

Now suppose $X$ is compact.
By \fullRef{T.2.2.5}, $X$ is a Spivak $n${--}ad,
and hence a Poincar\'e $n${--}ad.
Since $g\times f$ is a proper homotopy equivalence, it induces
isomorphisms on the proper homotopy groups.
We claim $\Delta(N\times M\Colon \pi_k)=\pi_k(N)\oplus\Delta(M\Colon \pi_k)$
for $N$ compact.
This is easily seen by using the cofinal collection of compact subsets
of $N\times M$ of the form $N\times C$, $C\subseteq M$ compact.
A similar result computes $\Delta\bigl(D(X\times Y)\Colon \pi_k\bigr)$.
Since $g\times f$ and $g$ induce isomorphisms,
$f_\ast\Delta(M\Colon \pi_k)\ \to\ \Delta(D(Y)\Colon \pi_k)$ is an isomorphism.
By inducting this argument over the various subspaces of $D(Y)$,
$f$ is seen to be a proper homotopy equivalence of $(m+1)${--}ads.
Hence $Y$ is a Poincar\'e duality $m${--}ad.
\end{proof}

\bigskip
\BEGIN{T.2.2.10}
$X$ a Poincar\'e duality $n${--}ad implies $\coverFA{X}$ is a Poincar\'e duality 
$n${--}ad for any cover of $X$.
If $X$ is compact or if $\coverFA{X}$ is a finite sheeted cover,
then the converse is true.
\end{Theorem}

\medskip\begin{proof}
Let $N$ be an s-r neighborhood for $X$.
Then $\coverFA{N}$ is an s-r neighborhood for $\coverFA{X}$, so
$D(\coverFA{X})=\coverFA{D(X)}$.
$X$ a Poincar\'e duality $n${--}ad implies $N\ \to\ D(X)$ is a proper
homotopy equivalence of $n${--}ads.
But then so is
$\coverFA{N}\ \to\ \coverFA{D(X)}$, so $\coverFA{X}$ is a Poincar\'e 
duality $n${--}ad.

If $\coverFA{X}$ is a Poincar\'e duality $n${--}ad, $X$ is a Spivak $n${--}ad
by \fullRef{T.2.2.5}.
Hence if $X$ is compact, it is a Poincar\'e duality $n${--}ad.

Now if $\coverFA{X}\ \to\ X$ is finite sheeted and we know 
$\coverFA{N}\ \to\ \coverFA{D(X)}$ is a proper homotopy
equivalence of $(n+1)${--}ads, we must show $N\ \to\ D(X)$ is a proper 
homotopy equivalence of $(n+1)${--}ads.
But if $\dim N - \dim X\geq 3$ (which we may freely assume), this map
is properly $2${--}connected.
Since $\Delta(\coverFA{N}\Colon \pi_k)\ \to\ \Delta(N\Colon \pi_k)$
is an isomorphism for $k\geq2$ when $\coverFA{N}$ is a finite sheeted
cover, $N\to D(X)$ is seen to be a proper homotopy equivalence.%
\footnote[1]{It is easy to be too naive about covers 
verses covering functors: indeed
this last sentence is not right.}
To proceed correctly, we repair the error in this last sentence.
The problem begins with the trees: the tree for $\coverFA{N}$
is a cover of the tree for $N$ and so may have a different end
structure. 
The groups at corresponding vertices are isomorphic but there tend
to be several vertices in $\coverFA{N}$ for each one in $N$.
What we do see is that the proper map $\coverFA{N}\ \to\ N$
induces epimorphisms $\Delta(\coverFA{N}\Colon \pi_k)\ \to\ \Delta(N\Colon \pi_k)$,
$k\geq 2$.
Now apply this remark to the relative homotopy groups.
The map is still epic and the domain is $0$, hence so is the range.
This argument can be applied to any piece of the $(n+1)${--}ad structure,
so $X$ is a Poincar\'e duality $n${--}ad.
\end{proof}

\bigskip\begin{xRemarks}
The full converse to Theorems \shortFullRef{T.2.2.9} 
and \shortFullRef{T.2.2.10} are false.
Let $X$ be any Spivak space which is not a Poincar\'e duality space.
Then $X\times\R^2$ is a counterexample to the converse of 
\Ref{T.2.2.9} as it is a Poincar\'e duality space
by \Ref{C.2.2.8.1}.
$X\times T^2$ is a counterexample to 
\Ref{T.2.2.10}, since $X\times T^2$ is not a
Poincar\'e duality space by \Ref{T.2.2.9},
but its cover $X\times\R^2$ is.
\end{xRemarks}

\bigskip
\BEGIN{T.2.2.11}
Let $\xi$ be any spherical fibration of dimension $\geq 2$ over
a locally compact, finite dimensional CW $n${--}ad $X$.
Then $X$ is a Poincar\'e duality $n${--}ad \iff\ $D(\xi)$ is a Poincar\'e duality
$(n+1)${--}ad.
\end{Theorem}

\medskip\begin{proof}
By \fullRef{T.2.2.4} or \fullRef{C.2.2.4.1}, we may assume $X$ and
$D(\xi)$ are Spivak ads, and we have the formula $U_\xi\ \cap\ 
[\xi\>] = [X]$.
Since the Thom isomorphism is valid for the $\Delta$ theory 
(see the appendix page \pageref{PoincareAppendix}), $\capf{[X]}{}$ is an isomorphism \iff\ 
$\capf{[\xi\>]}{}$ is an isomorphism.
Since $\dim \xi\geq 2$, a universal covering functor for $X$
induces one for $D(\xi)$.
\fullRef{T.2.2.7} now gives the desired
conclusions.
\end{proof}

\bigskip\begin{xRemarks}
The torsions of the Poincar\'e spaces occurring in Theorems
\shortFullRef{T.2.2.9}, \shortFullRef{T.2.2.10} and \shortFullRef{T.2.2.11} can be ``computed''.
In particular, $\tau(X\times Y)=A\bigl(\tau(X),\tau(Y)\bigr)$
where $A$ is the pairing $\sieb(X)\times\sieb(Y)\ \to\ \sieb(X\times Y)$
(see \fullRef{L.1.5.23} and the preceding discussion).
$\tau(\coverFA{X})={\text{ tr}}\ \tau(X)$,
where ${\text{ tr}}\colon \sieb(X)\ \to\ \sieb(\coverFA{X})$.
$\tau\bigl(D(\xi)\bigr) = (-1)^n\tau\bigl(D(\xi)\bigr)^t$, where
$n$ is the dimension of the fundamental class of $X$, and $t$ is
the transpose operation on $\sieb(D(\xi))$.
These formulas are not very hard to deduce and will be left to the
reader.
\end{xRemarks}

\medskip
We conclude this section by investigating the ``uniqueness'' of the Spivak 
normal fibration.
We first prove

\bigskip
\BEGIN{L.2.2.3}
Let $D(\xi)$ be a \cwation\ for some spherical fibration $\xi$ over
a Poincar\'e duality $n${--}ad.
If there is a stably parallelizable manifold $(n+1)${--}ad $N$ and a 
proper, degree one, homotopy equivalence $N\ \to\ D(\xi)$,
then $\xi$ is stably equivalent to the Spivak normal fibration.
\end{Lemma}

\medskip\begin{xRemarks}
Given all spherical fibrations over a Poincar\'e duality $n${--}ad $X$, 
we wish to determine which of these could be the normal fibration
of some complex having the same proper homotopy type as $X$.
In the compact case, Spivak showed that there was only one,
the one with the reducible Thom space.
\fullRef{L.2.2.3} shows that if $D(\xi)$
has the degree one proper homotopy type of a stably 
parallelizable manifold, then $\xi$ is the Spivak normal fibration
for $X$.
If $\xi$ is the normal fibration for some complex $Y$, $D(\xi)$ 
has the degree one proper homotopy type of a parallelizable
manifold, so in the non{--}compact case there is one and only
one Spivak normal fibration.
\end{xRemarks}

\medskip\begin{proof}
If the equivalence were simple, $N$ would be an s-r neighborhood
and this would follow from \fullRef{T.2.2.1}.
Now by Siebenmann \cite{bthirtythree}, 
$N\times S^1\ \to D(\xi)\times S^1$ is a simple equivalence.
$D(\xi)\times S^1$ is a simple \cwation\ for $\xi\times S^1$ over
$X\times S^1$.
$N\times S^1$ is an s-r neighborhood for $X\times S^1$.
${\mathcal P}(N_1,N,N)\times S^1\to\ N\times S^1$ makes the map
$N_1\times S^1\subseteq N\times S^1$ into a fibration,
so $\nu_X\times S^1$ is fibre homotopy equivalent to $\nu_{X\times S^1}$.
But $\xi\times S^1$ is stably fibre homotopy equivalent 
to $\nu\vert_{X\times S^1}$ by \fullRef{T.2.2.1}.
Hence $\nu_X$ is stably $\xi$.
\end{proof}

\bigskip
\BEGIN{T.2.2.12}
If $f\colon X\ \to\ Y$ is a proper homotopy equivalence 
between Poincar\'e duality $n${--}ads, then
$f^\ast\nu_Y\cong \nu_X$.
\end{Theorem}

\medskip\begin{proof}
Let $\xi=f^\ast(\nu_Y)$.
Then
\[\begin{matrix}%
D(\xi)&\to&\xi&\to&\nu_Y&\to&D(\nu_Y)\cr
&\searrow&\downarrow&&\downarrow&\swarrow\cr
&&X&\to&Y\cr
\end{matrix}\]
commutes.
The top horizontal row is a proper homotopy equivalence, as one easily
checks by applying $\Delta(\quad\Colon \pi_k)$ to everything.
Since $D(\nu_Y)$ has the degree one proper homotopy type of
a parallelizable manifold, so does $D(\xi)$.
Hence by \fullRef{L.2.2.3}, $\xi\cong\nu_x$.
\end{proof}

\medskip
Spivak's identification of the normal fibration actually proves a stronger
theorem.
We can prove this result as

\bigskip
\BEGIN{T.2.2.13}
Let $f\colon X\ \to\ Y$ be a degree one map of Poincar\'e duality
$n${--}ads.
If there is a spherical fibration $\xi$ over $Y$ such that
$f^\ast(\xi)\cong \nu_X$, then $\xi\cong\nu_Y$.
\end{Theorem}

\medskip\begin{proof}
\[\begin{matrix}%
D(\nu_X)&\to&\nu_X&\to&\xi&\to&D(\xi)\cr
&\searrow&\downarrow&&\downarrow&\swarrow\cr
&&X&\to&Y\cr
\end{matrix}\]
commutes, so it is not hard to show that the top row is a degree one map.
$U_{\nu_X}\ \cap\ \big[D(\nu_X)\bigr]=[X]$;
$U_{\xi}\ \cap\ \big[D(\xi)\bigr]=[Y]$; and $f_\ast[X]=[Y]$.
Hence the top row must take $\bigl[D(\nu_X)\bigr]$ to $\bigl[D(\xi)\bigr]$.
$D(\nu_X)$ has the proper homotopy type of a parallelizable $(n+1)${--}ad,
$N$ so there is a degree $1$ map $g\colon N\ \to\ D(\xi)$.
Since $N$ is parallelizable, there is a topological microbundle over $D(\xi)$
which pulls back to the normal bundle of $N$ (namely the trivial bundle).
If $\dim\xi\geq 2$ (which we may always assume) then the pair
$\Big(D\bigl(\xi(Z)\bigr),C\bigl(\xi(Z)\bigr)\Bigr)$, for $Z\subseteq Y$
as part of the $n${--}ad structure on $Y$, is properly $2${--}connected.
Hence by the remarks following \fullRef{T.3.1.2},
we can find a parallelizable manifold $M$ and a degree one proper homotopy
equivalence $M\ \to\ D(\xi)$.
By \fullRef{L.2.2.3}, $\xi\cong\nu_Y$.
\end{proof}

\medskip\begin{xRemarks}
Logically \fullRef{T.2.2.13} should follow \fullRef{T.3.1.2} in chapter 3.
We do not use the result until we are past that point so it does no harm
to include it here.
\end{xRemarks}

The chief purpose of \fullRef{T.2.2.13} is to
severely limit the bundles which can occur in a surgery problem.

\newpage
\section{The normal form for Poincar\'e duality spaces}
\newHead{II.3}
In order to get a good theory of surgery, one needs to be able to do
surgery on Poincar\'e duality spaces; at least one must be able to
modify fundamental groups.
The results of this section show that Poincar\'e duality spaces look like
manifolds through codimension $1$.
These results are a direct generalization of Wall \cite{bthirtynine} 
Section 2, especially pages 220{--}221.

\medskip\begin{xDefinition}
Let $X$ be a Poincar\'e duality $n${--}ad.
Then if $[X]\in H^\locf_n$, $X$ is said to have \emph{formal dimension $n$}.
($X$ is often said to have dimension $n$.)
\end{xDefinition}

\bigskip
\BEGIN{T.2.3.1}
Let $X$ be a Poincar\'e space of dimension $n\geq 2$.
The $X$ satisfies $D n$.
If $X$ is a connected Poincar\'e duality $m${--}ad $m\geq2$, of
dimension $n\geq 3$, then $X$ satisfies $D(n-1)$.
\end{Theorem}

\medskip\begin{proof}
This follows from definitions and \fullRef{T.1.6.2}.
\end{proof}

\bigskip
\BEGIN{T.2.3.2}
Let $X$ be a Poincar\'e duality space of dimension $n$, $\geq 4$.
Then $X$ has the proper homotopy type of $Y$, where $Y$ is a 
Poincar\'e duality space which is the union of two Poincar\'e duality pairs
$(Z,\partial H)$ and $(H,\partial H)$ where $H$ is a smooth manifold
of dimension $n$ formed from a regular manifold in $\R^n$ of
a given tree for $Y$ by adding $1$ handles along the boundary,
and where $Z$ is a subcomplex satisfying $D(n-2)$.
The torsion of this equivalence may have any preassigned value.
The map induced by inclusion $\Delta(H\Colon \pi_1)\ \to\ \Delta(Y\Colon \pi_1)$
is surjective.
\end{Theorem}

\medskip\begin{proof}
Let $\hat C_\ast$ be the dual chain complex for $X$, reindexed so that
there is a chain map $\capf{[X]}{}\colon \hat C_\ast\ \to\ 
C_{\ast}(X)$.
By \fullRef{T.1.6.3}, we can find a chain complex $Y$ with 
$C_\ast(Y)=\hat C_\ast$ in dimensions greater than $3$.
$C_3(Y)=\hat C_3\oplus junk$, and the complex $Y^2\ \cup\ junk$
satisfies $D2$.

Now we could have arranged things so that the only vertices of $X$
were the vertices of the tree.
This is seen as follows.
First we claim we can find a subcomplex $V\subseteq X$ which contains
all the vertices and such that $T\subseteq V$ is a proper deformation
retract.
We do this as follows.\footnote[1]{The original argument here was wrong.}
\hide
Let ${\mathcal S}=\{ U\ \vert\ U$ is  a $1${--}dimensional subcomplex of $X$,
$T\subset U$ with the inclusion a proper $0${--}equivalence, and
$U$ contains all the vertices of $X\ \}$.
${\mathcal S}\neq\emptyset$ as $X^1\in{\mathcal S}$.
${\mathcal S}$ is ordered by inclusion.
Let $U_1\supseteq U_2\supseteq\cdots$ be a totally ordered sequence 
in ${\mathcal S}$.
Then $\cap\ U_i$ is also in ${\mathcal S}$.
Let $V$ be a minimal element of ${\mathcal S}$, which exists by Zorn.
We claim $H_1(V)=0$, so if not, look at a cycle in $V$.
At least one of the $1${--}simplexes of the cycle is not in $T$ for $T$ has
no $1${--}cycles.
Let $V_1\subseteq V$ be all of $V$ less one of the one simplexes in
the cycle which is not in $T$.
Then $V_1$ is a subcomplex, $T\subseteq V_1$, and $V_1$ contains all
the vertices.
This contradiction shows $H_1(V)=0$.
Put a metric on the vertices of $X^1$ as follows.
If $p$ and $q$ are vertices, the length of a cellular path joining $p$
to $q$ is obtained by counting the number of $1${--}simplices
in the domain (or equivalently the number of vertices minus $2$).
Define $d(p,q)$ to be the minimum of the length amongst all cellular
paths joining $p$ to $q$.
It is not hard to see $d$ is a metric.
\endhide
Choose an increasing sequence of compact subcomplexes of the $1$-skeleton, 
$X^1$, $K_0\subseteq K_1\subseteq \cdots$, whose union is $X^1$.
Let $V$ be a subcomplex of $X^1$ with $T\subseteq V$.
Let $\{v_0, \cdots \}$ denote the vertices of $X^1-V$.
By the definition of a tree, there is a locally finite set of paths 
$\{\lambda_i\}$, with $\lambda_i$ joining $v_i$ to some vertex of $V$.
It is no problem to assume the $\lambda_i$ are cellular. 
\hide
We can refine the collection of $K_j$ if necessary to insure that
any $\lambda_i$ beginning at a vertex in $K_j$ lies entirely in $K_{j+1}$.
\endhide

Construct a sequence of increasing subcomplexes, $V_r$, of $X^1$,
two sets of points ( $\{ x^{(r)}_i\}$ the vertices in $X^1-V_{r}$
and $\{y^{(r)}_i\}$ the vertices of $V_r$ )
and two sets of locally finite cellular paths $\lambda^{(r)}_i\subset X^1$ 
(beginning at $x^{(r)}_i$ and ending at a vertex of $V_{r}$) 
and $\Lambda^{(r)}\subset V_r$ (beginning at $y^{(r)}_i$ and
ending at a vertex of $T$ )
inductively as follows.
$V_0=T$, $\lambda^{(0)}=\lambda_i$, $\Lambda^{(0)}_i$ the constant 
path at the relevant vertex.

$V_r$ is obtained from $V_{r-1}$ in two steps.
First adjoin all vertices $x^{(r-1)}_i$ of $X^1$ 
for which there exists a $j$ such that the path 
$\lambda^{(r-1)}_j$ has its next to the last vertex $x^{(r-1)}_i$.
Next add some $1${--}cells.
First, for a fixed $x^{(r-1)}_i$, there may be several $j$ satisfying
our condition: pick one, say $i_j$.
Then add the $1${--}cell from $x^{(r-1)}_i$ to $V_{r-1}$
given by the last $1${--}cell in $\lambda^{(r-1)}_{i_j}$.
We need only define $\lambda^{(r)}_i$ for the $ x^{(r)}_i$:
just define it to be the sub{--}path of $\lambda^{(r-1)}_i$ which
starts at $x^{(r)}_i$ and continues until it encounters a vertex in $V_r$.
This path is definitely properly contained in $\lambda^{(r-1)}_i$ since
the next to the last vertex in $\lambda^{(r-1)}_i$ is certainly in $V_r$.
(Of course $\lambda^{(r)}_i$ may be shorter than this.)
The $\lambda^{(r)}_i$ are a locally finite collection 
since they form a subcomplex of such a collection.
Define $\Lambda^{(r)}_i$ as follows.
If the vertex in question lies in $V_{r-1}$, use $\Lambda^{(r-1)}_i$.
Otherwise, look at the $1${--}cell out of  $y^{(r)}_i$. Its other end
lies in $V_{r-1}$ by definition so is $y^{(r-1)}_j$ for some $j$.
$\Lambda^{(r)}_i$ is the path which
follows the $1${--}cell and then $\Lambda^{(r-1)}_j$.
Check that the collection $\Lambda^{(r)}_i$ is a locally finite collection.

Let $V_\infty=\cup\ V_r$.
Consider any vertex $x^{(0)}_i$ and its path $\lambda^{(0)}_i$.
The distance from this vertex to the tree along this path is finite,
say $R$. 
Then $x^{(0)}_i\in V_R$ (it may of course land in a smaller $V_R$).
This is easy to check by induction on $R$.
It follows that $V_\infty$ contains all the vertices of $X^1$ and we
let $y^{(\infty)}_i$ be an enumeration of them.
Define $\Lambda^{(\infty)}_i\subset V_\infty$ as follows.
$y^{(\infty)}_i$ lies in some $V_r$ so define 
$\Lambda^{(\infty)}_i=\Lambda^{(r)}_i$ and note that 
$\Lambda^{(r+1)}_i$ is the same path so this is well{--}defined.

There is a deformation retraction 
$d\colon V_r \times[0,1]\to\ V_{r}$ of $V_r$ to $V_{r-1}$
obtained by collapsing the new $1${--}cells to the end attached
to $V_{r-1}$, so $V_\infty$ is a $1${--}complex with $H_1(V_\infty)=0$.
The paths $\Lambda^{(\infty)}_i$ show the inclusion $T\subseteq X_\infty$
is a proper $0${--}equivalence. 
It follows that $T\subseteq V_\infty$ is a proper homotopy equivalence
as desired.

Set $K=\Bar{V-T}$ and look at $X/K$.
The collapse map $X\to X/K$ is a proper homotopy equivalence.
For a proof, see \cite{bsix} Proposition 2.11, page 220.
Note that all the maps there may be taken to be proper.
$X/K$ has only the vertices of the tree for $0${--}cells.

Now, to return to our proof, we may assume $\hat C_n=C_0(X)$
has a generator for each vertex of our tree.
$\hat C_{n-1}$ has a generator for each $1${--}cell of $X$.
As in Wall \cite{bthirtynine} Corollary 2.3.2, each $(n-1)${--}cell
is incident to either two $n${--}cells, or to the same $n${--}cell twice.
Look at an attaching map $S^{n-1}\ \to\ X^{n-1}$ for an $n${--}cell.
This can be normalized to take a finite, disjoint, collection of discs
onto the $(n-1)${--}cells homeomorphically and to take the rest of
$S^{n-1}$ into the $(n-2)${--}skeleton. 
Each $(n-1)${--}cell eventually gets just two such discs mapped into it.
The $n${--}discs together with the $(n-1)${--}cells corresponding to the
$1${--}cells of the tree are seen to form a regular neighborhood in
$\R^n$ of the tree, and $H$ is obtained from this 
by attaching $1${--}handles.

If $Z$ is the part of $Y$ in dimensions $\leq n-2$ (or is $Y^2\ \cup\ junk$
if $n=4$), $Y=Z\ \cup_{\partial H}\ H$ where $H$ is formed 
from $n${--}discs corresponding to the $n${--}cells by attaching
$1${--}cells as indicated by the $(n-1)${--}cells.
Actually, we want to form the mapping cylinder of $\partial H\ \to\ Z$
and then tale the union along $\partial H$.
Since $H$ is a manifold, the result is clearly homeomorphic to $Y$.
We denote the mapping cylinder by $Z$ so $Y=Z\ \cup_{\partial H}\ H$,
and $\partial H$ is a subcomplex of $Z$ and hence $Y$.
Note that $Z$ still satisfies $D(n-2)$.
\vfill

Now $Z\subseteq Y$ is at least properly $2${--}connected, for
$Z$ always contains the $2${--}skeleton of $Y$.
Since $(H,\partial H)$ is a Poincar\'e duality space, 
\fullRef{T.2.1.4} says $(Z,\partial H)$ satisfies Poincar\'e 
duality with respect to the covering functor induced from the universal
covering functor for $Y$. 
But this is just the universal covering functor for $Z$ as $Z\subseteq Y$
is properly $2${--}connected.
$\partial H$ is a Poincar\'e duality space, so 
\fullRef{T.2.2.7} says $(Z,\partial H)$ is a Poincar\'e 
duality pair.
The statement about the torsion is contained in \fullRef{T.1.6.3},
so we finish by showing $\Delta(H\Colon \pi_1)\ \to\ \Delta(Y\Colon \pi_1)$
is onto.
Our proof is basically Wall \cite{bthirtynine} Addendum 2.3.3, but
is more complicated.
We too will use the construction of $Z$ and $H$ via the dual cell
decomposition.
In our original complex, there were $0${--}cells, $e^0_p$,
one for each $p$ a vertex of $T$.
There were $1${--}cells $e^1_i$ satisfying 
$\partial e^1_i=g_i e^0_p - e^0_q$
where $g_i$ is a loop at $p$.
The $g_i$ which occur generate $\Delta(Y\Colon \pi_1)$.
In the dual complex we have $n${--}cells, $e^n_p$ and $(n-1)${--}cells 
$e^{n-1}_i$ with 
$\displaystyle e^n_p=\sum_i (\pm g_i e^{n-1}_i)-\sum_j e^{n-1}_j$,
where the sign is given by the local coefficients on $Y$, and where the
sum runs over all $(n-1)${--}cells incident to $e^n_p$.
The core $1${--}disc of the handle corresponding to $e^{n-1}_i$
followed by the unique minimal path in $T$ from 
the endpoint of the $1${--}disc to its initial point point has homotopy
class $g_i$.
Hence $\Delta(H\Colon \pi_1)$ is onto $\Delta(Y\Colon \pi_1)$.
\end{proof}

\medskip
\BEGIN{C.2.3.2.1}
Let $X$ be a Poincar\'e duality space of dimension $3$.
Then $X$ has the proper homotopy type of $Y$, where $Y$ 
is the union of two Poincar\'e duality pairs $(Z,\partial H)$ and $(H,\partial H)$,
where $H$ is a regular neighborhood in $\R^3$ of a given tree for $X$,
and $Z$ is a subcomplex of $Y$ satisfying $D2$.
The torsion of this equivalence can be arbitrary.
\end{Corollary}

\medskip\begin{proof}
Using the dual cell decomposition as before, let $Z$ be the subcomplex
of $Y$ such that $\hat C_3=C_3(Y,Z)$ and such that $Z$ satisfies $D2$.
$\hat C_3$ has one $3${--}cell for each vertex of the tree.
Now there is a locally finite collection of paths from each $n${--}cell
to the vertex of the tree it represents.

Given $H$, a regular neighborhood of the tree in $\R^3$, we
describe a map $\partial H\ \to\ Z$ which extends to a map
$H\ \to\ Y$ such that the induced map $C_3(H,\partial H)\ \to\ 
C_3(Y,Z)$ is an isomorphism.
Hence $Z\ \cup_{\partial H}\ H$ has the proper homotopy type of $Y$
and we will be done.
The map is the following.
$H$ can be viewed as the connected sum of a collection of $n${--}discs,
one for each vertex of the tree, by tubes corresponding to the
$1${--}cells of the tree.
$H$ can then be properly deformed to the subcomplex consisting of
$n${--}discs joined by the cores of the connecting tubes.
$\partial H$ under this deformation goes to a collection 
of $(n-1)${--}spheres joined by arcs.
Map the $(n-1)${--}sphere to $Z$ by the attaching map of the
corresponding $n${--}cell in $Y$.
Map an arc between two such spheres to the paths to the tree,
and then along the unique path in the tree between the two vertices.
This map clearly has the necessary properties.
\end{proof}

\bigskip
\BEGIN{T.2.3.3}
Let $(X,\partial X)$ be a Poincar\'e duality pair of dimension $n$, $n\geq 4$.
Then $(X,\partial X)$ has the proper homotopy type of a Poincar\'e duality
pair $(Y,\partial Y)$ which is the union of a Poincar\'e duality pair
$(Z, \partial H\ \cup\ \partial Y)$ and 
a Poincar\'e duality pair $(H,\partial H)$, where $H$ is a regular neighborhood 
in $\R^n$ of any given tree for $Y$ with $1${--}handles added 
along the boundary, and $Z$ is a subcomplex of $Y$ satisfying $D(n-1)$.
The torsion of this equivalence may be given any preassigned value.
$\Delta(H\Colon \pi_1)\ \to\ \Delta(Y\Colon \pi_1)$ is onto.
\end{Theorem}

\medskip\begin{proof}
By \fullRef{T.2.3.2} or \fullRef{C.2.3.2.1}, we may assume $\partial X$
already looks like $K\ \cup\ M$, where $M$ is a regular neighborhood
for a tree of $\partial X$ in $\R^{n-1}$, and $K$ satisfies $D(n-2)$.

Let $\hat C_\ast$ be the dual complex for $C_{n-\ast}(X)$.
Then there is a chain map $\capf{[X]}{}\colon 
\hat C_\ast\ \to\ C_\ast(X,\partial X)$.
We apply \fullRef{T.1.6.4} to find a complex $Y$ with 
$C_\ast(Y)=C_\ast$ in dimensions greater than $3$ and with 
$\partial X \subseteq Y$.
$C_3(Y)=C_3\ \cup\ junk$.
Set $L$ to be $Y^{(n-1)}$.
Then $M\subseteq L$.
Normalize the attaching maps for the $n${--}cells as before.
If $Z=Y^{(n-2)}\ \cup\ M$ ($Y^2\ \cup\ junk\ \cup\ M$ if $n=4$),
then $Y=Z\ \cup\ H$ where $H$ has the advertised description. 
Notice $\partial H\ \cap\ \partial X$ can be $M$ if one likes.
As before, $(H, \partial H)$ is a Poincar\'e duality pair.
$\partial H\ \cap\ \partial X =M$, so
$\partial H=(\Bar{\partial H - M\ },\partial M)\ \cup\ (M,\partial M)$
and $\partial X=(K,\partial M)\ \cup\ (M,\partial M)$.
The rest of the proof proceeds as in the proof of
\fullRef{T.2.3.2}.
\end{proof}

\newpage
\section*{Appendix. The \cwation\ of a spherical fibration}
\label{II.A}
\renewcommand\rightmark{\Roman{chapter}.\ \nameref*{II.A}}
\LRTpageLabel{PoincareAppendix}
\bigskip
We recall the definition.
Let $\xi$ be a spherical fibration over a finite dimensional, 
locally finite CW $n${--}ad.
Assume $\xi\geq 2$.
Let $S(\xi)$ be the total space.
We seek an $n${--}ad $Y$, a proper map $f\colon Y\ \to\ X$,
and maps
$\vrule width 0pt depth 10pt height 12pt
\displaystyle S(\xi)\ {\scriptscriptstyle
\Atop{\Limitsarrow{10}{10}{\scriptscriptstyle g}{}{\rightarrowfill}}
{\Limitsarrow{10}{10}{}{\scriptscriptstyle h}{\leftarrowfill}}
}\ Y$ 
which commute with the two projections.
We also require that $Y$ have the proper homotopy type of a
locally compact, finite dimensional CW $n${--}ad.  $g\ \circ\ h$ must
be properly homotopic to the identity, and $h\ \circ\ g$ must
be fibre homotopic to the identity. 
We give $Y$ a simple homotopy type by finding an equivalent CW
complex for which the Thom isomorphism is simple.

We digress briefly to include a discussion of the Thom isomorphism.
If $D(\xi)$ is the total space of the disc fibration associated to $\xi$, 
we define $\Delta(D(\xi)\Colon h,\coverFA{})$ and
$\Delta(D(\xi), S(\xi)\Colon h,\coverFA{})$ to be the groups one
gets by applying the $\Delta$ construction to the groups
$h\bigl(\pi^{-1}(\coverFA{X-C},\hat p)\bigr)$ for $D(\xi)$ and
$h\bigl(\pi^{-1}(\coverFA{X-C},
\rho^{-1}(\coverFA{X-C},\hat p)\bigr)$ for $\bigl(D(\xi), S(\xi)\bigr)$,
where $p$ is a vertex of $X$, $\hat p$ is a lift of $p$ into 
$\coverFA{D(\xi)}$, and 
$\pi\colon \coverFA{D(\xi)}\ \to\ \coverFA{X-C}$ and
$\rho\colon \coverFA{S(\xi)}\ \to\ \coverFA{X-C}$ are the
projections for the fibrations over $\coverFA{X-C}$ by restriction
and pullback from $D(\xi)$ and $S(\xi)$ respectively.

Now the Thom class for $\xi$, $U_\xi$ goes under
$\coverFA{X-C}\ \to\ X-C\ \to\ X$
to the Thom class for $\coverFA{S(\xi)}$.
If $h$ is cohomology we modify the $\Delta$ groups above 
in the obvious manner.
We will denote by $\Delta_\ast(D(\xi)\Colon \coverFA{})$
the $\ast${-}\thx\ homology group with covering functor $\coverFA{}$\ .
$\Delta^\ast$ is the cohomology theory.
Then we have maps
$\cupf{U_\xi}{}\colon \Delta^{\ast}(D(\xi)\Colon \coverFA{})\ \to\
\Delta^{\ast+k}(D(\xi),S(\xi)\Colon \coverFA{})$ and
$\capf{U_\xi}{}\colon \Delta_{\ast+k}(D(\xi),S(\xi)\Colon \coverFA{})\ \to\ 
\Delta_\ast(D(\xi)\Colon \coverFA{})$.
They are easily seen to be isomorphisms.
The maps $h$ and $g$ induce isomorphisms of 
$\Delta_{\ast}(S(\xi)\Colon \coverFA{})$ and $\Delta_\ast(Y\Colon \coverFA{})$
with a similar result for cohomology (the reader should have no
trouble defining $\Delta_\ast(S(\xi)\Colon \coverFA{})$ or its cohomological 
analogue).
We also get isomorphisms of $\Delta_\ast(D(\xi),S(\xi)\Colon \coverFA{})$
and $\Delta_\ast(M_f,Y\Colon \coverFA)$, again with a similar result in 
cohomology.
Hence we can speak of a Thom isomorphism for the \cwation.

We first prove that if we can find a \cwation, we can give it a unique simple
homotopy type.
Let $C$ be a CW $n${--}ad the proper homotopy type of the \cwation\ $Y$
( $C$ locally compact, finite dimensional).
\topD{0}{$\begin{matrix}%
C&\RA{\ \rho\ }&Y\cr
&\hbox to 0pt{\hss$\hbox to 0pt{\hss$\scriptstyle f\circ\rho$}\searrow\hskip 10pt%
\swarrow\hbox to 0pt{$\scriptstyle{f}$\hss}$\hss}\cr
&X&
\end{matrix}$}{4}
commutes, where $\rho$ is a proper homotopy equivalence with 
$f\ \circ\ \rho$ cellular.
(It is easy to find such a $\rho$.)
Let $\tau_\rho$ denote the torsion of the corresponding Thom homology
isomorphism.
If $\lambda\colon K\ \to\ C$ is a proper homotopy equivalence, the
Thom isomorphism associated to $\rho\ \circ \lambda$ has torsion
$\tau_\rho+\rho_\ast \tau(\lambda)$ by 
\fullRef{L.1.5.22}.
Since we may pick $\tau(\lambda)$ arbitrarily, we can find a $\rho$
with $\tau_\rho=0$.

Suppose now we have $\lambda\colon K\ \to\ Y$ with $\tau_\lambda=0$,
$f\ \circ\ \lambda$ cellular).
Let $a\colon Y\ \to\ K$ be a proper homotopy inverse to $\lambda$.
Then
\[\begin{matrix}%
C&\RA{\ a\ \circ\ \rho\ }&K\cr
&\hbox to 0pt{\hss\hbox to 0pt{\hss$\scriptstyle f\circ\rho$}$\searrow\hskip 24pt%
\swarrow$\hbox to 0pt{$\scriptstyle{f\circ\lambda}$\hss}\hss}\cr
&X&
\end{matrix}\]
properly homotopy commutes.
We get a proper homotopy equivalence of pairs\hfill\penalty-10000
$F\colon (M_{f\circ\rho},C)\ \to\ (M_{f\circ\lambda},K)$
such that
$F\vert_C= a\ \circ\ \rho$, and
\topD{12}{$\begin{matrix}%
M_{f\circ\rho}&\RA{\ F\ }&M_{f\circ\lambda}\cr
&\hbox to 0pt{\hss$\searrow$}\hskip 10pt%
\hbox to 0pt{$\swarrow$\hss}\cr
&X&
\end{matrix}$}{4}
commutes.
By \fullRef{L.1.5.19}, $M_{f\circ\rho}\ \to\ X$
and $M_{f\circ\lambda}\ \to\ X$ are simple, so $F$ is a simple
equivalence.
The torsion of $F$ from $(M_{f\circ\rho},C)$ to $(M_{f\circ\lambda},K)$
is $\tau_\rho-\tau_\lambda=0$, so by \fullRef{T.1.5.1},
the torsion of $a\ \circ\ \rho$ on the subspace groups is zero.
But as $\dim\xi\geq 2$, $f\ \circ\ \rho$ and $f\ \circ\ \lambda$
are at least properly $2${--}connected.
Hence the subspace groups with the induced covering functor are the
absolute groups with the universal covering functor.
Hence $a\ \circ\ \rho$ is a simple homotopy equivalence, so the simple
homotopy type of a \cwation\ is unique.
\medskip

We now construct the promised $Y$.
Notice first that we can replace $X$ by any locally compact, 
finite dimensional CW $n${--}ad of the same proper homotopy type.
Hence we may as well assume $X$ is a locally finite simplicial $n${--}ad
of finite dimension.
This is seen as follows.
By \cite{beleven} Theorem 4.1 and Lemma 5.1, $X$ is the union of $A$
and $B$ where $A$ and $B$ are the disjoint union of finite complexes.
Each finite complex has the homotopy type of a finite simplicial complex,
and if a subcomplex is already simplicial, we need not disturb it.
Hence we get a locally finite simplicial complex $Y$ and a map 
$f\colon X\ \to\ Y$ by making subcomplexes of the form $C\ \cap\ D$,
$C\subset A$ and $D\subset B$ simplicial and then making $C$ and $D$
simplicial.
Then $Y=A^\prime\ \cup\ B^\prime$ where $f\colon A\ \to\ A^\prime$
and $f\colon B\ \to\ B^\prime$ are proper homotopy equivalences.
Also $f\colon E\ \to\ E^\prime$ is a proper homotopy equivalence
where $E=\{ C\ \cap\ D\ \vert\ C\subset A,\ D\subset B \}$.
The proper Whitehead Theorem shows $f$ 
is a proper homotopy equivalence.
$X$ being what it is, we can find open sets $C_i$ such that $X-C_i$
and $\Bar{C}_i$ are subcomplexes, each $\Bar{C}_i$ is compact,
and $\xi\vert_{\Bar{C}_i}$ is trivial.
Furthermore, $\cup\ C_i =X$, the $C_i$ are locally finite, and the $C_i$
are indexed by the positive integers.
We set $\displaystyle V_i=\mathop{\cup}_{j\leq i} C_j$
We can also find an increasing collection of open sets $U_i$ such that 
$U_i\subseteq V_i-C_i$, $\Bar{U}_i$ is compact, and 
$\displaystyle \mathop{\cup}_{i} U_i=X$.

We first construct spaces $Y_i$ and maps $g_i$ and $f_i$ inductively
so that
\[\begin{matrix}%
\xi\vert_{\bar{V}_i}&\RA{\ g_i\ }& Y_i\cr\noalign{\vskip6pt}
&\hbox to 0pt{\hss$%
\hbox to 0pt{\hss$\scriptstyle\pi\vert_{\Bar{V}_i}$}\searrow\hskip 20pt 
\swarrow\hbox to 0pt{$\scriptstyle f_i$\hss}
$\hss}\cr
&\Bar{V}_i\cr
\end{matrix}\leqno{A)}\]
commutes.

Let $Y_1=\Bar{V}_1\times S^k$, $k=\dim\xi\geq 2$.
$g_1$ and $f_1$ exist since $\xi\vert_{V_1}$ is trivial.
$f_1$ is just projection.
We now induct; i.e. we have
\begin{enumerate}
\item[1)] A space $Y_i$ and maps $g_i$ and $f_i$ such that A) commutes.
\item[2)] $g_i$ is a homotopy equivalence.
\item[3)] $Y_i=Y_{i-1}\  \cup_\rho\ \Bar{C}_i\times S^k$ via some
homotopy equivalence 
\item[]$\rho\colon Y_{i-1}\ \cap\ f^{-1}_{i-1}\bigl(
\Bar{V}_{i-1}\ \cap \Bar{C}_i\bigr)\ \to\ 
\bigl(\Bar{V}_{i-1}\ \cap\ \Bar{C}_i\bigr)\times S^k$.
\item[4)] $g_{i-1}\vert_{f^{-1}_{i-1}(U_{i-1})}=
g_{i}\vert_{f^{-1}_{i-1}(U_{i-1})}$ and $f_{i-1}\vert_{Y_{i-1}}=
f_{i}\vert_{Y_{i-1}}$.
\item[5)] Let ${\mathcal S}_r=\{ \Bar{C}_{i_1}\ \cap\ \Bar{C}_{i_2}\ \cap\ 
\cdots \Bar{C}_{i_r}\ \vert\ i_1 < i_2 < \cdots < i_r\ \}$.
\item[]If $C\in{\mathcal S}_r$, $g_i$ restricted to $f^{-1}_i(C\ \cap\ \Bar{V}_i)$ 
is a homotopy equivalence.
\end{enumerate}

Notice that $Y_1$, $g_1$ and $f_1$ satisfy  1){--}5). 
(Let $Y_0=\emptyset$.)

If we can verify 1){--}5), we can construct $Y$ as the increasing union of $Y_i$ 
with identifications. 
$g$ and $f$ can be defined from the $g_i$ and $f_i$ respectively by 4).

Inductively, $Y$ has the proper homotopy type of a locally compact, 
finite dimensional complex, since it is covered by finite complexes,
$\Bar{C}_i\times S^k$, of bounded dimension in a locally finite fashion.
For a better proof, see \fullRef{P.2.4.1} below.

Now given $Y_{i-1}$, $f_{i-1}$, and $g_{i-1}$, we construct
$Y_i$, $f_i$ and $g_i$.

By Dold \cite{beight}, $\xi\vert_{\Bar{V}_i}$ can be gotten from
$\xi\vert_{\Bar{V}_{i-1}}$ and $\xi\vert_{\Bar{C}_i}$ as follows.
Over $\Bar{C_i\ \cap\ V_{i-1}\ }$, we have an equivalence 
$\varphi\colon (\xi\vert_{\Bar{V}_{i-1}})\vert_{\Bar{C_i\ \cap\ V_{i-1}\ }}
\ \to\ (\Bar{C_i\ \cap\ V_{i-1}\ })\times S^k$.
Let $H_1=\xi\vert_{\Bar{V}_i-1}$, $H_2=\Bar{C}_i\times S^k$, and let
$H_3=\bigl\{ (x,w)\ \vert\ x\in H_1\vert_{\Bar{C_i\cap V_{i-1}\ }},
\ w\in\bigl((\Bar{C_i\cap V_{i-1}\ })\bigr)^I,\ 
\pi(x)=\pi\bigl(w(t)\bigr)$ for all $t\in I,\ \varphi(x)=w(1)\ \bigr\}$.
Then $\xi\vert_{\Bar{V}_i} \cong H_1\cup H_3 \cup H_2$, where
$H_1\vert_{\Bar{C_i\cap V_{i-1}}}$ is embedded in $H_3$ via
$x\mapsto \bigl(x, $ constant path at $\varphi(x)\bigr)$.
The embedding of $H_2\vert_{\Bar{C_i\cap V_{i-1}}}$ is harder
to describe.
Let $\varphi^\prime$ be the inverse to $\varphi$.
The $\varphi\ \circ\ \varphi^\prime$ is fibre homotopic to the identity.
Let $\psi$ be a fibre homotopy between these two maps, with
$\psi(\quad,0)={\text{ id}}$.
Then $H_2\vert_{\Bar{C_i\cap V_{i-1}}}$ is embedded in $H_3$ via
$x\mapsto \bigl(\varphi^\prime(x), \psi(x,t)\bigr)$.

We must now define the $\rho$ in 3).
We are given
\[\begin{matrix}%
H_1\vert_{\Bar{C_i\cap V_{i-1}}}&\RA{\hskip 20pt \varphi\hskip20pt }&
\Bar{C_i\cap V_{i-1}}\times S^k\cr
\noalign{\vskip2pt}
\downlabeledarrow[\Big]{}{g_{i-1}}&&\downlabeledarrow[\Big]{}{{\text{ id}}}\cr
\noalign{\vskip6pt}
Y_{i-1}\big\vert_{f^{-1}_{i-1}(\Bar{C_i\cap V_{i-1}})}&
\dottedBar{\hskip50pt}\hskip-10pt\to&
\Bar{C_i\cap V_{i-1}}\times S^k\cr
\noalign{\vskip2pt}
\downlabeledarrow[\Big]{}{f_{i-1}}&&\downlabeledarrow[\Big]{}{{\text{ proj}}}\cr
\noalign{\vskip6pt}
\Bar{C_i\cap V_{i-1}}&\RA{\hskip 20pt{\text{ id}}\hskip20pt }
&\Bar{C_i\cap V_{i-1}}\cr
\end{matrix}\leqno{B)}\]
We would like to fill in the dotted arrow with $\rho$ 
so that the diagram actually commutes. 
To do this, we may have to alter $\varphi$ within its fibre homotopy
class, but this will not change our bundle.

Since $g_{i-1}$ is a homotopy equivalence, it has an inverse, $h$.
$h$ may be assumed to g=be a fibre map, so $h\ \circ\ g_{i-1}$ is
a fibre homotopy equivalence.
Let $G$ be its fibre homotopy inverse.
The $G\ \circ\ h\ \circ\ g_{i-1}$ is fibre homotopic to the identity.
$g_{i-1}\ \circ\ (G\ \circ\ h)$ is homotopic to the identity.

Set $\rho=({\text{ id}})\ \circ\ \varphi\ \circ\ (G\ \circ\ h)$.
Then $\rho$ is a fibre map so the bottom square commutes.
Set $\varphi_1=({\text{ id}})^{-1}\ \circ\ \rho\ \circ\ g_{i-1}$.
Then $\varphi_1$ is fibre homotopic to $\varphi$, and B) commutes
with $\varphi_1$ in place of $\varphi$.
$\rho$ is a homotopy equivalence, so 3) is satisfied.

From now on, we assume $\varphi$ chosen so that B) commutes with
the $\rho$ along the dotted arrow.
Set $Y_i= Y_{i-1}\ \cup_\rho\ \Bar{C}_i\times S^k$.
$f_i$ is defined by $f_i\vert_{Y_{i-1}}= f_{i-1}$ and
$f_i\vert_{\Bar{\scriptstyle C}_i\times S^k}={\text{ proj}}$.
B) insures that this is well{--}defined on the intersection.

$g_i$ is unfortunately harder to define.
$\xi\vert_{\Bar{V}_i}\cong H_1\ \cup\ H_2\ \cup\ H_3$, so let
$\alpha\colon \xi\vert_{\Bar{V}_i}\ \to\  H_1\ \cup\ H_2\ \cup\ H_3$
be an equivalence.
$\alpha$ may be chosen to be the identity on $\xi\vert_{U_{i-1}}$.
We define a map $h\colon H_1\ \cup\ H_2\ \cup\ H_3\ \to\ Y_i$
as follows.
$g\vert_{H_1}=g_{i-1}$.
To define $g_i$ on the other two pieces, look at $\psi$, the fibre
homotopy between $\varphi\ \circ\ \varphi^\prime$ and ${\text{ id}}$.
This can be extended to a fibre map of 
$(\Bar{C}_i\times S^k)\times I\  \to\ (\Bar{C}_i\times S^k)$ since
${\text{ id}}\colon (\Bar{C_i\cap V_{i-1}})\times S^k\ \to\ 
(\Bar{C_i\cap V_{i-1}})\times S^k$ can clearly be extended.

Now define $(g\vert_{H_3})(x,w)=g_{i-1}(x)$.
Note our two definitions agree on $H_1\cap H_3$.
We could have defined $(g\vert_{H_3})(x,w)=w(1)$ equally well.
We define $(g\vert_{H_2})(x)=F(x)$.
If $x\in H_2\cap H_3$, then $(g\vert_{H_3})(x)=(g\vert_{H_3})\bigl(
\varphi^\prime(x),\psi(x,t)\bigr)=\psi(x,1)=
\varphi\ \circ\ \varphi^\prime(x)$.
$(g\circ_{H_2})(x)=F(x)=\varphi\ \circ\ \varphi^\prime(x)$ by
the definition of $F$.
Hence $g$ is well{--}defined and we set $g_i= g\ \circ\ \alpha$.

Now 4) clearly holds since $\alpha\vert_{f^{-1}_{i-1}(U_{i-1})}$ 
is the identity.
1) holds as $g\colon H_1\cup H_2\ \cup H_3\ \to\ Y_i$
preserves fibres by construction.
Hence we are left with showing 2) and 5).

For $r$ sufficiently large, $C\in{\mathcal S}_r$ implies 
$C\cap\Bar{V}_i=\emptyset$, since the collection $\{C_i\}$is locally finite
\setcounter{footnote}{0}\footnote{${\mathcal S}_r$ is defined in 5), page \pageref{ItemFive}}.
We show 5) by downward induction on $r$, since if 
$C\cap\Bar{V}_i=\emptyset$, 5) is obvious.
Assume we have established the result for $r=k+1$.
Let $C\in{\mathcal S}_k$. 
If $C\cap\bar{C_i}=\emptyset$, then $C\cap\Bar{V}_{i-1}=C\cap\Bar{V}_i$ 
and we are done since 5) holds for $g_{i-1}$ and $\alpha$ is a fibre
homotopy equivalence.
If $C\cap\Bar{V}_i=\Bar{C}_i\cap\Bar{V}_i$ we are done since $F$
is a fibre map.
So let $L=C\cap\Bar{V}_{i-1}$, and let $K=C\cap\Bar{C}_i$ with both $K$
and $L$ non{--}empty.
$g_i\vert_{f^{-1}_i(L)}$ is a homotopy equivalence, again since $\alpha$
is a fibre homotopy equivalence and $g_{i-1}\vert_{f^{-1}_i(K)}$ is also
a homotopy equivalence.
$K\cap L\subseteq \Bar{V}_{i-1}$, and $K\cap L\in{\mathcal S}_{k+1}$.
Hence $g_i\vert_{f^{-1}_i(K\cap L)}$ is a homotopy equivalence.
Therefore $g_i\vert_{f^{-1}_i(C)}$ is a homotopy equivalence and we
are done with 5).

Therefore we have a space $Y$ and maps
\topD{12}{$\begin{matrix}%
S(\xi)&\RA{\ g\ }&Y\cr
&\hbox to 0pt{\hss$\hbox to 0pt{\hss$\scriptstyle{\pi}$\hss}
\searrow\hskip10pt\swarrow\hbox to 0pt{$\scriptstyle{f}$\hss}$\hss}\cr
&X\cr
\end{matrix}\ .$}{4}
We claim $g$ is a homotopy equivalence.
Since by Milnor \cite{btwentytwo}, $S(\xi)$ has the homotopy type of 
a CW complex, this is equivalent to showing $g$ induces isomorphisms
in homotopy.
But $\displaystyle\pi_k(g)=\mathop{\lim}_{\Atop{\longrightarrow}{i}}
\pi_k(g_i)$, and since $\pi_k(g_i)=0$, $\pi_k(g)=0$.

Let $h\colon Y\ \to\ S(\xi)$ be a homotopy inverse for $g$.
By an easy argument like the one after diagram B), we may assume $h$
preserves fibres and the $h\ \circ\ g$ is fibre homotopic to the identity.
Notice that by construction $f^{-1}(x)$ is homeomorphic to a sphere
of dimension $\dim\xi$.
$\pi^{-1}(x)$ has the homotopy type of such a sphere.
Since $h\ \circ\ g$ is fibre homotopic to the identity,
$g_x\colon \pi^{-1}(x)\ \to\ f^{-1}(x)$ has a left inverse.
As both spaces are spheres of dimension $2$ or more, $g_x$ is a
homotopy equivalence.

Now in the terminology of Bredon \cite{btwo}, $f$ is $\psi${--}closed, and
$f^{-1}(x)$ is $\psi${--}taut, where $\psi$ is the family of compact supports.
(Note $Y$ is locally compact, so $\psi$ is paracompactifying, and then
apply (d) on page 52 to show $f^{-1}(x)$ is $\psi${--}taut.
$f$ is $\psi${--}closed easily from the definition, which is on page 53,
since $X$ is Hausdorff.)
Hence we have a Leray spectral sequence for the map 
$f\colon X\ \to\ Y$.
We have the Serre spectral sequence for $\pi\colon S(\xi)\ \to\ X$,
and $g$ induces a map between these two.
$g$ induces an isomorphism on the $E_2$ terms since it is a homotopy
equivalence on each fibre.
Hence $g\colon H^\ast_\cmpsup(Y)\ \to\ 
H^\ast_{\varphi}\bigl(S(\xi)\bigr)$ is an isomorphism, where $\varphi$ is
the set of supports whose image in $X$ is compact.

As $\dim\xi\geq2$, $\pi^\ast\colon H^\ast_\cmpsup(X)\ \to\ 
H^\ast_\varphi\bigl(S(\xi)\bigr)$ is an isomorphism for $\ast<2$.
Hence $f^\ast\colon H^\ast_\cmpsup(X)\ \to\ H^\ast_\cmpsup(Y)$ is 
an isomorphism for $\ast<2$, so
$f^\ast\colon H^0_{{\text{ end}}}(X)\ \to\ H^0_{{\text{ end}}}(Y)$ is an
isomorphism, so $f$ is a proper $0${--}equivalence.

We claim $f$ is a proper $1${--}equivalence.
To see this, note $f\vert_C$ is a $1${--}equivalence for $C\in{\mathcal C}_r$
all $r\geq1$.
Now an easy van{--}Kampen induction shows $f$ is a $1${--}equivalence
when restricted to any union of $\Bar{C}_i$'s.
Hence $f$ is a proper $1${--}equivalence.

Thus $g_\#\colon \Delta(S(\xi)\Colon \pi_1)\ \to\ \Delta(Y\Colon \pi_1)$ is an
isomorphism as both groups are isomorphic, via $\pi_\#$ and $f_\#$,
to $\Delta(X\Colon \pi_1)$.

Now we still have maps 

\[\begin{matrix}%
Y- f^{-1}(K_i)&{\scriptscriptstyle
\Atop{\Limitsarrow{10}{10}{\scriptscriptstyle h}{}{\rightarrowfill}}
{\Limitsarrow{10}{10}{}{\scriptscriptstyle g}{\leftarrowfill}}
}&
S(\xi\vert_{X-K_i})\cr
&\hbox to0pt{\hss$%
\hbox to 0pt{\hss$\scriptstyle{f}$}\searrow\hskip40pt\swarrow\hbox to 0pt{$\scriptstyle{\pi}$\hss}
$\hss}\cr
&X-K_i\cr
\end{matrix}\]
\newpage
\noindent
where $\displaystyle K_i=X-\ \mathop{\cup}_{j\geq i} \Bar{C}_i$.
$g$ restricted to each fibre is still a homotopy equivalence with
inverse induced from $h$.
For any cover $\coverFC{}$, of $X-K_i$, we get
\[\begin{matrix}%
\coverFC{Y- f^{-1}(K_i)\ }&{\scriptscriptstyle
\Atop{\Limitsarrow{10}{10}{\coverFC{\scriptscriptstyle \ h\ }}{}{\rightarrowfill}}
{\Limitsarrow{10}{10}{}{\coverFC{\ \scriptscriptstyle g\ }}{\leftarrowfill}}
}&
\coverFC{S(\xi\vert_{X-K_i})\ }\cr
&\hbox to0pt{\hss$%
\hbox to 0pt{\hss$\scriptstyle{\coverFC{\ f\ }}$}
\searrow\hskip40pt\swarrow\hbox to 0pt{$\scriptstyle{\coverFC{ \pi\ \ }}$\hss}
$\hss}\cr
&\coverFC{X-K_i\ }\cr
\end{matrix}\]
where the covers on the top row are induced covers from $\coverFC{}$ over
$X-K_i$.
$\coverFC{S(\xi\vert_{X-K_i})}$ is the same as
$S(\coverFC{\ \xi\ }\vert_{\coverFC{X-K_i}})$, 
the spherical fibration induced from $\xi\vert_{X-K_i}$ 
over $\coverFC{X-K_i}$.
$\coverFC{\ g\ }$ likewise induces a homotopy equivalence of fibres,
so as before we get
\[h^\ast\colon H^\ast_\varphi
\Bigl(\coverFC{S\bigl(\xi\vert_{X-K_i}\bigr)\ },
\coverFC{S\bigl(\xi\vert_{\partial(X-K_i)}\bigr)\ }\Bigr)\ \to\ 
H^\ast_\cmpsup\Bigl(\coverFC{Y-f^{-1}(K_i)\ },
\partial\bigl(\coverFC{Y-f^{-1}(K_i)\ }\bigr)\Bigr)
\]
is an isomorphism.
A word about the existence of these covers is in order.
Since $X-K_i$ is a CW complex, its cover exists.
The cover for $S(\xi\vert_{X-K_i})$ then also clearly exists.
We claim $Y-f^{-1}(K_i)$ is semi{--}locally $1${--}connected, from
which it follows that its cover also exists.
To see our claim, observe $f\colon Y-f^{-1}(K_i)\ \to\ X-K_i$ is
a $1${--}equivalence.
Given any point $y\in Y-f^{-1}(K_i)$, let $N\subseteq X-K_i$ be a
neighborhood of $f(y)$ such that $\pi_1(N)\ \to\ \pi_1(X-K_i)$
is the zero map.
Since $X-K_i$ is semi{--}locally $1${--}connected, such an $N$ exists.
Now $f^{-1}(N)$ is a neighborhood for $y$, and
$\pi_1\bigl(f^{-1}(N)\bigr)\ \to\ \pi_1\bigl(Y-f^{-1}(K_i)\bigr)$
is also zero.
Hence $Y-f^{-1}(K_i)$ is semi{--}locally $1${--}connected.
A similar argument shows $Y-f^{-1}(K_i)$ is locally path connected.

Therefore, $h^\ast\colon\Delta(S(\xi)\Colon \coverFC{})\ \to\ 
\Delta(Y\Colon \coverFC{})$ is an isomorphism for any covering functor
induced from one over $X$.
Since $f$ is a proper $1${--}equivalence, if we take a universal
covering functor for $X$, we get one for $Y$.
(The actual covering functor on $Y$ is the following. Any $A\in{\mathcal C}(Y)$
is contained in a unique minimal $f^{-1}(X-K_i)$ so let the cover over
$A$ be induced from the cover over this space.)

$g^\ast\colon \Delta(Y\Colon \coverFA{})\ \to\ \Delta(S(\xi)\Colon \coverFC{})$
is defined where $\coverFC{}$ is the covering functor induced by $g$ from 
$\coverFA{}$ over $Y$.
$g^\ast\ \circ\ h^\ast=(h\ \circ\ g)^\ast\colon 
\Delta^\ast(S(\xi)\Colon \coverFA{})\ \to\ \Delta^\ast(S(\xi)\Colon \coverFC{})$
is an isomorphism as $\coverFA{}$ and $\coverFC{}$ are equivalent
covering functors.
Hence $h^\ast\ \circ\ g^\ast=(g\ \circ\ h)^\ast\colon
\Delta^\ast(Y\Colon \coverFA)\ \to\ \Delta^\ast(Y\Colon \coverFA)$
is an isomorphism, so $h\ \circ\ g$ is a proper homotopy equivalence.
$g\ \circ\ h$ is already a fibre homotopy equivalence, and it
is not hard to change $h$ until $h\ \circ\ g$ is properly homotopic
to the identity and $g\ \circ\ h$ is fibre homotopic to the identity.

To finish, we need only show \fullRef{P.2.4.1} below.
We first need

\bigskip
\BEGIN{T.2.4.1}
Let $Y$ be a locally compact, separable ANR.
Then $Y$ is properly dominated by a locally{--}finite simplicial complex.
\end{Theorem}

\medskip\begin{proof}
Let $\alpha$ be an open covering of $Y$ 
by sets whose closure is compact.
Since $Y$ is metrizable , $Y$ is paracompact, so we can assume $\alpha$
is locally finite.

We now apply Hu \cite{bfifteen}, Theorem 6.1, page 138, to get a
locally finite simplicial complex $X$ and maps $\varphi\colon X\ \to\ Y$
and $\psi\colon Y\ \to\ X$ with $\varphi\ \circ\ \psi$
$\alpha${--}homotopic to the identity, i.e. if $H$ is the homotopy,
for each $y\in Y$, there exists $U\in\alpha$ such that $H(y,t)\in U$
for all $t\in[0,1]$.
By our choice of $\alpha$, $\varphi\ \circ\ \psi$ is properly
homotopic to the identity.

Now $X$ is actually the nerve of some cover $\delta$ in the proof of
Hu, Theorem 6.1.
In the proof, we may take $\delta$ to be star{--}finite and locally finite.
Then the nerve $X$ is a locally finite simplicial complex, and the map
$\varphi\colon X\ \to\ Y$ is proper.
To see this last statement, it is enough to show $\varphi^{-1}(U)$ is
contained in a compact subset of $X$ for any $U\in\alpha$.
Recall $\varphi$ is defined by picking a point in each $V\in\delta$ 
and sending vertex of the nerve which corresponds to $V$ to our
chosen point and then extending.
Our extension satisfies the property that any simplex lies entirely in
some element of $\alpha$.
So let $U_1$ be the union of all elements of $\alpha$ intersecting $U$.
$\Bar{U}_1$ is compact as $\alpha$ is locally finite, so let $U_2$
be the union of all elements of $\alpha$ intersecting $\Bar{U}_1$.
$\Bar{U_2}$ is again compact, so there are only finitely many elements
of $\delta$ which intersect $U_2$.
Let $K\subseteq X$ be the subcomplex generated by these elements
of $\delta$. 
$K$ is finite, hence compact, and $\varphi^{-1}(U)\subseteq K$.
\end{proof}

\medskip
\BEGIN{C.2.4.1.1}
Let $Y$ be a locally compact, separable ANR, and suppose the covering
dimension of $Y$, $\dim Y$, is finite 
(see Hurewicz and Wallman \cite{bsixteen} for a definition).
Then $Y$ is properly dominated by a locally finite simplicial complex of
dimension $\dim Y$.
\end{Corollary}

\medskip\begin{proof}
Make the same changes in Hu \cite{bfifteen} Theorem 6.1, page 164 
that we made to the proof of Theorem 6.1, page 138.
We get a simplicial complex $P$ and a proper map $\varphi\colon P\ \to\ Y$
such that for any map $f\colon X\ \to\ Y$ with $X$ a metric space 
of dimension $\leq\dim Y$, there exists a map $\psi\colon X\ \to\ P$
with $\varphi\ \circ\ \psi$ $\alpha${--}homotopic to $f$.
Moreover, $P$ has no simplices of dimension $>\dim Y$.
Apply this for $X=Y$, $f={\text{ id}}$.
\end{proof}

\medskip
\BEGIN{C.2.4.1.2}
A locally compact, separable ANR of dimension $\leq n$ satisfies $D n$.
\end{Corollary}

\medskip\begin{proof}
By \fullRef{C.2.4.1.1} and nonsense, it remains to show $Y$
is homogamous.
But an ANR is locally contractible (Hu \cite{bfifteen}, Theorem 7.1,
page 96), and any metric space is paracompact so 
\fullRef{C.1.1.2.1} applies.
\end{proof}

\bigskip
\BEGIN{P.2.4.1}
The space $Y$ which we constructed has the proper homotopy type of
a locally compact, finite dimensional CW complex.
\end{Proposition}

\medskip\begin{proof}
We first show $Y$ is a finite dimensional, locally compact, separable  ANR.
We then find a finite dimensional simplicial complex $Z$ and a proper
map $h\colon Z\ \to\ Y$ which is properly $n${--}connected for
any finite $n$.
Since both $Y$ and $Z$ satisfy $D n$ for some finite $n$, $h$ is a proper
homotopy equivalence.

\medskip
\insetitem{Step 1}
$Y$ is a finite dimensional, locally compact, separable ANR.
\smallskip
By Hu \cite{bfifteen} Lemma 1.1, page 177, Theorem 1.2, page 178, and
induction, each $Y_i$ is an ANR.
The induction is complicated by the necessity of showing 
$Y_{i-1}\ \cap\ f^{-1}_{i-1}(\Bar{V}_{i-1}\cap\Bar{C}_i)$ is an ANR.
Hence our induction hypothesis must be
\begin{enumerate}
\item[a$)_{k\hphantom{,r}}$] $Y_k$ is an ANR
\item[b$)_{k, r}$] $Y_k \cap f^{-1}_k(\Bar{V}_k\cap C)$ is an ANR
for all $C\in{\mathcal S}_r$.
\end{enumerate}

One then shows that for some finite $r$, b$)_{k, r}$ holds vacuously.
b$)_{k, s}$, $s>r$, and b$)_{k-1 , r}$ imply b$)_{k, r}$, so we get
b$)_{k, r}$ for all $r$.
b$)_{k, 1}$ and a$)_{k-1}$ imply a$)_k$, so we are done.

Since each $Y_i$ is an ANR, each $Y_i$ is a local ANR (Hu, Proposition 7.9,
page 97).
If $Y$ is metrizable, $Y$ is an ANR by Hu, Theorem 8.1, page 98.
Now $Y$ is $T_1$ and regular. 
To see this observe each $Y_i$ is $T_1$ and regular since it is metrizable.
Now if $U\subseteq Y$ is any compact set, there is a $Y_i$ with
$V\subseteq Y_i$ and $V$ homeomorphic to $U$.
With this result and the observation that $Y$ is locally compact, it is easy to
show $Y$ is $T_1$ and regular. $Y$ is locally compact because 
it has a proper map to a locally compact space $X$.
$Y$ is $\sigma${--}compact since $X$ is, so $Y$ is second countable.
Hence $Y$ is metrizable (see Kelly \cite{bseventeen} page 125)
and separable.

We are left with showing $Y$ has finite covering dimension.
By Nagami \cite{btwentyseven} (36-15 Corollary, page 206),
we need only show the small cohomological dimension with respect to
the integers (Nagami, page 199) is finite ($Y$ is paracompact since
it is $\sigma${--}compact and regular (see Kelly \cite{bseventeen}, 
page 172, exercise Y a) and b))).

To compute $d(Y\Colon \Z)$, look at the map $f\colon Y\ \to\ X$.
$f$ is a closed, onto map. $f$ is onto by construction, and $f$ is
closed since $Y$ is the increasing union of compact sets $\{D_i\}$,
so $F\subseteq Y$ is closed \iff\ $F\cap D_i$ is closed for all $i$, and
$f(F\cap D_i)$ is closed since $F\cap D_i$ is compact and $X$ is Hausdorff. 
We can find an increasing sequence of compact sets $\Bar{V}_i$ such
that $E\subseteq X$ is closed \iff\ $E \cap \Bar{V}_i$ is closed.
Since $f$ is proper, $D_i=f^{-1}(\Bar{V}_i$ has the expected properties.
But $f(F\cap D_i)=f(F)\cap\Bar{V}_i$ if $D_i=f^{-1}(\Bar{V}_i)$, so
$f$ is closed.
Hence by Nagami \cite{btwentyseven} (38-4 Theorem, page 216),
$d(Y\Colon \Z)\leq \dim X+k$ where $k$ is the dimension of the fibration $\xi$.
To see this, note $f^{-1}(x)$ is homeomorphic to $S^k$ for all$x\in X$,
so $d(F^{-1}(x)\Colon \Z)=k$.
Since $X$ is paracompact and metrizable, ${\text{ Ind}} X=\dim X = d(X\Colon \Z)=$
dimension of $X$ as a CW complex (see Nagami 8-2 Theorem for the
first equality; Nagami 36-15 Corollary shows the second; Nagami 37-12
Theorem and subdivision show the third
[this uses the fact that $X$ is a regular complex]).

\medskip\insetitem{Step 2}
\LRTpageLabel{ItemFive}
There is a locally compact, finite dimensional CW complex $Z$ and a proper
map $h\colon Z\  \to\ Y$ which is properly $n${--}connected for all $n$.
\smallskip
We define $Z$ and $h$ by induction; i.e. we have
\begin{enumerate}
\item[1)] a finite CW complex $Z_i$ and a map $h_i\colon Z_i\ \to\ Y_i$
\item[2)] $h_i$ is a homotopy equivalence
\item[3)] $h_i$ restricted to $(f_i\ \circ\ h_i)^{-1}(C\cap\Bar{V}_i)$
is a homotopy equivalence for all $C\in {\mathcal S}_r$, $r\geq 1$
\item[4)] $h_i\vert_{(h_{i-1}\circ f_{i-1})^{-1}(U_{i-1})}=
h_{i-1}\vert_{(h_{i-1}\circ f_{i-1})^{-1}(U_{i-1})}$
\item[5)] $Z_i=Z_{i-1}\ \cup_\lambda \Bar{C}_i\times S^k$ where
\item[] $\lambda\colon Z_{i-1}\ \cap\ (h_{i-1}\circ f_{i-1})^{-1}(
\Bar{V}_{i-1}\cap \Bar{C}_i)\ \to\ 
(\Bar{V}_{i-1}\cap\Bar{C}_i)\times S^k$
\item[] is a cellular homotopy equivalence.
\end{enumerate}

\medskip
If we can find such $Z_i$ and $h_i$, we can find $Z$ and $h\colon Z\ \to\ Y$. 
$h$ is clearly proper.
$h\vert_{(f\circ h)^{-1}(C)}\colon (f\circ h)^{-1}(C)\ \to \ f^{-1}(C)$
is a homotopy equivalence by 3) for all $C\in{\mathcal S}_r$, $r\geq 1$,
so $h\vert_{(f\circ h)^{-1}(D_i)}$ is a homotopy equivalence where
$\displaystyle D_i=\mathop{\cup}_{j\geq i} \Bar{C}_i$.
Thus $h$ induces isomorphisms on $H^0$ and $H^0_{{\text{ end}}}$,
and $\Delta(h\Colon \pi_s)=0$ for $s\geq 1$.
Hence we are done if we can produce $Z_i$ and $h_i$.

We proceed by induction on $i$. $Z_1=\Bar{V}_1\times S^k$ and
$h_1={\text{ id}}$.
1){--}5) are trivial, so suppose we have $Z_{i-1}$ and $h_{i-1}$.
We have
\[
Z_{i-1}\cap (f_{i-1}\circ h_{i-1})^{-1}(\Bar{V}_{i-1}\cap\Bar{C}_i)
\ \RA{\hskip10pt}\ 
Y_{i-1}\cap f^{-1}_{i-1}(\Bar{V}_{i-1}\cap\Bar{C}_i)
\ \RA{\ \rho\ }\ 
(\Bar{V}_{i-1}\cap\Bar{C}_i)\times S^k
\]
Let $\rho^\prime$ be this composition.
Deform $\rho^\prime$ to a cellular map as follows.
For some $r\geq1$, $C\in{\mathcal S}_r$ implies $C\cap\Bar{C}_i=\emptyset$.
Now deform $\rho^\prime$ to a cellular map over each 
$C \cap \Bar{C}_i \cap \Bar{V}_{i-1}$ for $C\in{\mathcal S}_r$, all $r\geq1$
and finally to a cellular map over $\Bar{C}_i\cap \Bar{V}_{i-1}$.
Denote this map by $\lambda$.

Let $Z_i=Z_{i-1}\ \cup_{\lambda}\ (\Bar{C}_i\times S^k)$.
We extend $h_{i-1}$ to a homotopy equivalence $h_i\colon Z_i\ \to\ Y_i$
which leaves $h_{i-1}$ fixed on 
$(f_{i-1}\circ h_{i-1})^{-1}(U_{i-1})$.
$h_i$ in fact can be chosen to be a homotopy equivalence on each
$(f_i\circ h_i)^{-1}(C\cap\Bar{V}_i)$ be extending inductively over the
various $C\in{\mathcal S}_r$.
1){--}5) hold and we are done.
\end{proof}

\chapter{The Geometric Surgery Groups}
\section{The fundamental theorems of surgery}
\newHead{III.1}

In this section we will prove three results which may be called the
fundamental theorems of surgery.
They constitute all the geometry needed to define surgery groups 
and to prove these groups depend only on the proper $1${--}type
of the spaces in question.
These results together with the $s${--}cobordism theorem constitute 
the geometry necessary to give a classification of paracompact 
manifolds in a given proper homotopy class \`a la Wall \cite{bfortyone},
Chapter 10.

Let \CAT\ denote either TOP, PL or DIFF.
If $X$ is a locally finite, finite dimensional CW $n${--}ad, and if $\nu$ is a
\CAT\ bundle over $X$, then
$\bor_M(X,\nu)$ is the space of cobordism classes of the following triples:
a \CAT\ manifold $n${--}ad $M$, $\dim M=m$; a proper map of $n${--}ads
$f\colon M\ \to\ X$; a stable bundle map $F\colon \nu_M\ \to\ \nu$,
where $\nu_M$ is the normal bundle of $M$ and $F$ covers $f$.
Such a triple is called a {\sl normal map}, and the cobordisms 
are called {\sl normal cobordisms}.

\bigskip
\BEGIN{T.3.1.1}
Given $\alpha\in\bor_m(X,\nu)$, there is a representative $(M,f,F)$ of 
$\alpha$ with $f$ properly $\bigl[\frac{m}{2}\bigr]${--}connected
if $X$ is a space. (\, $[\quad ]=$ greatest integer.)

For a pair $(X,\partial X)$, we have a representative 
$\bigl((M,\partial M), f,F)$ with
$f\colon M\ \to\ X$ properly $\bigl[\frac m2\bigr]${--}connected;
$f\colon \partial M\ \to\ \partial X$ properly 
$\bigl[\frac {m-1}2\bigr]${--}connected; and the pair map
$f\colon (M,\partial M)\ \to\ (X,\partial X)$ properly 
$\bigl[\frac m2\bigr]${--}connected.
If $\partial X\subseteq X$ is properly $0${--}connected, then the map
of pairs may be made properly homologically 
$\bigl[\frac{m+1}{2}\bigr]${--}connected provided $m\geq 3$.
\end{Theorem}

\medskip\begin{proof}
The proof follows Wall \cite{bforty}, Theorem 1.4. 
(See the remark following his proof.) 
We first remark that his Lemma 1.1 is equally valid in our case.

\medskip
\BEGIN{L.3.1.1}
Suppose $M$ and $X$ locally compact, finite dimensional CW complexes, 
$\psi\colon M\ \to\ X$ a map.
Then we can attach cells of dimension $\leq k$ to $M$ so that the
resulting complex is locally finite and so that the map is
properly $k${--}connected.
\end{Lemma}

\smallskip\begin{proof}
We may assume $\psi$ cellular by the cellular approximation theorem.
Then the mapping cylinder of $\psi$ is a locally compact, finite dimensional
complex, and $(M_\psi,M)$ is a CW pair.
Set $M^\prime =M^k_\psi\ \cup M$.
Note then that $M^\prime$is obtained from $M$ by adding cells of dimension
$\leq k$ and that $M^\prime\ \to\ M_\psi$ 
is properly $k${--}connected.
\end{proof}

Now given a representative $(N,g,G)$ for $\alpha$, attach handles of 
dimension $\leq k$ to $N$ to get $\psi\colon W\ \to\ X$ with 
$\partial W=N\ \cup\ M$, $\psi\vert_N=g$, and with $\psi$ covered
by a bundle map which is $G$ over $N$, and $\psi$ 
is properly $k${--}connected. 
The argument that we can do this is the same as for the compact case.
Wall \cite{bfortyone} Theorem 1.1 generalizes immediately to

\medskip
\BEGIN{L.3.1.2}
Given $\alpha\in\bor_m(X,\nu)$ with any representative $(M,f,F)$,
any element of $\Delta(f:\pi_k)$ determines proper regular homotopy
class of immersions of a disjoint collection of $S^k\times D^{m-k}$'s into
$M$ for $k\leq m-2=\dim M-2$.
\end{Lemma}

\smallskip\begin{proof}
Precisely as in Wall, \cite{bfortyone} Theorem 1.1, we get a stable
trivialization of the tangent bundle of $M$ over each sphere $S^k$
in our collection.
Given any sphere $S^k$, we see in fact that there is an open submanifold
$U\subseteq M$ such that we get a trivialization of the tangent bundle
of $U$ restricted to $S^k$ which agrees with the one for $\tau_M$.
In fact $U=f^{-1}\bigr($ a small neighborhood of the disc 
bounding $f(S^k)\bigl)$ will do (we have momentarily confused $S^k$
with its image in $M$).
Notice that we can pick such a collection of $U$'s to  be locally finite.
Now apply Hirsch \cite{bfourteen}, Haefliger \cite{btwelve}, or Lees 
\cite{bnineteen} to immerse each $S^k$ in its $U$ 
with trivial normal bundle.
This is a proper homotopy, so each $\alpha$ determines a proper map
which immerses each sphere.

It is not hard to show any two such immersions which are 
properly homotopic  are regularly properly homotopic.
\end{proof}

\medskip
If there is an embedding in the proper regular homotopy class of
$\alpha$, we can attach a collection of handles by $\alpha$
and extend our map and bundle map over the resulting trace of
the surgeries.
Notice that in an embedding, all the spheres have disjoint images, so
we can certainly do the surgery.
The map can be extended properly by construction, and one shows the
bundle map extends precisely as in the compact case 
(Wall \cite{bfortyone} Theorem 1.1).

\medskip
\BEGIN{L.3.1.3}
Given $\alpha\in\bor_m(X,\nu)$ with any representative $(M,f,F)$,
we can do surgery on any element $\beta\in\Delta(f:\pi_k)$
for $m>2k$.
\end{Lemma}

\smallskip\begin{proof} 
General position supplies us with an embedding.
\end{proof}

\medskip
We now return to the proof of \fullRef{T.3.1.1}.
By our lemmas, we see that if $m>2k$, we can get $W$ as advertised.
Now $W$ is obtained from $M$ by adding handles of dimension
$\geq (m+1)-k>k+1$, so $M\subseteq W$ is properly $k${--}connected.

In the pairs case, given a representative, we first fix up the boundary 
as above.
Then we can attach handles away from the boundary to get the
absolute map fixed up.
The long exact homotopy sequence shows that the pair map is properly
$\bigl[\frac m2\bigr]${--}connected.
If $m$ is even, we are done.
The case for $m=2k+1$ follows Wall \cite{bfortyone} Theorem 1.4.

We may assume that we have $f\colon (M,\partial M)\ \to\ (X,\partial X)$
properly connected up to the middle dimension on each piece.
Let $E$ be the disjoint union of the $(k+1)${--}cells of $M_f - M$.
Then we have a proper map $\partial E\ \to\ M_f$.
Since $\partial E$ is $k${--}dimensional, and since $(M_f,M)$ is properly
$k${--}connected, there is a proper homotopy of the attaching maps into
$M$.
$\partial E=\displaystyle\mathop{\disjointunion}_p S^k_p$, so embed these
spheres in $M$ with trivial normal bundle by Lemmas 
\shortFullRef{L.3.1.1} and \shortFullRef{L.3.1.2}.
Join each sphere to $\partial M$ by a locally finite collection of tubes,
one for each sphere.
(Since $H^0_{{\text{ end}}}(X,\partial X)=0$ by hypothesis, and since 
$M\ \to\ X$ is properly $1${--}connected (at least), and since
$\partial M\ \to\ \partial X$ is properly $0${--}connected, 
$H^0_{{\text{ end}}}(M,\partial M)=0$ so we can do this.)
Note in fact that we need only disturb $\partial M$ in a (\prex{--}assigned)
neighborhood of a set of base points.) 
By general position  we may assume all these tubes disjoint ($m\geq3$).
Hence we get framed embeddings of a collection of disjoint $D^k$'s.
We may assume (by adding trivial discs if necessary) that the centers
of our discs form a set of base points for $M$.

We claim that if we do these relative surgeries we will have killed
$K_k(M,\partial M)$ without affecting our other conditions.
Our proof of this claim is the same as Wall's.
Let $H$ denote the union of the handles, $N_0$ the constructed manifold,
$f_0\colon (N_0,\partial N_0)\ \to\ (X,\partial X)$ the resulting map.
Note that $(N_0,\partial N_0)\ \to\ (M, H\ \cup\ \partial N_0)$
is a proper excision map.
We can pick a set of base points for $\partial M$ away from 
$\partial M\ \cap\ H$.
As usual we can pick them so that they are a set of base points for
$f\colon \partial M\ \to\ \partial X$.
They are then also seen to be a set of base points for
$M$, $N_0$, $\partial N_0$, and $H\  \cup \partial N_0$.
With these base points and the above excision map we get an exact
sequence\hskip 10pt
$\Delta_k(H\cup\partial N_0,\partial M:M:\coverFA{})\ \to 
\Delta_{k+1}(f:\coverFA{})\ \to\ \Delta_{k+1}(f_0:\coverFA{})\ \to\ 
\Delta_{k-1}(H\cup\partial N_0,\partial M:M:\coverFA{})$.

Clearly the lower relative proper homotopy groups of $f_0$ vanish.
Notice $(H, H\cap\partial M)\ \to\ (H\cup\partial N_0,\partial M)$
is also a proper excision map.
Since $(H, H\cap\partial M)$ is a collection of copies of
$(D^k\times D^{k+1}, S^{k-1}\times D^{k+1})$,
$\Delta_\ast(H, H\cap\partial M)$ is $0$ except in dimension $k$.
If we pick base points in $H$, 
$\Delta_\ast(H \cup\partial N_0,\partial M:\coverFB{})=0$ 
also except in dimension $k$ ($\coverFB{}$ is any covering functor).
Hence $\Delta_{k-1}(H\cup\partial N_0,\partial M:M:\coverFA{})=0$.

Let $g\colon M\ \to\ X$ denote $f$ on $M$ to distinguish it from
$f$ on $(M,\partial M)$.
The collection of elements above generates $\Delta(g:\pi_{k+1})$.
Clearly the composite map \qquad
$\Delta(g:\pi_{k+1}) \ \to\ 
\Delta_{k+1}(g:\coverFA{})\ \to\ 
\Delta_{k+1}(f:\coverFA{})\ \to\ 
\Delta_{k+1}(f_0:\coverFA{})$ is the zero map.
But by Hurewicz, the first map is an isomorphism, and the second map is
onto since $\partial M\ \to\ \partial X$ is properly $k${--}connected.
Hence $\Delta_{k+1}(f:\coverFA{})\ \to\ 
\Delta_{k+1}(f_0:\coverFA{})$ is the zero map.

Now the last two paragraphs and our exact sequence show
$\Delta_{k+1}(f_0:\coverFA{})=0$ as claimed.
\end{proof}

\bigskip\begin{xRemarks}
Note throughout the proof that should 
$\partial X=\partial_1X\ \cup\ \partial_2X$ and 
$\partial M=\partial_1M\ \cup\ \partial_2M$, and if
$\partial_2M\ \to\ \partial_2X$ is already properly $r${--}connected,
then we need attach no cells of dimension less than $r$ to $\partial_2M$
in our construction
(provided $H^0_{{\text{ end}}}(X,\partial_1X)=0$, otherwise to get this
last part of the result we must attach some $k${--}cells in $\partial_2M$).
In particular, if $\partial_2M\ \to\ \partial_2X$ is a proper homotopy
equivalence, we can do our construction away from $\partial_2M$
(except possibly for the last step).
\end{xRemarks}

\bigskip
\BEGIN{T.3.1.2}
Let $f\colon (M,\partial M)\ \to\ (X,\partial X)$ be a degree one
normal map; i.e. a bundle over $X$ and 
a bundle map over $f$ are understood.
Let $(X,\partial X)$ be a Poincar\'e duality pair of formal 
dimension at least $6$.
Suppose $\partial X\subseteq X$ is a proper $1${--}equivalence.
Then $f$ is normally cobordant to 
$g\colon (N,\partial N)\ \to\ (X,\partial X)$ with $g$ a proper homotopy
equivalence of pairs.
The torsion of $g$ may have any \prex{--}assigned value.
The torsions of $g\colon \partial N\ \to\ \partial X$ and of $g$ as
a map of pairs is then determined. 
\end{Theorem}

\medskip\begin{proof}
The proof of the theorem divides into two cases.

\medskip\insetitem{\hskip\parindent Case 1} $\dim(X)=2k$.
\smallskip
By \fullRef{T.3.1.1}, we can do surgery on $f$ to make 
the map $f\colon M\ \to\ X$ $k${--}connected, and to make the
map $\partial f\colon \partial M\ \to\ \partial X$
$(k-1)${--}connected (properly connected actually, but we shall be sloppy).
Since $k\geq 3$, $f$, $\partial f$ and $\partial M\subseteq M$ are
all (proper) $1${--}equivalences.

Now subdivide $(M,\partial M)$ until the chain map 
$C_\ast(M,\partial M)\ \to\ C_\ast(X,\partial X)$ is onto.
$C_\ast(X,\partial X)$ is $C_\ast(X,\partial X:\Lambda,F)$ for a
collection of paths $\Lambda$ and a lift functor $F$.
The tree for $X$ should come from a tree for $\partial M$,
which we can clearly assume.
$C_\ast(M,\partial M)$ is defined in the same way only with lift
functor $f^{-1}F$.
Let $D_\ast(f)$ be the kernel complex.

Then $H_r\bigl(D_\ast(f)\bigr)=0$ for $r<k$ and 
$H^r\bigl(D_\ast(f)\bigr)=0$ for $r>k$.
Now \fullRef{T.1.5.5} shows $H_k\bigl(D_\ast(f)\bigr)$
is an s{--}free\ tree module.
Doing surgery on trivial $(k-1)${--}spheres in $\partial M$ replaces $M$
by its boundary connected sum with a collection of $(S^k\times D^k)$'s.
Hence we may as well assume $H_k\bigl(D_\ast(f)\bigr)$ is free and based.
Let $\{ e_i\}$ be a preferred basis for this module.

By the Namioka Theorem,
$\Delta(f:\pi_{k+1})\ \to\ H_k\bigl(D_\ast(f)\bigr)$ is an
isomorphism.
Thus the $e_i$ determine classes in $\Delta(f:\pi_{k+1})$.
These in turn determine a proper regular homotopy class of immersions
$e_i\colon (D^k\times D^k,\partial D^k\times D^k)\ \to\ (M,\partial M)$.
We claim the $e_i$ are properly regularly homotopic to disjoint embeddings.
It is clearly enough to show this for the restricted immersions
$\bar e_i\colon (D^k,\partial D^k)\ \to\ (M,\partial M)$, for then we
just use small neighborhoods of the $\bar e_i$ to get the $e_i$.

The proof for the $\bar e_i$ proceeds as follows.
Let $C_j$ be an increasing sequence of compact subsets of $M$ 
with $\displaystyle\mathop{\cap}_j C_j=M$.
Let $C_j$ be such that any element of $\pi_1(M-C_j)$, when pushed into
$\pi_1(M-C_{j-1})$, lies in the image of $\pi_1\bigl(\partial M\ \cap\ 
(M- C_{j-1})\bigr)$ (compatible base points are understood).
We can do this as $\partial M\subseteq M$ 
is a proper $1${--}equivalence.

We now proceed.
Only a finite number of the $\bar e_i$ do not lie in $M-C_2$.
Embed these disjointly by the standard piping argument.

Again only finitely many $\bar e_i$ which do lie in $M-C_2$ do not
lie in $M-C_3$.
Put these in general position.
The intersections and self{--}intersections can be piped into
$\partial M \cap(M-C_1)$ without disturbing the $\bar e_i$ we
embedded in the previous step.
This follows from Milnor \cite{btwentyfour}, Theorem 6.6, where we
see that, to do the Whitney trick, 
we need only move one of the protagonists.
Hence we can always leave the $\bar e_i$ from the previous steps fixed.

Continuing in this fashion, we can always embed an $\bar e_i$ which
lies in $M-C_{j}$ but not in $M-C_{j+1}$, in $M-C_{j-1}$.
This gives us a proper regular homotopy and establishes our claim.

We next perform handle subtraction.
Let $N$ be obtained from $M$ by deleting the interiors 
of the images of the $e_i$.
Let $U$ be the union of the images of the $e_i$.
Let $\partial N=N\cap\partial M$.

By our construction, there is a chain map 
$C_\ast\colon U\cup\partial M,\partial M)\ \to\ D_\ast(f)$ such that
\[\begin{matrix}%
0\to&C_\ast(U\cup\partial M,\partial M)&\to&
C_\ast(M,\partial M)&\to& C_\ast(M, U\cup\partial M)&\to0\cr
&\downarrow&&\parallel&&\downarrow\cr
0\to&D_\ast(f)&\to&C_\ast(M,\partial M)&\to&
C_\ast(X,\partial X)&\to0\cr
\end{matrix}\]
chain homotopy commutes.
$C_\ast(U\cup\partial M,\partial M)$ has homology only in
dimension $k$ where it is $H_k\bigl(D_\ast(f)\bigr)$.
The map $C_\ast(U\cup\partial M,\partial M)\ \to\ D_\ast(f)$
gives this isomorphism in homology by construction.

Hence $C_\ast(M, U\cup\partial M)\ \to\ C_\ast(X,\partial X)$
is a chain equivalence.
Now $(N,\partial N)\subseteq (M, U\cup\partial M)$ is a proper
excision map, so $g\colon (N,\partial N)\ \to\ (X,\partial X)$ is
a proper homotopy equivalence from $n$ to $X$.
It induces proper homology isomorphisms on $\partial N\ \to\ \partial X$
and is thus a proper homotopy equivalence there since 
$\partial X\subseteq X$ is $1${--}connected.
Hence $g$ is a proper homotopy equivalence of pairs.
By adding an $h${--}cobordism to $\partial N$, we can achieve any
torsion we like for the map $g\colon N\ \to\ X$.

Notice we have not assumed $X$ is a simple Poincar\'e duality $n${--}ad,
but even so, the torsion of $g$ determines the torsions of the
remaining maps. 
We leave it to the reader to derive the standard formulas and
remark that if $g$ is simple and if $X$ is simple, 
then $g$ is a simple proper homotopy
equivalence of $n${--}ads.\footnote[1]{This bit replaces the original
discussion which had assumed $X$ was simple.}

\medskip\insetitem{\hskip\parindent Case2} $\dim X=2k+1$.
\smallskip
This time, \fullRef{T.3.1.1} permits us to suppose that $f$
induces $k${--}connected maps $M\ \to\ X$ and 
$\partial M\ \to\ \partial X$, and moreover we may assume 
$K_k(M,\partial M)=0$.
Hence we get a short exact sequence of modules
$0\to\ K_{k+1}(M,\partial M)\ \to\ K_k(\partial M)\ \to\ 
K_k(M)\ \to0$.
( $K_\ast(M)$ is the tree of modules which is the kernel of the map
$H_\ast\bigl(C(M:\Lambda, f^{-1}F)\bigr)\ \to\ 
H_\ast\bigl(C(X:\Lambda^\prime, F)\bigr)$.
The other $K${--}groups are defined similarly.)
\fullRef{T.1.5.5} now tells us that each of these modules is s{--}free.
As before we can perform surgery on trivial $(k+1)${--}spheres in
$\partial M$ to convert all of the above modules to free modules.
Again we get a locally finite collection of immersions
$\bar e_i\colon (D^{k+1},\partial D^{k+1})\ \to\ (M,\partial M)$ 
representing a basis of $K_{k+1}(M,\partial M)$.

We can no longer modify the $\bar e_i$ by a proper regular homotopy
to get disjoint embeddings (we could do this if $\partial M
\subseteq M$ were properly $2${--}connected) but by the same
sort of argument as in the first part, we can modify the $\bar e_i$
until $\bar e_i\vert_{\partial D^{k+1}}$ is a collection of disjoint
embeddings.

The rest of the proof is the same as Wall's.
We have represented a basis of $K_{k}(\partial M)$ by framed,
disjoint embeddings $S^k\ \to\ \partial M$.
Attach corresponding $k+1${--}handles to $M$, thus performing surgery.
Let $U$ be the union of the added handles, and let $(N,\partial N)$
be the new pair.
Since our spheres are null homotopic in $M$, $M$ is just replaced 
(up to proper homotopy type) by $M$ with $(k+1)${--}spheres
wedged on in a locally finite fashion.
Hence $K_k(N)$ is free, with a basis given by these spheres.

Dually, the exact sequence of the triple
$\partial N\subseteq \partial N\cup U\subseteq N$, reduces, 
using excision, to
\[0\to\ K_{k+1}(N,\partial N)\ \to\ 
K_{k+1}(M,\partial M)\ \to\ 
K_{k}(U, U\cap\partial N:M)\ \to\ 
K_{k}(N,\partial N)\ \to0\ .\]
The map $K_{k+1}(M,\partial M)\ \to\ 
K_{k}(U, U\cap \partial N:M)$ is seen to be zero since it factors as
$K_{k+1}(M,\partial M)\ \to\ K_k(\partial M)\ \to\ 
K_k(U:M)\ \to\ K_k(U, U\cap\partial N:M)$ and $K_k(U:M)$ is zero.
(Note that in this composition, $K_k(\partial M)$ should be a subspace
group, but such a group is isomorphic to the absolute group in our case.)
Since $K_k(U, U\cap\partial N:M)$ is free, so is $K_k(N,\partial N)$
and $K_{k+1}(N,\partial N)\cong K_{k+1}(M,\partial M)$.

The attached handles correspond to a basis of $K_{k+1}(M, \partial M)$,
so the map \\$K_{k+1}(N)\ \to\ K_{k+1}(M,\partial M)$ is an
epimorphism, since $K_{k+1}(N)$ is free and based on a set of generators
for $K_{k+1}(M,\partial M)$ and the map takes each basis element
to the corresponding generator.
But $K_{k+1}(M, \partial M)$ is free on these generators, so this map
is an isomorphism.
Hence $K_{k+1}(N)\ \to\ K_{k+1}(N, \partial N)$ is an isomorphism. 

Now, by Poincar\'e duality, $K^k(N, \partial N)\ \to\ K^k(N)$ is an
isomorphism.
The natural maps 
$K^k(N,\partial N)\ \to\ \bigl(K_k(N,\partial N)\bigr)^\ast$ and 
$K^k(N)\ \to\ \bigl(K_k(N)\bigr)^\ast$
are isomorphisms by \fullRef{C.1.5.4.2} since all the modules
are free.
Hence the map $K_k(N)\ \to\ K_k(N,\partial N)$ is an isomorphism.
Thus $K_k(\partial N)=0$, so $f$ restricted to $\partial N$ is a
proper homotopy equivalence.

Next choose a basis for $K_k(N)$ and perform surgery on it.
Write $P$ for the cobordism so obtained of $N$ to $N^\prime$ say.
Consider the induced map degree $1$ and Poincar\'e triads
$(P:N\cup(\partial N\times I), N^\prime)\ \to\ 
(X\times I:X\times 0\cup\partial X\times I,X\times 1)$.
We will identify $N\cup(\partial N\times I)$ with $N$.
In the exact sequence
\[0\to\ 
K_{k+1}(N)\ \to\ 
K_{k+1}(P)\ \to\ 
K_{k+1}(P,N))\ \RA{\ d\ }\ 
K_{k}(N)\ \to\ 
K_{k}(P)\ \to0\]
the map $d$ is by construction an isomorphism.
Hence $K_k(P)=0$ and $K_{k+1}(N)\ \to\ K_{k+1}(P)$
is an isomorphism.

The dual of $d$ is $K_{k+1}(N,\partial N)\ \to\ K_{k+1}(P,N^\prime)$,
so it is an isomorphism (the map is the map induced by the inclusion).
Now, since $f$ on $\partial N$ is a proper homotopy equivalence,
$K_{k+1}(N)\ \to\ K_{k+1}(N,\partial N)$ is an isomorphism.
$K_{k+1}(N)\ \to\ K_{k+1}(P)$ is an isomorphism, so
$K_{k+1}(P)\ \to\ K_{k+1}(P,N^\prime)$ is an isomorphism.

Thus in the sequence
\[0\to\ 
K_{k+1}(N^\prime)\ \to\ 
K_{k+1}(P)\ \to\ 
K_{k+1}(P,N^\prime)\ \to\ 
K_{k}(N^\prime)\ \to0\]
we have $K_{k+1}(N^\prime)=K_k(N^\prime)=0$, so
$N^\prime\ \to\ X$ is a proper homotopy equivalence.
$\partial N^\prime\ \to\ \partial X$ is the same as
$\partial N\ \to\ \partial X$ (i.e. we did nothing to $\partial N$ as all
our additions were in the interior of $N$) and therefore is a proper
homotopy equivalence.
Hence we have an equivalence of pairs.
The statement about torsions is proved the same way as for Case 1.
\end{proof}

\bigskip\begin{xRemarks}
Note that our proof is still valid in the case $\partial X=\partial_1X\ \cup\
\partial_2X$ provided $\partial_1M\ \to\ \partial_1X$ is a
proper homotopy equivalence (of pairs if $\partial_1X\ \cap\ \partial_2X
\neq\emptyset$) and $\partial_2X\subseteq X$ is a proper 
$1${--}equivalence ( $(X:\partial_1X,\partial_2X)$ should be a Poincar\'e 
triad).
The proof is word for word the same after we note that $K_i(\partial_2N)
\ \to\ K_i(\partial M)$ is always an isomorphism and that 
we may attach all our handles away from $\partial_1M$.
By induction we can prove a similar theorem for $n${--}ads, which is
the result we needed to prove \fullRef{T.2.2.13}.
\end{xRemarks}

\medskip
Our approach to surgery is to consider the surgery groups as bordism 
groups of surgery maps.
To make this approach work well, one needs a theorem like
\fullRef{T.3.1.3} below.

\bigskip\begin{xDefinition}
Given a Poincar\'e duality $n${--}ad, 
a {\sl surgery map\/} is a map $f\colon M\ \to\ X$
where $M$ is a \CAT{--}manifold $n${--}ad, 
$f$ is a degree $1$ map of $n${--}ads, 
and there is a \CAT{--}bundle $\nu$ over $X$ and a bundle map
$F\colon \nu_M\ \to\ \nu$ which covers $f$.
\end{xDefinition}

Given a locally finite CW $n${--}ad $K$ with a class $\wone\in H^1(K;\cy2)$, 
we say \hfill\penalty-10000
$M\ \RA{\ f\ }\ X\ \RA{\ g\ }\ K$ is a {\sl surgery map over
$(K,\wone)$\/} provided $g$ is a map of $n${--}ads with $g^\ast\wone$
equal to the first Stiefel{--}Whitney class of $X$, and provided
$f$ is a surgery map.

Two surgery maps over $(K,\wone)$ are said to be {\sl bordant\/}
(over $(K,\wone)$) if there is a surgery $(n+1)${--}ad
$W\ \to\ \RA{\ F\ }\ Y\ \RA{\ G\ }\ (K\times I,\wone)$
which is one of the surgery maps on $K\times 0$ and the other
on $K\times 1$.

\bigskip
\BEGIN{T.3.1.3}
Let $M\ \RA{\ f\ }\ X\ \RA{\ g\ }\ K$ be a surgery map over $(K,\wone)$,\
a $3${--}ad.
Suppose the formal dimension of $X$ is at least $6$.
Then, if $f\vert_{\partial_1M}$ is a proper homotopy equivalence,
and if $\partial_2K\subseteq K$ is a proper $1${--}equivalence, 
we can find another surgery map
$N\ \RA{\ h\ }\ Z\ \RA{\ i\ }\ K$ over $(K,\wone)$ with $h$
a proper homotopy equivalence of $3${--}ads, and with $i$
bordant over $(K,\wone)$ to $g$ so that over $\partial_1K\times I$
the bordism map is $\partial_1M\ \to\ \partial_1X$ crossed with $I$.
\end{Theorem}

\medskip\begin{proof}
If $\partial_2X\subseteq X$ were a proper $1${--}equivalence we could
finish easily using \fullRef{T.3.1.2}.
The proof then consists of modifying $X$ and $\partial_2X$ 
to get this condition.
The idea is to do surgery, first on $\partial_2X$ (and then on $X$)
to get $\partial_2X\ \to\ \partial_2K$ a proper $1${--}equivalence
(similarly for $X\ \to\ K$) and then show that we can cover these
surgeries on $\partial_2M$ and $M$.

Look at the map $g\colon \partial_2X\ \to\ \partial_2K$.
By \fullRef{T.2.3.2}, $\partial_2X$ can be replaced by
$L\ \cup\ H$, where $H$ is a manifold and $L$ satisfies $D(n-3)$,
where $n$ is the formal dimension of $X$.
This replacement does not alter the bordism class in which we are working.
Let $\wone$ also denote the restriction of $\wone\in H^1(K;\cy2)$
to $\partial_2K$.
Let $\nu$ be the line bundle over $\partial_2K$ classified by $\wone$.
Let $g\colon H\ \to\ \partial_2K$ denote the induced map.

Then $\tau_H\oplus g^\ast\nu$ is trivial, for $H$ has the homotopy type
of a $1${--}complex so the bundle is trivial \iff\ its first Stiefel{--}Whitney
class vanishes (and $\wone(\tau_H\oplus g^\ast\nu)=0$ by
construction).
Hence we can find a bundle map $F\colon \nu_H\ \to\ \nu$.

By \fullRef{T.3.1.1} we can add $1$ and $2$ handles to
to $H$ to get $W$ with 
$\partial W=H\ \cup\ H^\prime\ \cup\ \partial H\times I$ and a map
$G\colon W\ \to\ \partial_2K$ with $G\vert_H=h$ and 
$G\vert_{H^\prime}$ a proper $1${--}equivalence.
Let $Y= L\times I\ \cup\ W$ by gluing $\partial H\times I$ to $L\times I$
via the map $\partial H\ \to\ L$ crossed with $I$.
$(Y: L\ \cup_{\partial H}\ H, 
L\ \cup_{\partial H}\ H^\prime\ \cup\ L\times I)$ is a Poincar\'e duality triad.
This follows since $(L,\partial H)$ is a Poincar\'e duality pair and $Y$ is
$(L,\partial H)\times I$ glued to the manifold triad 
$(W: H, \partial H\times I, H^\prime)$ along $\partial H\times I$.
$(L,\partial H)\times I$ is a Poincar\'e triad by \fullRef{T.2.2.9},
and we can glue by \fullRef{T.2.1.3} and \fullRef{T.2.2.7}.

Let $Z=L\ \cup_{\partial H}\ H^\prime$.
We have a map of $Y\ \to\ K\times I$ given by $L\ \to\ K$ 
crossed with $I$ on $L\times I$ and by $W\ \to\ K\times I$ on $W$.
We claim the restriction $Z\ \to\ \partial_2K\times 1$ is a proper
$1${--}equivalence.

To see this, note first that 
\topD{0}{$\begin{matrix}%
\partial H&\subseteq&H^\prime\cr
\cap\hskip .2pt\vrule width .3pt height 6pt depth .1pt&&
\cap\hskip .2pt\vrule width .3pt height 6pt depth .1pt\cr
L&\subseteq&Z\cr
\end{matrix}$}{5}
is a pushout.
$\partial H\subseteq L$ is properly $1${--}connected by construction
(see \fullRef{T.2.3.2}).
It follows from a Mayer{--}Vietoris argument that $H^\prime \subseteq Z$ 
induces isomorphisms on $H^0_{{\text{ end}}}$ and $H^0$.
Since $\Delta(\partial H:\pi_1)\ \to\ \Delta(L:\pi_1)$ is onto, it
follows from a van{--}Kampen argument that
$\Delta(H^\prime:\pi_1)\ \to\ \Delta(Z:\pi_1)$ is onto.

Now consider $H^\prime\subseteq Z\ \to\ \partial_2K$.
The composite is a proper $1${--}equivalence by construction.
The first map is properly $1${--}connected, as we saw in the last
paragraph.
It then follows that $Z\ \to\ \partial_2K$ is a proper $1${--}equivalence.

It is easy to extend our bundle $\nu$ over all of $Y$.
Wall \cite{bfortyone} pages 89{--}90 shows how to cover our
surgeries back in $\partial_2M$.
One changes $f\colon\partial_2M\ \to\ \partial_2X$ through a proper
homotopy until it is transverse regular to all our core spheres in
$H\subseteq\partial_2X$.
The inverse image of a core sphere back in $\partial_2M$ will be a collection
of disjoint spheres, and Wall shows that, if we do our surgery
correctly on these spheres, then we can extend all our maps and
bundles.
Hence we get $F\colon P\ \to\ Y$ and a bundle map $\nu_P\ \to\ \nu$,
where $\nu$ is the extended $\nu$ over $Y$.

Thus our original problem $M\ \to\ X\ \to\ K$ is normally cobordant
over $(K,\wone)$ to a problem for which $\partial_2X\ \to\ \partial_2K$
is a proper $1${--}equivalence.
We have not touched $\partial_1M\ \to\ \partial_1X$ so we still have
that this map is a proper homotopy equivalence.
In fact, the part of $\partial P$ over $\partial_1M$ is just a product.

Now use \fullRef{T.2.3.3} on $X$ and proceed as above to
get a problem for which $X\ \to\ K$ is a proper $1${--}equivalence.
Note that we need never touch $\partial X$ so
$\partial_1M\ \to\ \partial_1X$ is still a proper homotopy equivalence
and $\partial_2X\ \to\ \partial_2K$ 
is still a proper $1${--}equivalence.
\end{proof}

\bigskip
\section{Paracompact surgery{--}patterns of application}
\newHead{III.2}

It has been noted by several people (see especially Quinn 
\cite{btwentynine} or \cite{bthirty}) that the theorems in section 1,
the $s${--}cobordism theorem, and transverse regularity are all the
geometry one needs to develop a great deal of the theory of surgery.

We define surgery groups as in Wall \cite{bfortyone} Chapter 9.%
\footnote[1]{%
As pointed out
by Farrell and Hsiang, \cite{boneone} (pages 102-103) there is
a problem with Wall's definition of geometric surgery groups
which will affect us as well.
Many people have used variants of this definition.
Ranicki's annotation of Wall's book
\cite{bonefour}, page 92, mentions the problem.
In \cite{bonethreeA} I show that Wall's original definition is almost right. 
Wall's problem is a cavalier treatment of local coefficients and the 
treatment presented in \cite{bonethreeA} works here as well.}
Let $K$ be a locally compact CW $n${--}ad, and let $w\in H^1(K;\cy2)$
be an orientation.
An {\sl object of type $n$ over $(K,\wone)$\/} is a surgery map $f$
(see section 1) over $s_n K$ for which, if
$M\ \to\ X\ \to\ s_n K$ is the surgery map,
$\varphi\colon \partial_n M\ \to\ \partial_n X$ is a proper homotopy
equivalence of $n${--}ads.%

We write $(\varphi,f)\sim 0$ to denote the existence of a surgery map
over $(s_n K,\wone)$ such that $\partial_{n+1}$ is $(\varphi,f)$;
i.e. if $W\ \to\ Z\ \to\ s_{n+1}s_n K$ is the surgery map,
$\partial_{n+1}W\ \to\ \partial_{n+1}Z\ \to\ s_n K$ is our
original surgery problem; and such that
$\partial_n$ is a proper homotopy equivalence of $(n+1)${--}ads.
$(\varphi,f)\sim(\varphi_1,f_1)$ provided
$(\varphi,f) + -(\varphi_1,f_1)\sim 0$, where $+$ denotes disjoint union
and $-(\varphi,f)$ denotes the same object but with the reverse
orientation.
Write $L^h_m(K,\wone)$ for the group of objects of type $n$ and
dimension $m$ (i.e. $m$ is the dimension of $M$) modulo the
relation $\sim$.
One checks that $\sim$ is an equivalence relation and that disjoint
union makes these sets into abelian groups.

If we require all the torsions of all the proper homotopy equivalences 
in the above definitions to be $0$ (including the torsions for the
Poincar\'e duality $n${--}ads), we get groups $L^s_m(K,\wone)$.
If $c\subseteq\sieb(K)$ is a subgroup closed under the involution
induced by the orientation $\wone$, then we get groups
$L^c_m(K,\wone)$ by requiring all the torsions (including those
for the Poincar\'e duality $n${--}ads) to lie in $c$.
( $\sieb(K)$ is Siebenmann's group of proper simple homotopy types;
see Chapter 1, section 5, or \cite{bthirtythree}).

\bigskip
\BEGIN{T.3.2.1}
Let $\alpha\in L^c_m(K,\wone)$, $n+m\geq 6$.
Then if $M\ \RA{\ \varphi\ }\ X\ \RA{ f\ }\ K$ is a representative
of $\alpha$ with $f$ a proper $1${--}equivalence, $\alpha=0$ \iff\
there is a normal cobordism $W\ \to\ X\times I$ with
$\partial_-W\ \to\ X\times0$ our original map $\varphi$, and
$\partial_+W\ \to\ X\times 1$ a proper homotopy equivalence of
$n${--}ads with torsions lying in $c$.
\end{Theorem}

\medskip\begin{proof}
Standard from \fullRef{T.3.1.2}, by doing surgery on the
boundary object.
\end{proof}

\bigskip
\BEGIN{T.3.2.2}
\[\cdots \to\ L^c_m(\partial_n K,\wone)\ \to\ 
L^c_m(\delta_n K,\wone)\ \to\ 
L^c_m(K,\wone)\ \to\ 
L^c_{m-1}(\partial_n K,\wone)\ \to\ \cdots\]
is exact.
\end{Theorem}

\medskip\begin{proof}
A standard argument.
\end{proof}

\bigskip
\BEGIN{T.3.2.3}
If $f\colon K_1\ \to\ K_2$ is a proper map of $n${--}ads, we get an 
induced map 
$L^c_m(K_1,f^\ast\wone)\ \to\ L^{c^\prime}_m(K_2,\wone)$ where
$f_\#(c)\subseteq c^\prime$, $f_\#\colon\sieb(K_1)\ \to\ \sieb(K_2)$
If $f$ is a proper $1${--}equivalence, the induced map is an isomorphism
for $c=f^{-1}_\#(c^\prime)$.
\end{Theorem}

\medskip\begin{proof}
The induced map is easily defined by
$M\ \to\ X\ \to\ K_1$ goes to
$M\ \to\ X\ \to\ K_1\ \RA{\ f\ }\ K_2$.
For the last statement, if $m\geq 5$ this is just \fullRef{T.3.1.3},
if $K_1$ and $K_2$ are $1${--}ads.
If $K_1$ and $K_2$ are $n${--}ads, an induction argument shows the 
result for $n+m\geq6$.
The result is actually true in all dimensions and a proof can be given
following Quinn's proof in the compact case ( see \cite{btwentynine} or
\cite{bthirty}).
We will not carry it out here.
\end{proof}

\bigskip
\BEGIN{T.3.2.4}
Let $K$ be a $1${--}ad, and let $M^m\ \RA{\ \varphi\ }\ 
X\ \RA{\ f\ }\ K$ be a surgery map over $(K,\wone)$ with
$\varphi$ a proper homotopy equivalence and with $f$ a proper
$1${--}equivalence.
Suppose given $\alpha\in L^c_{m+1}(K,\wone)$, $m\geq5$,
and suppose all the torsions for $\varphi$ lie in $c$.
Then there is an object of type $1$, $W\ \to\ X\times I\ \to\ K$,
over $(K,\wone)$ with $\partial W= M\ \cup\ N$, $N\to\ X\times 1$
a proper homotopy equivalence whose torsion also lies in $c$, 
and such that the surgery obstruction for this problem is $\alpha$.
\end{Theorem}

\medskip\begin{proof}
The proof is basically Quinn's (see \cite{btwentynine}).
Given $\alpha$, there is always an object of type $1$,
$P\ \to\ Z\ \to\ K$, whose obstruction is $-\alpha$.
(We may always assume $\partial P$ and $\partial Z$ are non{--}empty
by removing a disc from $Z$ and its inverse image in $P$, which we
can modify to be a disc.)
$M\times I\ \to\ X\times I\ \to\ K$ is also an object of type $1$
over $K$.

Take the boundary connected sum of $Z$ and $X\times I$ by extending
$\partial Z\ \#\ X\times 0$ (we may always assume $X$ and $Z$ are
in normal form so we may take this sum in their discs).
Similarly we may extend $\partial P\ \#\ M\times 0$.
We get a new object of type $1$,
$P\ \#_{M\times 0}\ M\times I\ \to\ 
Z\ \#_{X\times 0}\ X\times I\ \to\ K$.

By the proof of \fullRef{T.3.1.3}, we may do surgery on
$Z\ \#_{X\times0}\ X\times I$ until the map of it to $K$ is a proper
$1${--}equivalence, and we may cover this by a normal cobordism
of $P\ \#_{M\times0}\ M\times I$.
In doing this, we need never touch $M\times 1$ or $X\times 1$.
Let $P^\prime\ \to\ Z^\prime\ \to\ K$ denote this new
object of type $1$.
Note that it still has surgery obstruction $-\alpha$.

Now using \fullRef{T.3.1.2}, we can do surgery on
$\varphi\colon P^\prime\ \to\ Z^\prime$ where $Z^\prime$ is
considered to be the triad $(Z^\prime: X\times1,$ any other
boundary components$)$.
$\varphi$ restricted to the other boundary components is a
proper homotopy equivalence, so we may do surgery leaving them fixed
($X\times 1\subseteq Z^\prime$ is a proper $1${--}equivalence).
Let $W$ be the normal cobordism obtained over $M\times 1$.
Then $W\ \to\ X\times 1\times I$ is a surgery map,
$\partial_{-}W\ \to\ X\times 1\times 0$ is our old map, and
$\partial_+W\ \to\ X\times 1\times 1$ is a proper homotopy equivalence.
We can make all our torsions lie in $c$, and then the surgery obstruction
for $W\ \to\ X\times I\ \to\ K$ must be $\alpha$.
\end{proof}

\bigskip\begin{xDefinition}
Let $\sts_{\CAT}(X)$, for $X$ a Poincar\'e duality space of dimension $n$, be
the set of all simple, degree $1$, homotopy equivalences 
$\varphi\colon N^n\ \to\ X$ ($N$ a \CAT{--}manifold) modulo the
relation $\varphi\sim\psi$ \iff\ there is a \CAT{--}homeomorphism $h$
such that
\topD{0}{$\begin{matrix}%
N&\lower 10pt\hbox{$\searrow$}\hbox to 0pt{$\scriptstyle \varphi$\hss}\cr
\downlabeledarrow{}{h}&&X\cr
M&\raise 10pt\hbox{$\nearrow$}\hbox to 0pt{$\scriptstyle {\psi}$\hss}\cr
\end{matrix}$}{5}
properly homotopy commutes.
\end{xDefinition}

A similar definition holds for $X$ a Poincar\'e $n${--}ad.

\bigskip
\BEGIN{T.3.2.5}
There is an exact structure sequence 
\[\cdots\ \to\ \big[\Sigma X, F/\CAT\bigr]\ \to\ 
L^s_{m+1}(X,\wone)\ \to\ \sts_{\CAT}(X)\ \to\ 
\bigl[X, F/\CAT\bigr]\ \RA{\ \theta\ }\ L^s_{m}(X,\wone)\]
where $\wone$ is the first Stiefel{--}Whitney class of the Poincar\'e duality 
space $X$ with dimension of $X$ being $m\geq5$.
We also insist that the Spivak normal fibration of $X$ lift to a
\CAT{--}bundle.
By exactness we mean the following.
First of all $\sts_{\CAT}(X)$ may be empty, but in any case,
$\theta^{-1}(0)$ is the image of $\sts_{\CAT}(X)$.
If $\sts_{\CAT}(X)$ is not empty, then $L^s_{m+1}(X,\wone)$
acts on it, and two elements of $\sts_{\CAT}(X)$ which agree in 
$\bigl[ X, F/\CAT\bigr]$ differ by an element of this action.
The sequence continues infinitely to the left.
($\Sigma X$ is the ordinary suspension of $X$.)
\end{Theorem}

\medskip\begin{proof} 
See Wall \cite{bfortyone}, Chapter 10.
\end{proof}

\bigskip
\BEGIN{T.3.2.6}
Let $\bar{\hskip10pt}$ be the involution defined on $\sieb(K)$
in Chapter 1, section 5. Define 
$A_m(K,\wone) = H^m\bigl(\cy2, \sieb(K)\bigr)$ 
where $\sieb(K)$ is made into a $\cy2${--}module by the involution
$\bar{\hskip10pt}$ (which depends on $\wone$).
If $K$ is an $n${--}ad, then
\[\cdots \to\ 
A_{m+1}(K,\wone)\ \to\ 
L^s_m(K,\wone)\ \to\ 
L^h_m(K,\wone)\ \to\ 
A_m(K,\wone)\ \to\ \cdots\]
is exact for $m+n\geq 6$.
\end{Theorem}

\medskip\begin{proof}
The map $L^s\ \to\ L^h$ is just the forgetful map.
The map $L^h\ \to\ A$ just takes the torsion of the part of the boundary
that was a proper homotopy equivalence and maps it into $A$ 
(if the proper homotopy equivalence is over more than one
component, sum the torsions).
The map $A\ \to\ L^s$ takes a proper homotopy equivalence 
$M^m\ \to\ X$ whose torsion hits an element in $A_{m+1}$, and
maps it to the obstruction to surgering the map to a simple
homotopy equivalence.
See Shaneson \cite{bthirtyone} for the details of proving these
maps well{--}defined and the sequence exact.
\end{proof}

\medskip
\BEGIN{C.3.2.6.1}
If $A^c_m(K,w)= H^m(\cy2, c)$,
\[\cdots \to\ 
A^c_{m+1}(K,\wone)\ \to\ 
L^s_m(K,\wone)\ \to\ 
L^c_m(K,\wone)\ \to\ 
A^c_m(K,\wone)\ \to\ \cdots\]
is exact for $m+n\geq 6$.
\end{Corollary}

\medskip
We now produce our major computation.
\bigskip
\BEGIN{T.3.2.7}
\footnote[1]{\text{
Note added in proof\ $^2$:
Compare Maumary, \emph{The open surgery
obstruction in odd dimensions}}. Notices Amer. Math. Soc.
17 (number 5) p.848.\hfill\par $^2$Footnote $1$ is an original footnote. A better
reference is Maumary \cite{bonetwo}.}
Let $(L,\partial L)$ be a finite CW pair.
Form a new CW $n${--}ad $K$ by
$K=L\ \cup\ \partial L\times[0,\infty)$.
Suppose $\partial L$ is the disjoint union of subcomplexes
$\partial_i L$, $i=1$, \dots, $n$. 
Let $L^c_m(L,\wone)$ denote the Wall group for $L$
for homotopy equivalences which are simple over $L$ and
which, over $\partial_i L$, have torsions in $c_i$, where
$c_i=\ker\Bigl(\wh\bigl(\pi_1(\partial_i L)\bigr)\ \to\ 
\wh\bigl(\pi_1(L)\bigr)\Bigr)$.
Then there is an isomorphism
$L^c_m(L,\wone)\ \to\ L^s_m(K,\wone)$.
Combining this with Wall's long exact sequence we see
\[\cdots\zt
\mathop{\oplus}_{i=1}^n 
L^{c_i}_m\bigl(\pi_1(\partial_i L),\wone\bigr)
\to L^s_m\bigl(\pi_1(L),\wone\bigr)\to
L^s_m\bigl(K,\wone\bigr)\to
\mathop{\oplus}_{i=1}^n 
L^{c_i}_{m-1}\bigl(\pi_1(\partial_i L),\wone\bigr)\zt\cdots
\]
for $m\geq 7$.
\end{Theorem}

\medskip\begin{proof}
The map $L^s_m\bigl(\pi_1(L),\wone\bigr)\ \to\ 
L^s_m\bigl(K,\wone\bigr)$ is given by
$M\ \to\ X\ \to\ L$ goes to 
$M\ \cup\ \partial M\times[0,\infty)
\ \to\ 
X\ \cup\ \partial X\times[0,\infty)\ \to\ 
L\ \cup\ \partial L\times[0,\infty)$.

Siebenmann's thesis \cite{bthirtytwo} shows this map is a monomorphism.
To show that the map is onto we can assume
$W\ \RA{\ \varphi\ }\ Z\ \to\ K$ is a surgery map and that
$Z$ is a manifold using \fullRef{T.3.2.4} (this
representation theorem is also needed to show injectivity).
By Siebenmann \cite{bthirtytwo}, we can assume $Z$ is collared; i.e.
$\displaystyle Z=N\ \cup\ \bigl(\mathop{\cup}_{i=1}^n
\partial_i N\times[0,\infty)\bigr)$.
By making $\varphi$ transverse regular to the $\partial_i N$, we
get a problem over $L$, say $V\ \to\ N\ \to\ L$.
We claim
$\displaystyle 
V\ \cup\ \bigl(\mathop{\cup}_{i=1}^n \partial_i V\times[0,\infty)\bigr)
\ \to\ 
N\ \cup\ \bigl(\mathop{\cup}_{i=1}^n \partial_i N\times[0,\infty)\bigr)
$
has the same surgery obstruction as $W\ \to\ Z$.
But this is seen by actually constructing the normal cobordism 
using Siebenmann's concept of a $1${--}neighborhood and
some compact surgery.
\end{proof}

\medskip
\BEGIN{C.3.2.7.1}
We can improve $m\geq7$ to $m\geq6$.
\end{Corollary}

\smallskip\begin{proof}
Using recent work of Cappell{--}Shaneson \cite{bfive},
one can get a modified version of Siebenmann's main theorem.
One can not collar a $5${--}manifold, but one can at least get an
increasing sequence of cobordisms whose ends are
$\partial_i N\ \#\ S^2\times S^2\ \#\ \cdots\ \#\ S^2\times S^2$.
This is sufficient.
\end{proof}

\bigskip
Actually, one would hope that these surgery groups would be
periodic, just as the compact ones are.
This is actually the case, but the only proof I know involves describing
surgery in terms of algebra.
This can be done, but the result is long and will be omitted.

We briefly consider splitting theorems.
The two{--}sided codimension $1$ splitting theorem holds; i.e. if
$W$ has the simple homotopy type of 
$Z=(X,\partial X)\ \cup\ (Y,\partial X)$ with 
$\partial X\subseteq X$ a proper $1${--}equivalence, then the map
$W\ \to\ Z$ can be split.
The proof is the same as for the compact case.
Hence we also get codimension greater than or equal to $3$
splitting theorems for proper submanifolds. 
In fact, most of Wall \cite{bfortyone} Chapter 11
goes over with minor modifications. 

We are unable to obtain a one{--}sided splitting theorem in general,
due to the lack of a Farrell fibering theorem in the non{--}compact case.

We also note in passing that one could define surgery spaces as in
\cite{btwentynine} and \cite{bthirty}.
We than get the same basic geometric constructions; e.g. assembly maps
and pullback maps.
We have nothing new to add to the theory, so we leave the reader the 
exercise of restating \cite{btwentynine} so that it is valid
for paracompact surgery spaces\setcounter{footnote}{0}\footnote{We presumably 
could define proper algebraic bordism groups following Ranicki.
See for example \cite{Ranicki}.}.
\vfill

\newpage
\renewcommand{\rightmark}{Bibliography}

\vfill

\begin{bibdiv}
\begin{biblist}

\bib{bone}{book}{
   author={Bass, Hyman},
   title={Algebraic $K$-theory},
   publisher={W. A. Benjamin, Inc., New York-Amsterdam},
   date={1968},
   pages={xx+762},
   review={\MR{0249491}},
}


\bib{btwo}{book}{
   author={Bredon, Glen E.},
   title={Sheaf theory},
   publisher={McGraw-Hill Book Co., New York-Toronto, Ont.-London},
   date={1967},
   pages={xi+272},
   review={\MR{0221500}},
 }

\bib{bthree}{misc}{
      author={Browder, W.},
       title={Cap products and Poincar\'e duality},
        date={1964},
        note={Cambridge University notes \BR{Browder}},
}

\bib{bfour}{misc}{
      author={Browder, W.},
       title={Surgery on simply{--}connected manifolds},
        date={1969},
        note={Princeton University notes \BR{Browder}},
}

\bib{bfive}{unpublished}{
      author={Cappell, S.},
      author={Shaneson, J.},
       title={On four dimensional surgery and applications},
       note={preprint \BR{CS}},
}

\bib{bsix}{book}{
    author={Cooke, George E.},
   author={Finney, Ross L.},
   title={Homology of cell complexes},
   series={Based on lectures by Norman E. Steenrod},
   publisher={Princeton University Press, Princeton, N.J.; University of
   Tokyo Press, Tokyo},
   date={1967},
   pages={xv+256},
   review={\MR{0219059}},
}

\bib{bseven}{article}{
   author={Dold, Albrecht},
   title={Partitions of unity in the theory of fibrations},
   journal={Ann. of Math. (2)},
   volume={78},
   date={1963},
   pages={223--255},
   issn={0003-486X},
   review={\MR{0155330}},
}

\bib{beight}{book}{
   author={Dold, Albrecht},
   title={Halbexakte Homotopiefunktoren},
   language={German},
   series={Lecture Notes in Mathematics},
   volume={12},
   publisher={Springer-Verlag, Berlin-New York},
   date={1966},
   pages={156},
   review={\MR{0198464}},
}


\bib{bnine}{article}{
      author={Farrell, T.},
      author={Wagoner, J.},
       title={A torsion invariant for proper $h${--}cobordisms/ an algebraic criterion for a map to be a proper homotopy equivalence},
        date={1969},
     journal={Notices Amer. Math. Soc.},
      volume={16},
       pages={988 and 1090},
      note={\BR{FWB}},
}
\bib{bten}{unpublished}{
      author={Farrell, T.},
      author={Wagoner, J.},
       title={Infinite matrices in algebraic $K${--}theory and topology},
       note={To appear: \BR{FWB}},
}

\bib{beleven}{unpublished}{
      author={Farrell, T.},
      author={Taylor, L.},
      author={Wagoner, J.},
       title={The {W}hitehead theorem in the proper category},
       note={To appear: \BR{FWA}},
}


\bib{btwelve}{article}{
   author={Haefliger, Andr{\'e}},
   title={Lissage des immersions. I},
   language={French},
   journal={Topology},
   volume={6},
   date={1967},
   pages={221--239},
   issn={0040-9383},
   review={\MR{0208607}},
 }

\bib{bthirteen}{book}{
   author={Hilton, P. J.},
   title={An introduction to homotopy theory},
   series={Cambridge Tracts in Mathematics and Mathematical Physics, no. 43},
   publisher={Cambridge, at the University Press},
   date={1953},
   pages={viii+142},
   review={\MR{0056289}},
 }

\bib{bfourteen}{article}{
   author={Hirsch, Morris W.},
   title={Immersions of manifolds},
   journal={Trans. Amer. Math. Soc.},
   volume={93},
   date={1959},
   pages={242--276},
   issn={0002-9947},
   review={\MR{0119214}},
}

\bib{bfifteen}{book}{
   author={Hu, Sze-tsen},
   title={Theory of retracts},
   publisher={Wayne State University Press, Detroit},
   date={1965},
   pages={234},
   review={\MR{0181977}},
 }

\bib{bsixteen}{book}{
   author={Hurewicz, Witold},
   author={Wallman, Henry},
   title={Dimension Theory},
   series={Princeton Mathematical Series, v. 4},
   publisher={Princeton University Press, Princeton, N. J.},
   date={1941},
   pages={vii+165},
   review={\MR{0006493}},
}

\bib{bseventeen}{book}{
  author={Kelley, John L.},
   title={General topology},
   publisher={D. Van Nostrand Company, Inc., Toronto-New York-London},
   date={1955},
   pages={xiv+298},
   review={\MR{0070144}},
}



\bib{beighteen}{article}{
      author={Kirby, R.},
      author={Siebenmann, L.},
       title={Foundations of topology},
        date={1969},
     journal={Notices Amer. Math. Soc.},
      volume={16},
      note={\BR{KS}},
}

\bib{bnineteen}{article}{
   author={Lees, J. Alexander},
   title={Immersions and surgeries of topological manifolds},
   journal={Bull. Amer. Math. Soc.},
   volume={75},
   date={1969},
   pages={529--534},
   issn={0002-9904},
   review={\MR{0239602}},
}

\bib{btwenty}{book}{
   author={Lefschetz, Solomon},
   title={Introduction to Topology},
   series={Princeton Mathematical Series, vol. 11},
   publisher={Princeton University Press, Princeton, N. J.},
   date={1949},
   pages={viii+218},
   review={\MR{0031708}},
}

\bib{btwentyone}{book}{
      author={Lundell, A.},
      author={Weingram, S.},
       title={The topology of CW complexes},
   publisher={Van Nostrand Reinhold Co.},
     address={New York},
        date={1969},
}

\bib{btwentytwo}{article}{
   author={Milnor, John},
   title={On spaces having the homotopy type of a ${\rm CW}$-complex},
   journal={Trans. Amer. Math. Soc.},
   volume={90},
   date={1959},
   pages={272--280},
   issn={0002-9947},
   review={\MR{0100267}},
}

\bib{btwentythree}{article}{
   author={Milnor, J.},
   title={Whitehead torsion},
   journal={Bull. Amer. Math. Soc.},
   volume={72},
   date={1966},
   pages={358--426},
   issn={0002-9904},
   review={\MR{0196736}},
}

\bib{btwentyfour}{book}{
   author={Milnor, John},
   title={Lectures on the $h$-cobordism theorem},
   series={Notes by L. Siebenmann and J. Sondow},
   publisher={Princeton University Press, Princeton, N.J.},
   date={1965},
   pages={v+116},
   review={\MR{0190942}},
}

\bib{btwentyfive}{book}{
   author={Mitchell, Barry},
   title={Theory of categories},
   series={Pure and Applied Mathematics, Vol. XVII},
   publisher={Academic Press, New York-London},
   date={1965},
   pages={xi+273},
   review={\MR{0202787}},
}

\bib{btwentysix}{article}{
   author={Moore, John C.},
   title={Some applications of homology theory to homotopy problems},
   journal={Ann. of Math. (2)},
   volume={58},
   date={1953},
   pages={325--350},
   issn={0003-486X},
   review={\MR{0059549}},
}
\vfill
\bib{btwentyseven}{book}{
   author={Nagami, Kei{\^o}},
   title={Dimension theory},
   series={With an appendix by Yukihiro Kodama. Pure and Applied
   Mathematics, Vol. 37},
   publisher={Academic Press, New York-London},
   date={1970},
   pages={xi+256},
   review={\MR{0271918}},
}

\bib{btwentyeight}{article}{
    author={Namioka, I.},
   title={Maps of pairs in homotopy theory},
   journal={Proc. London Math. Soc. (3)},
   volume={12},
   date={1962},
   pages={725--738},
   issn={0024-6115},
   review={\MR{0144345}},
}

\bib{btwentynine}{thesis}{
      author={Quinn, F.},
       title={A geometric formulation of surgery},
        type={Ph.D. Thesis},
        date={1969},
}

\bib{bthirty}{inproceedings}{
   author={Quinn, Frank},
   title={A geometric formulation of surgery},
   conference={
      title={Topology of Manifolds (Proc. Inst., Univ. of Georgia, Athens,
      Ga., 1969)},
   },
   book={
      publisher={Markham, Chicago, Ill.},
   },
   date={1970},
   pages={500--511},
   review={\MR{0282375}},
}

\bib{bthirtyone}{article}{
    author={Shaneson, Julius L.},
   title={Wall's surgery obstruction groups for $G\times Z$},
   journal={Ann. of Math. (2)},
   volume={90},
   date={1969},
   pages={296--334},
   issn={0003-486X},
   review={\MR{0246310}},
}


\bib{bthirtytwo}{thesis}{
      author={Siebenmann, L.},
       title={The obstruction to finding a boundary for an open manifold of dimension greater than five},
        type={Ph.D. Thesis},
        date={1965},
}

\bib{bthirtythree}{article}{
   author={Siebenmann, L. C.},
   title={Infinite simple homotopy types},
   journal={Nederl. Akad. Wetensch. Proc. Ser. A 73 = Indag. Math.},
   volume={32},
   date={1970},
   pages={479--495},
   review={\MR{0287542}},
}

\bib{bthirtyfour}{unpublished}{
      author={Siebenmann, L.},
       title={A report on topological manifolds},
      note={To appear in {\sl Proc. ICM} (Nice 1970), \BR{Sieb}},
}

\bib{bthirtyfive}{book}{
   author={Spanier, Edwin H.},
   title={Algebraic topology},
   publisher={McGraw-Hill Book Co., New York-Toronto, Ont.-London},
   date={1966},
   pages={xiv+528},
   review={\MR{0210112}},
}

\bib{bthirtysix}{article}{
   author={Spivak, Michael},
   title={Spaces satisfying Poincar\'e duality},
   journal={Topology},
   volume={6},
   date={1967},
   pages={77--101},
   issn={0040-9383},
   review={\MR{0214071}},
}

\bib{bthirtyseven}{article}{
   author={Wall, C. T. C.},
   title={Finiteness conditions for ${\rm CW}$-complexes},
   journal={Ann. of Math. (2)},
   volume={81},
   date={1965},
   pages={56--69},
   issn={0003-486X},
   review={\MR{0171284}},
}

\bib{bthirtyeight}{article}{
   author={Wall, C. T. C.},
   title={Finiteness conditions for ${\rm CW}$ complexes. II},
   journal={Proc. Roy. Soc. Ser. A},
   volume={295},
   date={1966},
   pages={129--139},
   review={\MR{0211402}},
}

\bib{bforty}{article}{
   author={Wall, C. T. C.},
   title={Surgery of non-simply-connected manifolds},
   journal={Ann. of Math. (2)},
   volume={84},
   date={1966},
   pages={217--276},
   issn={0003-486X},
   review={\MR{0212827}},
}

\bib{bthirtynine}{article}{
   author={Wall, C. T. C.},
   title={Poincar\'e complexes. I},
   journal={Ann. of Math. (2)},
   volume={86},
   date={1967},
   pages={213--245},
   issn={0003-486X},
   review={\MR{0217791}},
}

\bib{bfortyone}{book}{
   author={Wall, C. T. C.},
   title={Surgery on compact manifolds},
   note={London Mathematical Society Monographs, No. 1},
   publisher={Academic Press, London-New York},
   date={1970},
   pages={x+280},
   review={\MR{0431216}},
}

\bib{bfortytwo}{article}{
   author={Whitehead, George W.},
   title={On spaces with vanishing low-dimensional homotopy groups},
   journal={Proc. Nat. Acad. Sci. U. S. A.},
   volume={34},
   date={1948},
   pages={207--211},
   issn={0027-8424},
   review={\MR{0028027}},
}

\bib{bfortythree}{article}{
   author={Whitehead, J. H. C.},
   title={Combinatorial homotopy. I},
   journal={Bull. Amer. Math. Soc.},
   volume={55},
   date={1949},
   pages={213--245},
   issn={0002-9904},
   review={\MR{0030759}},
}

\bib{bfortyfour}{book}{
   author={Wilder, Raymond Louis},
   title={Topology of Manifolds},
   series={American Mathematical Society Colloquium Publications, vol. 32},
   publisher={American Mathematical Society, New York, N. Y.},
   date={1949},
   pages={ix+402},
   review={\MR{0029491}},
 }

\newpage
\centerline{\bf Additional References}
\bib{Browder}{book}{
   author={Browder, William},
   title={Surgery on simply-connected manifolds},
   note={Ergebnisse der Mathematik und ihrer Grenzgebiete, Band 65},
   publisher={Springer-Verlag},
   place={New York},
   date={1972},
   pages={ix+132},
   review={\MR{0358813} (50 \#11272)},
}

\bib{CS}{article}{
   author={Cappell, Sylvain E.},
   author={Shaneson, Julius L.},
   title={On four dimensional surgery and applications},
   journal={Comment. Math. Helv.},
   volume={46},
   date={1971},
   pages={500--528},
   issn={0010-2571},
   review={\MR{0301750} (46 \#905)},
}

\bib{boneone}{article}{
    author={Farrell, F. T.},
   author={Hsiang, W. C.},
   title={Rational $L$-groups of Bieberbach groups},
   journal={Comment. Math. Helv.},
   volume={52},
   date={1977},
   number={1},
   pages={89--109},
   issn={0010-2571},
   review={\MR{0448372}},
 }

\bib{FTW}{article}{
   author={Farrell, F. T.},
   author={Taylor, L. R.},
   author={Wagoner, J. B.},
   title={The Whitehead theorem in the proper category},
   journal={Compositio Math.},
   volume={27},
   date={1973},
   pages={1--23},
   issn={0010-437X},
   review={\MR{0334226} (48 \#12545)},
}

\bib{FWA}{article}{
   author={Farrell, F. T.},
   author={Wagoner, J. B.},
   title={Infinite matrices in algebraic $K$-theory and topology},
   journal={Comment. Math. Helv.},
   volume={47},
   date={1972},
   pages={474--501},
   issn={0010-2571},
   review={\MR{0319185} (47 \#7731a)},
}

\bib{FWB}{article}{
   author={Farrell, F. T.},
   author={Wagoner, J. B.},
   title={Algebraic torsion for infinite simple homotopy types},
   journal={Comment. Math. Helv.},
   volume={47},
   date={1972},
   pages={502--513},
   review={\MR{0319185} (47 \#7731b)},
}
	
\bib{KS}{book}{
   author={Kirby, Robion C.},
   author={Siebenmann, Laurence C.},
   title={Foundational essays on topological manifolds, smoothings, and triangulations},
   note={With notes by John Milnor and Michael Atiyah;
   Annals of Mathematics Studies, No. 88},
   publisher={Princeton University Press},
   place={Princeton, N.J.},
   date={1977},
   pages={vii+355},
   review={\MR{0645390} (58 \#31082)},
}

\bib{bonetwo}{inproceedings}{
   author={Maumary, Serge},
   title={Proper surgery groups and Wall-Novikov groups},
   conference={
      title={Algebraic $K$-theory, III: Hermitian $K$-theory and geometric
      applications (Proc. Conf., Battelle Memorial Inst., Seatlle, Wash.,
      1972)},
   },
   book={
      publisher={Springer, Berlin},
   },
   date={1973},
   pages={526--539. Lecture Notes in Math, Vol. 343},
   review={\MR{0377938}},
 }

\bib{Ranicki}{article}{
   author={Ranicki, Andrew},
   title={Algebraic Poincar\'e cobordism},
   conference={
      title={Topology, geometry, and algebra: interactions and new
      directions },
      address={Stanford, CA},
      date={1999},
   },
   book={
      series={Contemp. Math.},
      volume={279},
      publisher={Amer. Math. Soc., Providence, RI},
   },
   date={2001},
   pages={213--255},
   review={\MR{1850750}},
   doi={10.1090/conm/279/04563},
}
\bib{Sieb}{article}{
   author={Siebenmann, L. C.},
   title={Topological manifolds},
   conference={
      title={Actes du Congr\`es International des Math\'ematiciens},
      address={Nice},
      date={1970},
   },
   book={
      publisher={Gauthier-Villars},
      place={Paris},
   },
   date={1971},
   pages={133--163},
   review={\MR{0423356} (54 \#11335)},
}


\bib{bonethree}{inproceedings}{
      author={Taylor, L.~R.},
       title={Surgery groups and inner automorphisms},
        date={1973},
   booktitle={{A}lgebraic $K$-theory {III}: Hermitian $K$-theory and
 Geometric Applications (proc. conf., Battelle Memorial Inst., Seattle, Wash.,
  1972)},
      series={{L}ecture {N}otes in {M}ath.},
      volume={343},
   publisher={{S}pringer{--}{V}erlag},
     address={{B}erlin},
       pages={471\ndash 477},
}


\bib{bonethreeA}{article}{
   author={Taylor, Laurence R.},
   title={Unoriented geometric functors},
   journal={Forum Math.},
   volume={20},
   date={2008},
   number={3},
   pages={457--467},
   issn={0933-7741},
   review={\MR{2418201}},
   doi={10.1515/FORUM.2008.023},
   eprint={{https://arxiv.org/abs/math/0606651}{arXiv:math/0606651v1}},
}

\bib{bonefour}{book}{
   author={Wall, C. T. C.},
   title={Surgery on compact manifolds},
   series={Mathematical Surveys and Monographs},
   volume={69},
   edition={2},
   note={Edited and with a foreword by A. A. Ranicki},
   publisher={American Mathematical Society, Providence, RI},
   date={1999},
   pages={xvi+302},
   isbn={0-8218-0942-3},
   review={\MR{1687388}},
   doi={10.1090/surv/069},
}
\newpage
\centerline{\bf Some Papers Citing This Work}

\bib{MR2014912}{article}{
   author={Ayala, R.},
   author={C{\'a}rdenas, M.},
   author={Muro, F.},
   author={Quintero, A.},
   title={An elementary approach to the projective dimension in proper
   homotopy theory},
   journal={Comm. Algebra},
   volume={31},
   date={2003},
   number={12},
   pages={5995--6017},
   issn={0092-7872},
   review={\MR{2014912}},
   doi={10.1081/AGB-120024863},
}
\bib{MR1848146}{book}{
   author={Baues, Hans-Joachim},
   author={Quintero, Antonio},
   title={Infinite homotopy theory},
   series={$K$-Monographs in Mathematics},
   volume={6},
   publisher={Kluwer Academic Publishers, Dordrecht},
   date={2001},
   pages={viii+296},
   isbn={0-7923-6982-3},
   review={\MR{1848146}},
   doi={10.1007/978-94-009-0007-3},
}
\bib{MR1758300}{article}{
   author={Block, Jonathan},
   author={Weinberger, Shmuel},
   title={Arithmetic manifolds of positive scalar curvature},
   journal={J. Differential Geom.},
   volume={52},
   date={1999},
   number={2},
   pages={375--406},
   issn={0022-040X},
   review={\MR{1758300}},
}
\bib{MR671656}{article}{
   author={Farrell, F. T.},
   author={Hsiang, W. C.},
   title={The stable topological-hyperbolic space form problem for complete
   manifolds of finite volume},
   journal={Invent. Math.},
   volume={69},
   date={1982},
   number={1},
   pages={155--170},
   issn={0020-9910},
   review={\MR{671656}},
   doi={10.1007/BF01389189},
}
\bib{MR973309}{article}{
   author={Farrell, F. T.},
   author={Jones, L. E.},
   title={A topological analogue of Mostow's rigidity theorem},
   journal={J. Amer. Math. Soc.},
   volume={2},
   date={1989},
   number={2},
   pages={257--370},
   issn={0894-0347},
   review={\MR{973309}},
   doi={10.2307/1990978},
}
\bib{MR1388311}{article}{
   author={Ferry, Steven C.},
   author={Pedersen, Erik K.},
   title={Epsilon surgery theory},
   conference={
      title={Novikov conjectures, index theorems and rigidity, Vol.\ 2 },
      address={Oberwolfach},
      date={1993},
   },
   book={
      series={London Math. Soc. Lecture Note Ser.},
      volume={227},
      publisher={Cambridge Univ. Press, Cambridge},
   },
   date={1995},
   pages={167--226},
   review={\MR{1388311}},
   doi={10.1017/CBO9780511629365.007},
}
\bib{MR666159}{article}{
   author={Freedman, Michael H.},
   title={A surgery sequence in dimension four;\ the relations with knot
   concordance},
   journal={Invent. Math.},
   volume={68},
   date={1982},
   number={2},
   pages={195--226},
   issn={0020-9910},
   review={\MR{666159}},
   doi={10.1007/BF01394055},
}
\bib{MR872483}{article}{
   author={Freedman, Michael H.},
   title={A geometric reformulation of $4$-dimensional surgery},
   note={Special volume in honor of R. H. Bing (1914--1986)},
   journal={Topology Appl.},
   volume={24},
   date={1986},
   number={1-3},
   pages={133--141},
   issn={0166-8641},
   review={\MR{872483}},
   doi={10.1016/0166-8641(86)90054-4},
}
\bib{MR689390}{article}{
   author={Hambleton, Ian},
   title={Projective surgery obstructions on closed manifolds},
   conference={
      title={Algebraic $K$-theory, Part II},
      address={Oberwolfach},
      date={1980},
   },
   book={
      series={Lecture Notes in Math.},
      volume={967},
      publisher={Springer, Berlin-New York},
   },
   date={1982},
   pages={101--131},
   review={\MR{689390}},
}
\bib{MR1268592}{article}{
   author={Hambleton, Ian},
   author={Madsen, Ib},
   title={On the computation of the projective surgery obstruction groups},
   journal={$K$-Theory},
   volume={7},
   date={1993},
   number={6},
   pages={537--574},
   issn={0920-3036},
   review={\MR{1268592}},
   doi={10.1007/BF00961217},
}
\bib{MR2212277}{article}{
   author={Jahren, Bj{\o}rn},
   author={Kwasik, S{\l}awomir},
   title={Manifolds homotopy equivalent to $\bold R{\rm P}^4\#\bold R{\rm
   P}^4$},
   journal={Math. Proc. Cambridge Philos. Soc.},
   volume={140},
   date={2006},
   number={2},
   pages={245--252},
   issn={0305-0041},
   review={\MR{2212277}},
   doi={10.1017/S0305004105008893},
}
\bib{MR520500}{article}{
   author={Jones, Lowell},
   title={The nonsimply connected characteristic variety theorem},
   conference={
      title={Algebraic and geometric topology},
      address={Proc. Sympos. Pure Math., Stanford Univ., Stanford, Calif.},
      date={1976},
   },
   book={
      series={Proc. Sympos. Pure Math., XXXII},
      publisher={Amer. Math. Soc., Providence, R.I.},
   },
   date={1978},
   pages={131--140},
   review={\MR{520500}},
}
\bib{MR869712}{article}{
   author={Kwasik, S{\l}awomir},
   title={On periodicity in topological surgery},
   journal={Canad. J. Math.},
   volume={38},
   date={1986},
   number={5},
   pages={1053--1064},
   issn={0008-414X},
   review={\MR{869712}},
   doi={10.4153/CJM-1986-051-6},
}
\bib{MR1191376}{article}{
   author={Kwasik, S{\l}awomir},
   author={Schultz, Reinhard},
   title={Vanishing of Whitehead torsion in dimension four},
   journal={Topology},
   volume={31},
   date={1992},
   number={4},
   pages={735--756},
   issn={0040-9383},
   review={\MR{1191376}},
   doi={10.1016/0040-9383(92)90005-3},
}
\bib{MR1451758}{article}{
   author={Kwasik, S{\l}awomir},
   author={Schultz, Reinhard},
   title={Inductive detection for homotopy equivalences of manifolds},
   journal={$K$-Theory},
   volume={11},
   date={1997},
   number={3},
   pages={287--306},
   issn={0920-3036},
   review={\MR{1451758}},
   doi={10.1023/A:1007702426750},
}
\bib{MR3349793}{article}{
   author={Kwasik, S{\l}awomir},
   author={Schultz, Reinhard},
   title={Tangential thickness of manifolds},
   journal={Proc. Lond. Math. Soc. (3)},
   volume={110},
   date={2015},
   number={5},
   pages={1281--1313},
   issn={0024-6115},
   review={\MR{3349793}},
   doi={10.1112/plms/pdv009},
}
\bib{MR579574}{article}{
   author={Pedersen, Erik Kjaer},
   author={Ranicki, Andrew},
   title={Projective surgery theory},
   journal={Topology},
   volume={19},
   date={1980},
   number={3},
   pages={239--254},
   issn={0040-9383},
   review={\MR{579574}},
   doi={10.1016/0040-9383(80)90010-5},
}
\bib{MR1747539}{article}{
   author={Pedersen, Erik Kj{\ae}r},
   title={Continuously controlled surgery theory},
   conference={
      title={Surveys on surgery theory, Vol. 1},
   },
   book={
      series={Ann. of Math. Stud.},
      volume={145},
      publisher={Princeton Univ. Press, Princeton, NJ},
   },
   date={2000},
   pages={307--321},
   review={\MR{1747539}},
}
\bib{MR1944688}{article}{
   author={Pedersen, Erik Kj{\ae}r},
   title={Controlled methods in equivariant topology, a survey},
   conference={
      title={Current trends in transformation groups},
   },
   book={
      series={$K$-Monogr. Math.},
      volume={7},
      publisher={Kluwer Acad. Publ., Dordrecht},
   },
   date={2002},
   pages={231--245},
   review={\MR{1944688}},
   doi={10.1007/978-94-009-0003-5-15},
}
\bib{MR1818774}{article}{
   author={Weiss, Michael},
   author={Williams, Bruce},
   title={Automorphisms of manifolds},
   conference={
      title={Surveys on surgery theory, Vol. 2},
   },
   book={
      series={Ann. of Math. Stud.},
      volume={149},
      publisher={Princeton Univ. Press, Princeton, NJ},
   },
   date={2001},
   pages={165--220},
   review={\MR{1818774}},
}

\end{biblist}
\end{bibdiv}
\end{document}